\RequirePackage{fix-cm}
\documentclass[]{article}
\usepackage{amsmath,amsthm,doi,url,hyperref}
\hypersetup{colorlinks=true, linkcolor=blue!80!black, urlcolor=blue!80!black, citecolor=blue!80!black}
\usepackage{fullpage}
\usepackage{authblk}
\newcommand{\keywords}[1]{Keywords: \textit{#1}}
\usepackage[caption=false]{subfig}%
\usepackage[numbers,sort&compress]{natbib}%
\bibpunct[, ]{[}{]}{,}{n}{,}{,}%

\theoremstyle{plain}%

\theoremstyle{definition}

\theoremstyle{remark}

\usepackage{graphicx}

\usepackage{setspace}

\usepackage{thumbpdf,lmodern}
\usepackage{amsfonts,amssymb}
\usepackage{amsmath}
\usepackage{framed}
\usepackage{xcolor}
\usepackage{hyperref}
\hypersetup{breaklinks=true, colorlinks=true, linkcolor={blue!90!black}, citecolor={blue!90!black}, urlcolor={blue!90!black}}

\usepackage[ruled]{algorithm2e}
\usepackage{algpseudocode}
\usepackage{algorithmicx}
\usepackage{thumbpdf} %

\newcommand{\vecx}{\mathbf{x}}
\newcommand{\s}{\mathbf{s}}
\newcommand{\R}{\mathbb{R}}
\newcommand{\veck}{\mathbf{k}}
\newcommand{\vecY}{\mathbf{y}}
\newcommand{\vecYu}{\bar{\mathbf{y}}}

\newcommand{\Lan}{\boldsymbol{\Lambda}_n}

\newcommand{\LaN}{\boldsymbol{\Lambda}_N}
\newcommand{\An}{\mathbf{A}_n}

\newcommand{\Kn}{\mathbf{K}_n}
\newcommand{\KN}{\mathbf{K}_N}

\newcommand{\Nmax}{N_\textrm{max}}
\newcommand{\pmax}{p_\textrm{max}}

\newcommand{\gammai}{\gamma_{ \rm inc} }
\newcommand{\gammad}{\gamma_{ \rm dec} }
\newcommand{\Var}[1]{\mathbb{V}\left[{#1}\right]}
\newcommand{\Esp}[1]{\mathbb{E}\left[{#1}\right]}
\newcommand{\Cov}[1]{\mathbb{C}ov\left[{#1}\right]}

\usepackage{thmtools}
\usepackage{nameref}
\usepackage[nameinlink]{cleveref}

\crefname{defi}{definition}{definitions}
\Crefname{defi}{Definition}{Definitions}
\crefname{lemma}{lemma}{lemmas}
\Crefname{lemma}{Lemma}{Lemmas}
\crefname{assumption}{assumption}{assumptions}
\Crefname{assumption}{Assumption}{Assumptions}

\newcommand{\defined}{\triangleq}

\begin{document}

\title{Adaptive Replication Strategies in Trust-Region-Based Bayesian Optimization of Stochastic Functions}

\author[a]{Micka\"el Binois}
\author[b]{Jeffrey Larson}

\affil[a]{Inria, Université Côte d’Azur, CNRS, LJAD, Sophia Antipolis, France. {\tt mickael.binois@inria.fr}}
\affil[b]{Mathematics and Computer Science Division, Argonne National Laboratory.  {\tt jmlarson@anl.gov}}
\date{}

\maketitle

\begin{abstract} 

We develop and analyze a method for stochastic simulation optimization based
on Gaussian process models within a trust-region framework. We focus on
settings where the variance of the objective function is large, making accurate
estimation challenging and often requiring many evaluations. To address this
regime, we combine local modeling with adaptive replication, allowing the
method to allocate repeated evaluations where they are most beneficial.
We introduce several mechanisms to promote and adapt replication, including
modifications to the acquisition function and cost-aware evaluation strategies.
These components enable our approach to scale effectively when high levels of
sampling are required to reduce noise.
We refer to the resulting method as
OGPIT, for Optimization by Gaussian Processes In Trust regions.
Numerical experiments show that adaptive
replication can substantially improve computational efficiency while preserving solution accuracy
compared to baseline methods, in particular
when evaluation costs are taken into account.
\end{abstract}


\section{Introduction}

We consider the optimization of a stochastic oracle or function $y: D \subset \R^d
\mapsto \R$, where evaluations of $y$ are accessible only through noise-corrupted evaluations, $y(\vecx) \defined
f(\vecx) + \epsilon(\vecx)$, 
where $\epsilon(\vecx)$ is a random variable with 
mean zero and finite variance $r^2(\vecx)$. In practice, $r^2(\vecx)$ may be unknown, and we assume that the noise is Gaussian.
Because $y(\vecx)$ is stochastic, relying on single samples 
is typically insufficient for determining the underlying function value.
Instead, we will have to sample $y(\vecx)$ multiple times in order to solve the 
problem: 
\begin{equation}
\text{find~} \vecx^* \in \arg \min_{\vecx \in D} \Esp{y(\vecx)}.
\label{eq:optpb}
\end{equation}

We assume that the cost of
 evaluating $y(\vecx)$ is nontrivial and therefore seek methods that require
 the fewest possible oracle queries to solve~\eqref{eq:optpb}. Given the
 probabilistic nature of this optimization problem, we seek to solve it using
 probabilistic models, namely, Gaussian processes (GPs). Among surrogate
 models, GPs are used for their prediction accuracy and analytic
 tractability, plus their ability to handle multi-modality. More importantly, they offer robust uncertainty quantification
 and effective noise filtering, both of which are critical for stochastic
 optimization.

A possible strategy for separating signal from noise in~\eqref{eq:optpb} is to
fix a number of replications/evaluations of $y$ to be performed at any $\vecx$
requested by the optimization method~\cite{Jalali2016,menhorn2017trust}. This
simplification is useful when the signal-to-noise ratio is low, but it may be
wasteful away from portions of the domain $D$ that are critical to
optimization. Limited research has been conducted on methods that adaptively
control the number of replicates. This control can be achieved either by allowing the
method to determine the number of replicates when selecting new points
(single stage) or by including a replication phase (in multistage
approaches); see, for example,~\cite{Jalali2016,pedrielli2016g}. In this work we focus on
adaptively managing replication when proposing new points, the single-stage
case, so that the overall evaluation budget is more effectively distributed
across the parameter space. Such a replication approach is also beneficial
when using GP surrogates because it reduces the associated computational cost of
updating the GP with new observations. These updates can quickly become a
bottleneck when the total number of oracle evaluations grows (something that
is needed to converge precisely). Using relatively larger evaluation budgets  is
referred to as the high-throughput regime in BO, which also involves parallel evaluation (e.g., of
replications)~\cite{moss2023inducing}.

Although our framework is general, this work is motivated by an application
where the cost of evaluating $y(\vecx)$ incurs a setup cost that is paid only 
once for all subsequent replications. 
Specifically, the cost of obtaining $p$
evaluations of $y$ at $\vecx$ is $c(p) \defined c_0 + c_1 \times p$. Note that
the cost $c_0$ is incurred any time replications are requested, so obtaining
an appropriate $p$ is important. That is, twice querying $y$ at $\vecx$ for a
single replication costs $2 (c_0 + c_1)$; but if it is known that two
replications are desired, the cost can be reduced to $c_0 + 2 c_1$.
Such a cost structure arises when optimizing variational circuit parameters for
quantum computing applications where preparing the circuit cost ($c_0$) is
approximately 100 to 1000 times more expensive than observing one sample
(or ``shot'') of the system~\cite{Shi2024}.

Such a model emphasizes the need to balance setup costs with the number of
replications; this is similar to the ``cost-aware Bayesian optimization''
community (see, e.g.,~\cite{klein2017fast,garnett2023bayesian}) where the goal
is to balance evaluation cost and the regret---the difference between the objective value at the
returned solution and the optimal value. 
This related (but distinct) problem occurs when the evaluation cost function depends
on $\vecx$ (see, e.g.,~\cite{ramesh2022movement, Pricopie_2024}). Such a 
cost/benefit trade-off is also present with multifidelity objectives, where
evaluations at different fidelity levels have different evaluation costs impacting the accuracy. In
these situations, some researchers have developed acquisition
functions that jointly select the next point and fidelity level;
see~\cite{huang2006sequential,kandasamy2017multi,Sacher2021}. This approach may be an  option if limiting to a few possible replication budgets per observation. Lifting this restriction to use only a few possible fidelity levels, Kandasamy et al.~\cite{kandasamy2017multi} proposed  a continuous
 version with a GP on both $D$ and fidelity level. In the presence of stochastic noise where fidelity corresponds to the number of replicates, however, this setup effectively reduces to a noisy GP model~\cite{Picheny2013b}. Most of these
 techniques have been explored in either deterministic or high
 signal-to-noise regimes, highlighting the need for novel methods tailored to
 the demanding scenario of high variance, rapid convergence requirements, and
 nontrivial setup costs.

To this end, we utilize a trust-region (TR) framework. 
TR methods provide a systematic way to explore the
domain by iteratively refining a localized region around promising
candidates, and they can be paired with different surrogate models.
The appeal of this combination has been demonstrated in recent works. For
instance, TRs have been employed with polynomials~\cite{Larson2016,Shashaani2018}, radial
basis functions~\cite{wild2008orbit,SMWCAS11} and,
more recently, Gaussian
processes~\cite{eriksson2019scalable,diouane2022trego,regis2016trust,amine2018efficient,nguyen2022local},
although most of these GP studies assume a deterministic $f$ or a high
signal-to-noise regime. The difficulty is that when reducing the TR radius, the signal-to-noise ratio decreases as well.
Lower signal-to-noise ratios have been explored by the
operations research community~\cite{menhorn2017trust,pedrielli2016g}, often
with stochastic kriging models~\cite{ankenman2010stochastic}, which can
identify solutions near a local optimum
but rarely achieve high precision using only GPs. That is, these methods 
have difficulty finding solutions with regret considerably less than the noise level. They also do not jointly find the next point and evaluation budget; in other words, they are not single-stage methods.
Some researchers propose switching to a different optimizer (e.g.,~\cite{mcleod2018optimization}) or adopting a
strategy that discards GPs in the local search~\cite{muller2022accelerating,Mathesen2017}.
Here we focus on improving the GP-based local search.

In order to accommodate noise, possibly with a heteroscedastic variance, and potential non-stationarity in the
objective, one option is to transform inputs and outputs, as in~\cite{cowen2022hebo} with learnable parametric transformations of the input and output spaces. Another option is to employ local models~\cite{gramacy2015local}
combined with replication, a path we pursue because it can significantly reduce
computational overhead. This approach is particularly well suited to the
high-noise regime in which repeated evaluations at promising points are
essential for accuracy, yet overall evaluation budgets remain limited.

Our main contributions are fourfold. First, we rely on replication and
adaptively control the number of replicates when selecting new points. This strategy may be
governed by the evaluation cost problem structure. Second, we make key modifications to the
standard GP trust-region framework to handle low signal-to-noise ratios.
Third, we propose a trust-region method built on local GP models
that scales efficiently with the large number of observations required for high
accuracy. Fourth, we provide software implementing OGPIT (Optimization by Gaussian Processes In Trust regions)
and demonstrate its good performance and improved scalability
against other methods on a range of test problems.

We do not seek to formally establish convergence guarantees or rates for the modified algorithm. 
The underlying GP and trust-region frameworks that we are building on have almost-sure convergence results (say, under suitable regularity assumptions for the GP, acquisition function and local search, see e.g.,~\cite{wang2019controlling,wu2024behavior}, or~\cite[Ch.~10]{garnett2023bayesian}). 
Our modifications are designed to preserve almost-sure convergence 
guarantees of these methods while improving their practical performance in noisy
and high-throughput regimes. Also, asymptotic convergence
behavior is rarely observed in practical finite-budget settings. This further
motivates our emphasis on performance and robustness rather than asymptotics.

\section{Formalization}
This section presents background information on Gaussian processes, Bayesian optimization, and trust-region methods. 
\subsection{Gaussian process regression}

Assuming that $\epsilon \sim \mathcal{N}(0, r^2(\vecx))$ for simplicity, we model
$y$ by a prior zero-mean Gaussian process $Y$ with a covariance function $k$.
Given $N$ observations $\vecY \defined (y_i)_{1 \leq i \leq N}$ at
$\vecx_1, \dots, \vecx_N$, following standard textbooks (e.g.,~\cite
{Rasmussen2006,gramacy2020surrogates}), the prediction at any given $\vecx$ is $Y
(\vecx)|\vecY \sim \mathcal{N}(m_N(\vecx), s^2_N(\vecx))$ with mean and
variance given by
\begin{align}
    m_N(\vecx) &\defined \mathbb{E} \left[Y(\vecx)|\vecY\right] =
\veck(\vecx)^\top (\KN + \LaN)^{-1} \vecY,\\
    s^2_N(\vecx) &\defined
    \mathbb{V}\left[Y(\vecx)|\vecY\right] = r^2(\vecx) + k(\vecx, \vecx) - \veck(\vecx)^\top (\KN
+ \LaN)^{-1} \veck(\vecx),
\end{align} where $\LaN \defined Diag(r^2(\vecx_1), \dots, r^2
 (\vecx_N))$, $\veck(\vecx) \defined(k(\vecx, \vecx_i))_{1\leq i \leq N}$ and
 $\KN \defined (k(\vecx_i, \vecx_j))_{1 \leq i,j \leq N}$. Since the noise
 has mean zero, the corresponding GP model for $f$ is such that $\mathbb{E}\left[f
 (\vecx)|\vecY\right] = m_N(\vecx)$ while $\mathbb{V}\left[f(\vecx)|\vecY\right] = k
 (\vecx, \vecx) - \veck(\vecx)^\top (\KN + \LaN)^{-1} \veck(\vecx)$. There
 exist extensions to non-Gaussian noise, and therefore non-Gaussian
 likelihood, requiring additional approximations; see, for example,~\cite
 {hensman2015mcmc}. In the remainder of this paper we focus on the Gaussian
 noise framework. 

A useful tool in this context is replication, that is, repeating experiments at
the same $\vecx_i$. In this Gaussian framework, $a_i$ evaluations at $\vecx_i$,
$y_i^{(1)}, \dots, y_i^{(a_i)}$, give an equivalent Gaussian observation with mean
$\bar{y}_i \defined \frac{1}{a_i} \sum_{j=1}^{a_i} y_i^{(j)}$ and variance $r^2(\vecx_i)/a_i$. If $a_i$ is large
enough, then the empirical mean and variance give good approximations of $f
(\vecx_i)$ and $r^2(\vecx_i)$. Hence, by controlling $a_i$, one can control the
variance of the equivalent aggregated observation. The benefits are as follows:
\begin{itemize}
 \item The output of the algorithm is one $(\vecx_i, \bar{y}_i)$ with variance smaller ($a_i > 1$) than or equal ($a_i=1$) to that of a single evaluation.
 \item Replicating brings speedups for GPs by using the equivalent observations, lowering the $\mathcal{O}(N^3)$ complexity to $\mathcal{O} (n^3)$ with $n$ the number of unique $\vecx_i$s.
 \item $r^2(\vecx_i)$ is easier to estimate, since the empirical variance may be computed.
\end{itemize}
In the heteroscedastic output noise variance case, by denoting the number
of replicates at observed locations with $\An \defined Diag(a_1, \dots, a_n)$, we can give the GP
prediction for $y$ equivalently by
\begin{align}
\label{eq:gpn}
m_n(\vecx) &\defined \veck(\vecx)^\top (\Kn + \Lan \An^{-1})^{-1} \vecYu,\\ 
s_n^2(\vecx) &\defined r^2(\vecx) + k(\vecx, \vecx) - \veck(\vecx)^\top (\Kn + \Lan \An^{-1})^{-1} \veck(\vecx).
\label{eq:gpvar}
\end{align}
In the remainder, depending on the context, we will use either the generic full-$N$ notation or the replicate-$n$ one (dropping the dependence on $\An$ in a slight abuse of notation).

The covariance function $k$ is typically selected from a parametric family, such as the Gaussian or Matérn covariance functions. Here we use the Matérn 5/2 covariance function:
\[
k(\vecx, \vecx') = \sigma^2 \prod \limits_{j=1}^{d} \left(1 + \frac{\sqrt{5} |x_j - x_j'|}{l_j} + \frac{5(x_j - x_j')^2}{3 l_j^2} \right) \exp\left(- \frac{\sqrt{5}|x_j - x_j'|}{l_j}\right)
\]
 with parameters $\sigma^2$ for the process variance and lengthscales $l_i$.
Most packages propose learning these hyperparameters based on the likelihood. When the $r^2(\vecx_i)$ are unknown, they need to
be estimated as well, which can be difficult. A naive way is to just estimate a
constant noise, which is sensible if the trust region is small enough.
Alternatives for learning $r^2$ include stochastic kriging~\cite{ankenman2010stochastic} and
likelihood-based approaches~\cite{binois2021hetgp}. Once fitted, the GP model can be employed for optimization.

\subsection{Bayesian optimization}

For a broad overview of derivative-free optimization, we refer the reader to~\cite{larson2019derivative}. Within this field, Bayesian optimization (BO) is a
powerful framework that uses a probabilistic surrogate---often a Gaussian
process---to model the objective function~\cite{garnett2023bayesian}. At each
iteration, an \emph{acquisition function} (or \emph{infill criterion})
determines where to sample next, balancing \emph{exploration} of regions of parameter space with
high predictive variance against \emph{exploitation} of regions of parameter space with low
predictive mean. By updating the surrogate model with each new observation, BO
adapts its search strategy and efficiently converges to high-quality solutions,
even in noisy or expensive-to-evaluate settings.

Common acquisition functions are
based on the notion of improvement, that is, the random variable 
$I(\vecx)|\vecY \defined \max(0, T_n - Y(\vecx)|\vecY)$, where $T_n$ is the best observed value
in the deterministic case (i.e., $T_n = \min_{1 \leq i \leq n} y_i$). This leads to the probability of improvement or
the expected improvement (EI): $EI(\vecx) \defined \Esp{I(\vecx)|\vecY}$. Both have closed-form expressions, but the latter is generally preferred, specifically:     
\begin{equation}
    \label{eq:EI}
    EI(\vecx) = (T_n - m_n(\vecx)) \Phi\left( \frac{T_n - m_n(\vecx)}{s_n(\vecx)} \right) + s_n(\vecx) \phi \left(\frac{T_n - m_n(\vecx)}{s_n(\vecx)} \right),
\end{equation}
where $\phi$ and $\Phi$ are the probability and cumulative density functions of the standard Gaussian
distribution, respectively. For notational simplicity, next we drop the implicit conditioning on the current data, $\vecY$.

In the stochastic setting, $T_n$ is typically the best predictive mean over all
observations~\cite{gramacy2020surrogates} (i.e., $T_n = \min_{1 \leq i \leq n} m_n(\vecx_i)$). We follow this plugin approach here, which simplifies computation of the acquisition function. More Bayesian versions
are included in~\cite{Letham2019}, and detailed discussions are provided by~\cite[7.2]{gramacy2020surrogates},~\cite[8.2]{garnett2023bayesian}. Accounting for
all the uncertainty in the improvement, both on the reference $T_n$ if taken as the
best solution and on GP prediction, computations link to those of a discretized
knowledge gradient (KG); see, for example,~\cite{picheny2013benchmark,garnett2023bayesian}. KG is a
less myopic option than EI: $KG(\vecx) \defined \Esp
{\min_\vecx m_{n}(\vecx) - \min_\vecx m_{n+1}(\vecx)}$. However, this criterion has no
analytical expression, which incurs a larger computational cost and the need
for approximations.

Other noisy adaptations have been proposed, such as with the augmented expected
improvement (AEI)~\cite{huang2006sequential,Jalali2016}: $AEI(\vecx) \defined EI
(\vecx) \left(1- \frac{r(\vecx)}{s_n(\vecx)} \right)$ with $T_n$ taken from a best
predictive quantile. The second term leads to the ratio of the reduction in
the posterior standard deviation after a new replicate is added, but only at
the current best point. It thus takes into account the effect of replication,
although indirectly through $s_n(\vecx)$, as noticed in~\cite
{picheny2013benchmark}. 

Yet another variant of EI, also taking into account the effect of a new
observation, is the conditional improvement (CI), that is, the improvement at
a reference location $\vecx$ after another location $\vecx_{n+1}$ is added~\cite
{gramacy2020surrogates}: $CI(\vecx, \vecx_{n+1}) \defined I(\vecx|\vecx_
{n+1})$. The expected conditional improvement (ECI) is expressed as 
    \begin{align}
    \label{eq:CI}
    ECI(\vecx) &\defined \Esp{CI(\vecx)|\vecY}\\
    \nonumber
    & = (T_n - m_n(\vecx)) \Phi\left( \frac{T_n - m_n(\vecx)}{s_{n+1}(\vecx)} \right) + s_{n+1}(\vecx) \phi \left(\frac{T_n - m_n(\vecx)}{s_{n+1}(\vecx)} \right).
    \end{align}

Note that this is almost EI but using the future variance as if $\vecx_{n+1}$ was
observed, since $s_{n+1}$ does not depend on the observations (see Eq.~\eqref{eq:gpvar}). 
ECI is generally integrated over the input domain, leading to the
integrated expected conditional improvement (IECI) that depends only on $\vecx_
{n+1}$. But it is computationally costly since the spatial integral is
not in closed form.

The main issue of EI is that it is myopic: it is optimal for a single new
evaluation, but it does not take into account the effect of future
evaluations nor of replication. So does another popular infill criterion,
based on an upper confidence bound (UCB, for maximizing), $UCB
(\vecx) \defined m_n(\vecx) - \beta s_n(\vecx)$ with a tunable parameter $\beta$. 
One might also use the entropy of the minimizer value or location, which is sample efficient but generally computationally costly.

Developed for noisy optimization, the expected quantile improvement (EQI) has
several appealing properties~\cite{Picheny2013a}. It looks at a given level
$\beta \in [0.5, 1)$ of GP quantile $q_n(\vecx) = m_n(\vecx) + \Phi^{-1}(\beta) s_n
(\vecx)$, and the improvement is taken over the decrease of the
 $\beta$-quantile brought by a new observation at $\vecx_{n+1}$. The
resulting expression is similar to that of EI, Eq.~\eqref{eq:EI}, where $T_n
= \min_{1\leq i \leq n} q_n(\vecx_i)$, while $m_n(\vecx)$ and $s_n(\vecx)$ are replaced
by their counterparts for $q_{n+1}$, also in closed form. Depending on the
level $\beta$ and the noise of the future observation, EQI can suggest to
replicate. Tuning the value of the future noise is also discussed by~\cite{Picheny2013a}; the
proposed strategy is to consider that in the search of $\vecx_{n+1}$, all the
remaining budget will be used there. In fact a fixed elementary budget is
used at each iteration. A more dynamical procedure is to keep replicating
until the value of the EQI becomes less than a fraction of its initial value,
before selecting a new point. A limitation is that the next point is chosen
based on a virtual replication number; hence it is not really single stage.

Extensions of these works to select batches of new points have been proposed;
see, for example,~\cite{chevalier2013fast,binois2021portfolio} and references
therein. Most, however, do not take into account the option to
replicate.

\begin{figure}
\includegraphics[width=0.5\textwidth,trim= 20 40 15 22,clip]{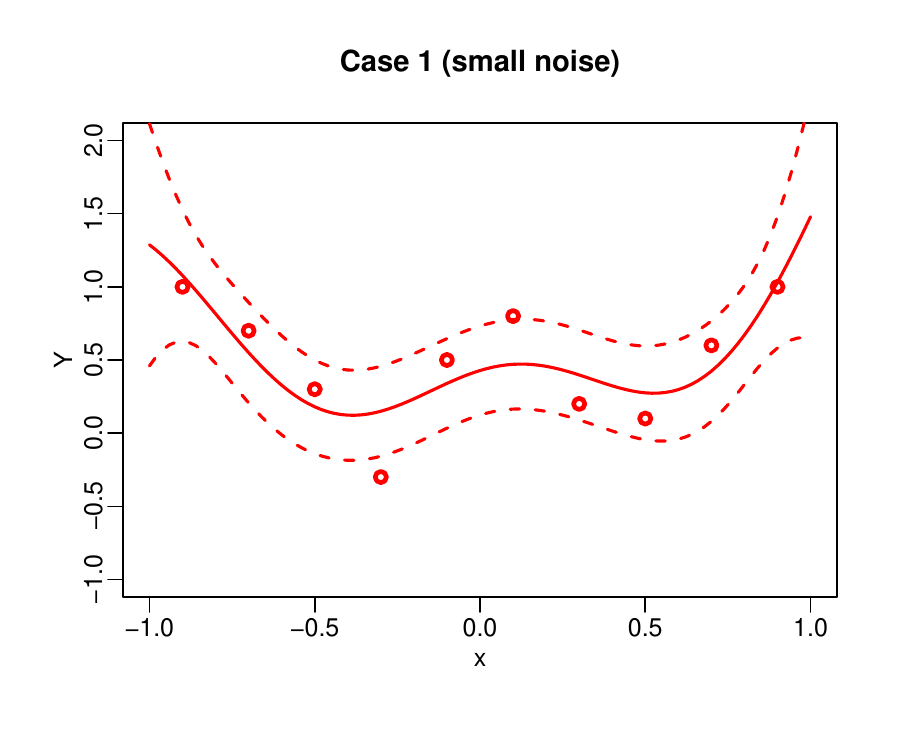}%
\includegraphics[width=0.5\textwidth,trim= 20 40 15 22,clip]{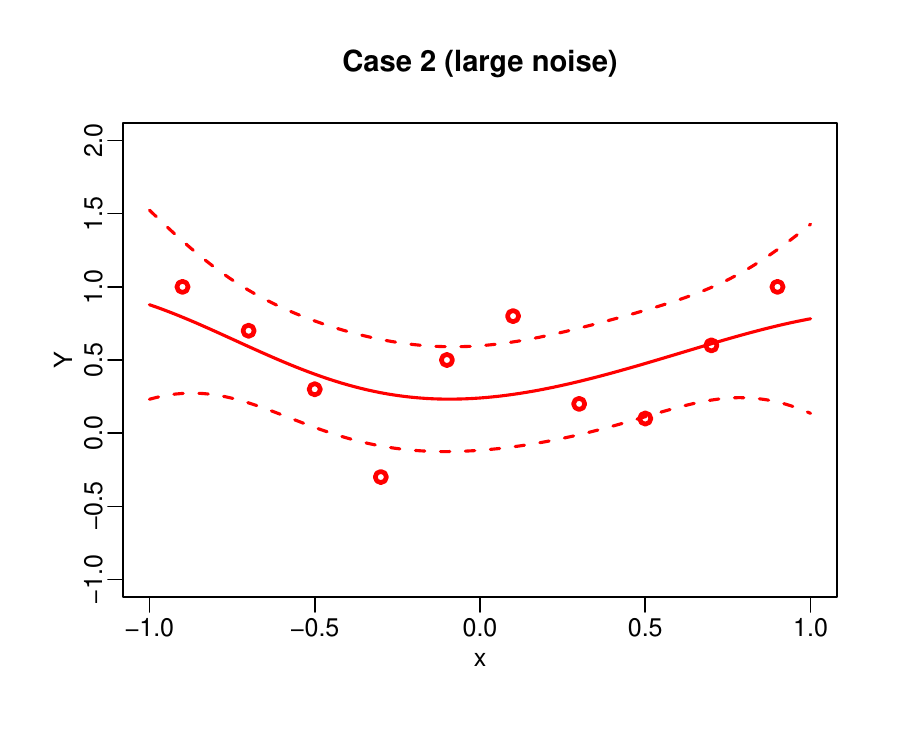}\\
\includegraphics[width=0.5\textwidth,trim= 20 40 15 40,clip]{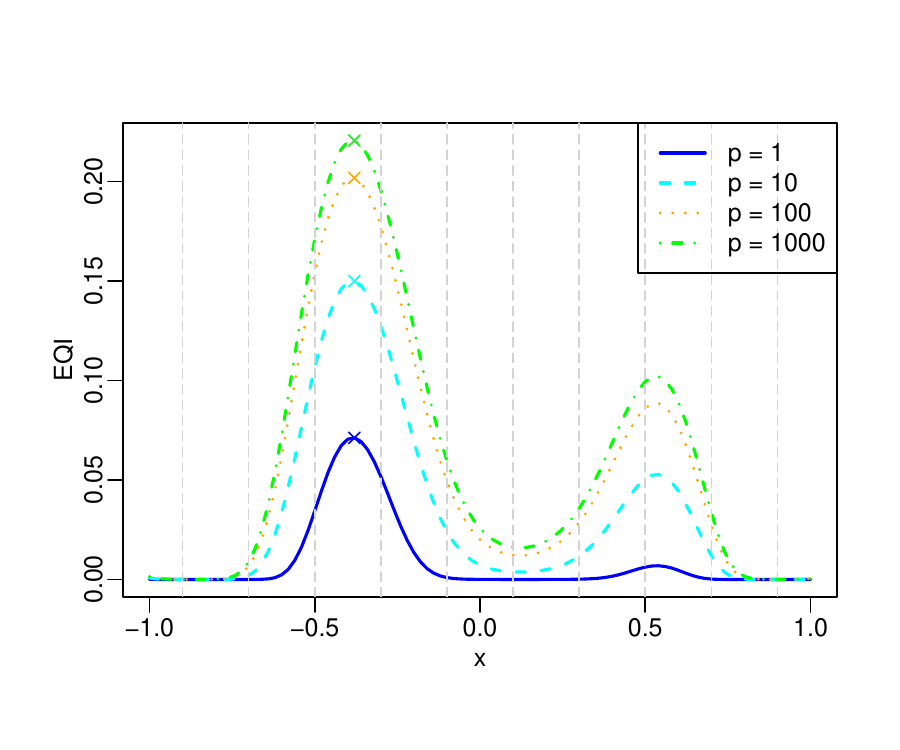}%
\includegraphics[width=0.5\textwidth,trim= 20 40 15 40,clip]{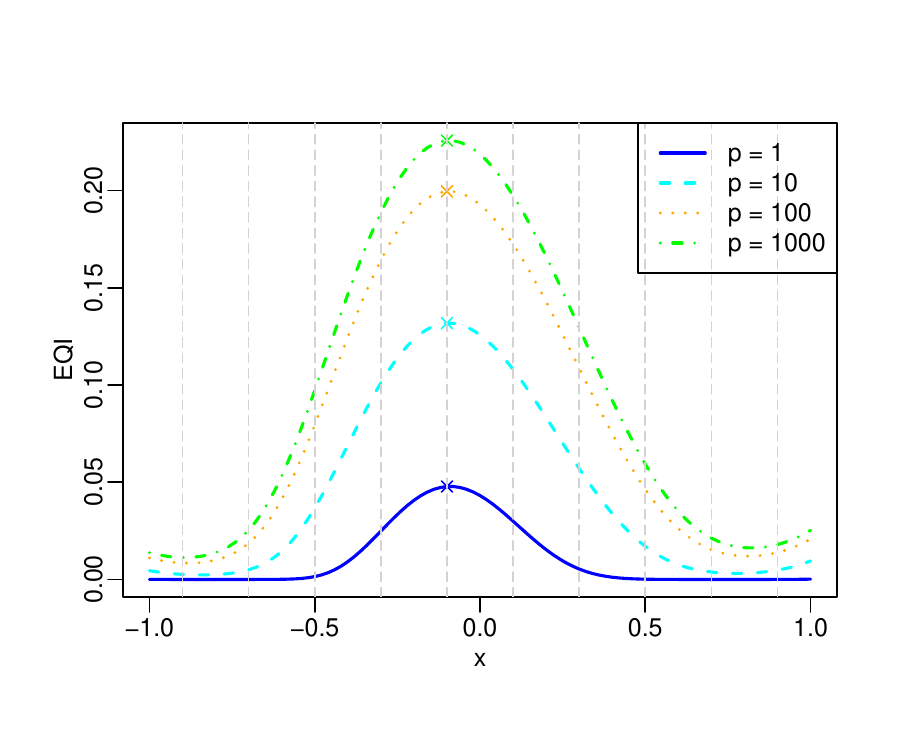}\\
\includegraphics[width=0.5\textwidth,trim= 20 40 15 40,clip]{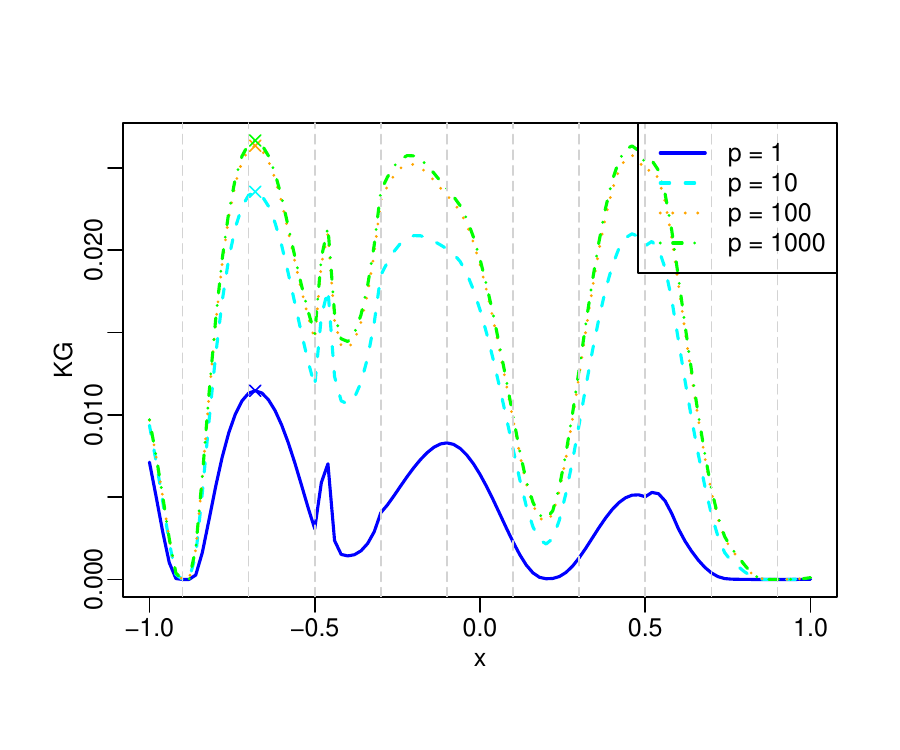}%
\includegraphics[width=0.5\textwidth,trim= 20 40 15 40,clip]{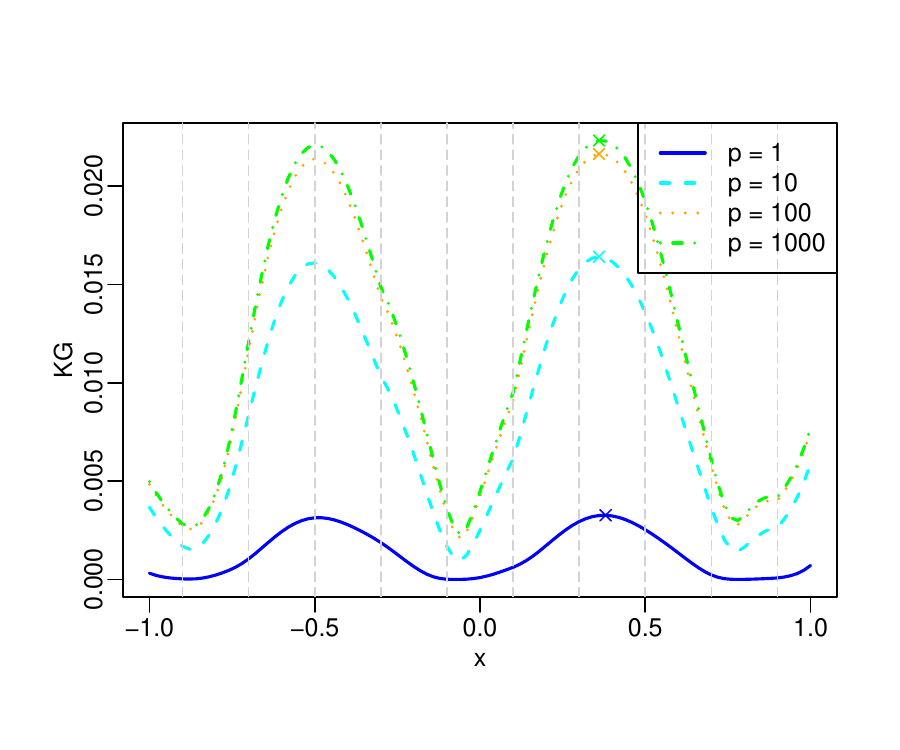}
\caption{GP models and the KG and EQI infill criteria, depending
 on the number of replicates ($p$) at the future points, for two cases (left and right). Vertical lines indicate design points. Crosses mark the respective maxima.}
\label{fig:GPandcrits}
\end{figure}

We illustrate in \Cref{fig:GPandcrits} two criteria on two different setups with different
signal-to-noise ratios. KG and EQI behave
differently: KG wants new points on either side of the current local optima,
while EQI wants to add points on them. The diminishing returns effect of adding
more replicates is also visible on both criteria: the difference between
curves for 1 and 10 replicates is much larger than for 100 and 1000
replicates.

\subsection{TR methods and combination with BO}

Trust regions have been used in various optimization strategies, for example with
mesh-adaptive direct search~\cite{Audet_2021}. The driving mechanism of TR
methods is the control of the trust-region radius over iterations. In case
progress is made, the radius can grow, whereas it is reduced otherwise. Several
combinations of BO with trust regions have been proposed in the literature
and have attracted attention in benchmarks; see, for example,~\cite{regis2016trust,eriksson2019scalable,diouane2022trego,santoni2024comparison}.
A possible motivation for such
attention is that TR framework promotes more precise local convergence while
still allowing for global exploration capabilities using a global Bayesian
model.

Derivative-free optimization (DFO)-TR methods maintain a
surrogate model $m_n(\vecx)$ of the objective in the neighborhood
$\mathcal{B}(\vecx_n,\Delta_n)$ centered at the current iterate $\vecx_n$ with radius
$\Delta_n$. The model is constructed 
using previously sampled points and is designed to satisfy geometric
(poisedness) conditions that ensure $m_{n}$ it is sufficiently close to $f$ on $\mathcal{B}(\vecx_n,\Delta_n)$. 
Each iteration, a trial step $\s_n$ is computed by (approximately) minimizing
\[
    \min_{\|\s\| \le \Delta_n} m_n(\vecx_n + \s).
\]
This step is accepted (and $\vecx_{n+1} \gets \vecx_n + \s$) or is rejected ($\vecx_{n+1}
\gets \vecx_n$) based on the $\rho_n$ ratio of actual decrease over predicted decrease. 
When the objective is stochastic, values in the numerator must be replaced by estimates.
If $\rho_n$ is sufficiently large, $\Delta_n$ can be increased; otherwise, the
step is rejected and $\Delta_n$ is reduced. Under suitable model accuracy
conditions (e.g., fully linear or fully quadratic models), this framework
yields global convergence to first- or second-order stationary points

When the objective is stochastic, the acceptance ratio $\rho_n$ becomes a
random quantity, and noise can lead to incorrect step acceptance or rejection.
Model construction is also more delicate: regression models must balance bias
and variance, and sample sizes may need to increase as $\Delta_n \to 0$ to
ensure sufficient model accuracy relative to the trust-region radius.
Estimating predicted reduction and enforcing sufficient decrease conditions
require variance-aware criteria. This can involve repeated sampling or adaptive
sampling rules. Without careful control of noise (e.g., increasing sample
sizes, using common random numbers, or employing probabilistic acceptance
tests) the algorithm may stall at suboptimal points. Stochastic DFO-TR methods
must couple model fidelity, sampling effort, and radius control to maintain
both convergence guarantees and practical efficiency.

The differences between these methods depend on the balance between the two
frameworks. The globabilization strategy is either via global BO steps with
local TR ones with TREGO~\cite{diouane2022trego}, entertaining several trust regions in TuRBO~\cite{eriksson2019scalable}, or performing
restarts with TRIKE~\cite{regis2016trust} or GSTAR~\cite{pedrielli2016g}. The other main difference is on the success condition,
which can be based on the actual improvement over EI, consecutive
successes/failures, or a sufficient decrease condition. The infill criterion
in the local steps also varies, from the predictive mean to EI or Thompson
sampling. More local GP methods have proposed criteria based on the
probability of descent; see, for example,~\cite{nguyen2022local}. The convergence rates
for local BO are derived in~\cite{wu2024behavior}.

Most of these do not deal specifically with the noisy setup. SNOWPAC~\cite
{menhorn2017trust} does, improving a polynomial surrogate by a GP one, also for constrained optimization, but it relies
on the computation of a robustness measure obtained by sampling.
Note that SNOWPAC requires a measure of the magnitude of the noise in the objective, something that is not available in many situations we consider. For noisy simulator, GSTAR~\cite{pedrielli2016g} relies on an ensemble of global and local GP models similar to stochastic kriging. But it remains in a two stage frameworks: new designs are added with a fixed number of replicates while more replicates are added to already sampled points in a dedicated phase.

\section{Exploration-exploitation-replication trade-off}

With this background information, we can now provide greater detail about our
approach to decide both the next iterate $\vecx_{n+1}$ and its associated number of
replicates $a_{n+1}$. 
We also propose a new infill criterion, borrowing
from links between batch and look-ahead formulations to handle the setup cost case. 

\subsection{Replication} 

First, as seen in~\eqref{eq:gpn}, using replicates can directly lower the
computational expense of model building. Without replication, to lower
computational cost with large $N$, GP modeling often relies on techniques
such as pseudo-inputs for sparse GPs~\cite{snelson2006variable}, specialized
linear algebra methods~\cite{wang2019exact}, or the Vecchia
approximation~\cite{cao2022scalable}. In
the low signal-to-noise regime, however, single evaluations offer limited
information, and the overhead of deciding on a new point for just one
evaluation may not be justified. An alternative is parallel/batch
acquisition functions to evaluate multiple points simultaneously, but this
introduces a more demanding optimization and approximation problem. In many
cases, choosing a single new point with an appropriate replication budget
proves both simpler and more effective.

Second, replication can be warranted by the
acquisition function itself. In deterministic settings, this possibility does
not arise; but in stochastic contexts, previous works have demonstrated that
replicating at a single location can be optimal for certain global accuracy
criteria (e.g., integrated mean square prediction error~\cite{binois2019replication}) or for the EQI criterion~\cite{Picheny2013a}.

Replication is still not the panacea. As the number of replicates increases,
the marginal benefit diminishes because the noise standard deviation decreases
in proportion to the square root of the number of replicates (see also
\Cref{fig:GPandcrits}). Thus, finding the right replication level is critical.
As detailed below, this decision can stem directly from the acquisition
function or incorporate user-defined incentives such as a setup cost function. At a
higher level, strategies include revisiting previously evaluated points to add
a fixed number of replicates~\cite{Jalali2016}, i.e., two-stage approaches, or enforcing a minimal distance
between points~\cite{binois2019replication}.
Discrete domains also naturally lend themselves to replication because they prevent adding points arbitrarily close to one another, as may occur in continuous spaces.

Handling low signal-to-noise ratios often requires further
adaptations, and incorporating these ideas into a trust-region framework brings
additional subtleties. We address these considerations in the subsequent
sections.

\subsection{Looking ahead}

The most common infill criteria---such as EI or UCB---do not naturally suggest a replication level for a given
point because they account only for existing observations, not the number of
future replicates. One might divide by an evaluation cost function that depends on $\vecx$~\cite{snoek2012practical}, but that approach is outside our current scope. A
more direct remedy is to adopt a look-ahead perspective, where policies
consider the effect of $p$ future observations at points
$\vecx_{N+1}, \dots, \vecx_ {N + p}$, shortened in $\vecx_{(N+1):(N+p)}$.
In the GP context, closed-form update
formulas~\cite{chevalier2013corrected,gramacy2020surrogates} facilitate such
approaches. For example, adding $a_{n+1} = p$ replicates at $\vecx_{n+1}$ updates the GP
mean and variance as follows:
\begin{align}
\label{eq:up2}
m_{N+p}(\vecx) = m_{n+1}(\vecx) &= m_n(\vecx) +  \frac{s_n^2(\vecx, \vecx_{n+1})}{s_n^2(\vecx_{n+1}) +r^2(\vecx_{n+1}) / p} (\bar{y}_{n+1} - m_n(\vecx_{n+1})),\\
s_{N+p}^2(\vecx) = s_{n+1}^2(\vecx) &= s_{n}^2(\vecx) - \frac{s_n^2(\vecx, \vecx_{n+1})^2}{s_n^2(\vecx_{n+1}) +r^2(\vecx_{n+1}) / p}.
\label{eq:up3}
\end{align}

These update formulas have largely been used in look-ahead infill
criteria---such as KG or entropy-based metrics---because the predictive
variance term remains independent of the actual future observations, even
though the predictive mean does not. However, fully propagating uncertainty
from potential future samples makes look-ahead methods computationally
challenging beyond a few steps. Various approximations therefore arise in
practice~\cite{Ginsbourger2010,garnett2023bayesian}. One strategy is to
restrict the set of future decisions, for example by allowing only one new
point in a future step and devoting all other steps to replication~\cite{binois2019replication}. In that scenario, if adding a new point is more
advantageous in a later iteration, then replication at $\vecx_{n+1}$ is
optimal for now---thus balancing exploration, exploitation, and replication.

But we can also rely on the ``strong connection between the optimal batch and
sequential policies''~\cite{garnett2023bayesian}, such that a batch policy is
a reasonable approximation of the sequential/look-ahead one. An example 
can be found in the GLASSES method~\cite{Gonzalez2015}.

\subsection{Blending all together: a new infill criterion geared to 
trust-region methods}

We primarily focus on \emph{improvement-based} infill criteria rather than, for
example, selecting iterates using only the predictive mean $m_n$. In our experience, the
predictive mean in high-noise settings becomes relatively flat, often placing
its minimum near the trust-region center $\vecx_c$. By contrast, improvement-based
methods foster more exploration.

Our approach here will take advantage of the closed-form expression for the future
predictive variance without explicitly modeling the noisy future
observations, thereby retaining computational efficiency. One criterion meeting
this requirement is IECI, although it involves a costly spatial integral.

Our novel infill criterion is qERCI, for parallel expected reduction in
conditional improvement. More precisely, we consider the reduction in
improvement at several reference locations brought by a batch of candidate
evaluations. The rationale is that in a TR framework, besides the new point(s), the main points of
interest are the current center $\vecx_c$, the estimated optimum $\vecx_n^*$ (which can coincide with the
center, but not necessarily). A key point is that we want to estimate it a
priori, among existing points, not a posteriori as with KG, which would become more computationally expensive. Plus in high-noise regimes, where the GP prediction is comparatively flat (see
\Cref{fig:GPandcrits}), preliminary estimates of the minimum value and location based on the predictive mean appearing in the acquisition function are unlikely to change drastically once the new
observations are incorporated. 

The expected reduction in improvement at $q$ designs $\vecx'_{1:q}$ when adding new points $\vecx_{(N+1):(N+p)}$ is given by the general qERCI criterion:
\begin{align}
    qERCI \left(\vecx'_{1:q} | \vecx_{(N+1):(N+p)} \right) &\defined  \Esp{\left(T_N - \min \limits_{i \in \left\{1, \dots, q \right\}} Y(\vecx'_i )\right)_+} - \\
    \nonumber
    & \Esp{\left(T_N - \min \limits_{i \in \left\{1, \dots, q \right\}} Y(\vecx'_i )| \vecx_{(N+1):(N+p)} \right)_+}\\
    &= qEI \left(\Esp{Y(\vecx'_{1:q})}, \Cov{Y(\vecx'_{1:q})}\right) - \\
    \nonumber
    & qEI\left(\Esp{Y(\vecx'_{1:q})}, \Cov{Y(\vecx'_{1:q})| \vecx_{(N+1):(N+p)}}\right).
\end{align}

The first part is the standard parallel expected improvement (qEI), based on $Y(\vecx)$, while the second term is the parallel expected improvement for $Y(\vecx) | \vecx_{(N+1):(N+p)}$.
Fortunately, qEI also has a closed-form expression~\cite{chevalier2013fast}, and its
gradient can be computed analytically~\cite{marmin2015differentiating}, making evaluations relatively
fast when $q$ is small. For larger $q$, an efficient approximation is
provided by~\cite{Binois2015c}. 

\begin{figure}
\includegraphics[width=0.5\textwidth,trim= 25 30 10 20,clip]{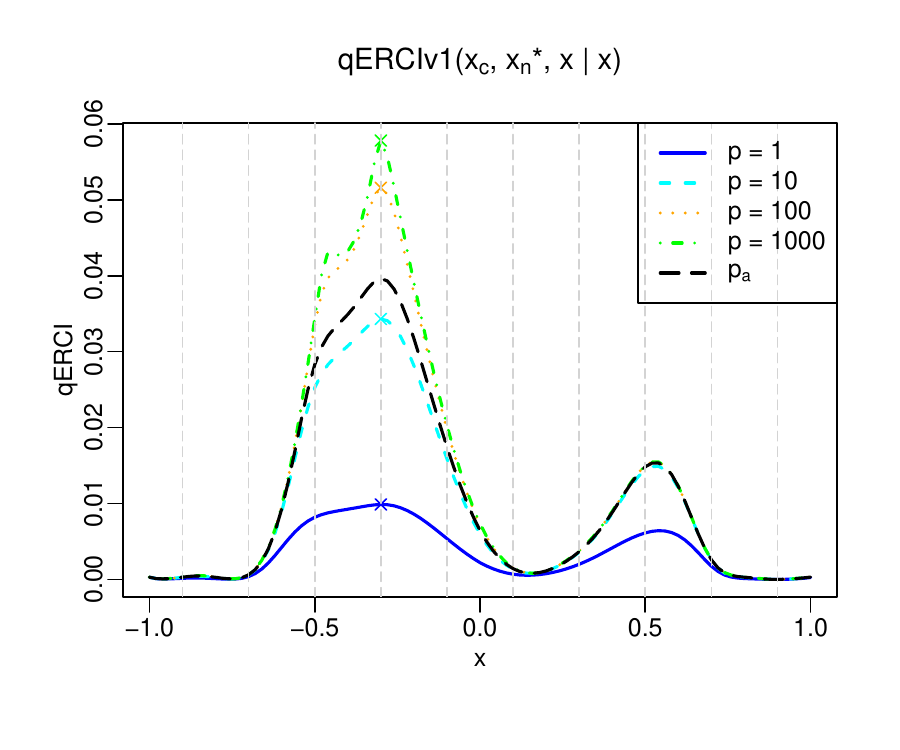}%
\includegraphics[width=0.5\textwidth,trim= 25 30 10 20,clip]{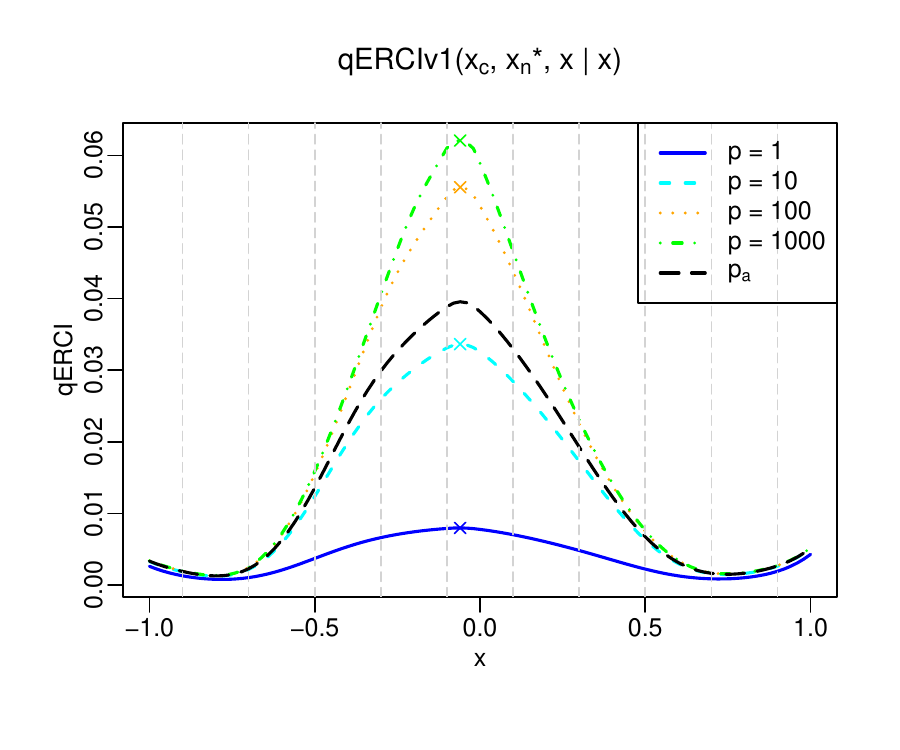}\\
\includegraphics[width=0.5\textwidth,trim= 25 40 20 20,clip]{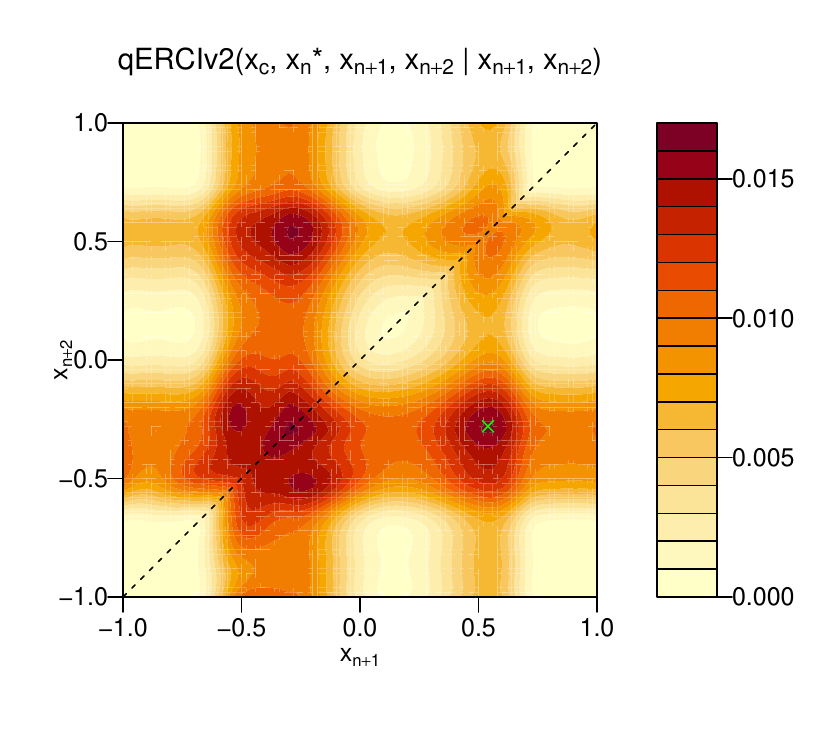}%
\includegraphics[width=0.5\textwidth,trim= 25 40 20 20,clip]{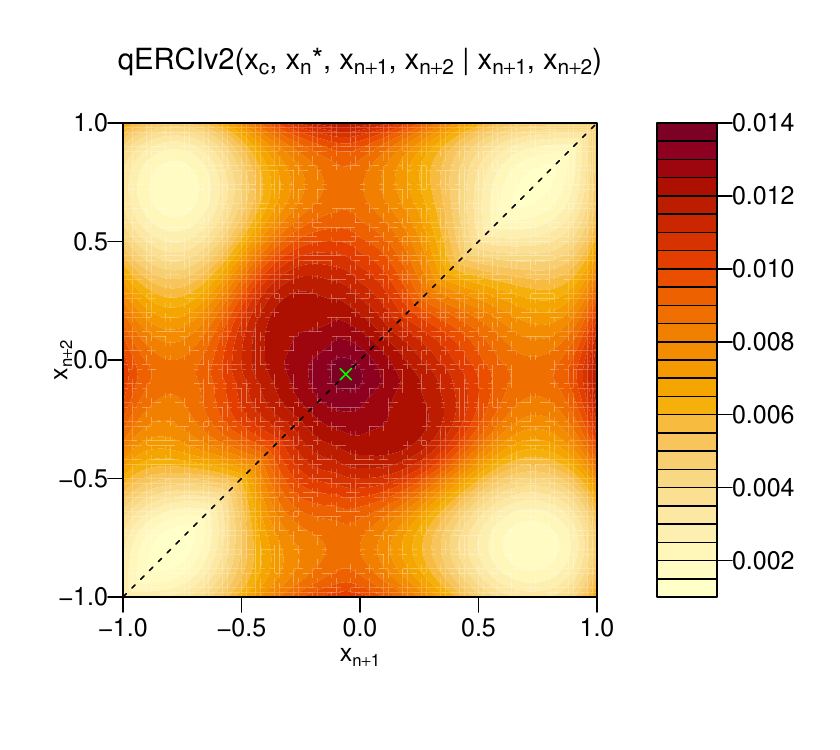}
\caption{Parallel reduction in improvement example for the cases in \Cref{fig:GPandcrits}, one case per column.
Crosses mark the respective maxima.
 Top: qERCI for one new point
 $\vecx_{n+1}$ allowing for various numbers of replicates ($a_{n+1} = p$).
 Bottom: qERCI for two new points with a single replicate. The dashed line marks where $\vecx_{n+1}=\vecx_{n+2}$.}
\label{fig:qRIt}
\end{figure}

In practice, $\vecx_{(N+1):(N+p)}$ are the points that will be tentatively evaluated, including replicates, while $\vecx'_{1:q} = \left(\vecx_c, \vecx^*_n,  \vecx_{(n+1):(N+p)} \right)$. The behavior of the qERCI
criterion when adding a single new design is depicted in \Cref{fig:qRIt}. As can be seen, depending on the
number of possible replicates and noise variance, the best location may be a
replicate of already evaluated designs.

\subsection{Joint selection of $\vecx_{n+1}$ and $a_{n+1}$}
\label{sub:joint_selection_of_vecx_n_1}

For non-myopic infill criteria, additional replications always provide more
information gain. Because evaluations are costly, however, controlling the
replication budget $a_{n+1}$ is essential in practice. 

Rather than deciding $a_{n+1}$ once $\vecx_{n+1}$ is found, it can be decided
beforehand based on a minimal variance reduction. Indeed, in a low
signal-to-noise setting, a single observation typically has limited impact on
the GP prediction. Consequently, we propose replicating enough to
achieve a predefined reduction in the predictive variance. That is, let this adaptive number of replicates
$p_a(\vecx)$ be the fewest number of replications $p$ at $\vecx_{n+1} = \vecx$ such that $\frac{s_{n}^2(\vecx) - s_{n+1}^2(\vecx)}
{s_{n}^2(\vecx)} \geq T_a$, where $T_a$ is a user-specified threshold (e.g., $T_a = 0.2$).
From~\eqref{eq:up3}, 
\begin{equation}
p_a(\vecx) =\lceil r^2(\vecx_{n+1}) (s_n^2(\vecx, \vecx)/T_a - s_n^2(\vecx, \vecx))\rceil.
\label{eq:pa}
\end{equation}

Then the acquisition function is computed with $a_{n+1} = p_a$ replicates at $\vecx$, which we call version 1 of qERCI:
\begin{equation}
\label{eq:qerciv1}
    qERCIv1(\vecx) \defined qERCI \left( \left(\vecx_c, \vecx^*_n,  \vecx \right) | \left(\vecx, \overset{p_a(\vecx)~\text{times}}{\cdots}, \vecx \right) \right).
\end{equation} 

This is
illustrated in \Cref{fig:qRIt}. Compared with options with a fixed number of
replicates, this adaptive version can suggest taking more or fewer
replicates and does not follow a specific fixed $p$ curve.

Next we focus on the case where there is a setup cost $c_0$
(say, meshing for finite element models, or an experimental setup) that 
needs to be paid only once for several replicates, each with a smaller cost of
$c_1$. The cost for $p$ evaluations of a given $\vecx$ is then $c_0 + p c_1$.
The benefit of introducing a cost function is to balance with the increased value from replicating of look-ahead acquisition functions (that show the diminishing return property when adding new replicates). Thus, the cost must increase sufficiently to penalize adding too many replicates.

Introducing a second possible new design is another option to assess when to stop replicating.
In this case it suffices to optimize qERCI with respect to $\vecx_{n+1}, a_{n+1}$ and $\vecx_{n+2}, a_{n+2}$, with up to $a_{n+1} + a_{n+2} \leq
 \pmax$ replicates. By dividing by the setup costs, it should naturally
limit the number of replicates; and with $c_0 > 0$, the interest of a second new
point may be diminished if not bringing enough gain compared to replicates on the first one. This is the second version of qERCI that we propose:
\begin{equation}
\label{eq:qerciv2}
    qERCIv2(\vecx, a, \vecx', a', c_0, c_1) \defined  \frac{qERCI\left( \left( \vecx_c, \vecx^*_n,  \vecx, \vecx' \right) | \left(\vecx, \overset{a~\text{times}}{\cdots}, \vecx, \vecx', \overset{a'~\text{times}}{\cdots}, \vecx' \right)\right)}{c_0 \times (\mathbf{1}_{a > 0} + \mathbf{1}_{a' > 0}) + c_1 \times (a + a')}.
\end{equation}

The replication budget may be optimized as a continuous variable for simplicity: Eq.~\eqref{eq:up3} does not need an integer valued $p$ but ultimately it must be rounded for the evaluation.
If splitting the replication budget between two candidates is preferred ($a_{n+1}> 0$ and $a_{n+2}>0$), the
one with the largest number of replicates is evaluated.

We entertain this strategy on the previous examples, optimizing $qERCIv2$ over $x,x',a,a'$ by running a particle swarm optimization for global search, before applying L-BFGS-B to the best candidate. Depending on the cost
function, as seen in \Cref{tab:qRI2cost1,tab:qRI2cost2}, all the
evaluation budget is used or split between two candidates, or not.

\begin{table}[htpb]
\caption{Results for the left case in \Cref{fig:qRIt}. The maximum
 number of replicates is $\pmax=10$. For zero costs, the value of the denominator in $qERCIv2$ is set to 1.}
\label{tab:qRI2cost1}
\centering
\begin{tabular}{c|c|c|c|c|c|c}
$c_0$ & $c_1$ & $x_{n+1}$ & $a_{n+1}$ &  $x_{n+2}$ &  $a_{n+2}$ & qERCIv2 \\ \hline
0     & 0     &   -0.29   &  7 & 0.54 & 3 & 0.041    \\ \hline
1     & 1     &   -0.29   &  2 & -& - & 0.016   \\ \hline
1     & 0.1   &   -0.29   &  6 & - & - & 0.028   \\ \hline
1     & 0.001 &   -0.29   &  10  & - & - & 0.034 
\end{tabular}
\end{table}

\begin{table}[htpb]
\caption{Results for the right case in \Cref{fig:qRIt}. The maximum number
 of replicates is 10. For zero costs, the value of qERCI is divided by 1 in the search.}
\label{tab:qRI2cost2}
\centering
\begin{tabular}{c|c|c|c|c|c|c}
$c_0$ & $c_1$ & $x_{n+1}$ & $a_{n+1}$ & $x_{n+2}$ &  $a_{n+2}$ & qERCIv2 \\ \hline
0     & 0     &   -0.061        &  9 & 1 & 1 & 0.036  \\ \hline
1     & 1     &   -0.061        &  2 & - & - & 0.014   \\ \hline
1     & 0.1   &   -0.059        &  8 & - & - & 0.031   \\ \hline
1     & 0.001 &   -0.059        & 10  & - & - &  0.034 \\ \hline
\end{tabular}
\end{table}

\section{Adaptations for TR-Based Bayesian Optimization}

Although trust-region-based BO methods have shown promise in the deterministic
case, several important considerations arise when adapting them to the low
signal-to-noise ratio regime.

\subsection{Handling of non-stationarity and replication}

To address non-stationarity in the objective function, one can adopt
non-stationary covariance kernels~\cite{paciorek2006spatial}, pseudo-inputs~\cite{snelson2005sparse}, or local models~\cite{gramacy2015local}. We choose
local models because they also reduce computational cost and naturally pair
with the TR framework, which restricts attention to observations
within the current trust region rather than the entire domain.

In practice, we scale inputs in the trust region to $[-1, 1]^d$ and scale 
outputs by subtracting their mean and scaling by their standard deviation,
which avoids numerical difficulties with very small or large output values and
simplifies the setting of hyperparameter bounds. 
We then apply simple kriging (i.e., assuming a known constant mean of zero).

For replication, we rely on the approach of~\cite{binois2018practical} to leverage the associated computational savings. 
This method can also either estimate a constant noise level or learn an input-dependent noise GP for $r^2(\vecx)$ by a joint likelihood approach on the mean and variance GP processes. The stochastic kriging~\cite{ankenman2010stochastic} approach could be used as well, but it requires a second GP fitted on empirical variance
estimates to predict the variance at arbitrary locations (hence a minimal degree of replication). When the noise variance depends on the input, learning it is beneficial, but it also comes with additional inference and computational complexity. Moreover, when the noise variance is smooth and the TR radius decreases, the benefit of modeling heteroscedasticity becomes less important. Consequently, unless the noise is assumed to be constant, we allow the heteroscedastic model to revert to an homoscedastic one as proposed by~\cite{binois2021hetgp}, based on the comparison between their likelihoods.
In both cases, the
model is trained on the $n_b$ closest points to the TR center. This helps
improve the local accuracy of the model and further reduce computational costs.

\subsection{Acceptance test}
\label{sec:accept}

From the benchmarking study in~\cite{Jalali2016}, it is identified that accurately
ranking and selecting the ``best'' evaluated point is crucial. In a trust-region
setting, this task is simpler: we need to determine only whether the new
candidate outperforms the current center.

A key quantity governing the TR method for deterministic objectives is the
acceptance ratio: $\rho_n =
\frac{f(\vecx_c) - f(\vecx_{n+1})}{m_n(\vecx_c) - m_n(\vecx_{n + 1})}$. When noise is
present, using the noisy function values is not reasonable, but their updated model values
could be used instead: $m_{n+1}(\vecx_c)$ and $m_{n+1}(\vecx_{n+1})$. 
Then we can actually compute the acceptance ratio, although the numerator and
denominator of $\rho_n$ are no longer independent quantities. 

In our framework we also compute the denominator in the acceptance ratio
using leave-one-out predictions under the same (updated) GP model used for
the numerator---in other words, the updated GP model after observing the data at $\vecx_
{n+1}$. The leave-one-out predictions at $\vecx_i$, that is, after removing the corresponding observation, have closed-form
expressions~\cite{dubrule1983cross,Bachoc2013}:
\begin{align}
\tilde{m}_{n}(\vecx_i) &= \bar{y}_i - [(\Kn + \Lan \An^{-1})^{-1} \bar{\vecY}]_i/[(\Kn + \Lan \An^{-1})^{-1}]_{i,i},\\
\tilde{s}_n^2(\vecx_i) &= [(\Kn + \Lan \An^{-1})^{-1}]_{i,i} - r^2(\vecx_i)/a_i.
\end{align}
This approach differs from simply relying on the previous GP model in that
the covariance hyperparameters can be more precisely estimated with the additional data at the newly sampled design.

To decide on the next center, the decision is generally based on the value of
the mean. To make this step more robust, we add a constraint on the relative
variance between $Y(\vecx_c)$ and $Y(\vecx_{n+1})$. Indeed, a deceptive case
can occur when the variance at $Y(\vecx_{n+1})$ is large and the
corresponding mean is lower than for $Y(\vecx_c)$, while this can simply be
an overoptimistic estimate due to the noise variance. To prevent this case,
the number of replicates $a_{n+1}$ at $\vecx_{n+1}$ is increased such that $\Var{Y
(\vecx_{n+1})} \leq 4 \Var{Y(\vecx_c)}$. If this inequality is not satisfied
after updating with the observations, the new point cannot be accepted as the
center until more data is available.

Another corner case that may occur is when the new point selected by the
acquisition function leads to a predicted decrease that is negative
(i.e., corresponding to exploration, with large predictive variance). There
the acceptance ratio should be computed differently: $\rho_n =
\frac{m_{n+1}(\vecx_c) - m_{n+1}(\vecx_{n+1}) - (\tilde{m}_{n+1}(\vecx_c) -
\tilde{m}_{n+1}(\vecx_{n+1}))}{|\tilde{m}_{n+1}(\vecx_{n+1}) - \tilde{m}_{n+1}
 (\vecx_{c})|}$. This modification accommodates \emph{good surprises} by
 allowing for the possibility that the actual observed improvement might
 exceed the model’s initial negative prediction.

\subsection{Reducing the trust-region radius}

\label{sec:trreduc}

As trust-region methods progress, they must decrease the TR radius in order to improve
model quality and find continued decrease. In noisy settings, however, this can
be problematic: as the radius shrinks by a coefficient $\gammad$, the local signal-to-noise ratio also
declines, causing the predictive surface to become flatter and more uncertain.
This, in turn, makes identifying points where the function is improved
increasingly difficult, which can then trigger further radius reductions and
exacerbate the issue. In such cases, the TR radius should \emph{not} be
decreased.

\begin{figure}
\centering
\includegraphics[width=0.3333\textwidth, trim= 25 40 20 22, clip]{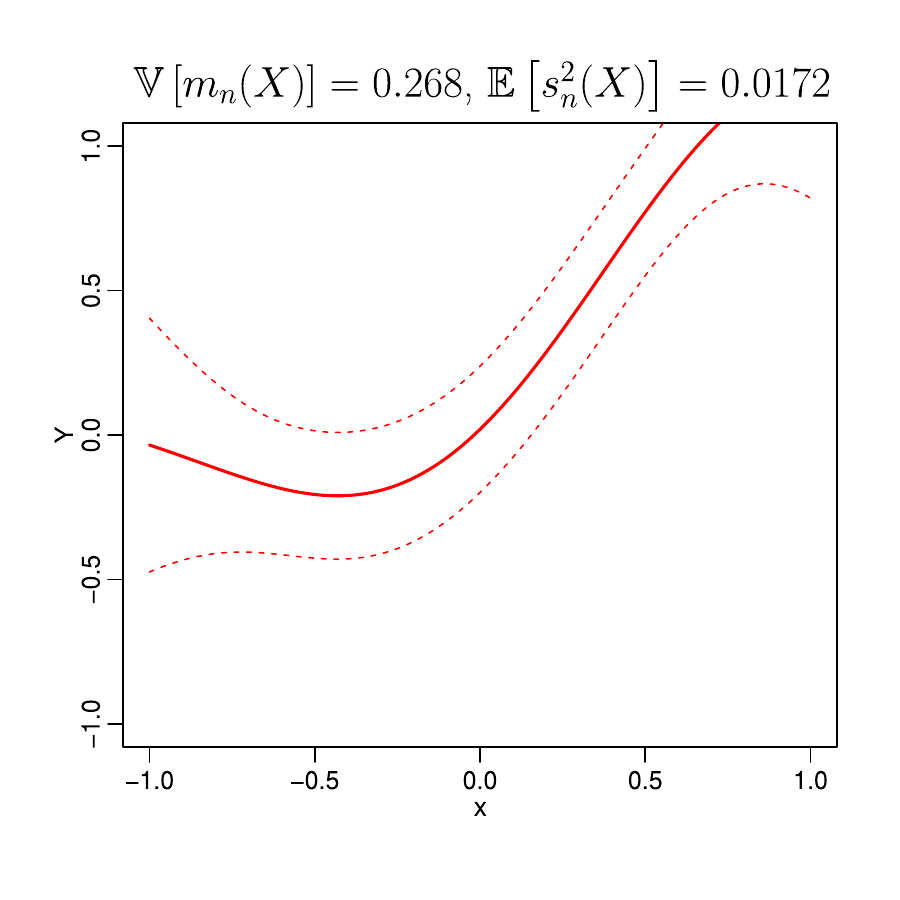}%
\includegraphics[width=0.3333\textwidth, trim= 25 40 20 22, clip]{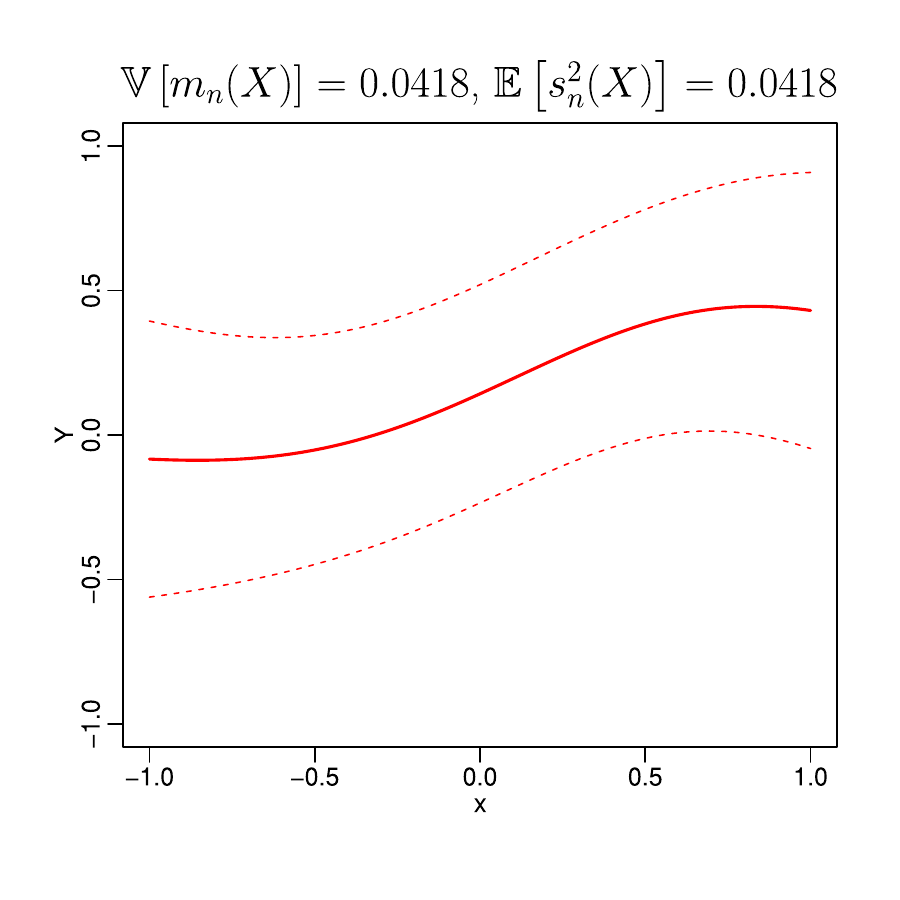}%
\includegraphics[width=0.3333\textwidth, trim= 25 40 20 22, clip]{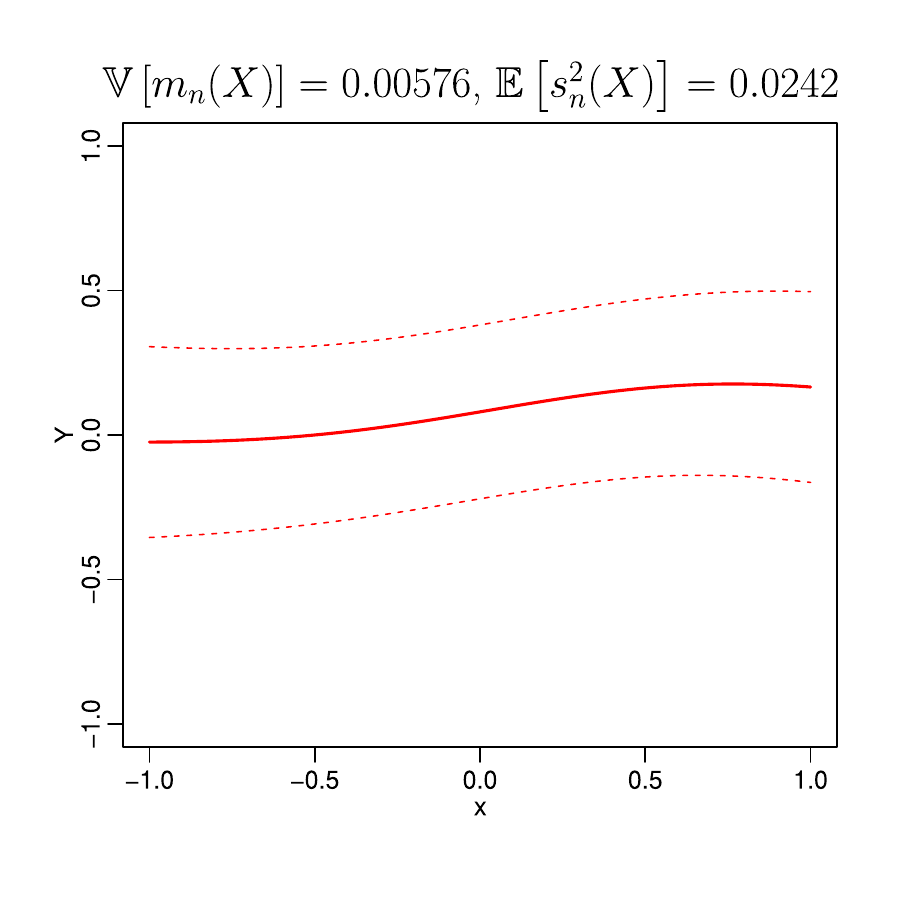}%
\caption{Gaussian processes with various variance of the mean vs. mean predictive variance ratios.}
\label{fig:mnsnv}
\end{figure}

To identify when this situation arises, we compare the
contributions to the variance of $y$ over the TR domain $\Omega$, that is, 
$\Var{y(X)| X \sim \mathcal{U}(\Omega)}$. For a GP model $Y$, the law of total variance gives 
\begin{align}
\mathbb{E}[Y] &= \mathbb{E}[\mathbb{E}[Y|X]] = \mathbb{E}\left[m_n(X)\right],\\
\Var{Y} &= \mathbb{E}[\Var{Y | X}] + \Var{\mathbb{E}[Y| X]} = \mathbb{E}\left[s_n^2(X)\right] + \Var{m_n(X)}.
\label{eq:totvar}
\end{align}
For $\Var{Y}$, the first term is the integrated mean squared prediction error (IMSE), 
and the second term can be written as 
$$\Var{m_n(X)} = \mathbb{E}\left[m_n(X)^2\right] - \mathbb{E}\left[m_n(X)\right]^2 = (\Kn^{-1}\vecY)^\top \Var{k(X, \vecx_{1:n})}\left(\Kn^{-1}\vecY\right).$$
Both have known expressions for the standard Gaussian or Mat\'{e}rn kernels; see, for example,~\cite{binois2019replication}. 
We illustrate the idea in \Cref{fig:mnsnv}. On the left, the leading
contribution is the variance of the mean, while on the right it is the predictive variance. In this case, reducing the TR radius would be detrimental. 

\subsection{Summary of the approach}

These elements are generic and may be integrated within existing TR frameworks. Next we detail one instantiation.

The full OGPIT pseudocode is summarized in \Cref{alg:OGPIT}. A GP prior can be defined with an arbitrary number of initial points (even no
points). In practice, Step 2 uses an initial design of experiments consisting of a Latin hypercube sample (LHS) of size $N_0 = \min(10, 2d)$, selected according to the maximin criterion. No replicates are needed at this stage.
While $N$ appears to represent the current number of observations, the actual number of iterations depends on the number of unique designs $n$.
The initial TR center $\vecx^0$ can be provided; otherwise, it corresponds to the initial design point with the best predictive mean.
In Step 6, if the number of designs in the TR is too small, new points are added while enforcing a minimal variance reduction based on the latest GP model (either the initial one from Step 2 or one built on the $n_b$ closest points to the TR center in Step 8). These new designs are uniformly sampled in the TR. This allows more variety of distances between pairs of design points, and has been shown to help learning the GP hyperparameters~\citep{zhang2022batch}. Using space filling points is possible, but this will likely add points on the boundary of the TR. As in TRIKE~\cite{regis2016trust} or ORBIT~\cite{wild2008orbit}, the minimal number of designs is set to $d+1$.
For computational speed, in Steps 8 and 13 it is possible to relearn the GP hyperparameters only every few iterations, following a fixed schedule as in~\cite{Lyu2021}.

In Step 9, the acquisition function is optimized over the TR. In our implementation, it starts from a uniformly sampled set of candidates of size $\min(100d, 5000)$ before performing local optimization with L-BFGS-B starting from the best candidate. For EI, since it does not take the replication effort into account, $a_{n+1}$ is simply determined based on the minimal variance reduction $T_a$. For qERCIv2, which involves more variables to optimize, the initial uniform sample is replaced by a particle swarm optimization with 50 particles over 40 iterations. As discussed in Section~\ref{sec:accept}, Step 10 imposes that the future predictive variance at the new point is not too large in order to stabilize the update of the TR center; $a_{n+1}$ may be increased accordingly in Step 11, before evaluating the new design with $a_{n+1}$ replicates (Step 12). The GP model is then updated in Step 13 to estimate whether the progress brought by the new point is sufficient (Step 14). 
If so, the acceptance ratio $\rho_n$ is computed based on the updated GP model (Step 15) by comparing the predictive means at the center and new point with the leave-one-out counterparts of these predictions (see Section~\ref{sec:accept}). If the ratio is sufficiently large, the TR radius is increased (Step 17); otherwise, the TR radius may be decreased (Step 19) depending on the estimated signal-to-noise ratio (see Section~\ref{sec:trreduc}). The same condition is checked if the iteration was unsuccessful (Step 22). This handling of the TR radius is closer in spirit to TRIKE~\cite{regis2016trust} than to TurBO~\cite{eriksson2019scalable}, where the number of consecutive successes or failures triggers TR radius changes.

\begin{algorithm}[htp]
\caption{OGPIT: Optimization by Gaussian Processes In Trust Regions}\label{alg:OGPIT}
\SetKwInOut{Require}{Require}
\Require{
    $0 < \gammad (0.8) < 1 < \gammai (1/\gammad)$,
    $0 < \eta (0.2),
    \beta (10^{-3}) < 1$,
    $n_{b}$ (number of neighbors),
    $N_0~(2\times d$, initial number of observations),
    acquisition function $\alpha$,
    $\Nmax$ (budget), $\pmax$ (max number of replication),
    $\vecx^0$ (initial TR center),    
    $0 < \Delta_0$ (initial TR radius),
    $T_a$ (0.2, minimal variance reduction). 
}
Set $\vecx_c = \vecx^0$, $\Delta = \Delta_0$\;
Construct initial design of experiments: e.g., maximin LHS of size $N_0$\;
$N=N_0$\;
Build initial GP model\;

\While{$N \leq N_{max}$ and $\Delta > \Delta_{min}$}{

    \If{Number of points in the TR $< d+1$}{Augment the DoE to have $d+1$ points\;}

    Update GP model $Y_N$ (mean $m_n$, variance $s^2_n$) conditioned on the $n_{b}$ nearest neighbors of $\vecx_c$ (normalize inputs in the TR, normalize outputs)\;
    
    Compute $\s, a_{n+1} \in \arg \underset{\s: \left\| s \right\| \le \Delta,~ p \leq \pmax}{\max} \alpha(\vecx_c + \s)$\;
    \If{$s_{n+1}^2(\vecx_c + \s) > 4 s_n^2(\vecx_c)$}{Increase $a_{n+1}$ (up to $\pmax$) based on~\eqref{eq:pa}\;}
    Evaluate $f(\vecx_c + \s)$ with replication number $a_{n+1}$ \;
    Update GP: $Y_{N+a_{n+1}}$: $m_{n+1}, s_{n+1}^2$\;
    \If{$m_{n+1}(\vecx_c) - m_{n+1}(\vecx_c+\s) \ge \beta \min \left\{ \Delta, \Delta^2 \right\}$}
    {
    Calculate $\displaystyle \rho_n = \frac{m_{n+1}(\vecx_c) - m_{n+1}(\vecx_c + \s)}{ \tilde{m}_{n+1}(\vecx_c)-\tilde{m}_{n+1}(\vecx_c+\s)}$\;
    \eIf{$\rho_n \ge \eta$ and $s_{n+1}^2(\vecx_c + \s) \leq 4 s_n^2(\vecx_c)$}
    {
      $\vecx_{c} = \vecx_c + \s$; $\Delta = \gammai \Delta$\;
    }
    {
      \If{$\Var{m_n(X)} \geq 10 \mathbb{E}[s_n^2(X)]$}{$\Delta = \gammad \Delta$\;}
    }
  }\Else
  {
    \If{$\Var{m_n(X)} \geq 10 \mathbb{E}[s_n^2(X)]$}{$\Delta = \gammad \Delta$\;}
  }
 }
  $N = N+a_{n+1}$ \;
  Return $\vecx_c$ and estimated value $m_n(\vecx_c)$

\end{algorithm}

The effect of the various parameters of the approach ($\gamma_{dec}, T_a$, the future predictive variance threshold in Step 10, or the IMSE criterion in Step 22) is studied in Appendix~\ref{ap:A}. Summarizing the results, $\gamma_{dec}$ has little influence compared to the IMSE parameter. In particular, reducing the IMSE constant, which makes the TR reduction criterion less strict, degrades performance. The effect of the minimal variance reduction $T_a$ shows that it is not necessary to reduce it too aggressively, e.g., by 50\%. The same conclusion applies to the variance ratio between the TR center and the new design: a reduction that is too large is detrimental, while reducing it excessively also starts degrading performance.
The remaining parameter values ($\beta, \eta$) take acceptable values from the literature.

\section{Empirical evaluation}

We now compare the capabilities of the proposed GP-based TR method for
stochastic optimization, dubbed OGPIT for optimization by Gaussian processes
in trust regions. In a first phase we test the ability to reach a local
optimum efficiently; then we evaluate the proposed cost-aware version.

We focus on local optimization to avoid confusion coming from the global
optimization level. Also, issues with existing methods become preponderant when
the signal-to-noise ratio decreases, as happens when the trust-region radius
becomes small. We consider two sets of benchmark functions. 
Benchmark 1 includes simple validation test functions plus a typical function
from the BO literature (i.e., Branin): sphere function (2d, 4d, 6d), squared sphere
function(2d, 4d, 6d), Branin (2d), and Rosenbrock (2d, 4d). 
Benchmark 2 is obtained from the BenDFO repository~\cite{bendforepo,more2009benchmarking}, which consists of nonlinear least squares problems of the form $\sum \limits_{i=1}^m f_i(\vecx)^2$. We restrict the analysis to the smooth problem type and instances of dimension lower than 10, resulting in 36 different test problems. 

As baselines, we compare with the state-of-the-art GP-based trust-region
optimizer TuRBO~\cite{eriksson2019scalable} and regular BO from the BoTorch
\texttt{Python} library~\cite{balandat2020BoTorch}. We remark that TuRBO, while it can run on stochastic problems, is
initially dedicated for deterministic problems and thus is not expected to
find precise solutions. BoTorch is a global optimization method, hence
not necessarily suited for precise convergence on local optimization
problems. Closer to polynomial-surrogate-based TR methods, we also run SNOWPAC on Benchmark 1, as it also relies on a GP model.
Since these methods do not exploit replication, the computational cost
increases quickly with $N$. To keep reasonable timings, we enforce a 10-minute stopping
time, which is too small for SNOWPAC on Benchmark 2.

TuRBO relies on the Thompson sampling acquisition function, BoTorch and OGPIT on EI, while SNOWPAC optimizes the mean prediction directly. TuRBO, BoTorch and SNOWPAC evaluate only a single new point without replicate at each iteration, while OGPIT uses $p_a$ replicates (and $\pmax = 500)$.

We let each method use its own initial design procedure, but with the same number of initial designs. For SNOWPAC, we provide a starting solution, as detailed in  \Cref{tab:startP}. For Benchmark 2, as it contains starting points, these are also passed to all methods.

We entertain versions of the test functions with no noise for validation, then with zero mean
Gaussian additive noise with standard deviation equal to 0.001, 0.01, and 0.1
(resp.\ 0.001, 0.1, and 10) for Benchmark 1 (resp.\ 2). We compute the regret
over iterations, that is, the difference to the reference solution. For OGPIT,
the regret is computed on the trust-region center value. For BoTorch, SNOWPAC and
TuRBO, it is computed on all evaluated points, hence favoring these methods.
The results are summarized through data profiles of the regret. For a
given amount of evaluations, a data profile shows the proportion of problems solved for a
given precision, based on the best solution found among the considered methods. Data profiles
are detailed in~\cite{more2009benchmarking}. Results are normalized by $d+1$ to compare problems of various dimensions.

Each test problem is repeated 30 times since results depend on the initial
designs or starting point. Two versions of the code, one in \texttt{R} and
one in \texttt{Python}, are available here: \url{https://github.com/POptUS/StochDFO}.

\subsection{Local optimization without setup cost}

\begin{figure}
\includegraphics[width=0.333\textwidth, trim= 0 0 40 15, clip]{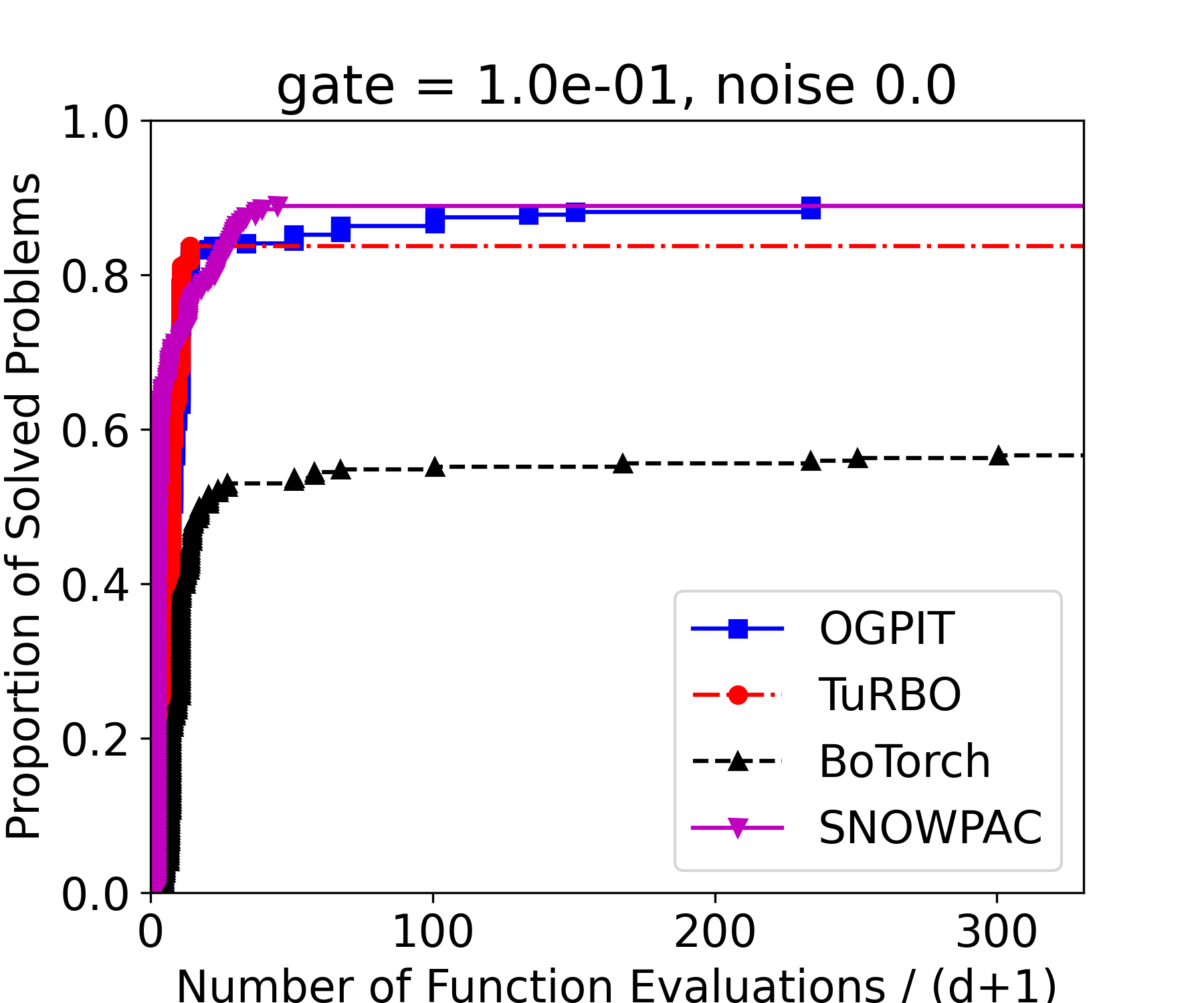}%
\includegraphics[width=0.333\textwidth, trim= 0 0 40 15, clip]{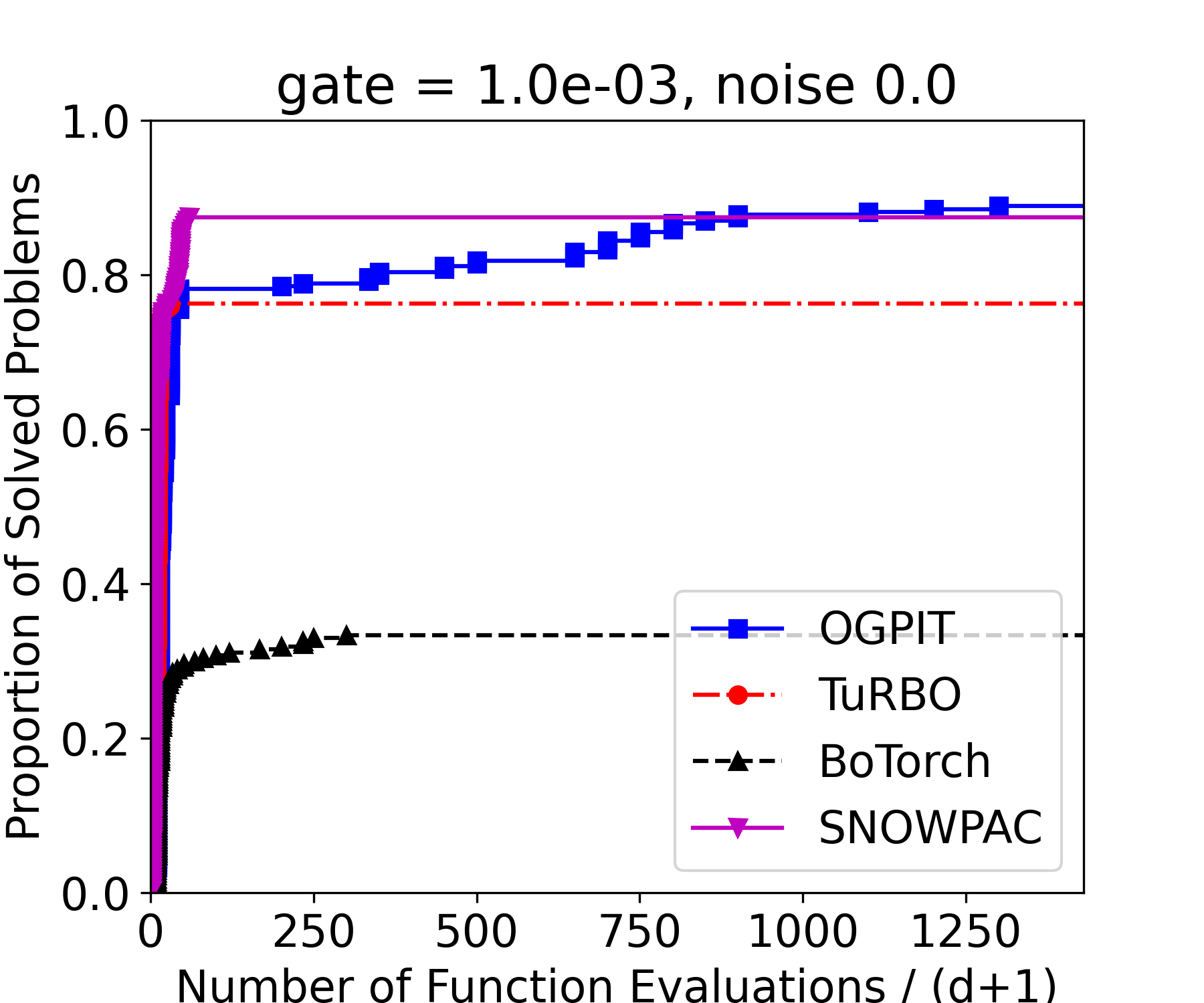}%
\includegraphics[width=0.333\textwidth, trim= 0 0 40 15, clip]{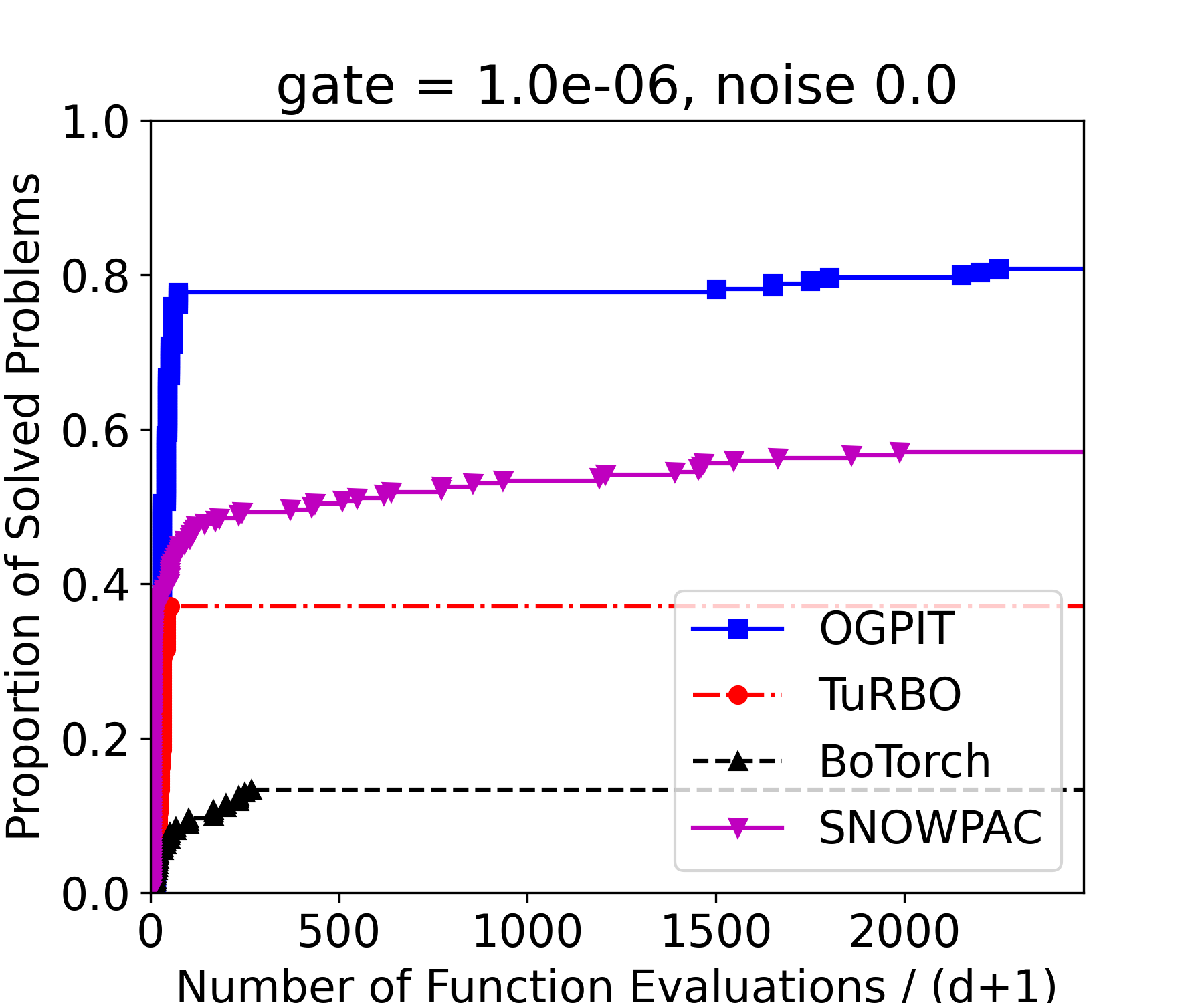}\\
\includegraphics[width=0.333\textwidth, trim= 0 0 40 15, clip]{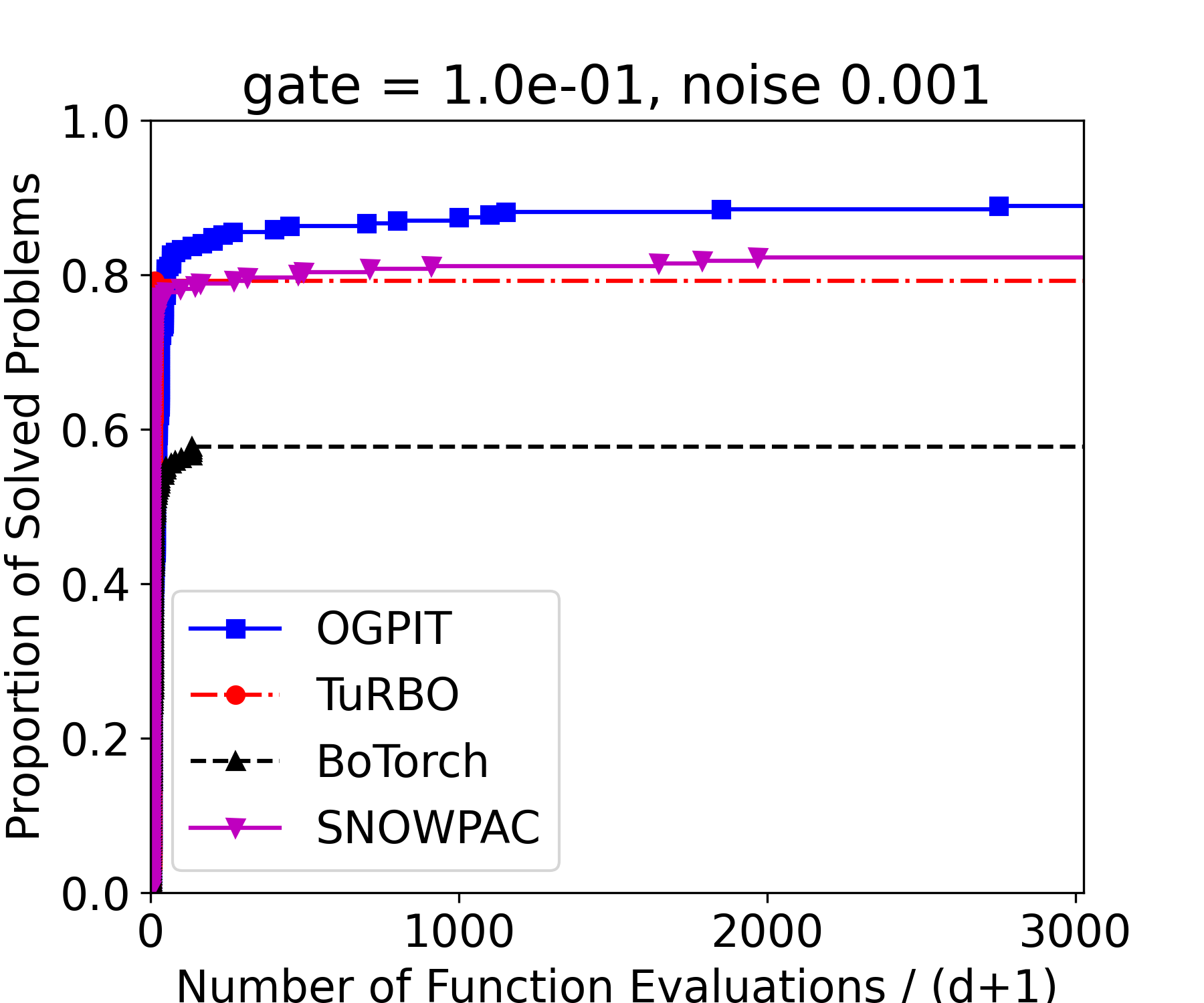}%
\includegraphics[width=0.333\textwidth, trim= 0 0 40 15, clip]{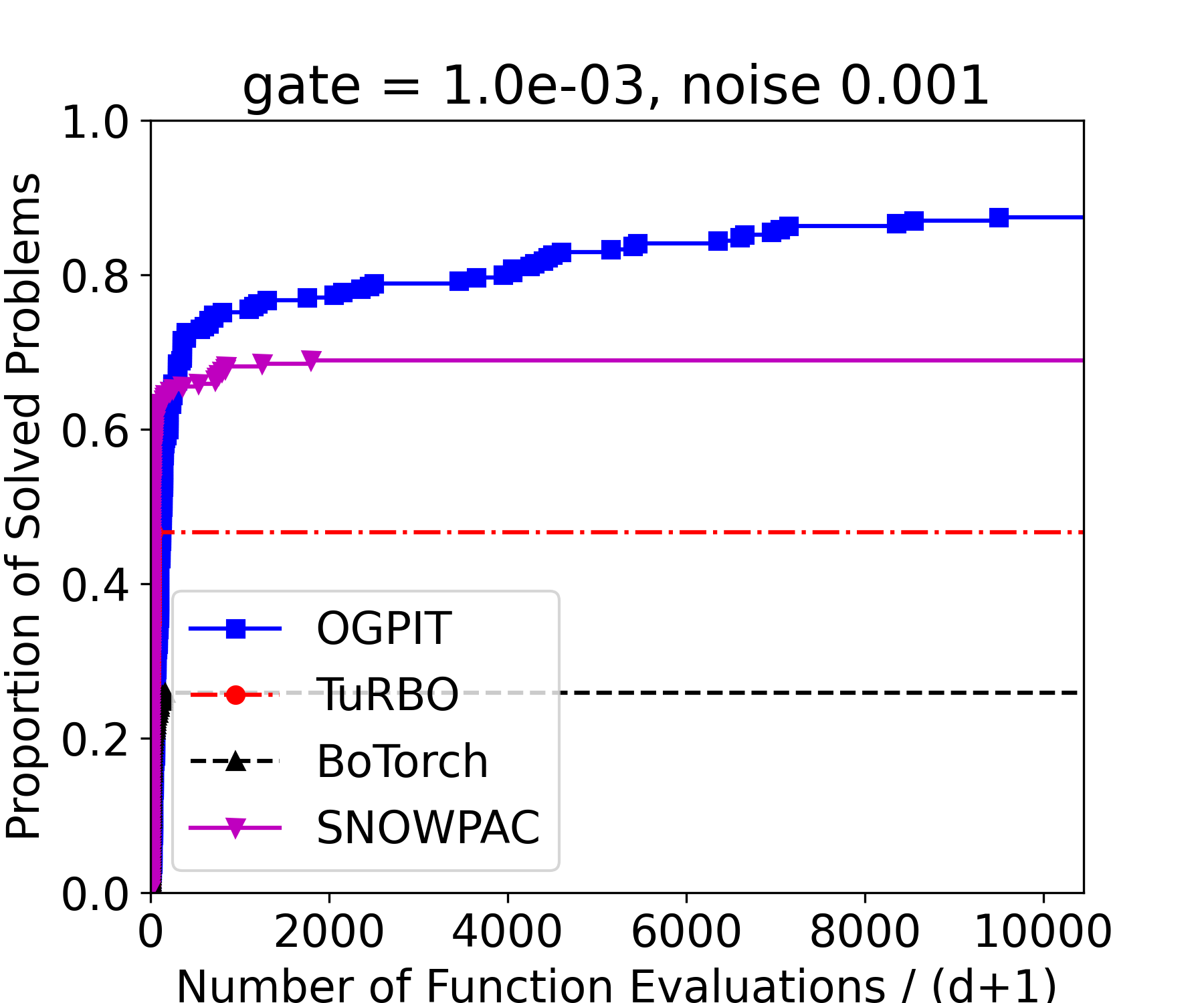}%
\includegraphics[width=0.333\textwidth, trim= 0 0 40 15, clip]{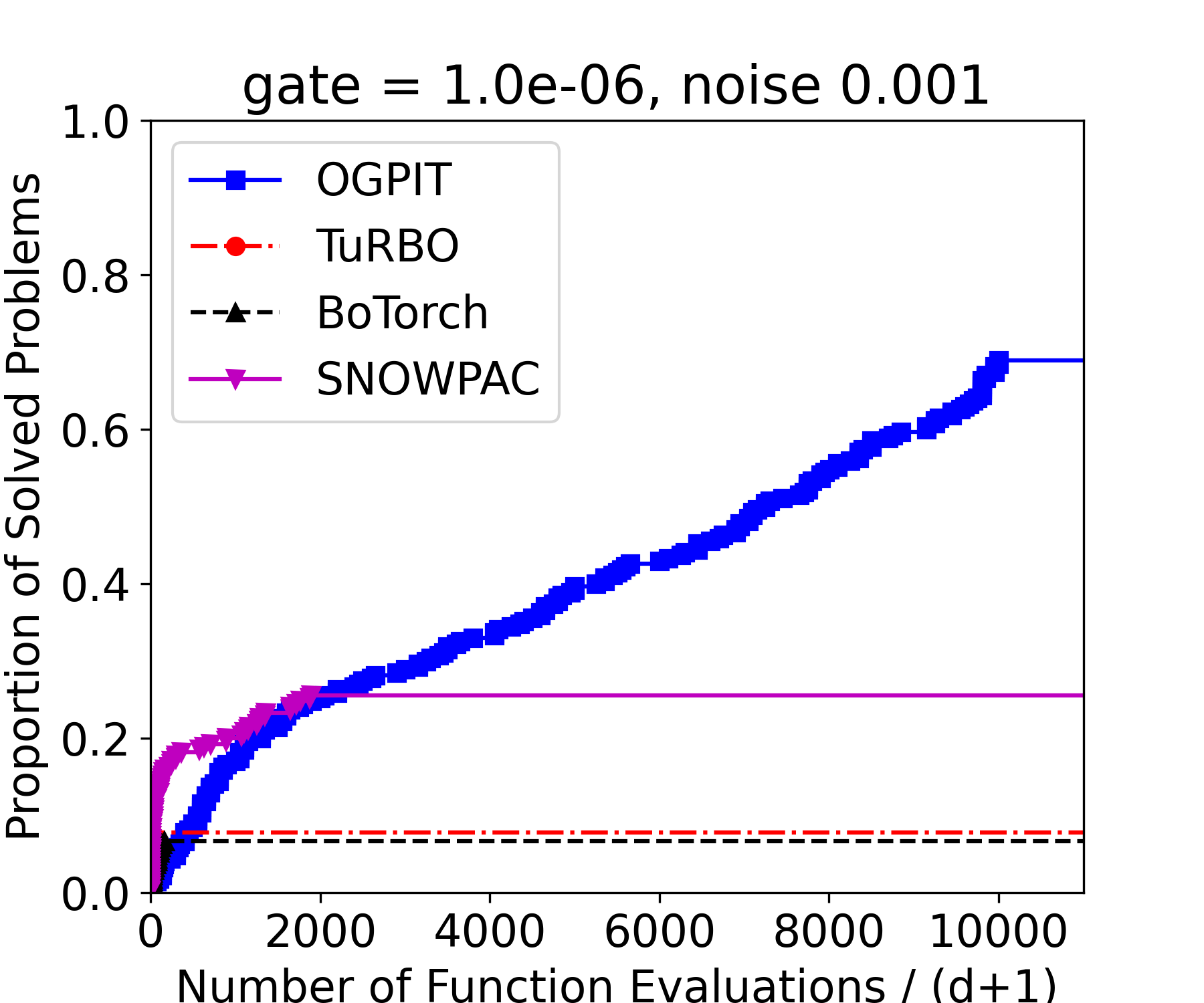}\\
\includegraphics[width=0.333\textwidth, trim= 0 0 40 15, clip]{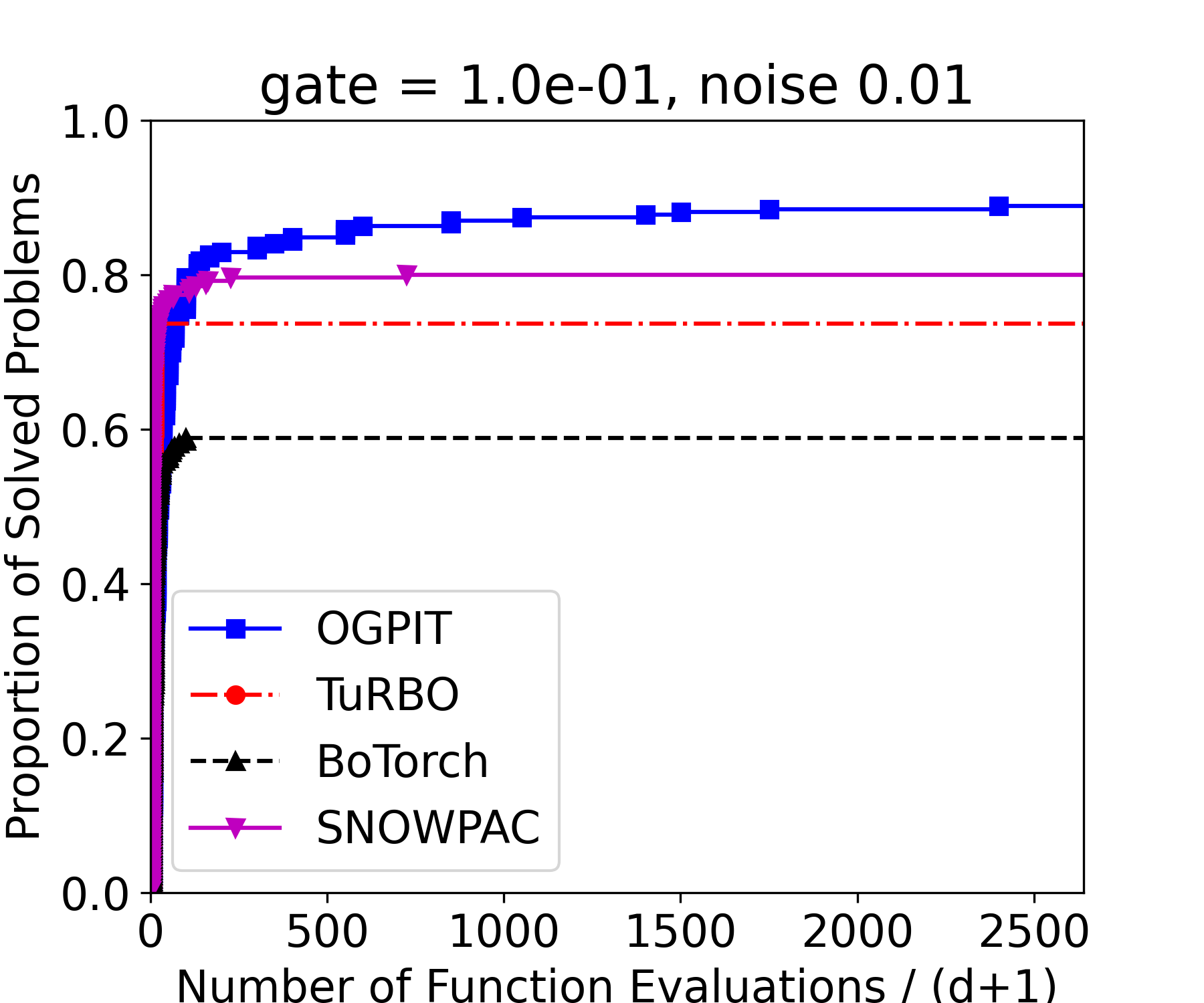}%
\includegraphics[width=0.333\textwidth, trim= 0 0 40 15, clip]{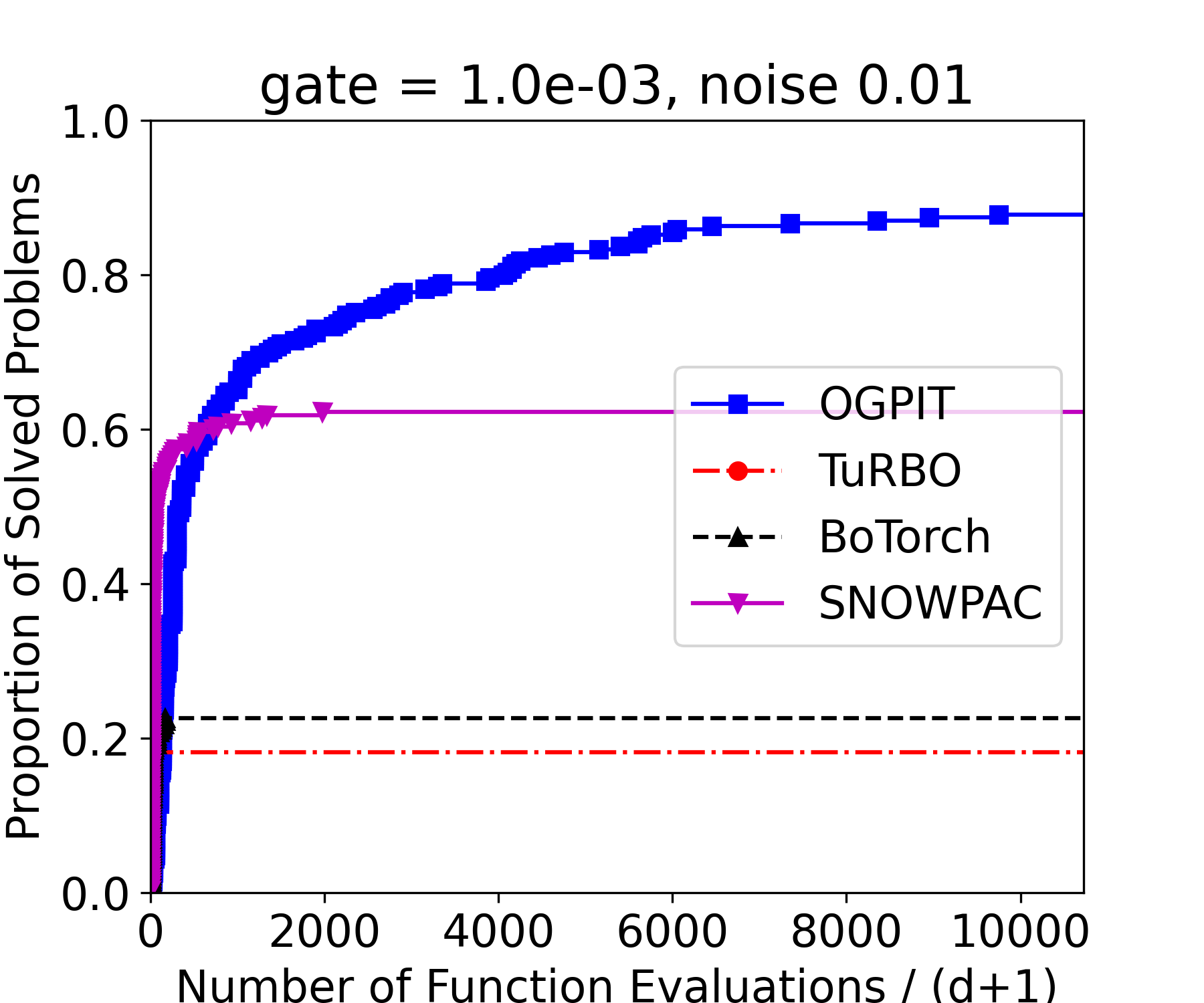}%
\includegraphics[width=0.333\textwidth, trim= 0 0 40 15, clip]{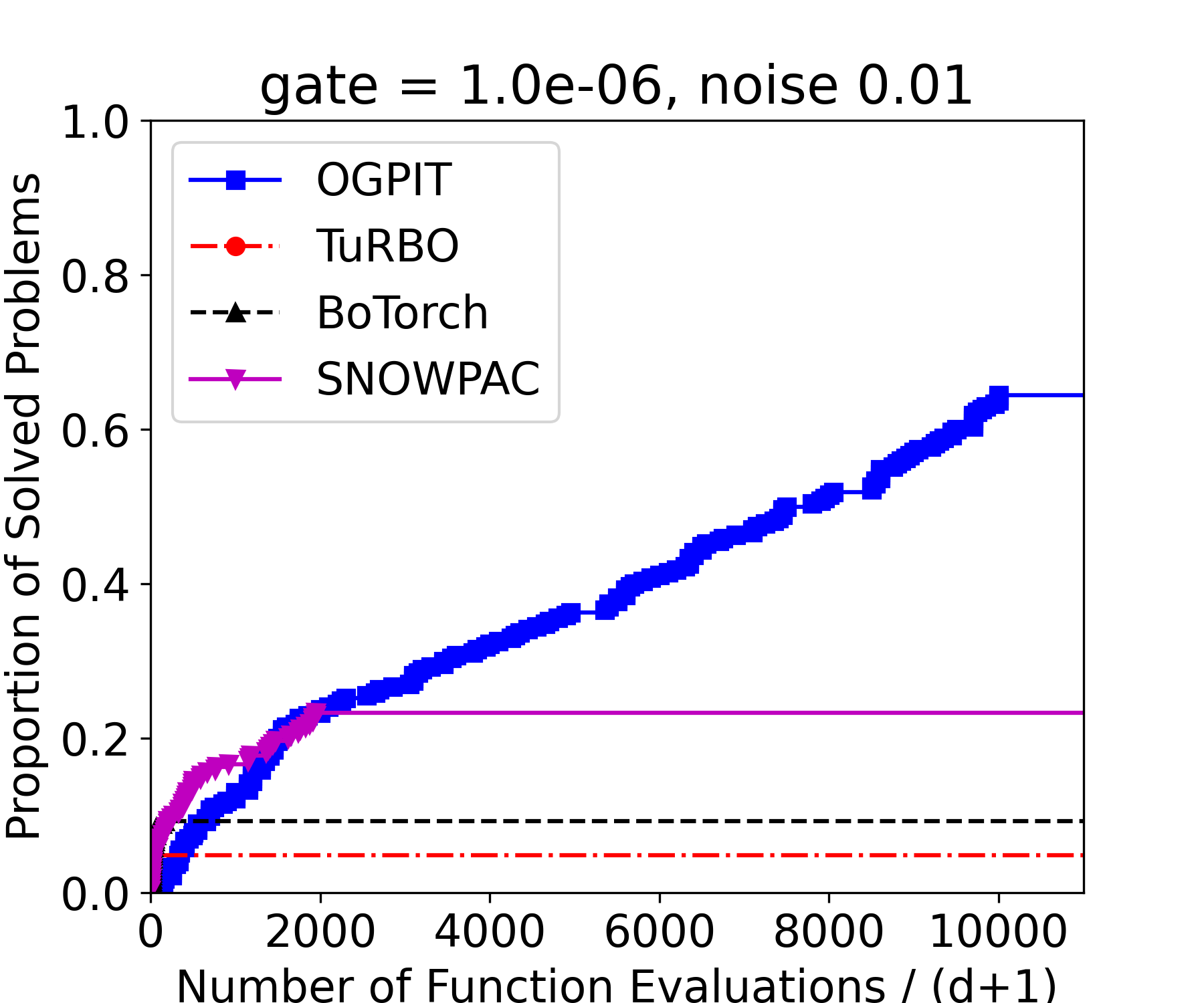}\\
\includegraphics[width=0.333\textwidth, trim= 0 0 40 15, clip]{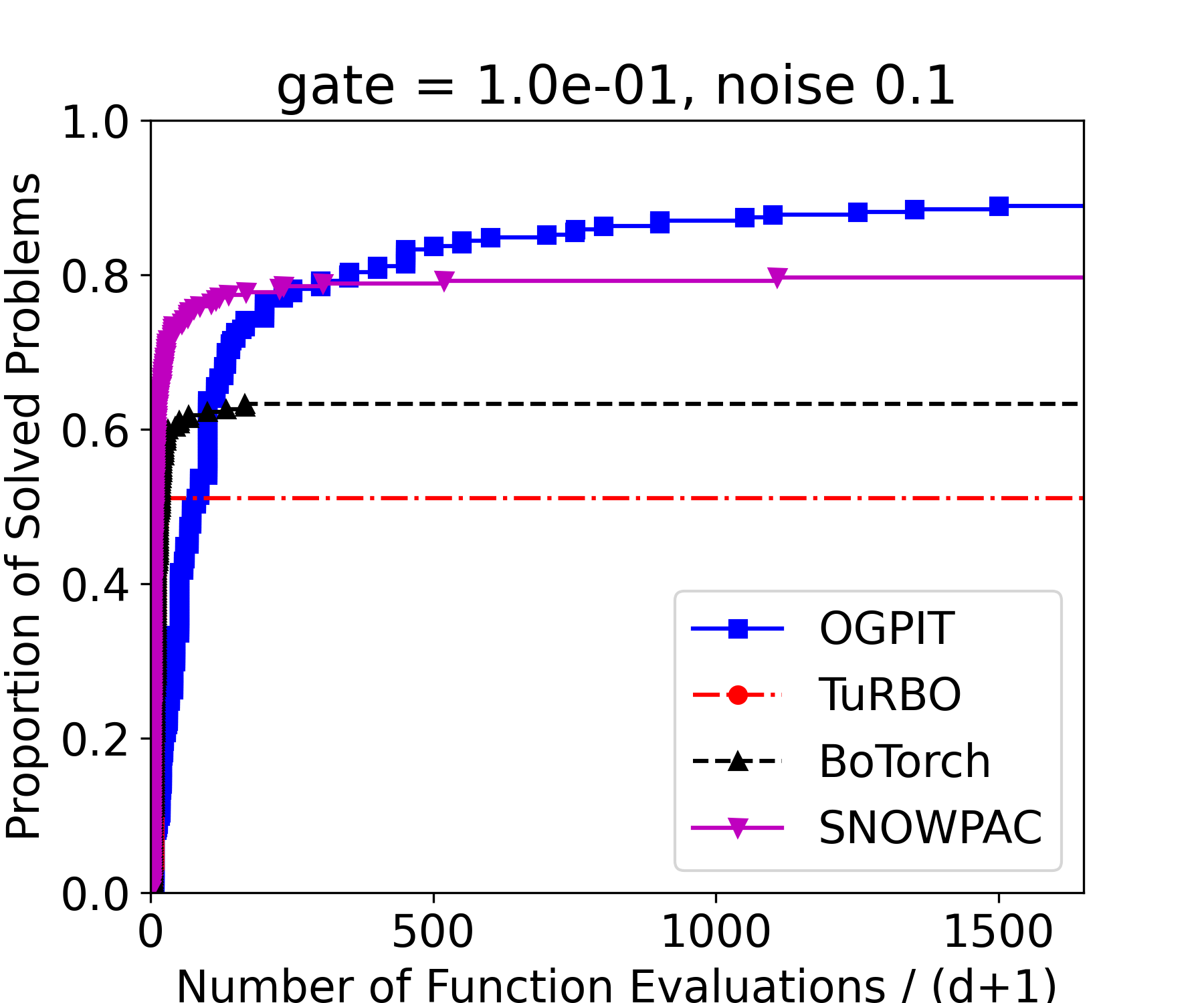}%
\includegraphics[width=0.333\textwidth, trim= 0 0 40 15, clip]{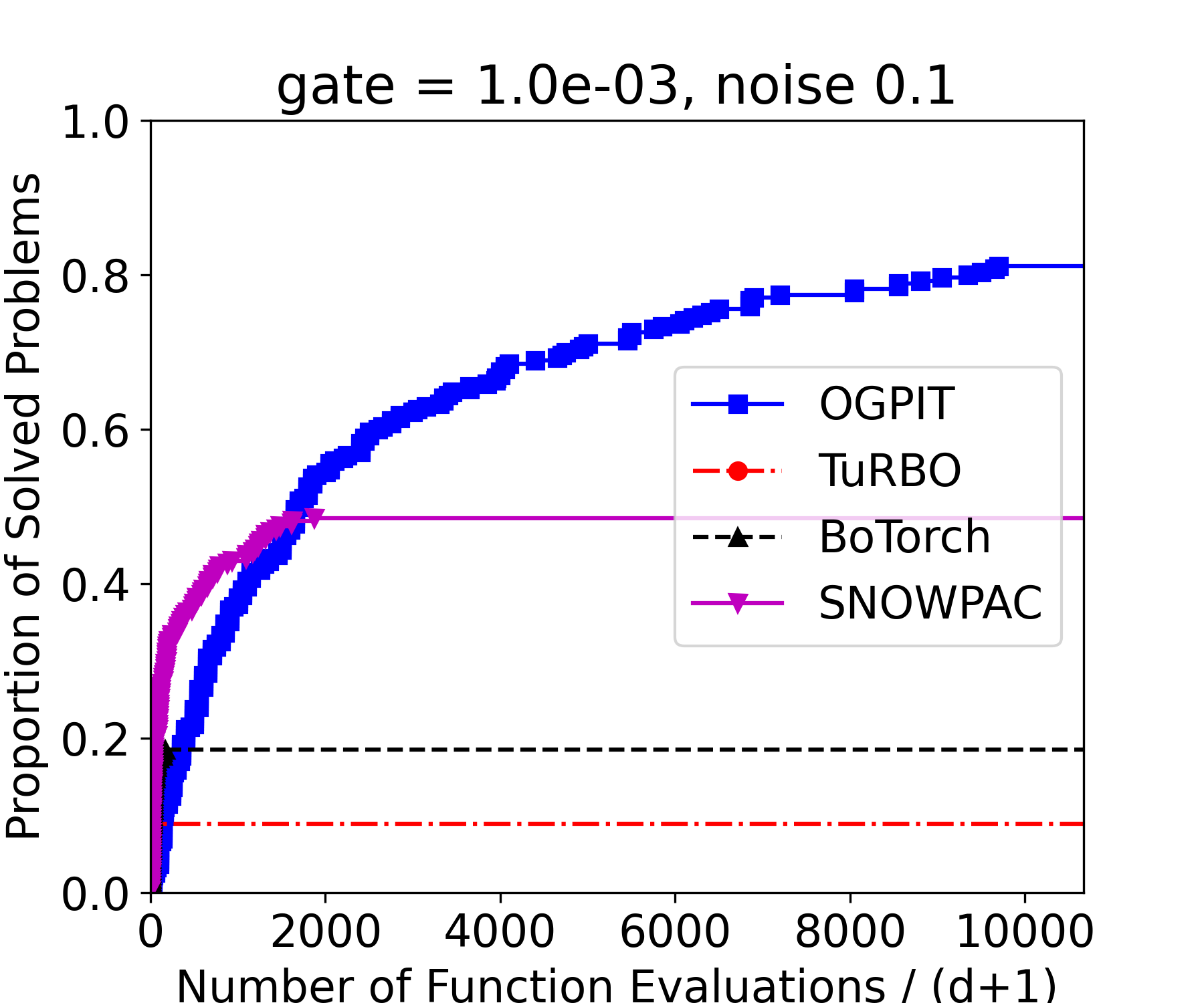}%
\includegraphics[width=0.333\textwidth, trim= 0 0 40 15, clip]{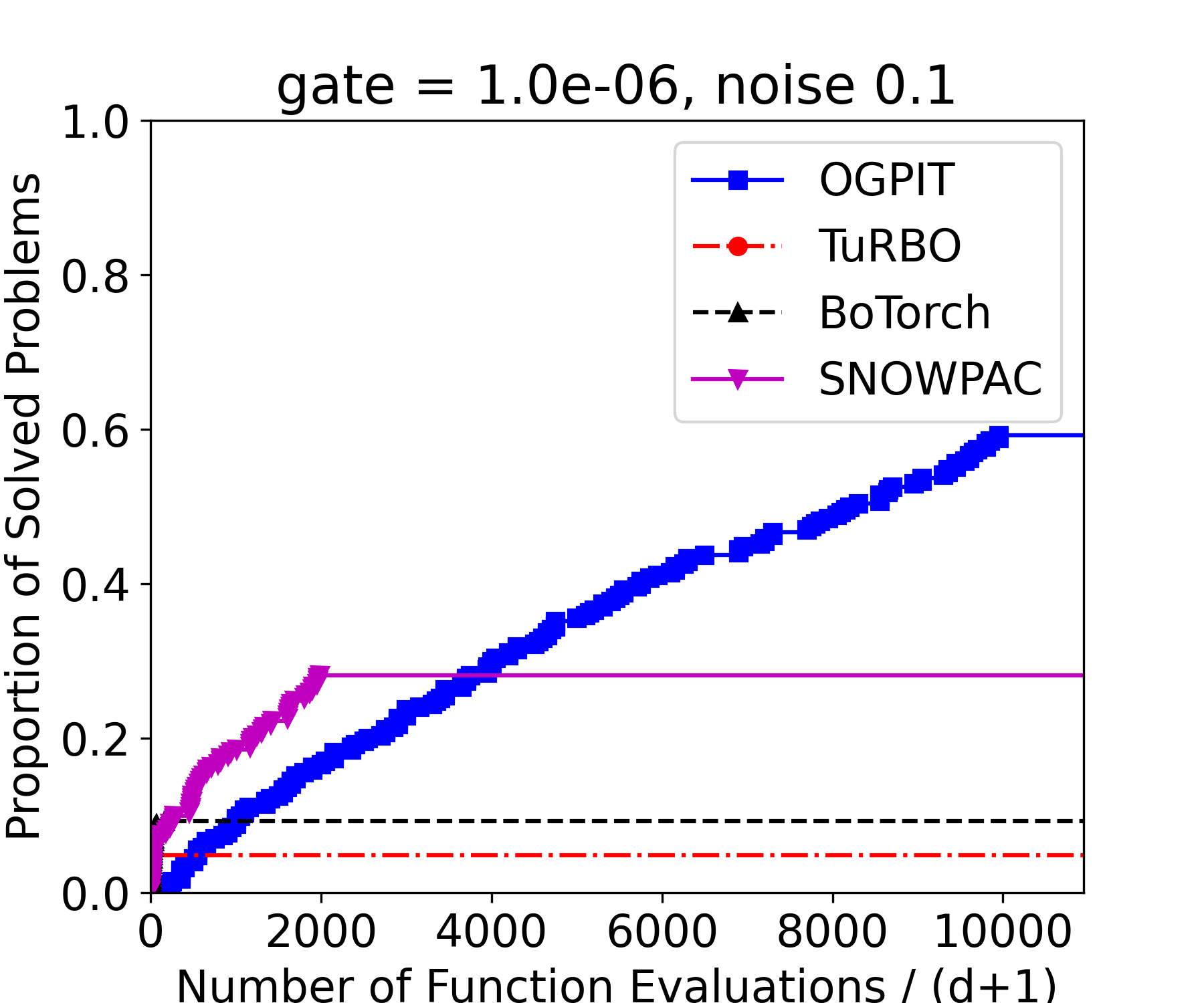}
\caption{Data profiles of the method against TuRBO, SNOWPAC and BoTorch on Benchmark 1.}
\label{fig:localbenchprofs1}
\end{figure}

We use a maximum budget of $\Nmax=10^4(d+1)$. The performance is compared
first on Benchmark 1, with results provided in \Cref
{fig:localbenchprofs1}. For noiseless problems, OGPIT and TuRBO 
perform similarly for the highest gates, but OGPIT and SNOWPAC outperforms on the
smaller ones. BoTorch underperforms for all cases, failing to converge
precisely on these simple problems (only $\sim$ 60\% of problems solved even
for the largest gate). As the noise variance increases, the performance of
TuRBO degrades compared with OGPIT, especially for the smallest gates. SNOWPAC shows the best results in the early stages, helped by its quadratic surrogate on these convex test problems. As the noise increases, the performance degrades compared to OGPIT. Then, as it cannot scale to larger number of iterations, its performance plateaus.

We illustrate the outcome of the adaptive replication strategy in \Cref{ap:B}, where we show in \Cref{fig:repanalysis} how the number of replication evolves over iterations for OGPIT. It depends both on the problem instance and its dimension, which makes fixing a replication number a priori difficult. Indeed, fixing the replication effort in OGPIT degrades the performance, see \Cref{fig:fixedrep}.

\begin{figure}
\includegraphics[width=0.333\textwidth, trim= 0 0 40 15, clip]{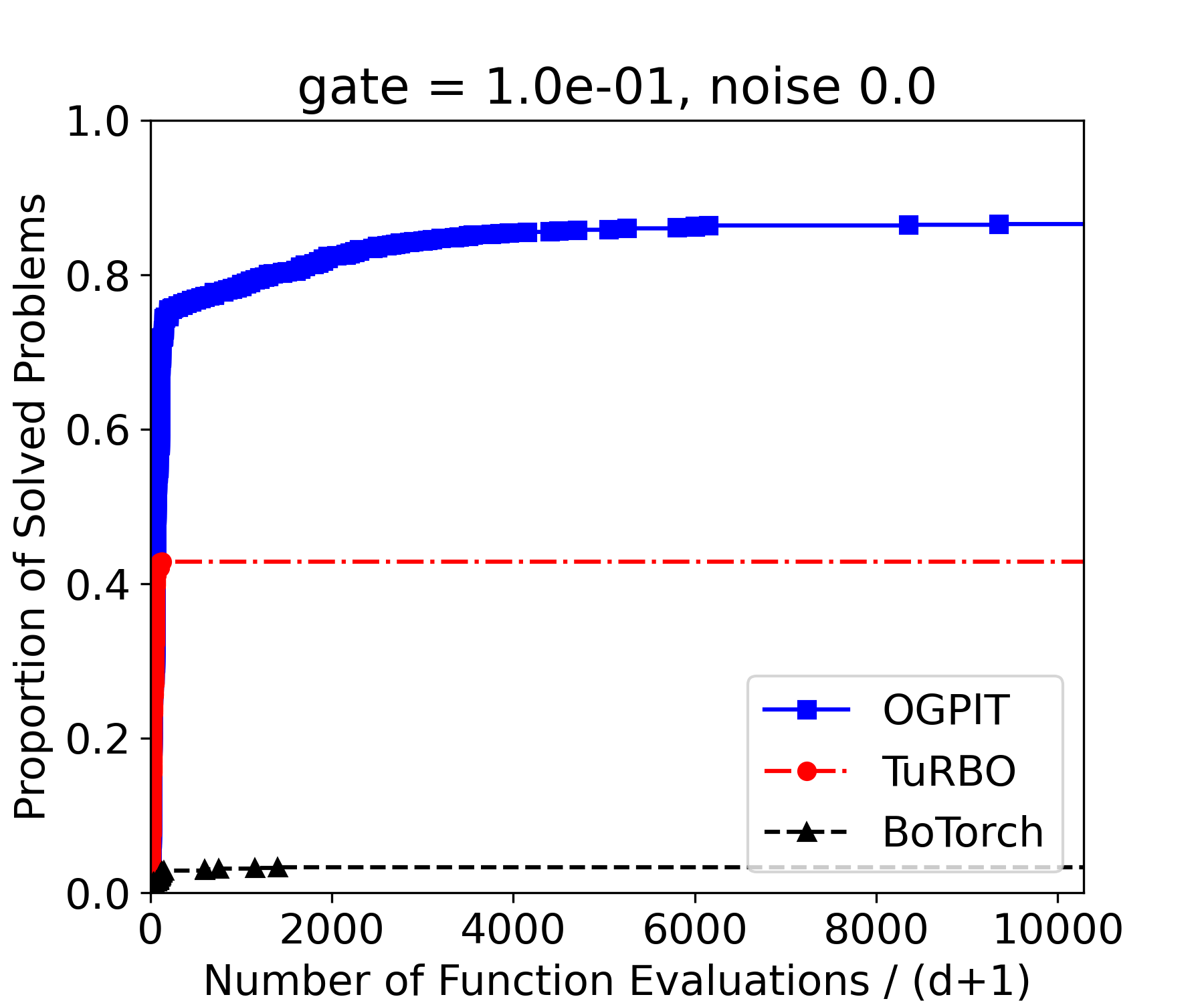}%
\includegraphics[width=0.333\textwidth, trim= 0 0 40 15, clip]{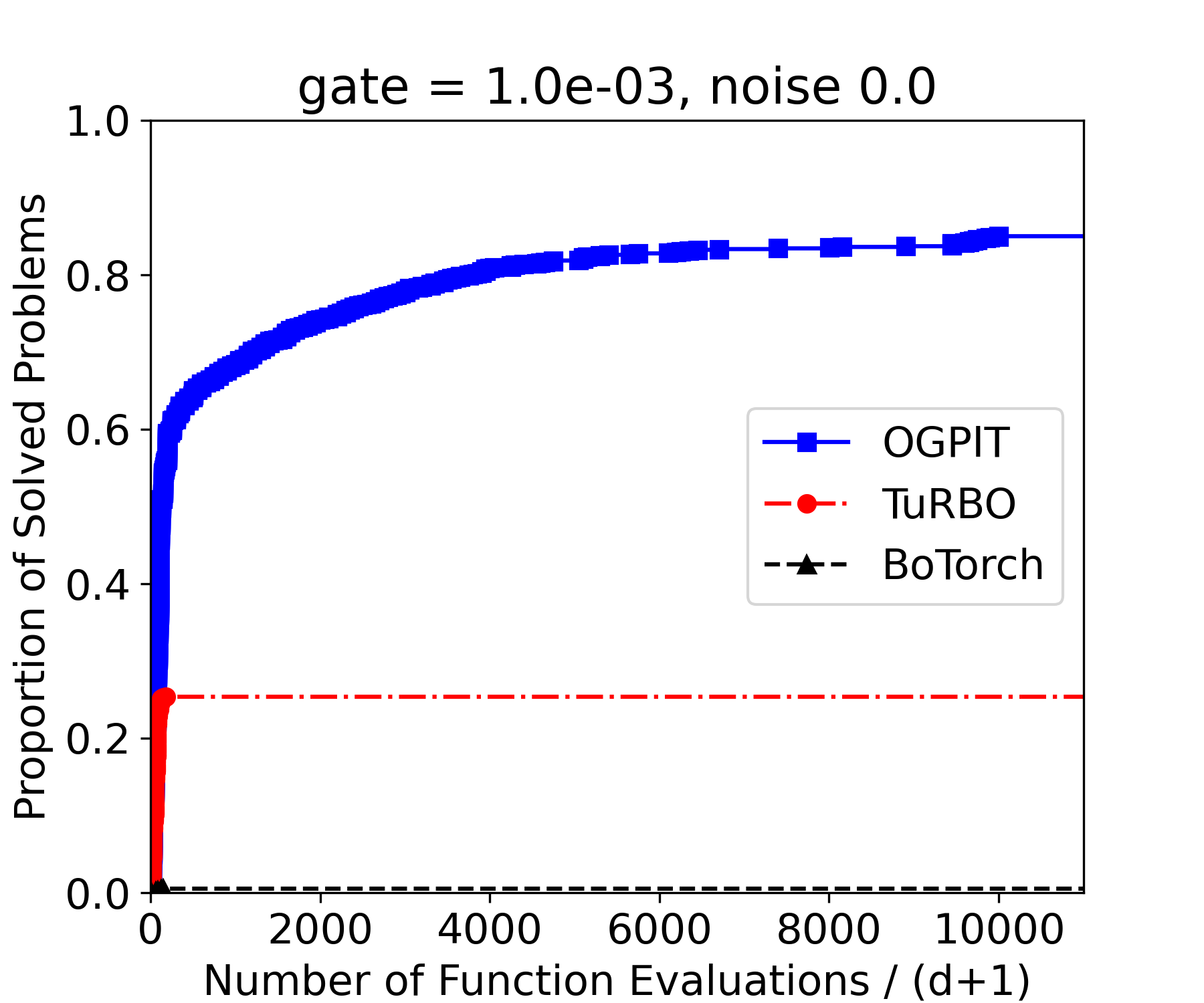}%
\includegraphics[width=0.333\textwidth, trim= 0 0 40 15, clip]{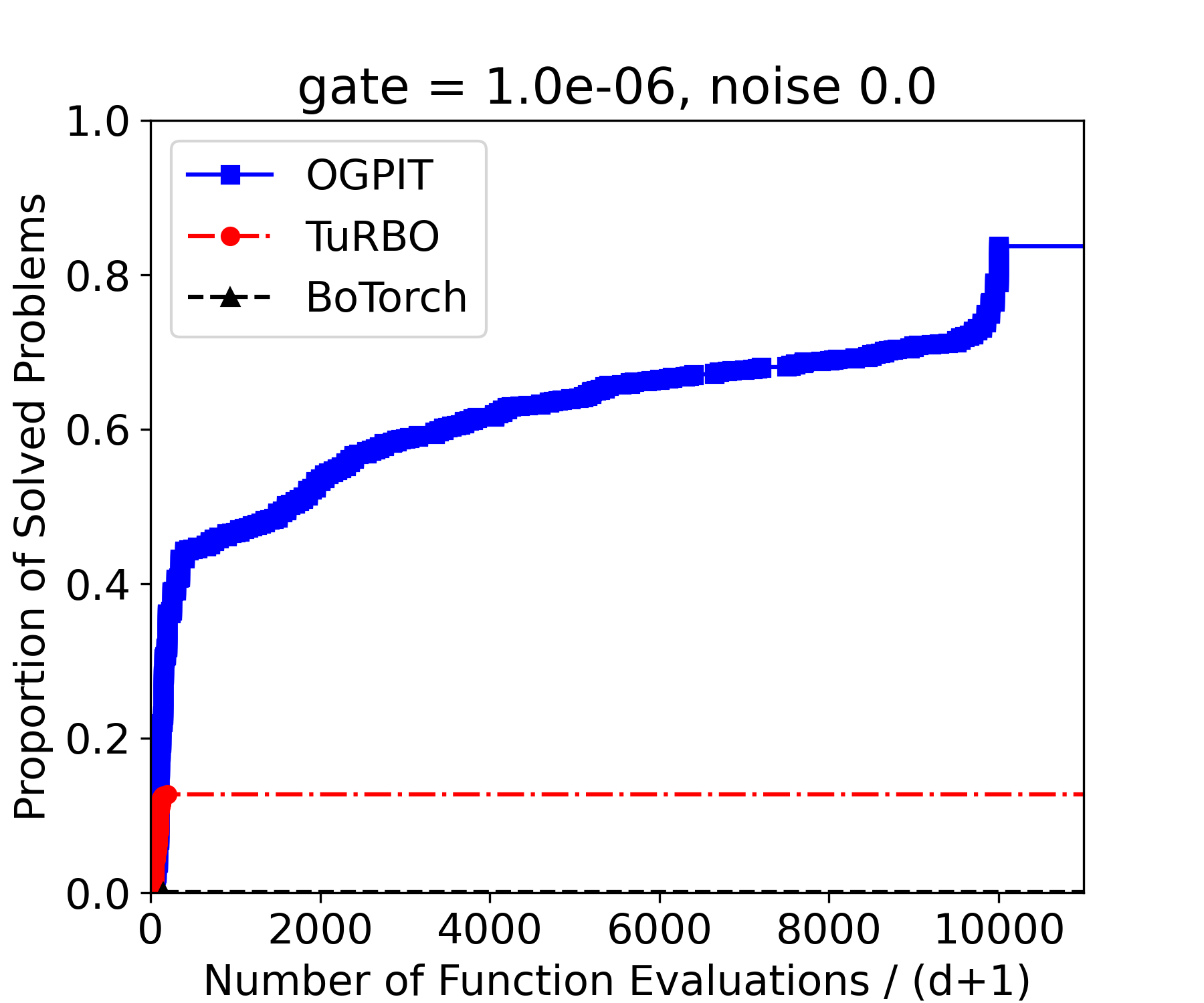}\\
\includegraphics[width=0.333\textwidth, trim= 0 0 40 15, clip]{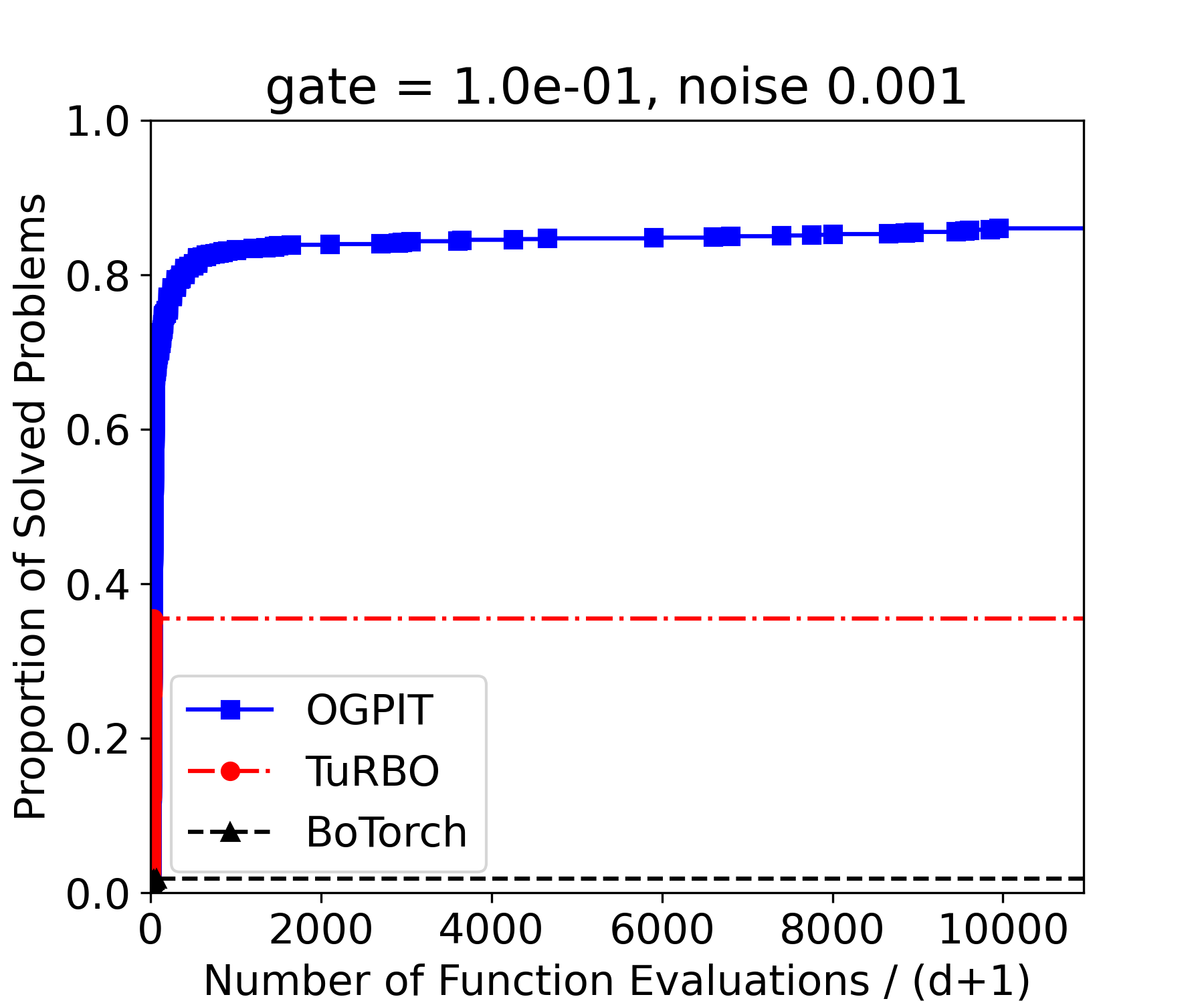}%
\includegraphics[width=0.333\textwidth, trim= 0 0 40 15, clip]{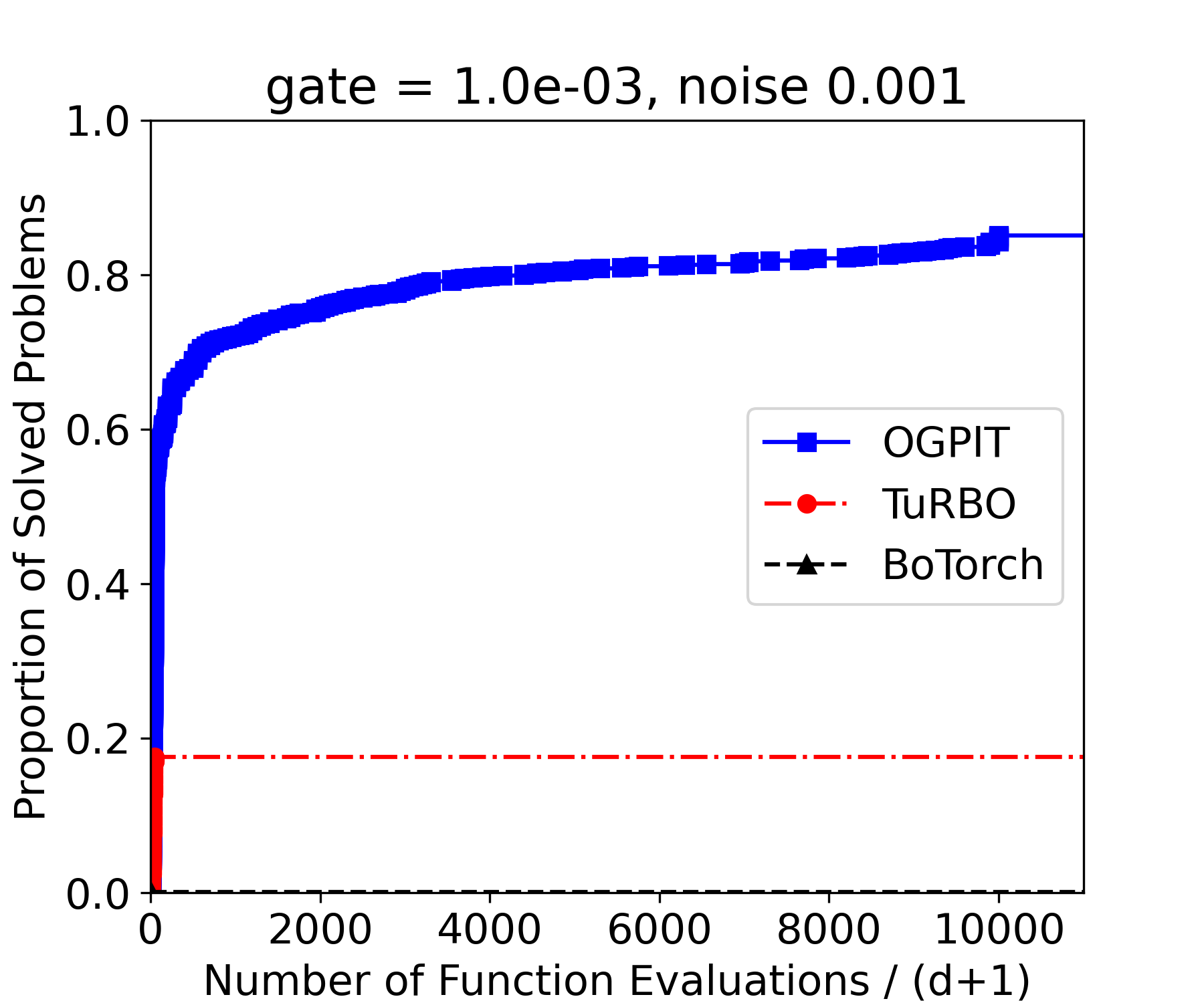}%
\includegraphics[width=0.333\textwidth, trim= 0 0 40 15, clip]{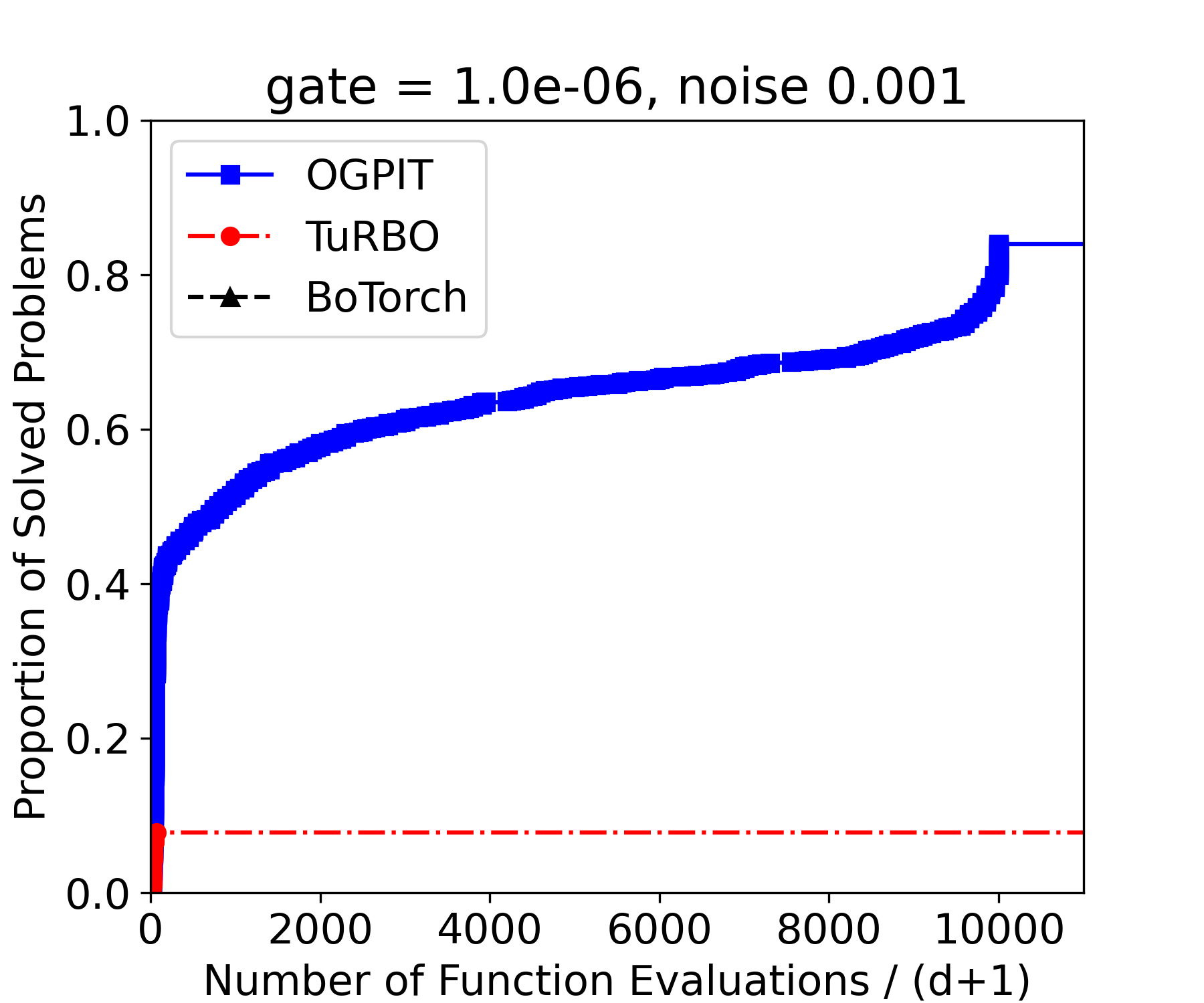}\\
\includegraphics[width=0.333\textwidth, trim= 0 0 40 15, clip]{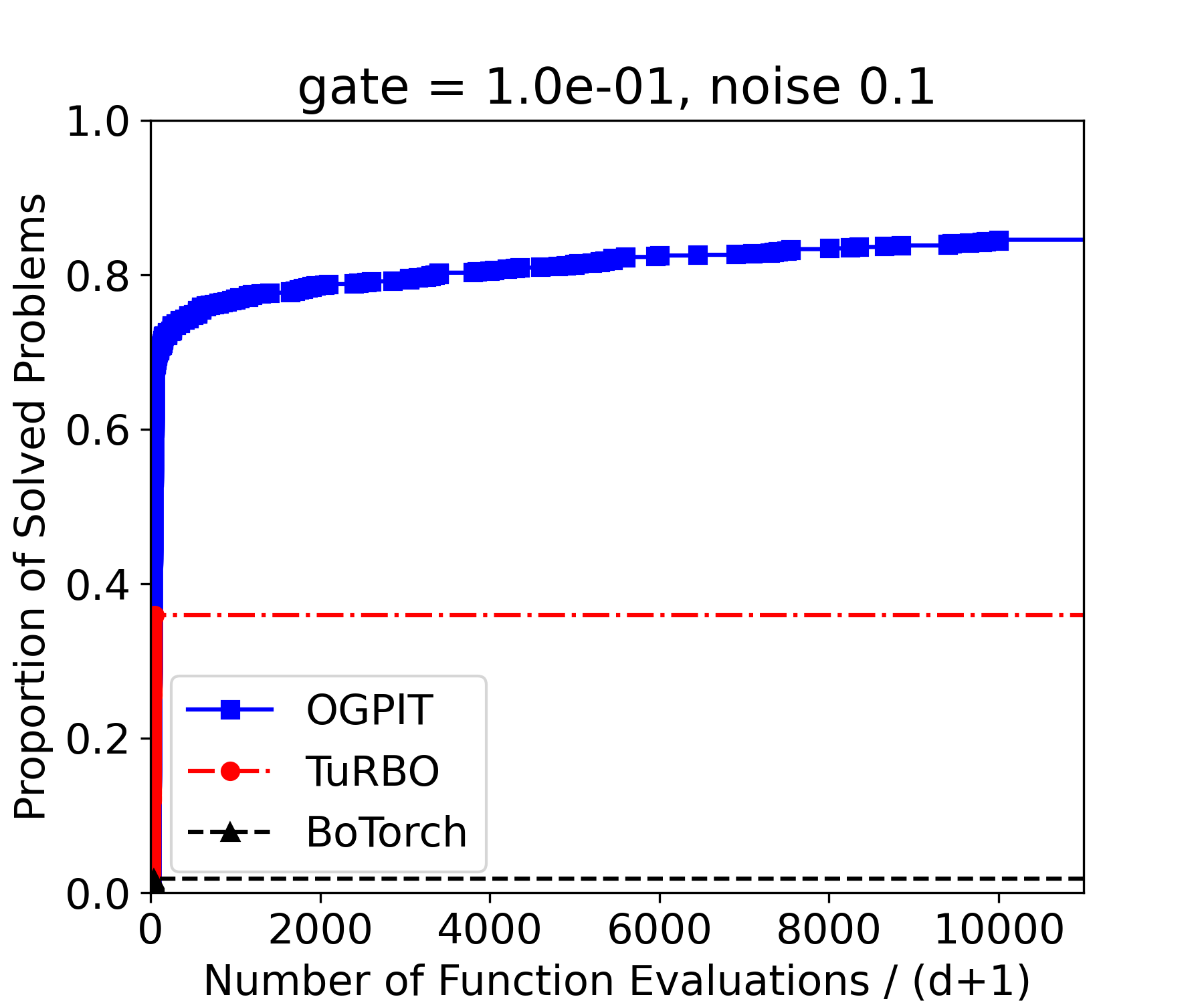}%
\includegraphics[width=0.333\textwidth, trim= 0 0 40 15, clip]{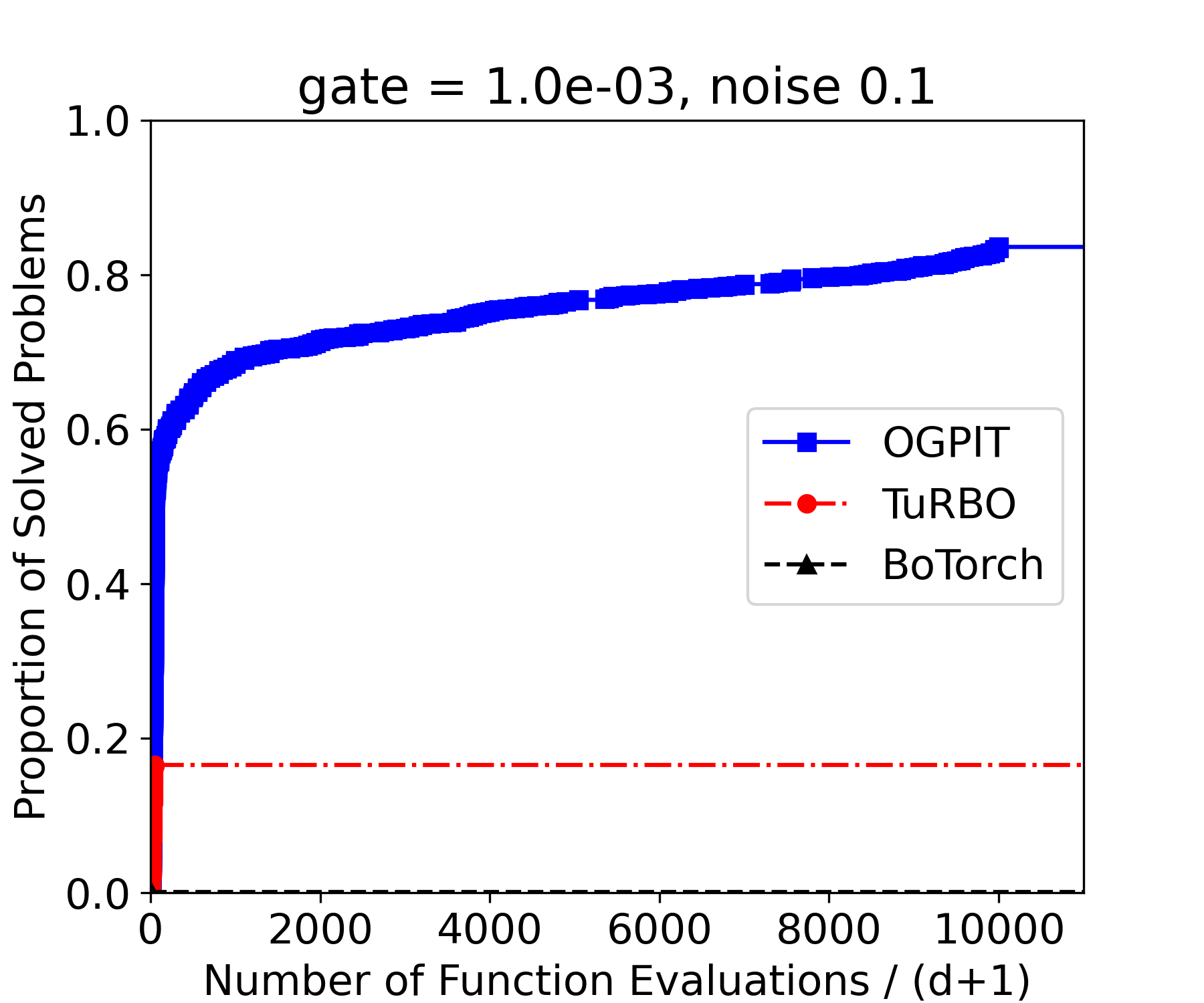}%
\includegraphics[width=0.333\textwidth, trim= 0 0 40 15, clip]{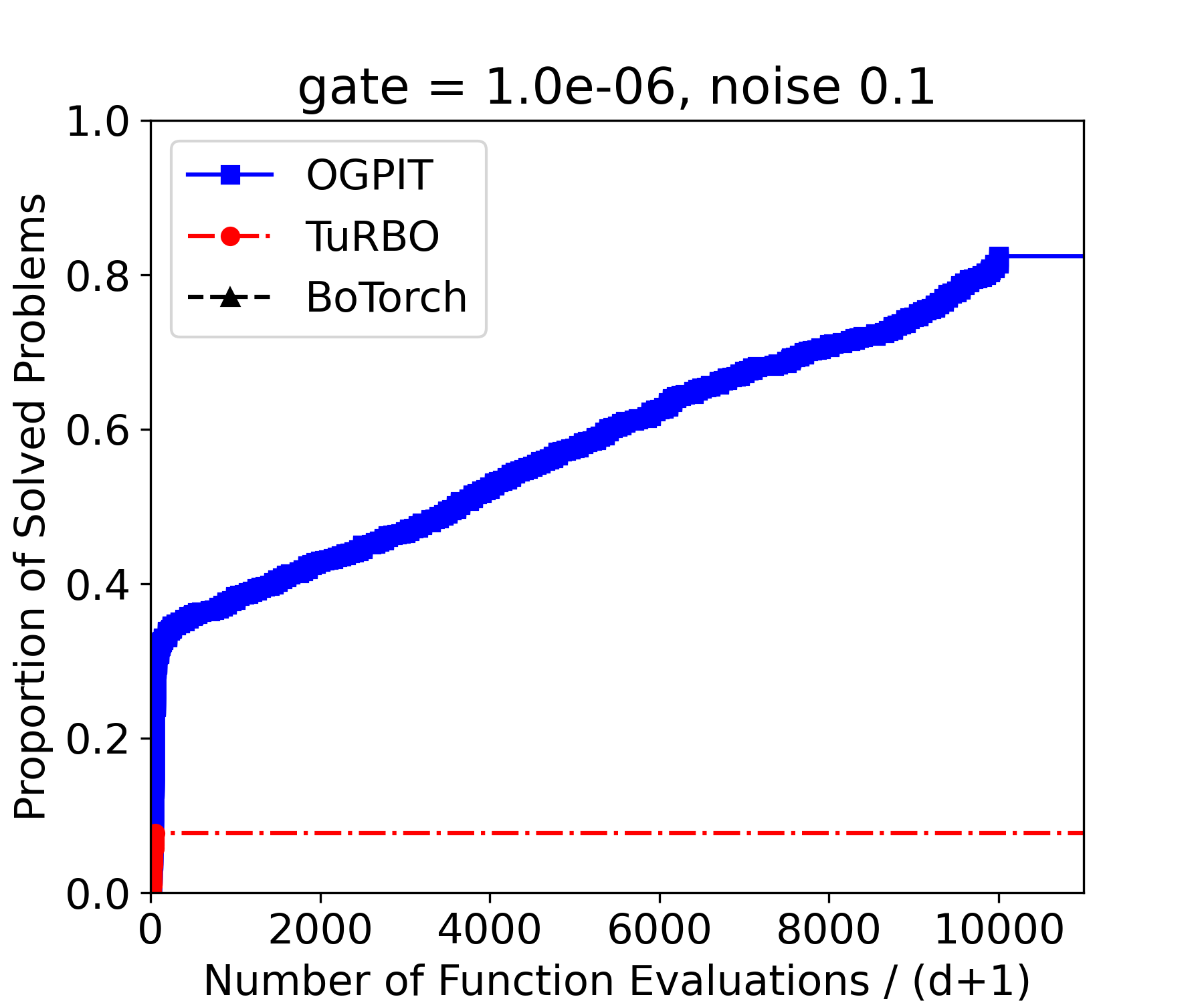}\\
\includegraphics[width=0.333\textwidth, trim= 0 0 40 15, clip]{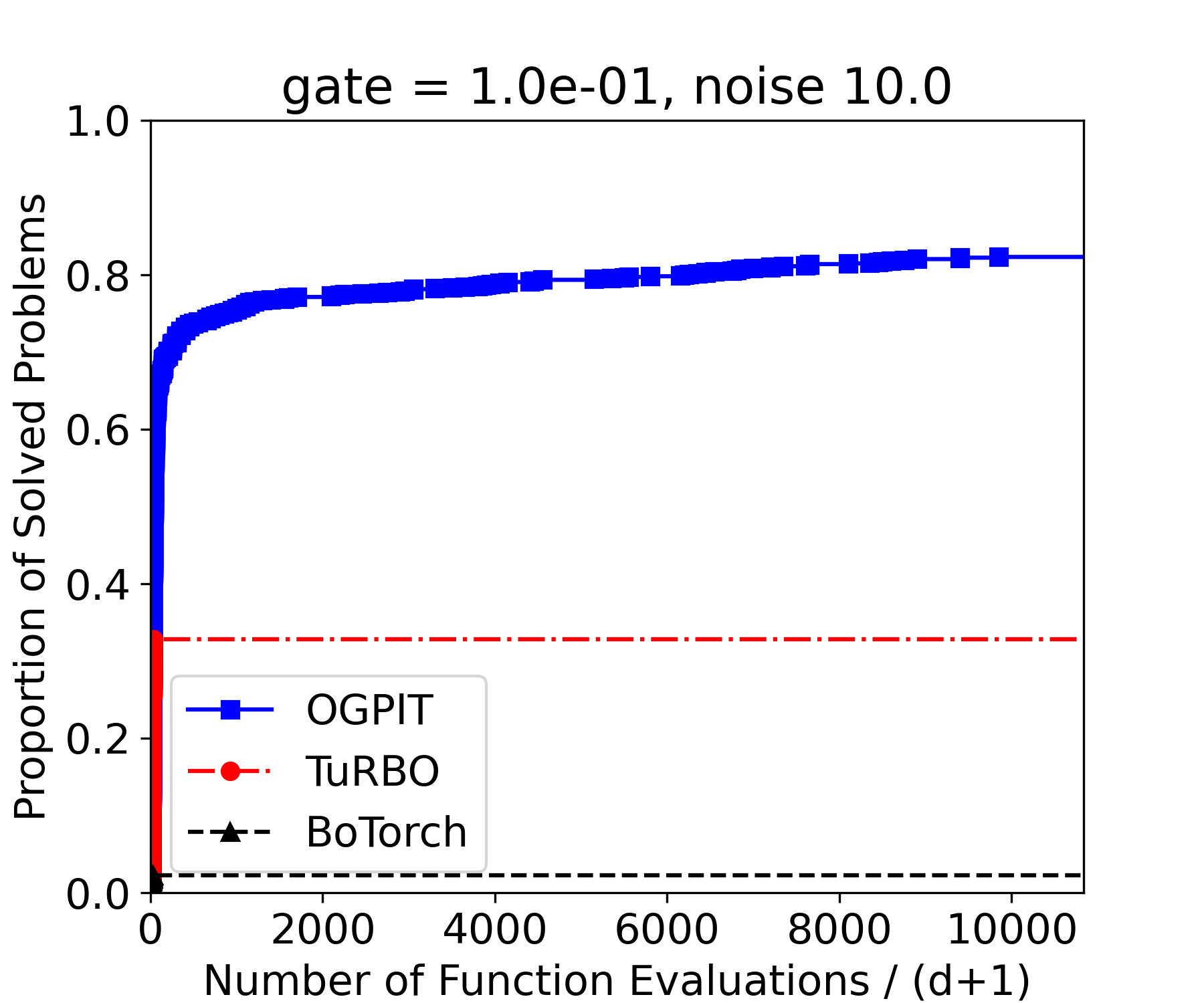}%
\includegraphics[width=0.333\textwidth, trim= 0 0 40 15, clip]{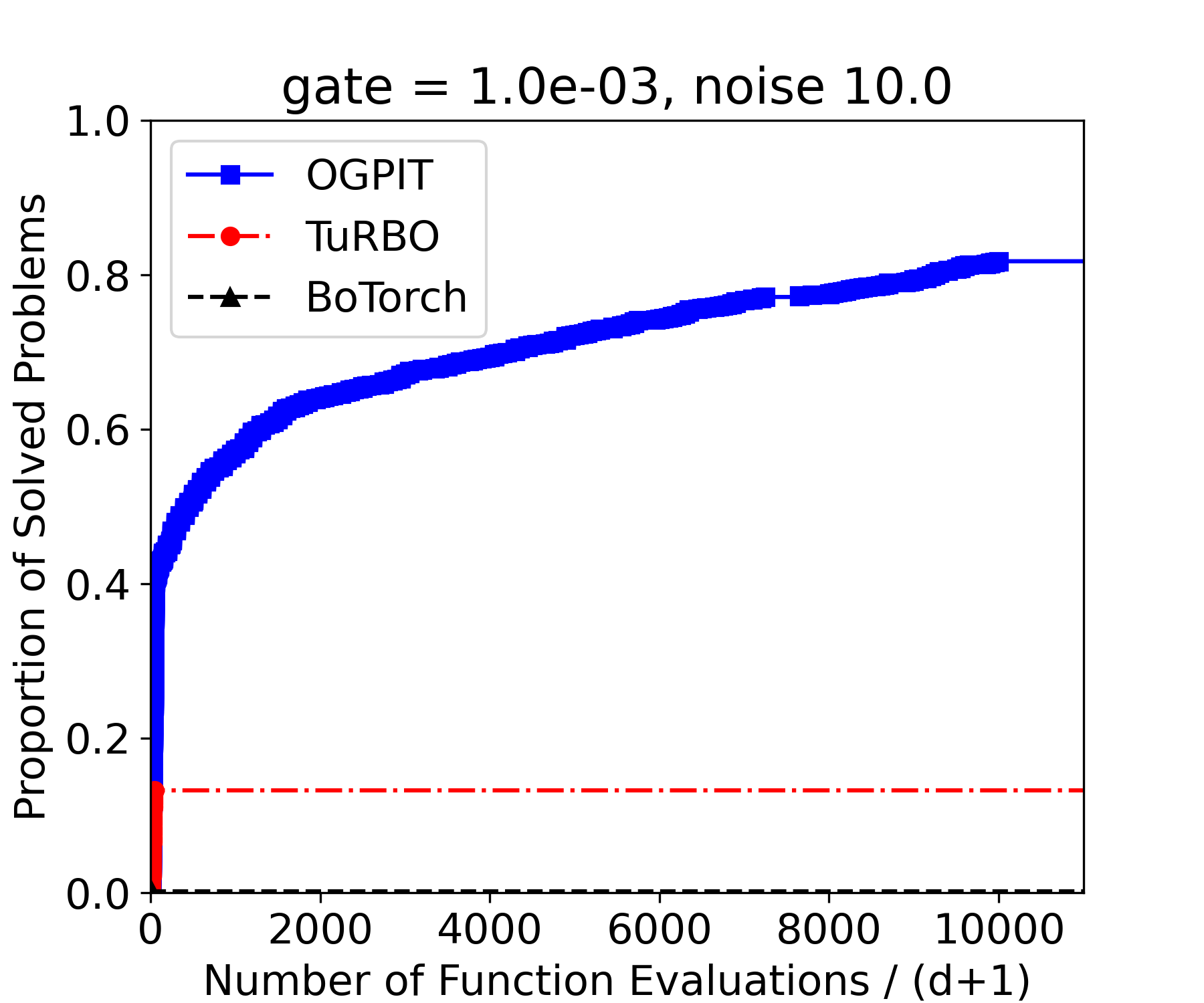}%
\includegraphics[width=0.333\textwidth, trim= 0 0 40 15, clip]{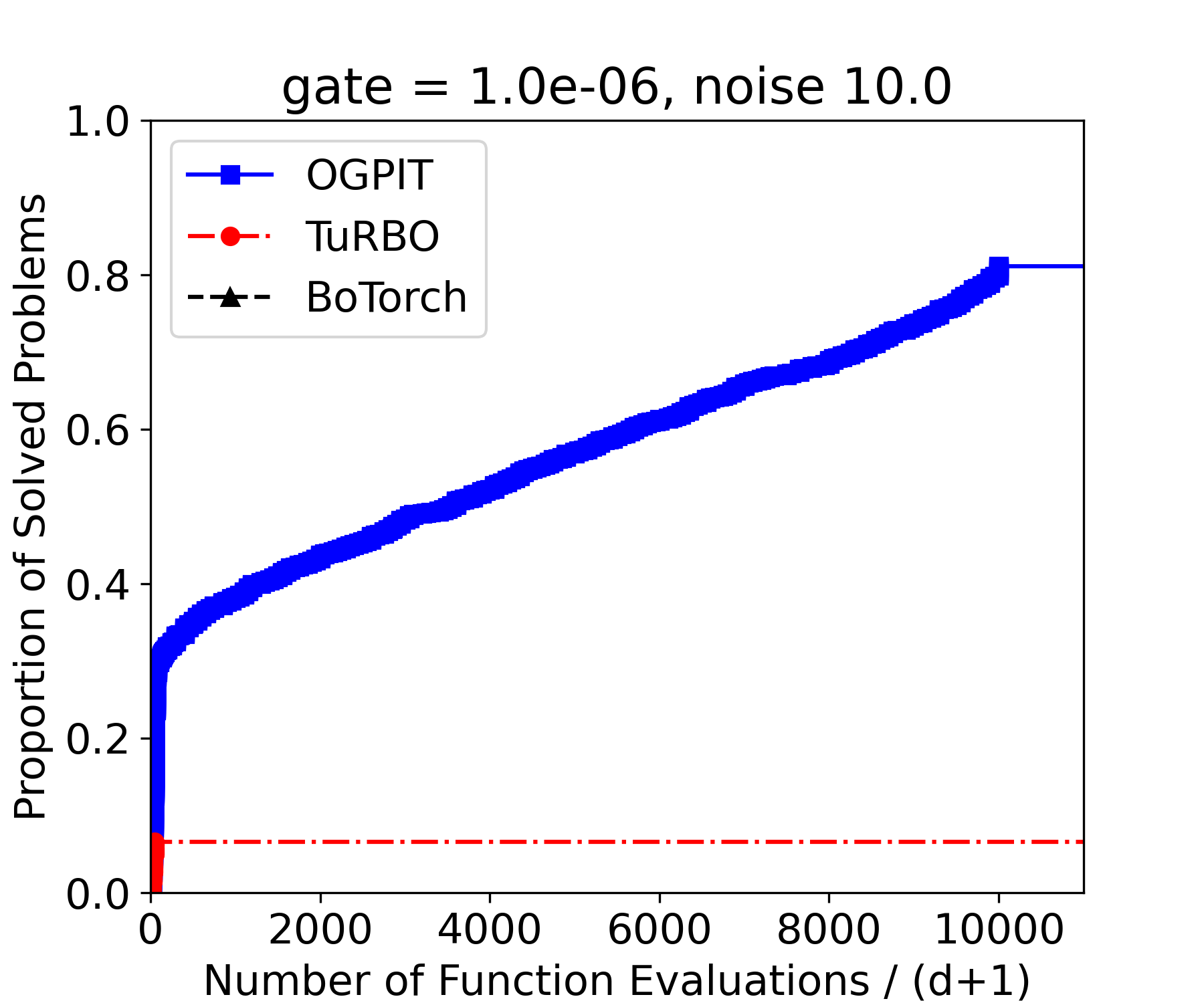}
\caption{Data profiles of the method against TuRBO and BoTorch on Benchmark 2.}
\label{fig:localbenchprofs2}
\end{figure}

Next we provide results for the more complex Benchmark 2 in \Cref
{fig:localbenchprofs2}. Again, OGPIT outperforms TuRBO, with larger gaps with
increasing noise. The results show the effect of the lower signal-to-noise ratio as
the TR radius reduces. Here BoTorch performs poorly for all cases, even the
noiseless one, because of the bad conditioning close to the solution and because
of the large
difference in objective values observed in single problems in this benchmark
set. This is generally difficult for surrogate
models. Local models are more suited in this context. For the first
iterations, most methods perform similarly. Then TuRBO and BoTorch tend to
stagnate while OGPIT keeps improving.

\subsection{Local optimization with setup cost}

For emulating the behavior with our motivating setup cost structure, we ran
OGPIT on Benchmark 1 and present average results over costs. We use a
maximum cost of 1000, with $c_0 =1$ and $c_1=10^{-4},10^{-3},10^{-2}$. 

As the default method, as in the preceding section we use EI as the acquisition
function. We then compare with the two variants of qERCI introduced in
\Cref{sub:joint_selection_of_vecx_n_1}. First this is OGPIT with qERCI v1, based on the
adaptive strategy using $p_a$ ($T_a = 0.2$) and dividing by the setup and replication cost. The
number of replicates is such that the variance reduction at the new design is
 at least 20\%. 

The second method, dedicated for this setup cost is the qERCI v2 acquisition function. There OGPIT is with qERCI optimized over $\vecx_{n+1}$ and $\vecx_{n+2}$. It
entertains adding up to two new points with replicates based on the reduction in
improvement divided by the cost. This emulates looking ahead by a batch
strategy. Then only the design with the largest number of replicates is
actually evaluated.

\begin{figure}
\includegraphics[width=0.5\textwidth, trim= 20 35 20 20, clip]{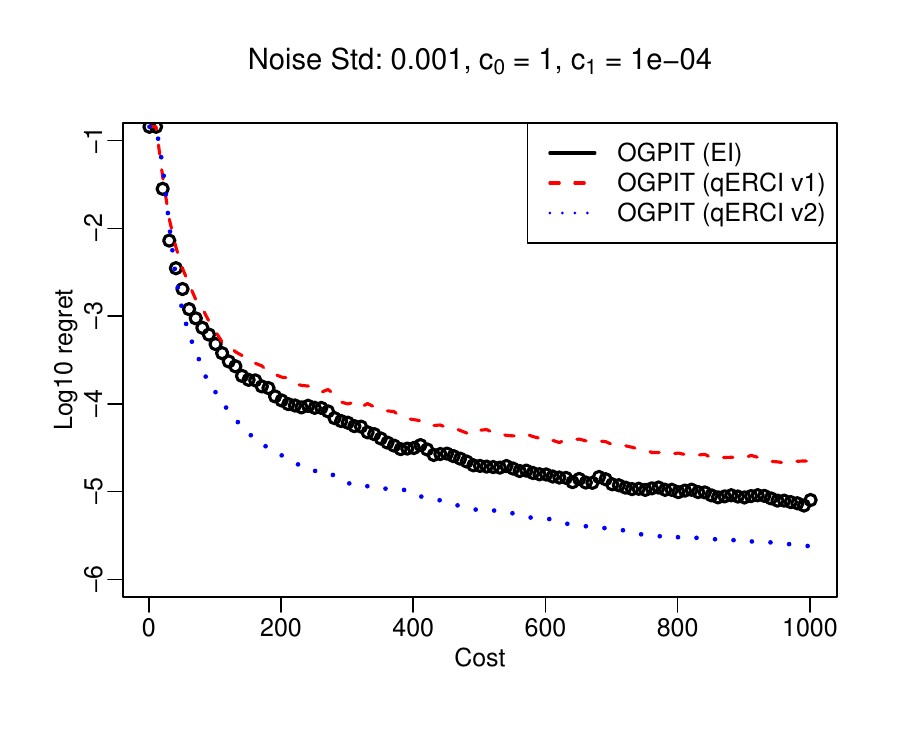}%
\includegraphics[width=0.5\textwidth, trim= 20 35 20 20, clip]{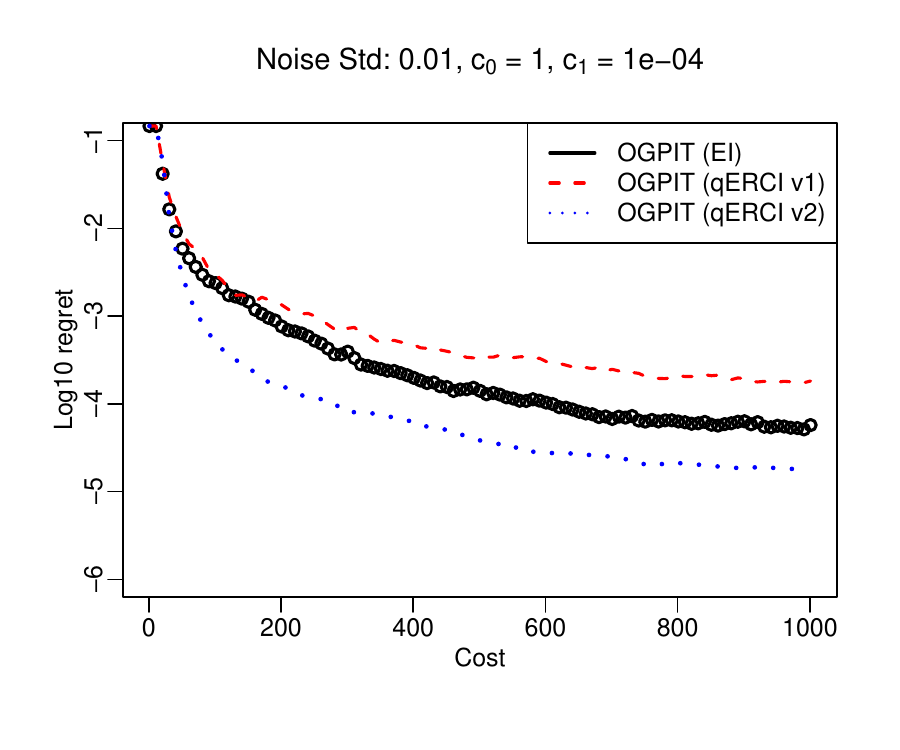}\\
\includegraphics[width=0.5\textwidth, trim= 20 35 20 20, clip]{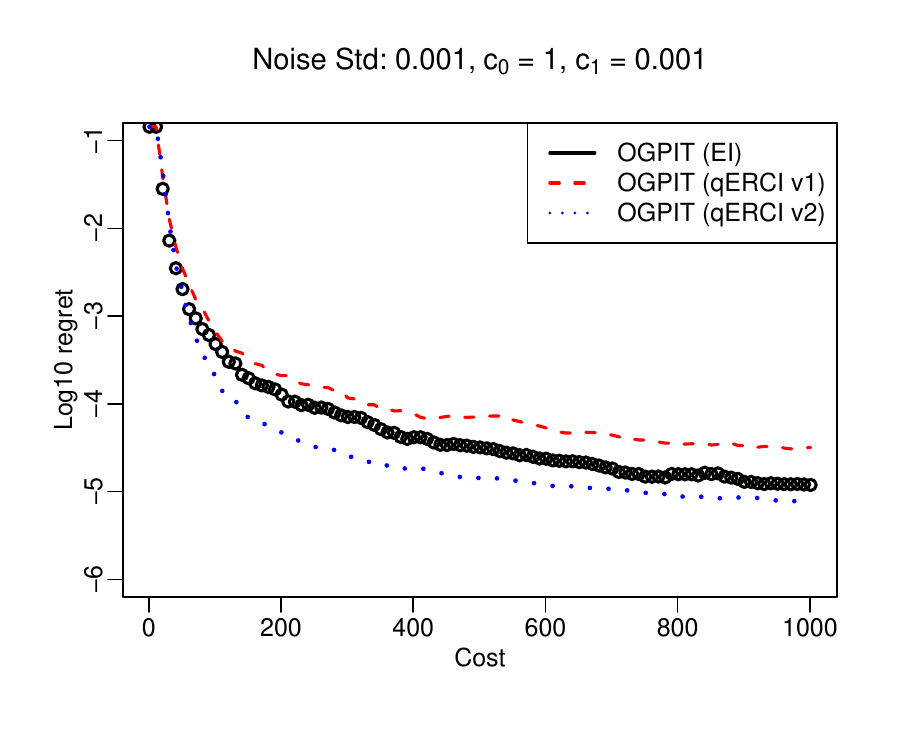}%
\includegraphics[width=0.5\textwidth, trim= 20 35 20 20, clip]{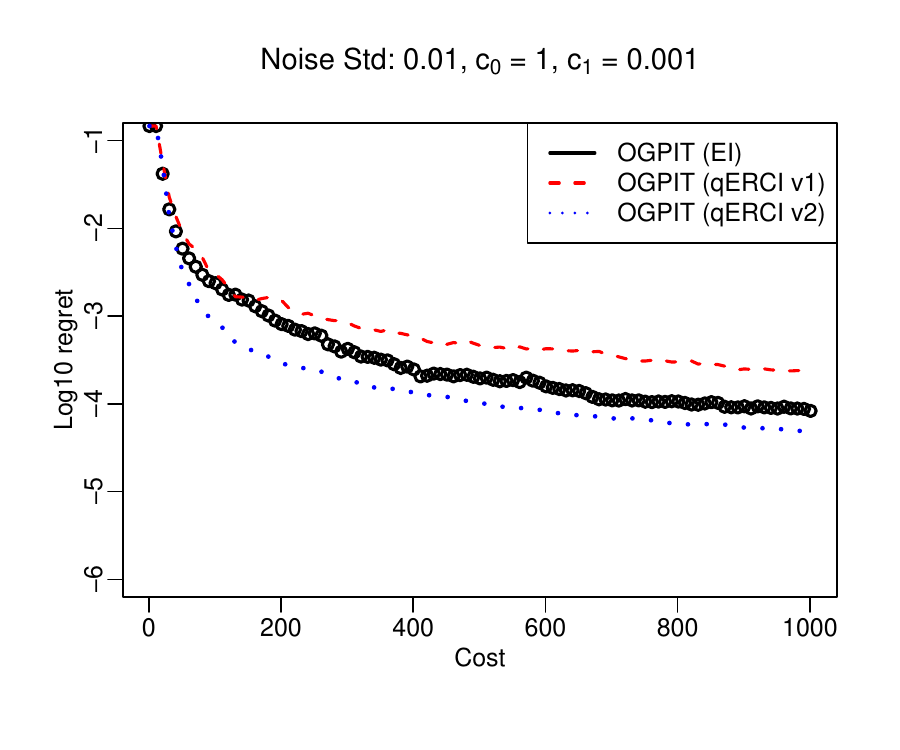}\\
\includegraphics[width=0.5\textwidth, trim= 20 40 20 20, clip]{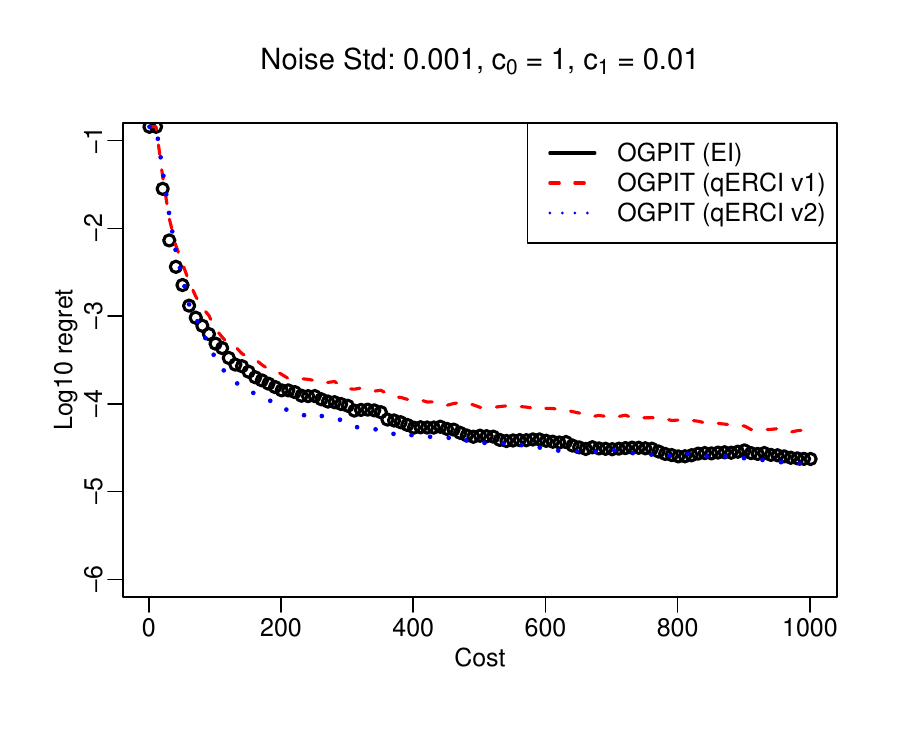}%
\includegraphics[width=0.5\textwidth, trim= 20 40 20 20, clip]{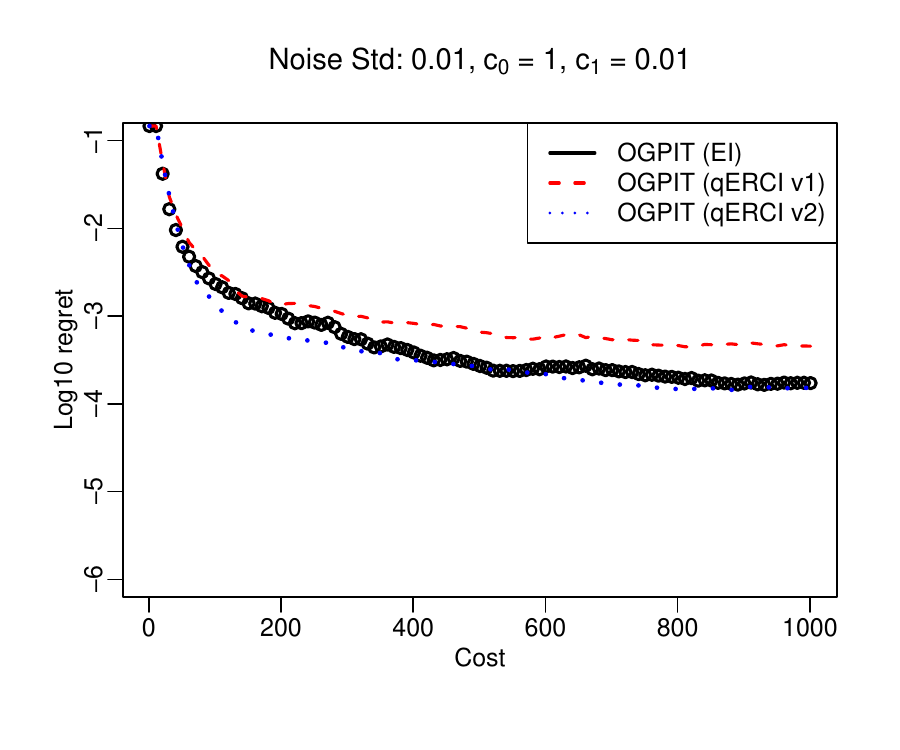}%
\caption{Results on Benchmark 1 with two noise levels and varying setup costs.}
\label{fig:costres}
\end{figure}

The results are provided in \Cref{fig:costres}. For the smallest replication
cost $c_1$, the cost-aware version (qERCIv2) gives the best results with
a large margin. For the medium replication cost scenario, it still
outperforms the others. When the replication cost increases, it becomes
equivalent to the default. The naive strategy of
reducing the variance by 20\% performs worse than the default
acquisition function for all cost setups. Looking ahead further is thus
preferable. The noise variance does not modify the ranking of the methods; it
just translates the curves. Indeed, more evaluations are needed to reach a
given regret level.

\subsection{Motivating quantum test case}

We consider a stochastic optimization problem arising from the quantum
approximate optimization algorithm (QAOA)~\cite{Farhi_2014}.
QAOA is a prominent variational quantum algorithm and a candidate for demonstrating
quantum advantage due to its relatively modest requirements on qubit counts
and circuit depth compared with fully fault-tolerant approaches~\cite{Preskill_2018,Cerezo_2021}. We consider using QAOA to solve
Max-Cut problem instances on a collection of graphs $G=(V,E)$.
Given a QAOA depth $\mathring{p}$, the algorithm requires $2\mathring{p}$ continuous parameters,
typically written as
\[
\theta = (\mathring{\gamma}_1,\ldots,\mathring{\gamma}_p,\mathring{\beta}_1,\ldots,\mathring{\beta}_p),
\]
which define alternating applications of a cost Hamiltonian derived from the
Max-Cut objective and a mixing Hamiltonian. The resulting objective value
corresponds to the expected Max-Cut energy obtained by measuring the final
quantum state.
Finding high-quality parameters is widely
recognized as a key bottleneck in the practical performance of QAOA~\cite{Zhou_2020,Brandao_2018,Hao_2025}.

The graphs we consider include the Chv\'atal graph and Erd\H{o}s--R\'enyi graphs
of increasing size, and the default experiments use $\mathring{p}=1$ or $\mathring{p}=2$ (so $d=2\mathring{p}=2$ or $d=4$ parameters) ranging from $0$ to $\pi/2$, rescaled to $[0,1]$.
We benchmark our approach on this stochastic objective, where the expectation
is estimated using a finite number of measurement shots. As with most
near-term quantum algorithms, QAOA is fundamentally stochastic in practice due
to measurement noise and finite sampling. For the modest graph sizes considered
here, however, the exact expectation value can be computed deterministically
via classical simulation; we use these true values solely for benchmarking and
diagnostic purposes. Each noisy function call returns multiple stochastic
realizations, and the scalar test function aggregates these by taking the
sample mean. This yields an unbiased estimator of the expected objective.

\begin{figure}
\includegraphics[width=0.5\textwidth, trim= 25 40 20 20, clip]{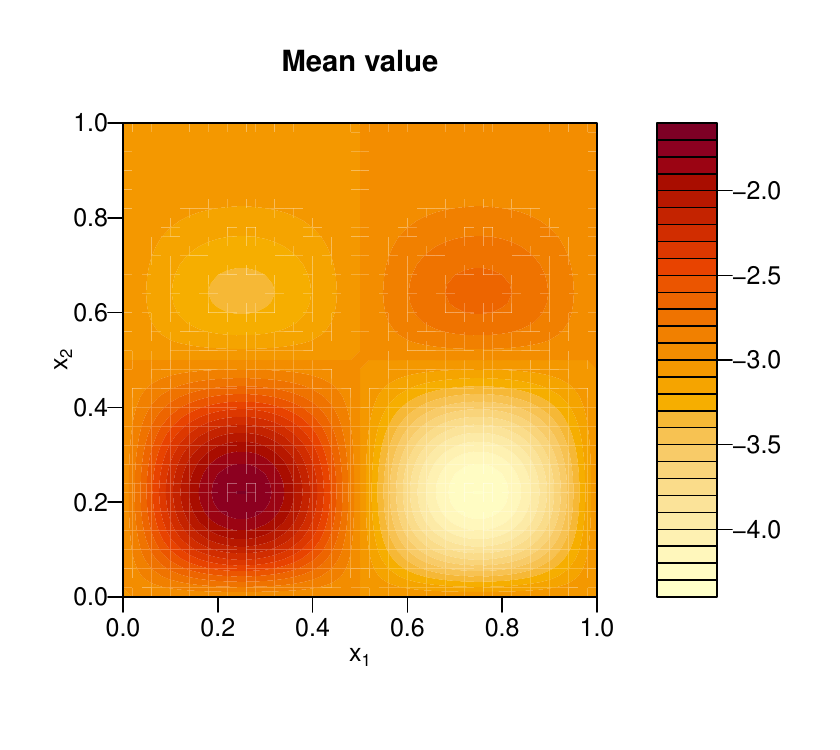}%
\includegraphics[width=0.5\textwidth, trim= 25 40 20 20, clip]{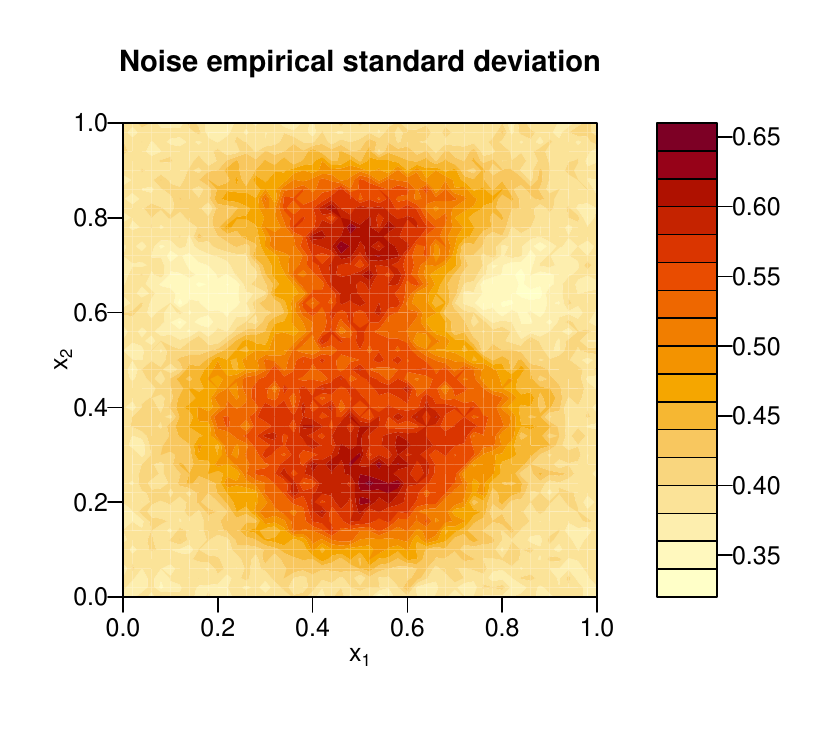}
\caption{Illustration of the quantum test case with ground truth (left) and estimated noise standard deviation estimated over 5000 replicates.}
\label{fig:quantum2d}
\end{figure}

This quantum test case exhibits several structural features that can be
exploited by a trust-region Bayesian optimization framework. First,
the noise variance is heteroscedastic (see
\Cref{fig:quantum2d}), varying significantly across the parameter
domain due to the underlying quantum state amplitudes and measurement
probabilities.
Second, the
problem is representative of settings with substantial setup costs: on
average, preparing a circuit instance (i.e., configuring the parameterized
quantum circuit and compiling it for execution) is approximately
$100$--$1000$ times more expensive than a single measurement shot.
This cost structure arises because circuit preparation involves classical
preprocessing, parameter loading, and potential compilation or transpilation
steps, whereas individual shots correspond only to repeated measurements of
the same prepared circuit. We henceforth use $c_0=1$ and $c_1=0.001$ in our experiments, with a maximum cost of $250$ for the $2d$ problem and $500$ for the $4d$ problem. The initial design of experiments are of size 10 and 200, respectively, to ensure that the global optima are identified over the 5 repetitions of the experiment.

\begin{figure}
\includegraphics[width=0.5\textwidth, trim= 20 35 20 20, clip]{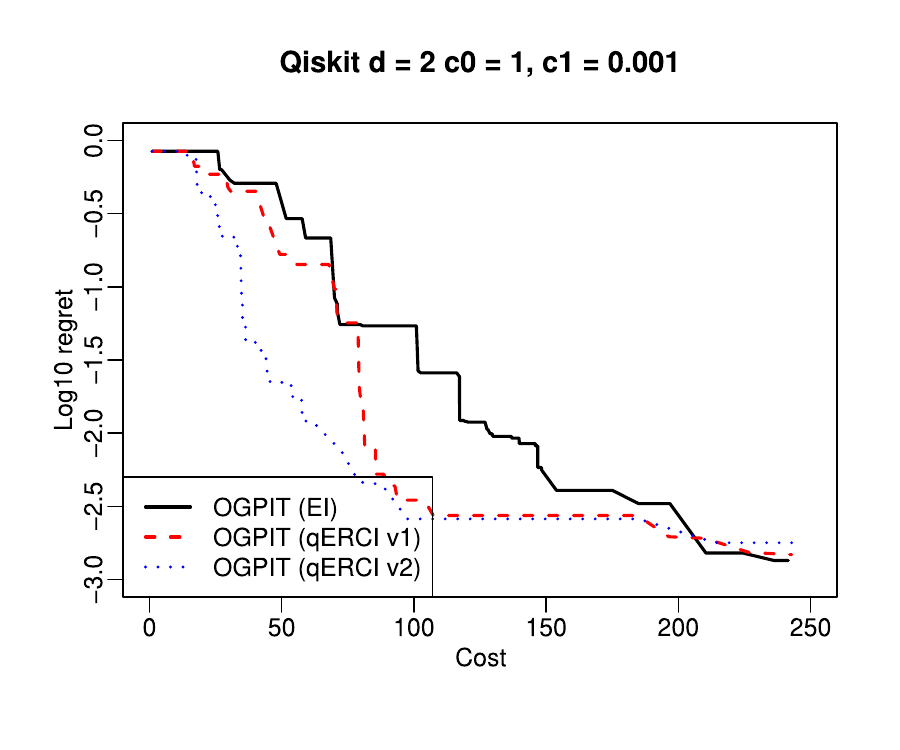}%
\includegraphics[width=0.5\textwidth, trim= 20 35 20 20, clip]{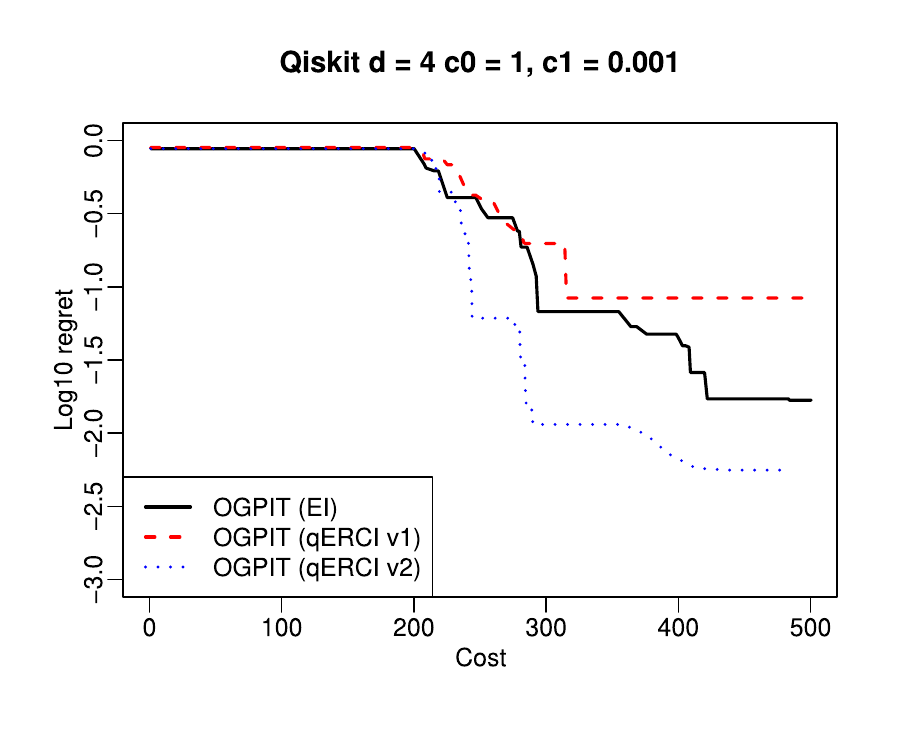}
\caption{Results on the quantum test case, in dimension 2 (left) and dimension 4 (right).}
\label{fig:costquant}
\end{figure}

The results are presented in \Cref{fig:costquant}, where the second version of qERCI again performs better at reducing the regret efficiently with respect to the cost. We note that the final regret is well below the noise variance.
In the 2-dimension case, while similar solution qualities are found, using knowledge of the cost structure uses approximately 46\% less computational costs than the default approach to reach a regret below $10^{-2}$. This is achieved using relatively more replications than unique evaluations: OGPIT with EI reaches this regret with an average of 10,832 evaluations over 119 unique designs while qERCI reaches it with around 13,145 evaluations over 56 unique designs.

\section{Conclusion and perspectives}

Standard Bayesian optimization provides efficient global optimization
capabilities but is generally not well suited for precise convergence. For this purpose,
trust-region-based Bayesian optimization strategies are more appropriate, and they 
perform well in benchmarks on noise-free problems; see, for
instance,~\cite{santoni2024comparison}. 

As we detail, the difficulty of current methods is exacerbated when noise is
present. We therefore propose adaptations for noisy problems, including cases with large
noise variance. In particular, selecting both the design to evaluate and the
replication level allows us to control the accuracy of the model while maintaining
computational scalability as the budget increases. Additional modifications
to the trust-region mechanism are proposed, relying on properties of Gaussian
processes. Targeting applications with significant setup costs but where
repeating evaluations is comparatively inexpensive, we develop a dedicated
acquisition function. We demonstrate the effectiveness of these
modifications, showing good results even on noiseless test problems while
yielding significant improvements when noise is present or when setup costs are substantial.

Although we focus on local optimization performance, possible future work
includes studying the interplay with the global search phase of TR-based BO
methods, since several options exist (e.g.,~\cite{regis2016trust,diouane2022trego,eriksson2019scalable,Mathesen2017,mathesen2021stochastic}). 
Results with the polynomial-based TR method SNOWPAC showed a strong initial performance. Integrating such a polynomial trend as the mean in GP modeling is then a promising direction to explore. In addition, it may help to obtain convergence guaranties, for instance relying on the ASTRO-DF framework \cite{Shashaani2018}.
Detecting when the
function is unimodal in the trust region, as done in~\cite{mcleod2018optimization} but in the presence of
noise, is a promising direction. Further work is also needed to study
the effect of the input dimension when noise is present; high dimensionality
is already a challenge even without noise and typically requires structural assumptions; see, for example,~\cite{Binois2022}.

\section*{Acknowledgments}
This work was supported by the U.S.~Department of Energy, Office of Science,
Office of Advanced Scientific Computing Research, Scientific Discovery through
Advanced Computing (SciDAC) Program through the FASTMath Institute under
Contract No.~DE-AC02-06CH11357.
The authors are grateful to the OPAL infrastructure from Université Côte d'Azur for providing resources and support.

\bibliographystyle{plainnat}      %
\bibliography{ogpit}

\vfill
\framebox{\parbox{.90\linewidth}{\scriptsize The submitted manuscript has been
created by UChicago Argonne, LLC, Operator of Argonne National Laboratory
(``Argonne''). Argonne, a U.S.\ Department of Energy Office of Science
laboratory, is operated under Contract No.\ DE-AC02-06CH11357. The U.S.\
Government retains for itself, and others acting on its behalf, a paid-up
nonexclusive, irrevocable worldwide license in said article to reproduce,
prepare derivative works, distribute copies to the public, and perform publicly
and display publicly, by or on behalf of the Government. The Department of
Energy will provide public access to these results of federally sponsored
research in accordance with the DOE Public Access Plan
\url{http://energy.gov/downloads/doe-public-access-plan}.}}

\appendix

\section{Sensitivity analysis of OGPIT's parameters}
\label{ap:A}
We show the effect of modifying the parameters of the method compared to the given default values on the problems of Benchmark 1.

\subsection{Effect of the trust region increase/decrease coefficient}

\begin{figure}
\includegraphics[width=0.333\textwidth, trim= 0 0 40 15, clip]{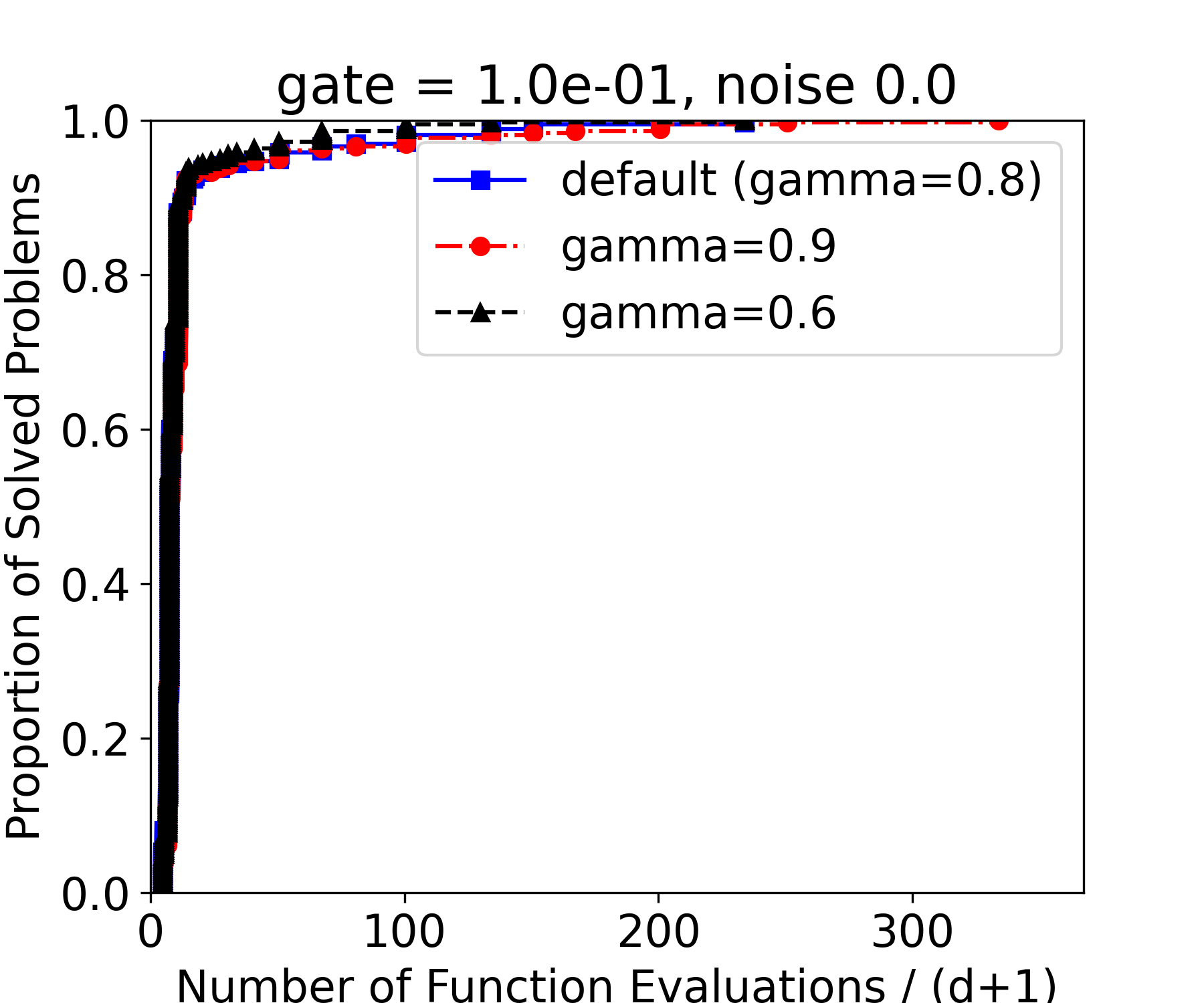}%
\includegraphics[width=0.333\textwidth, trim= 0 0 40 15, clip]{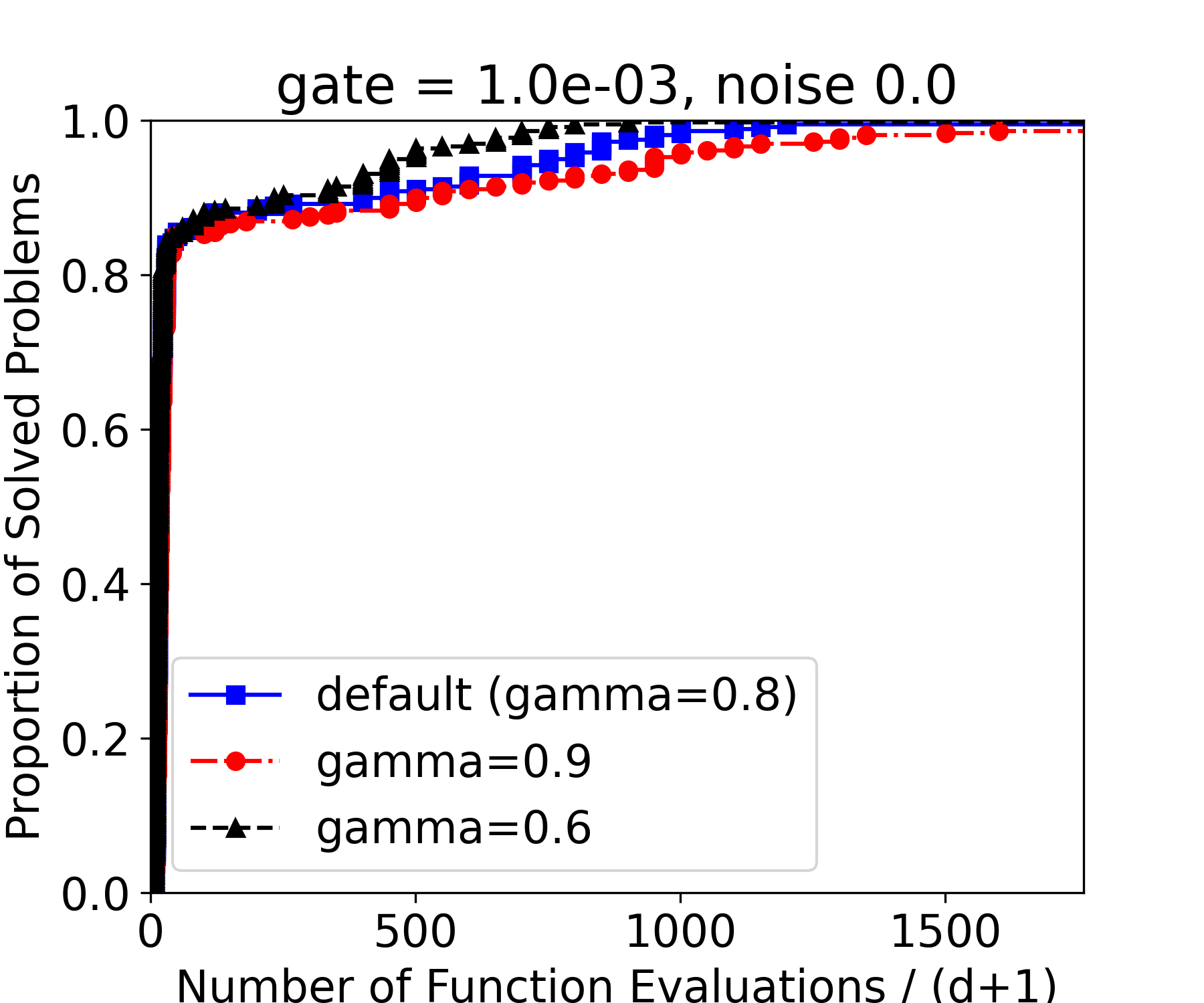}%
\includegraphics[width=0.333\textwidth, trim= 0 0 40 15, clip]{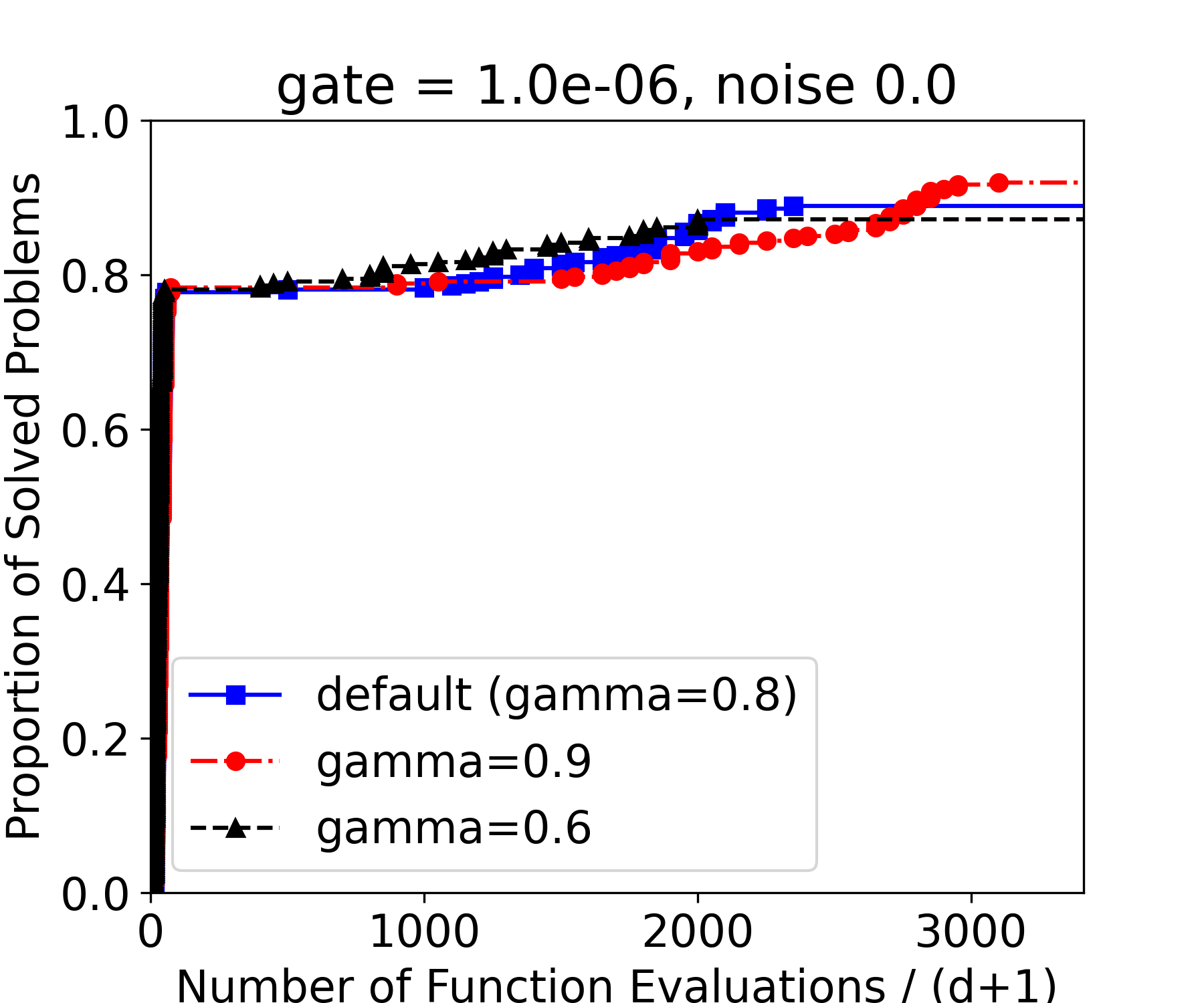}\\
\includegraphics[width=0.333\textwidth, trim= 0 0 40 15, clip]{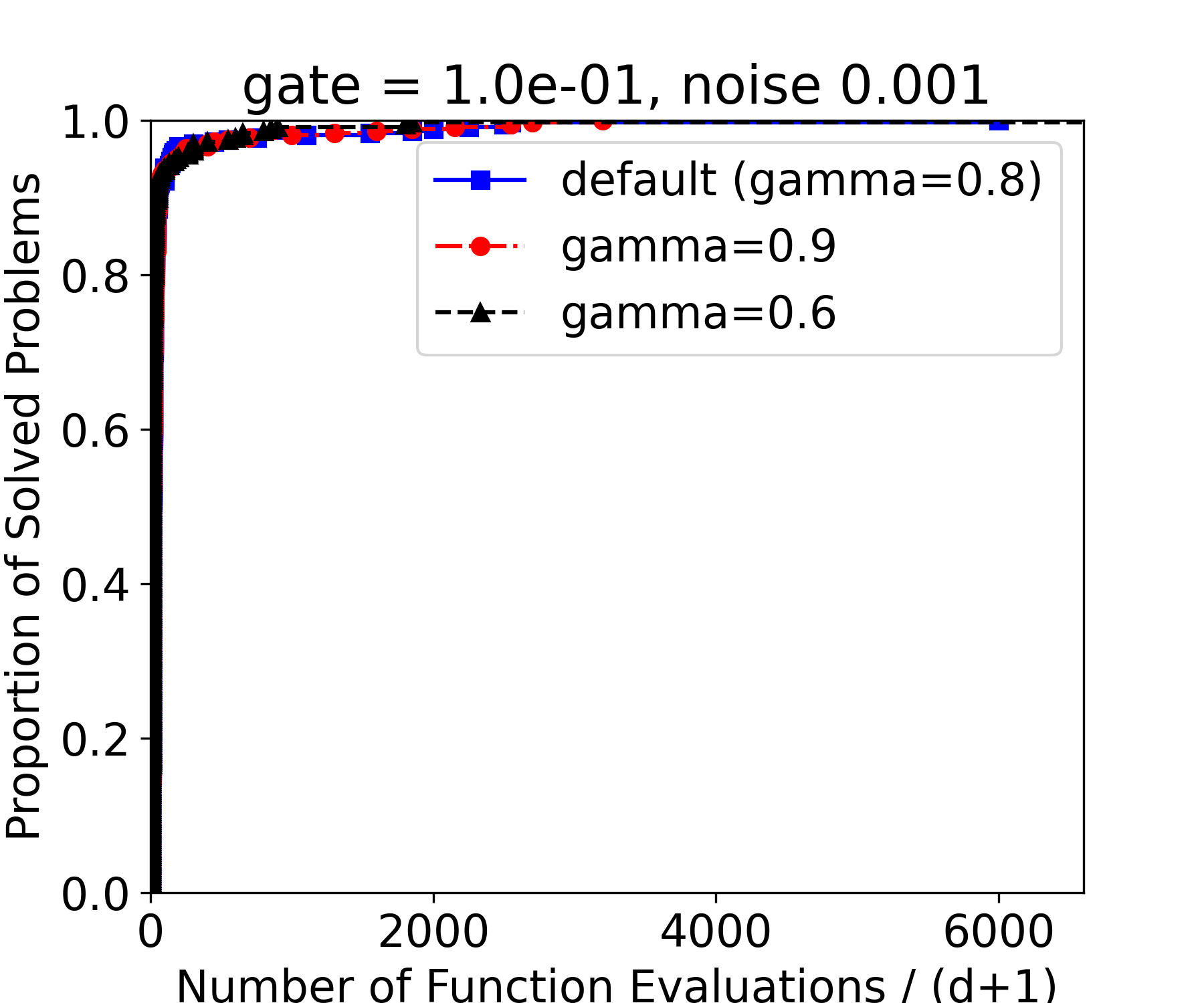}%
\includegraphics[width=0.333\textwidth, trim= 0 0 40 15, clip]{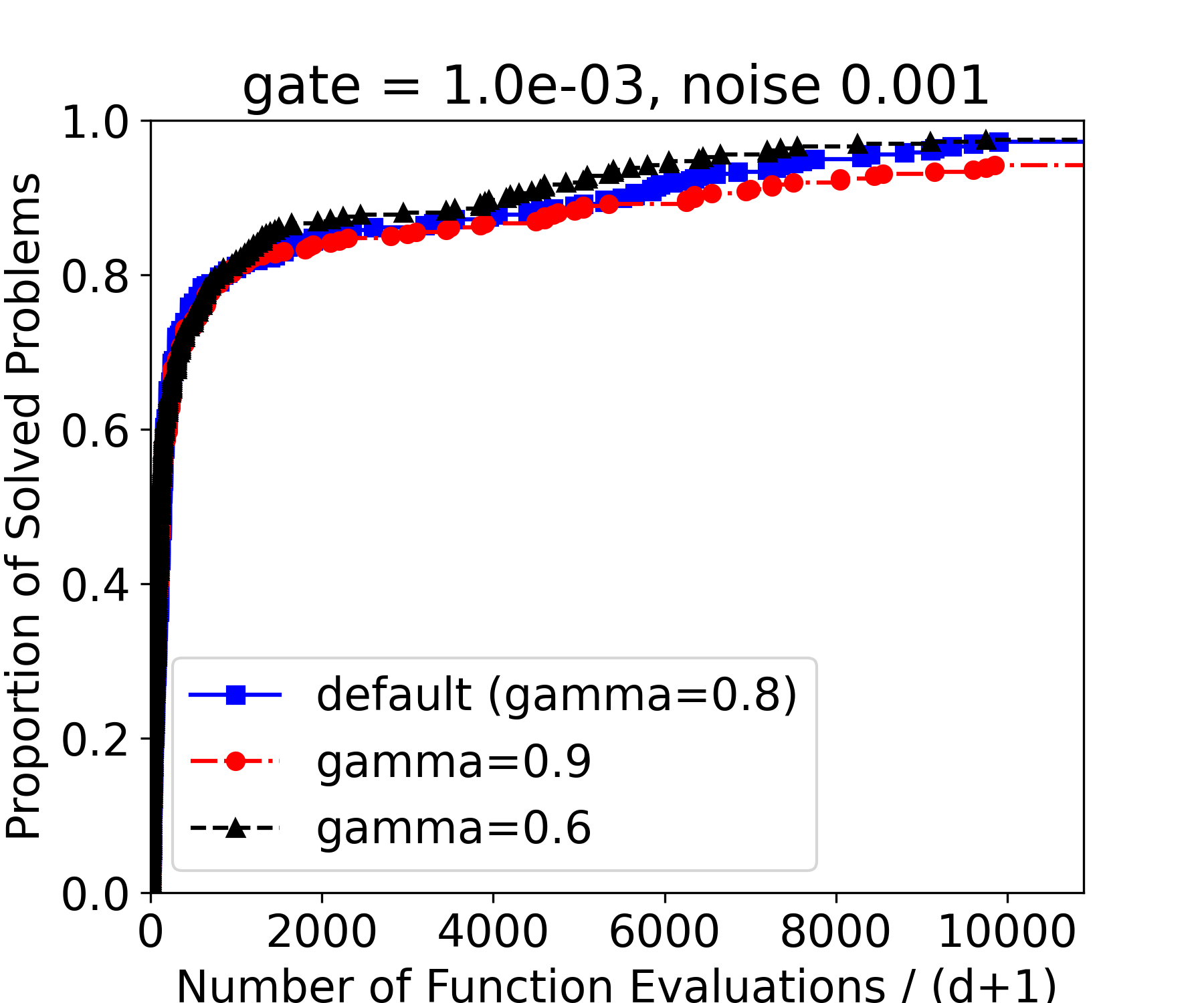}%
\includegraphics[width=0.333\textwidth, trim= 0 0 40 15, clip]{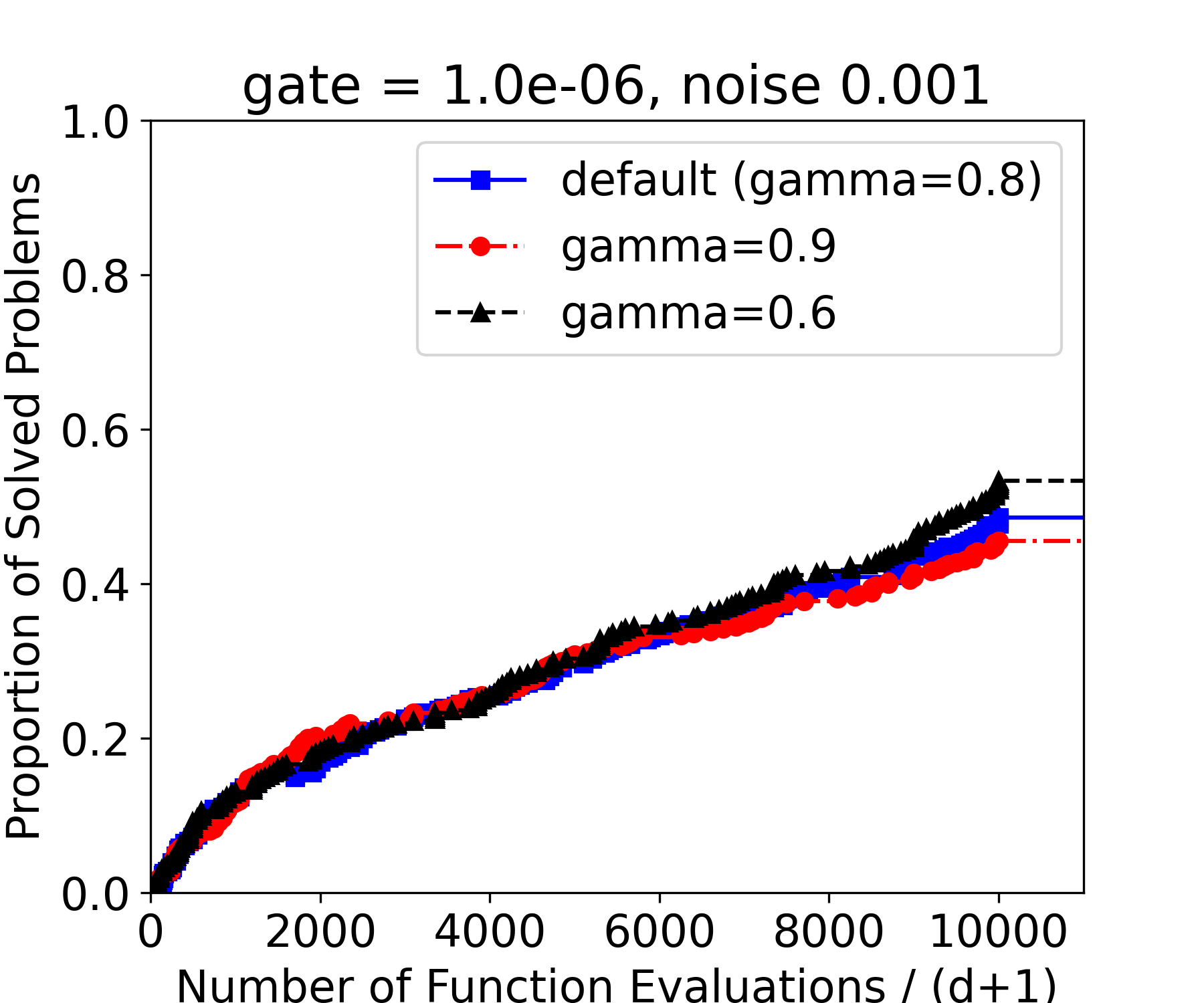}\\
\includegraphics[width=0.333\textwidth, trim= 0 0 40 15, clip]{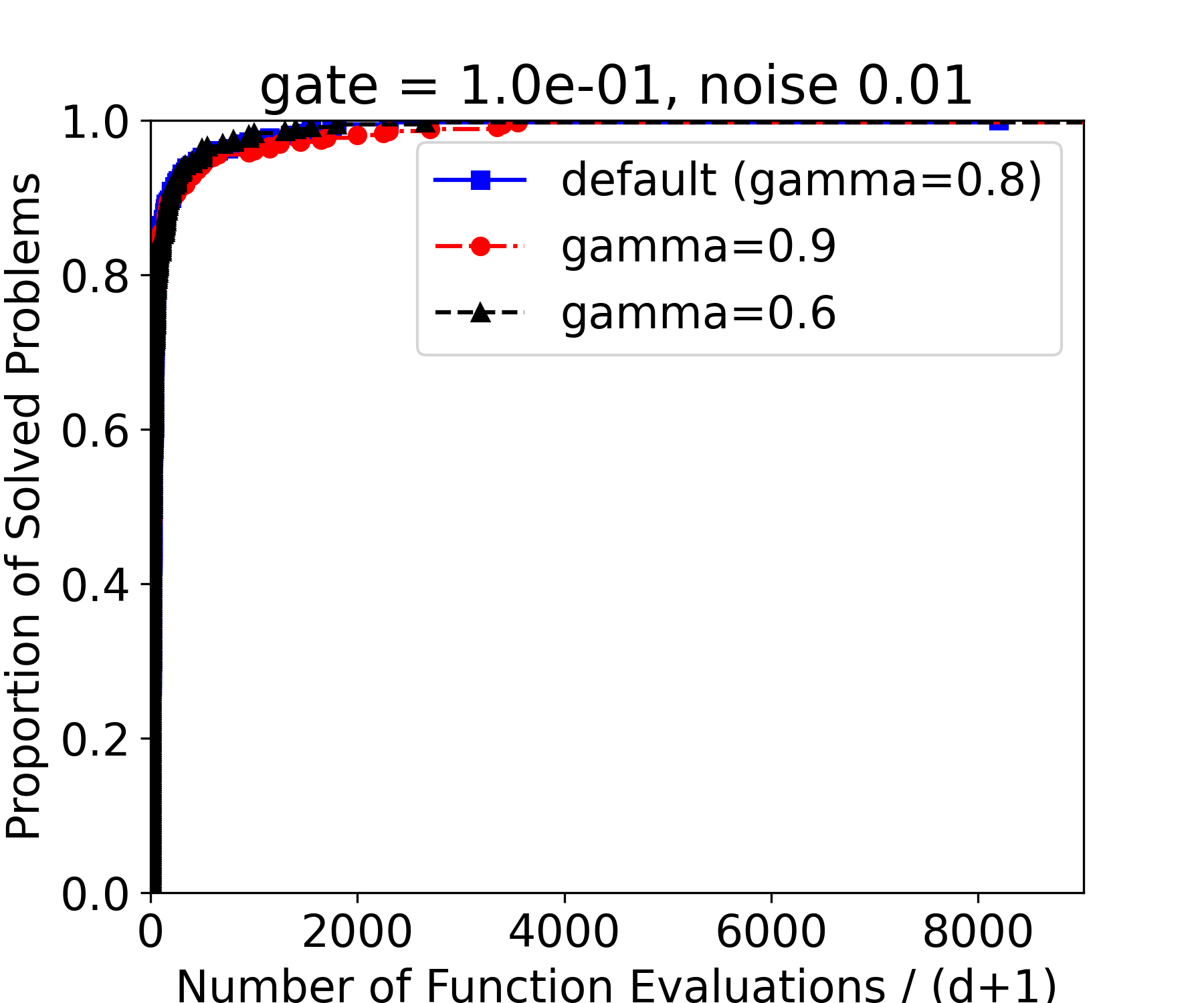}%
\includegraphics[width=0.333\textwidth, trim= 0 0 40 15, clip]{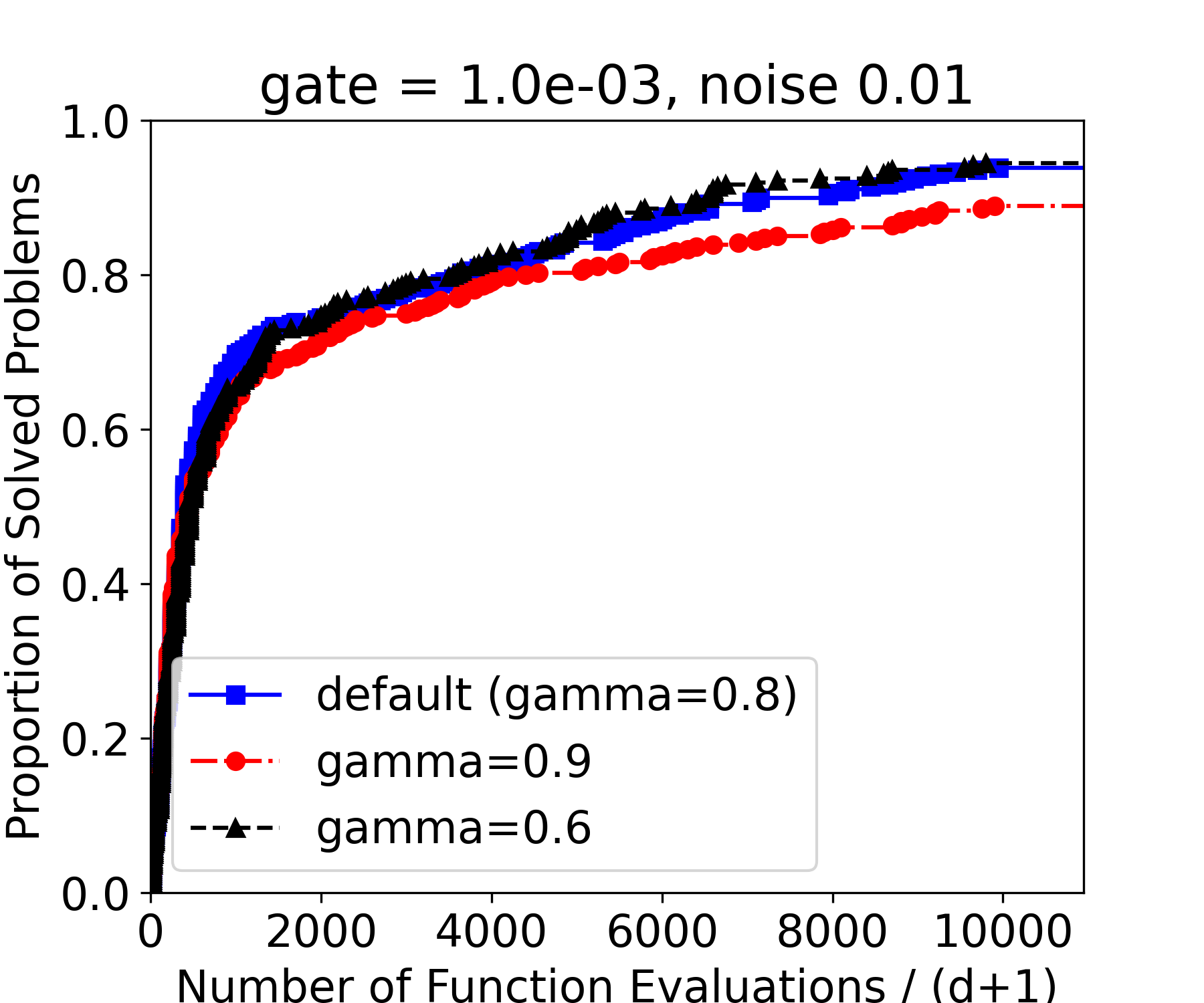}%
\includegraphics[width=0.333\textwidth, trim= 0 0 40 15, clip]{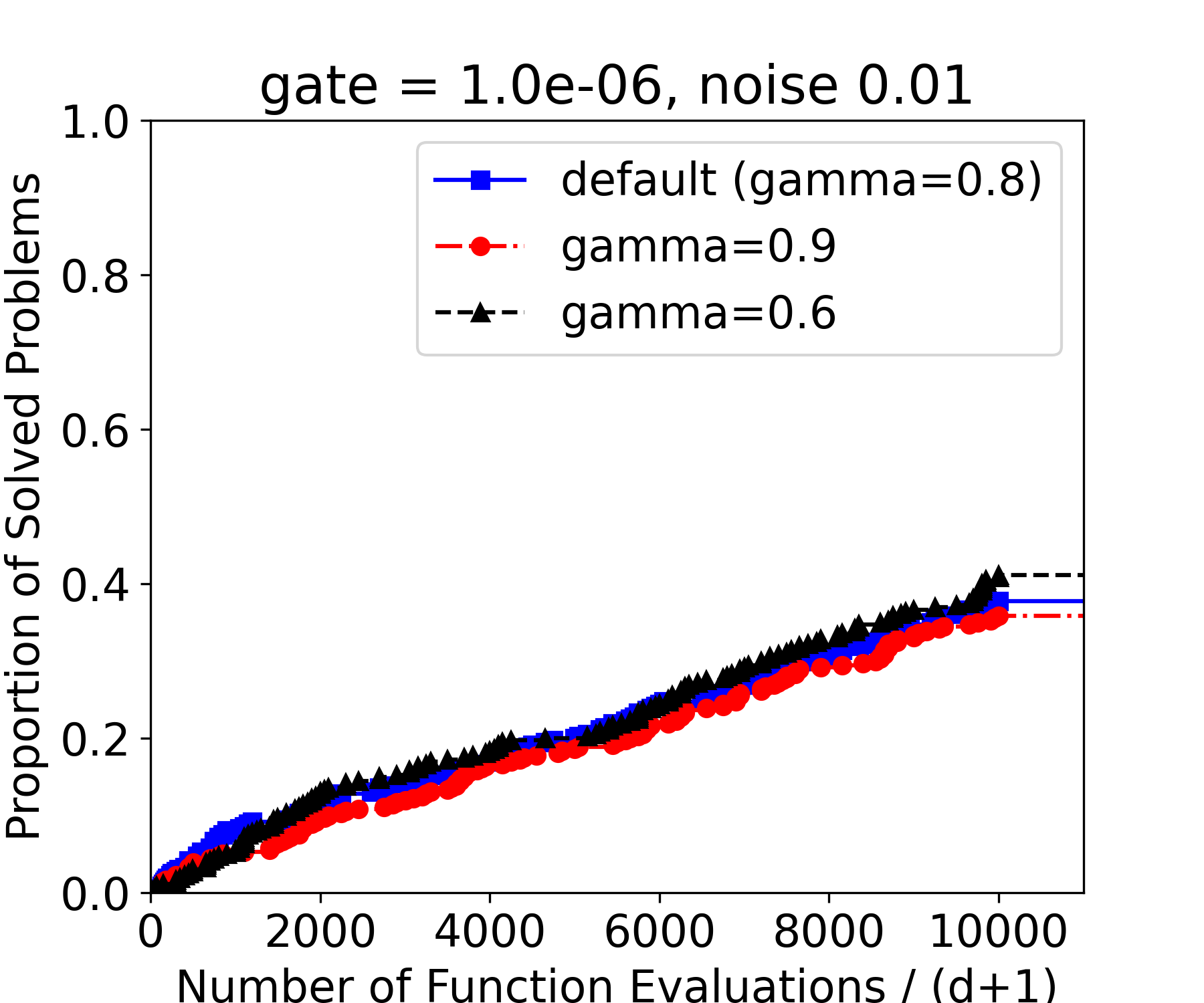}\\
\includegraphics[width=0.333\textwidth, trim= 0 0 40 15, clip]{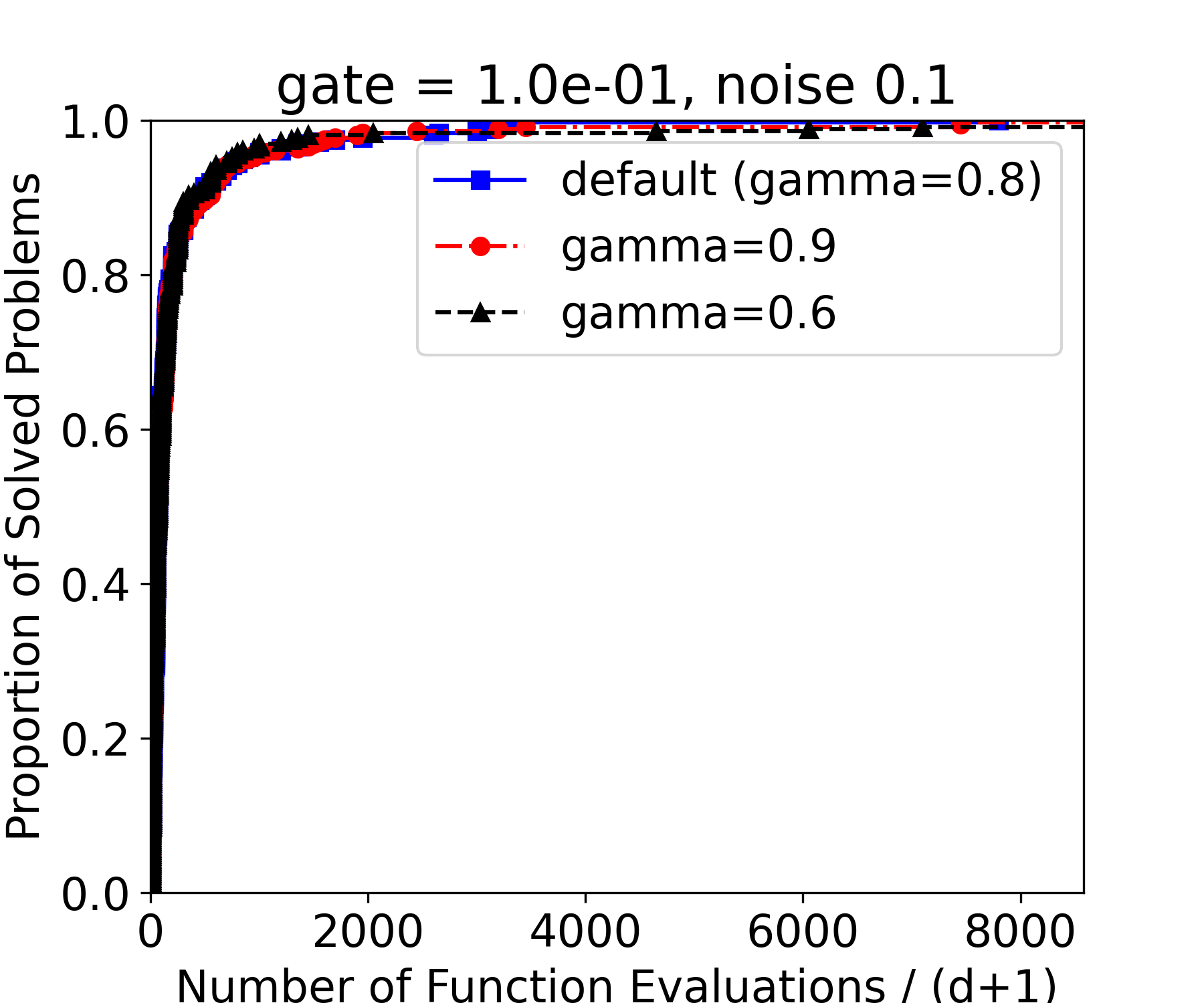}%
\includegraphics[width=0.333\textwidth, trim= 0 0 40 15, clip]{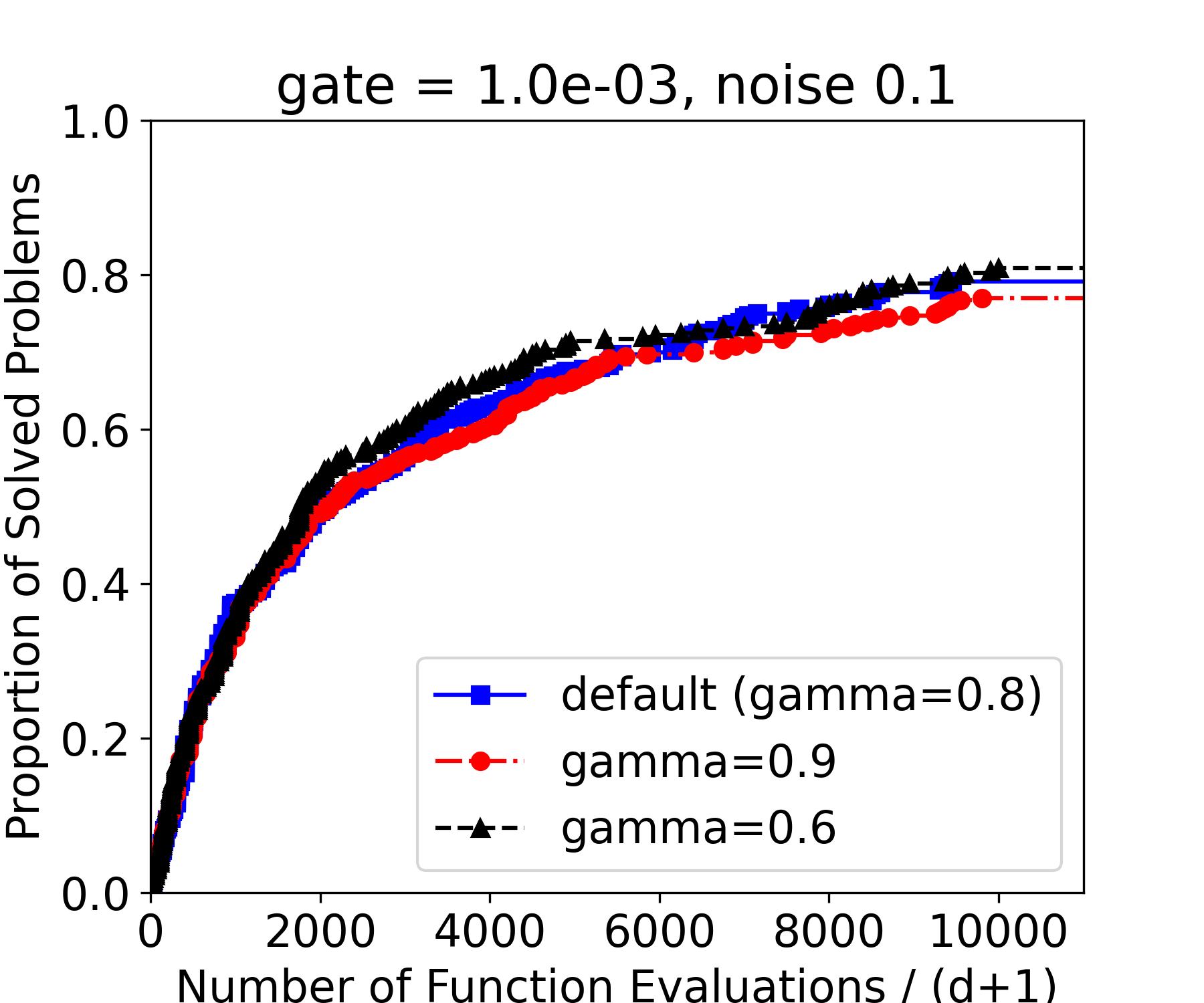}%
\includegraphics[width=0.333\textwidth, trim= 0 0 40 15, clip]{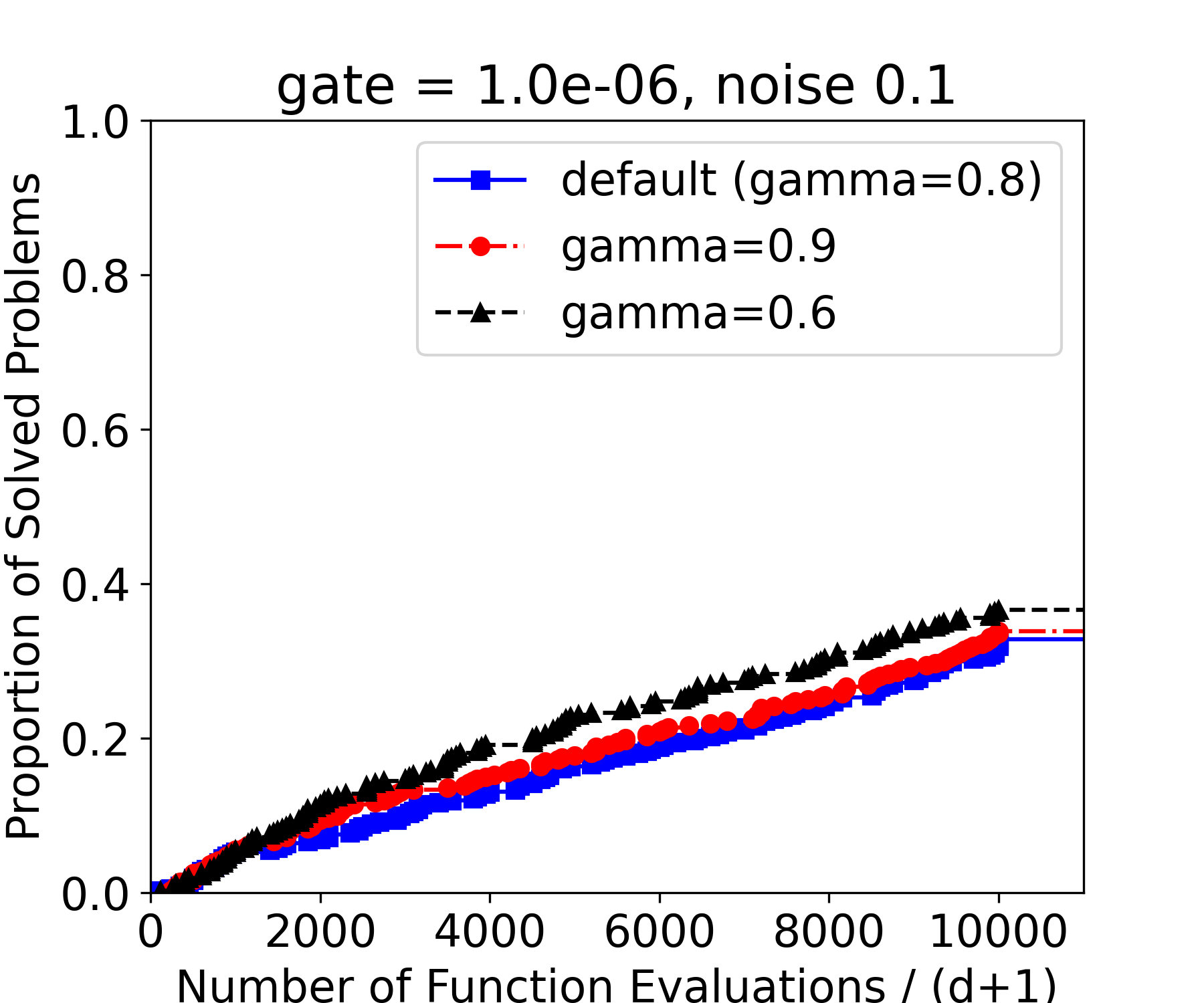}
\caption{Data profiles for varying $\gamma_{dec}$, over 40 repetitions.}
\label{fig:localbenchprofsgamma}
\end{figure}

The coefficient for decreasing the TR radius, $\gamma_{dec}$ is studied
in \Cref{fig:localbenchprofsgamma}. There is no clear best value, while the
smallest coefficient ($\gamma_{dec} = 0.6$) is the best only once, for the
deterministic version, for the smallest gate. A larger value ($\gamma_{dec} =
0.9$) is sometimes better but the difference with the default is slight.

\subsection{Parameter for the IMSE ratio}

\begin{figure}
\includegraphics[width=0.333\textwidth, trim= 0 0 40 15, clip]{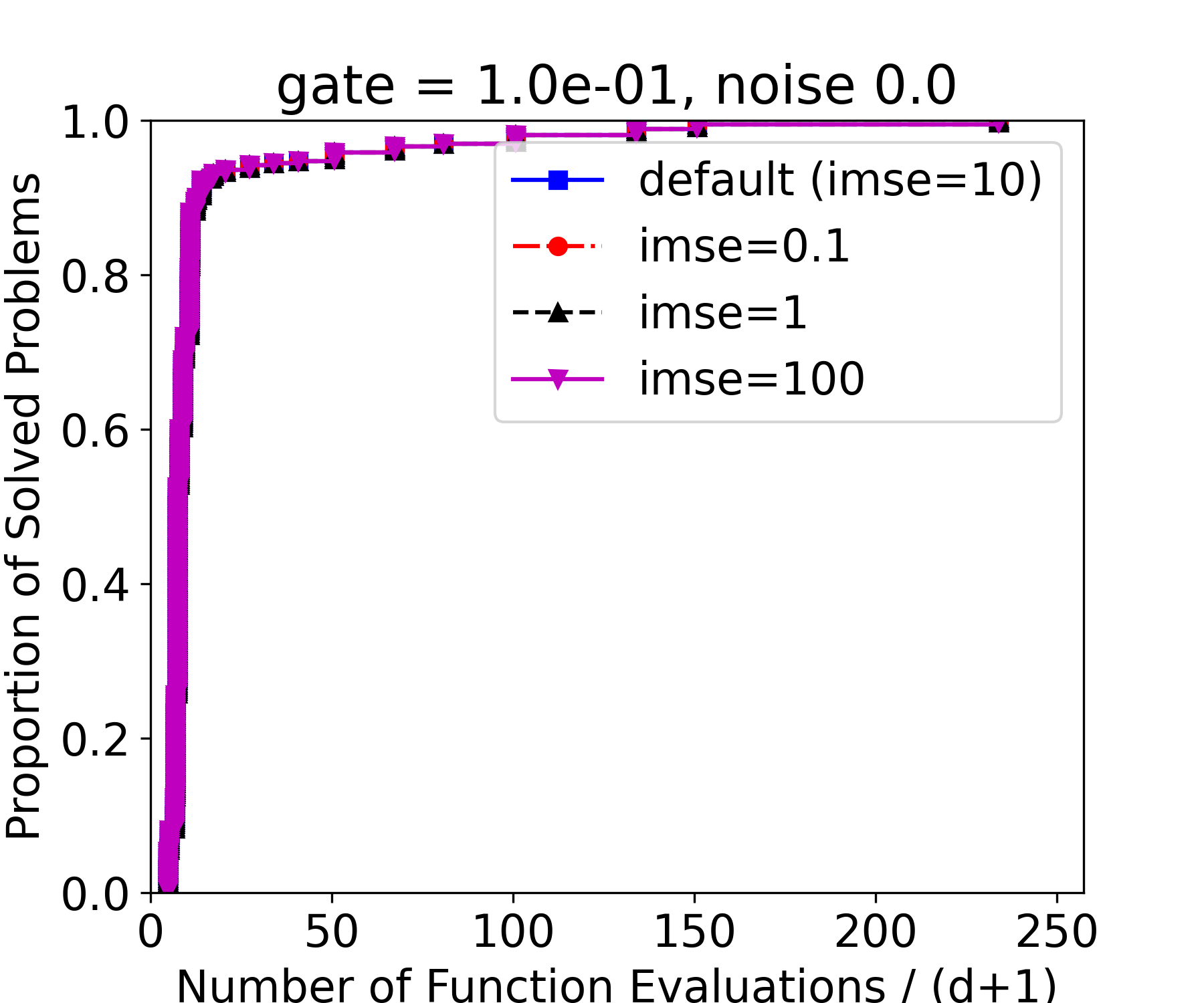}%
\includegraphics[width=0.333\textwidth, trim= 0 0 40 15, clip]{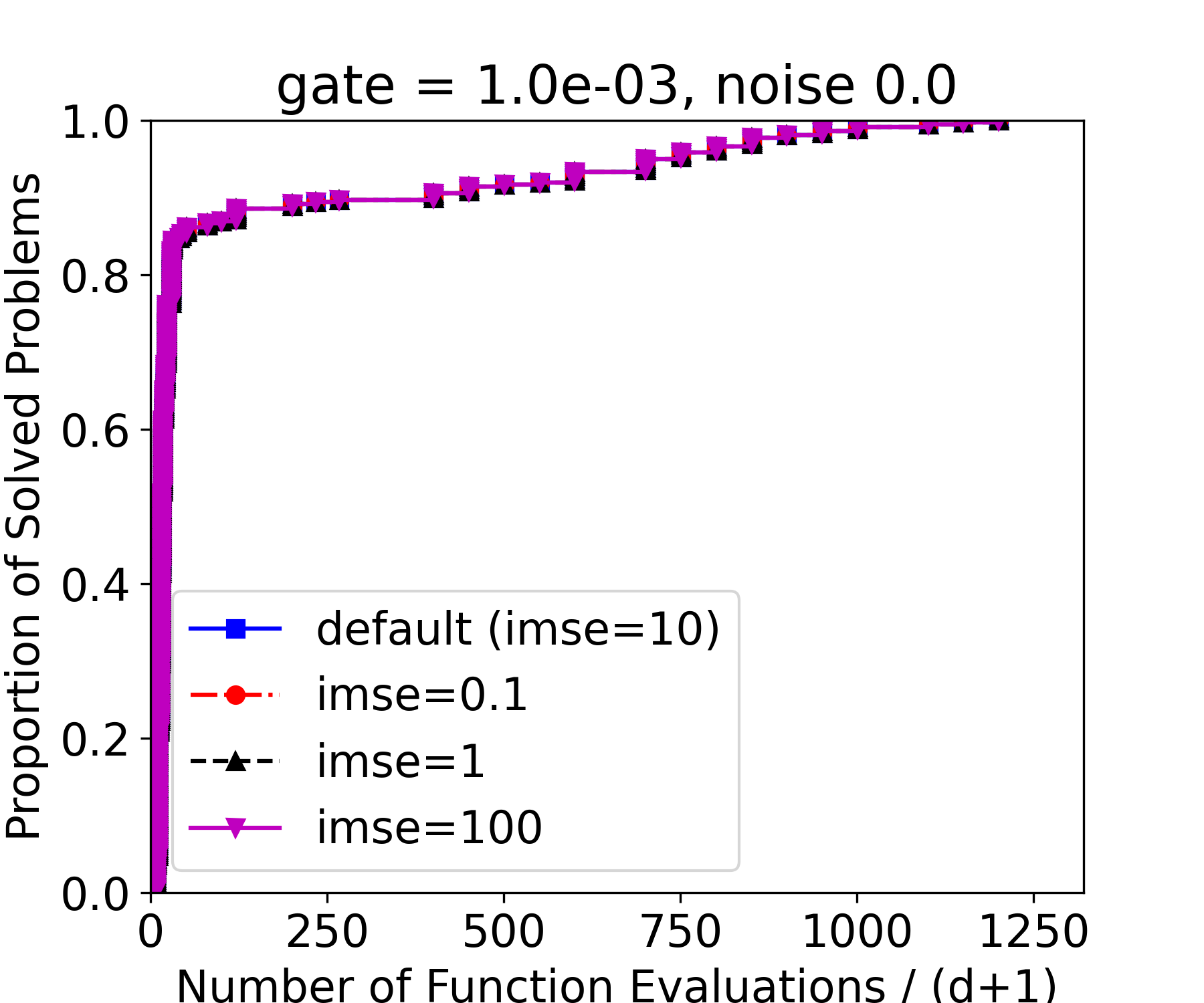}%
\includegraphics[width=0.333\textwidth, trim= 0 0 40 15, clip]{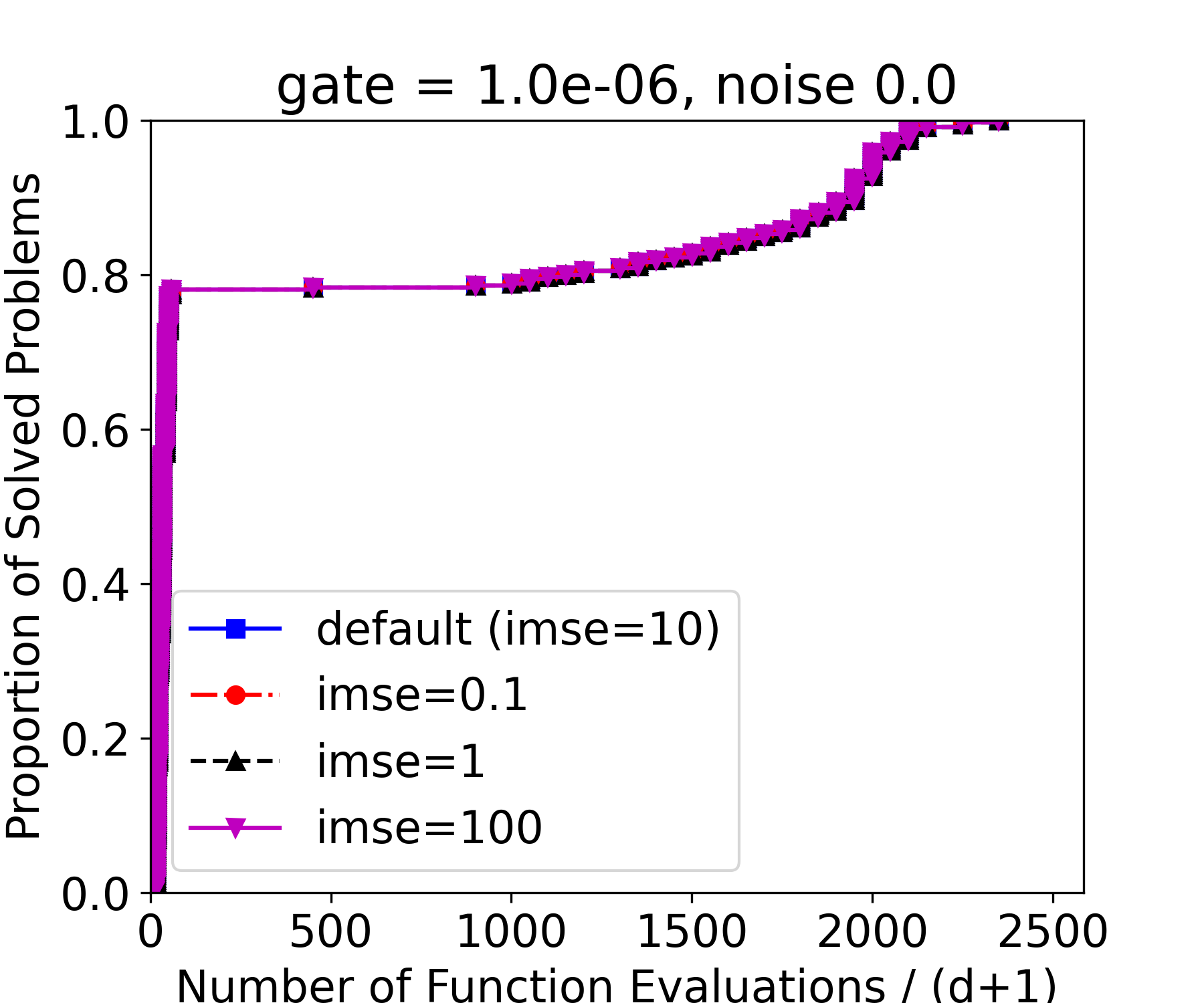}\\
\includegraphics[width=0.333\textwidth, trim= 0 0 40 15, clip]{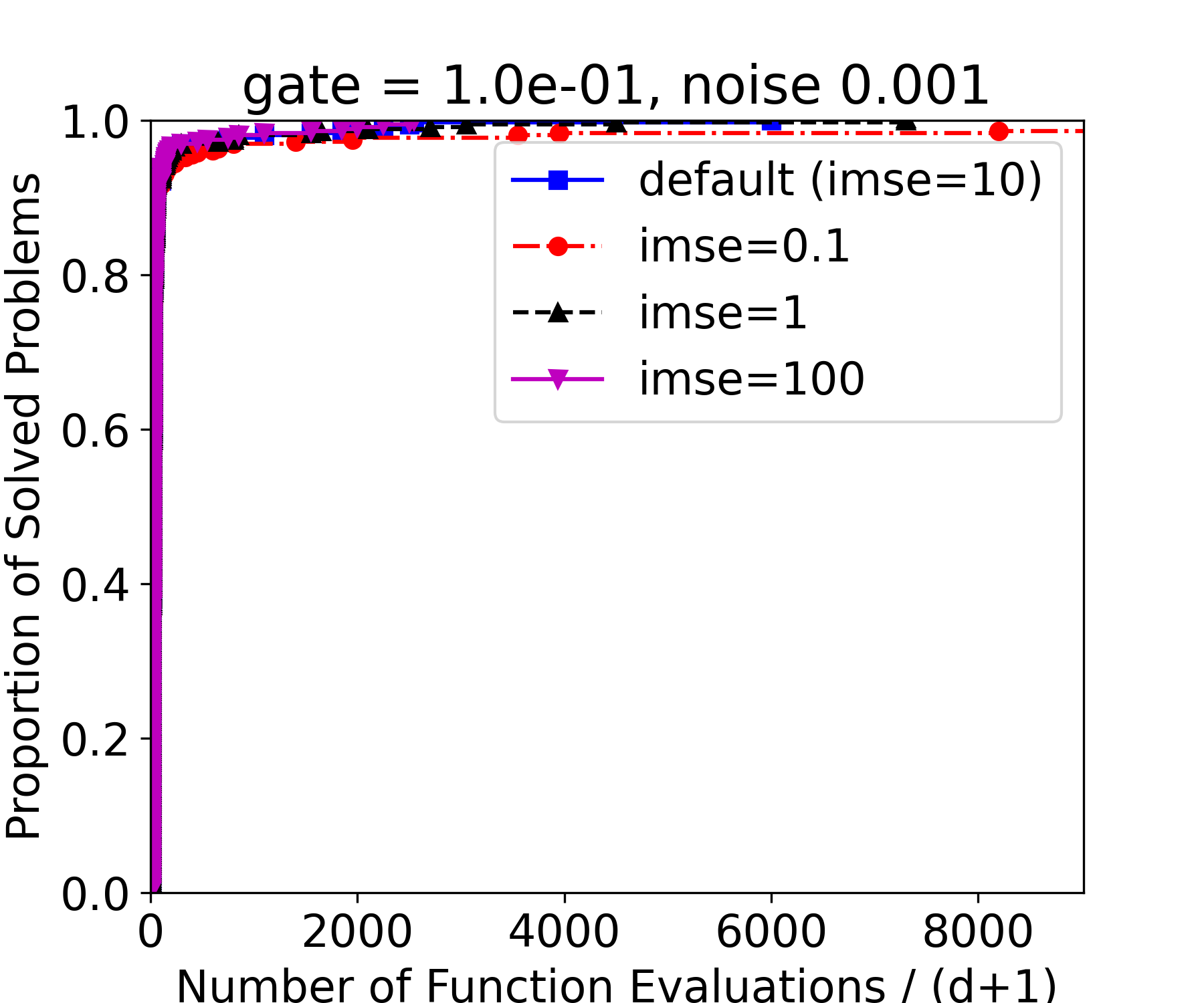}%
\includegraphics[width=0.333\textwidth, trim= 0 0 40 15, clip]{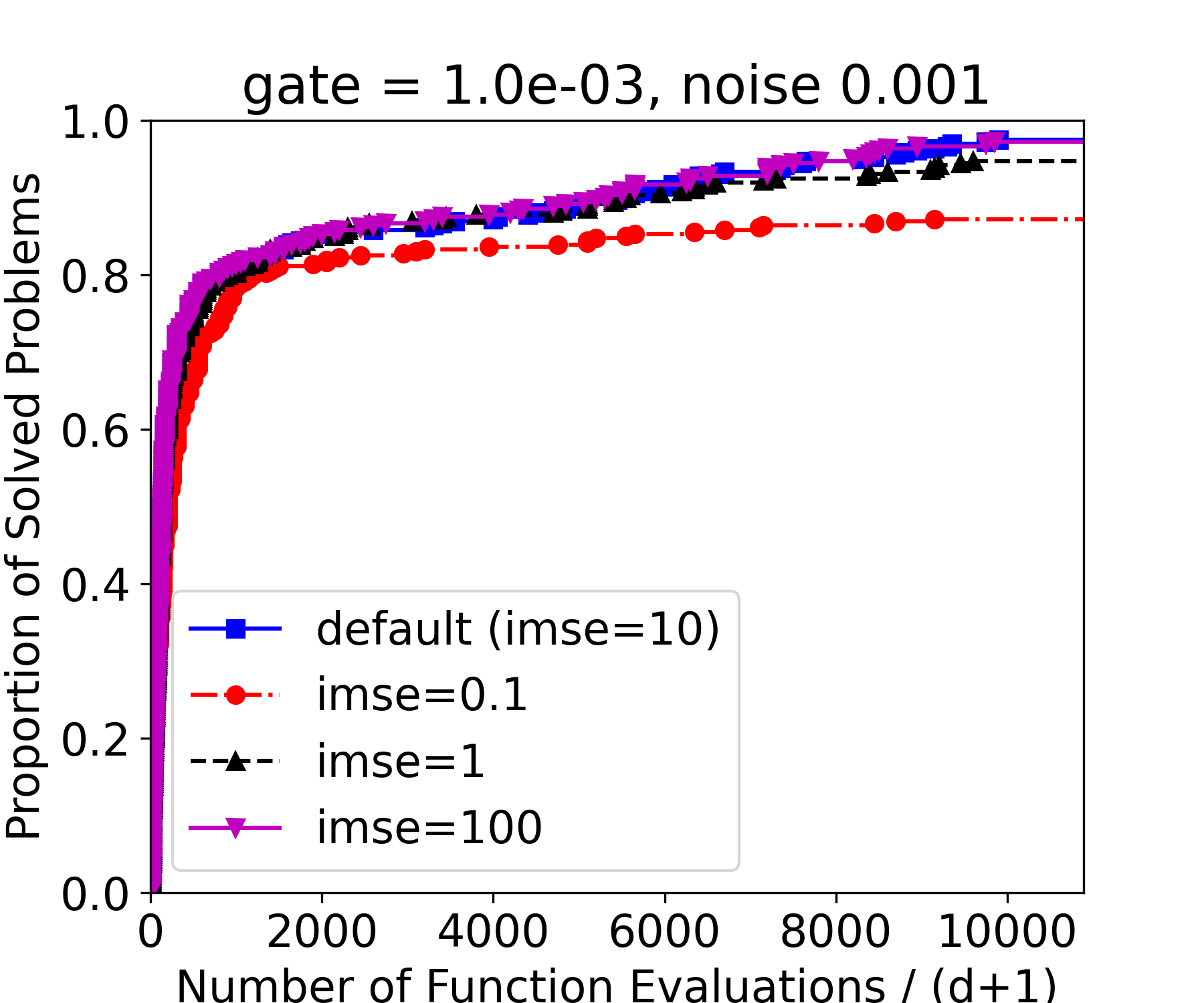}%
\includegraphics[width=0.333\textwidth, trim= 0 0 40 15, clip]{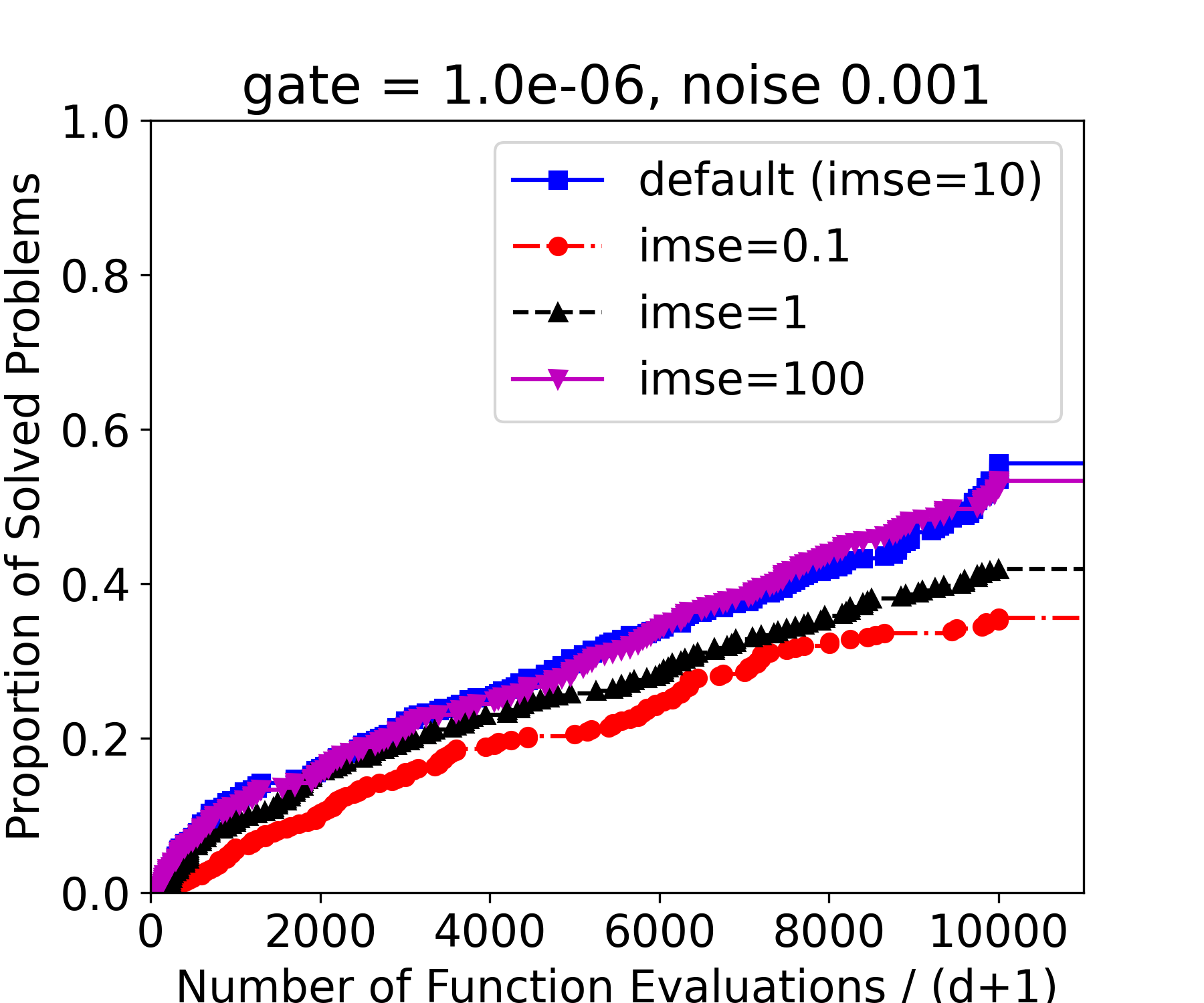}\\
\includegraphics[width=0.333\textwidth, trim= 0 0 40 15, clip]{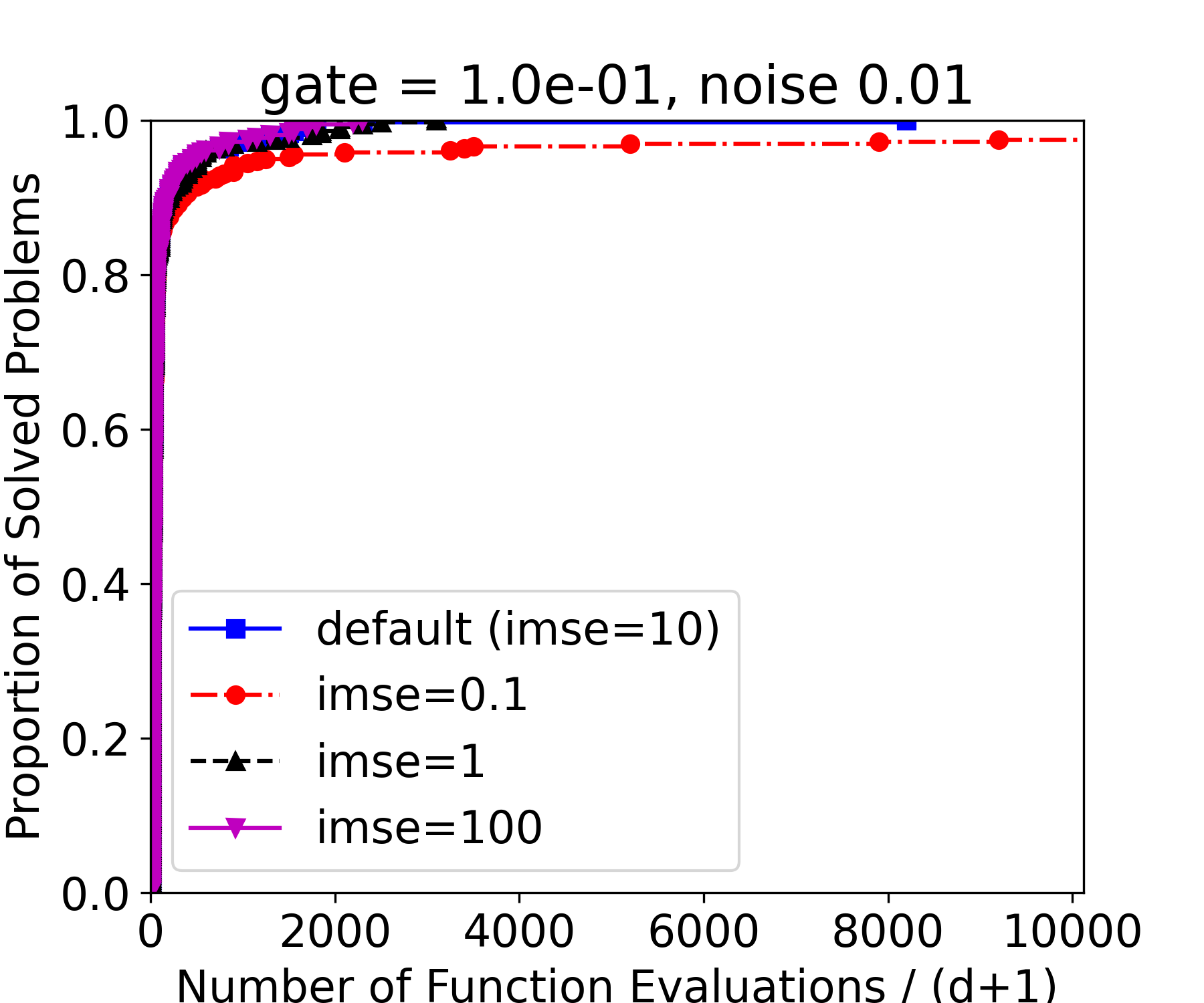}%
\includegraphics[width=0.333\textwidth, trim= 0 0 40 15, clip]{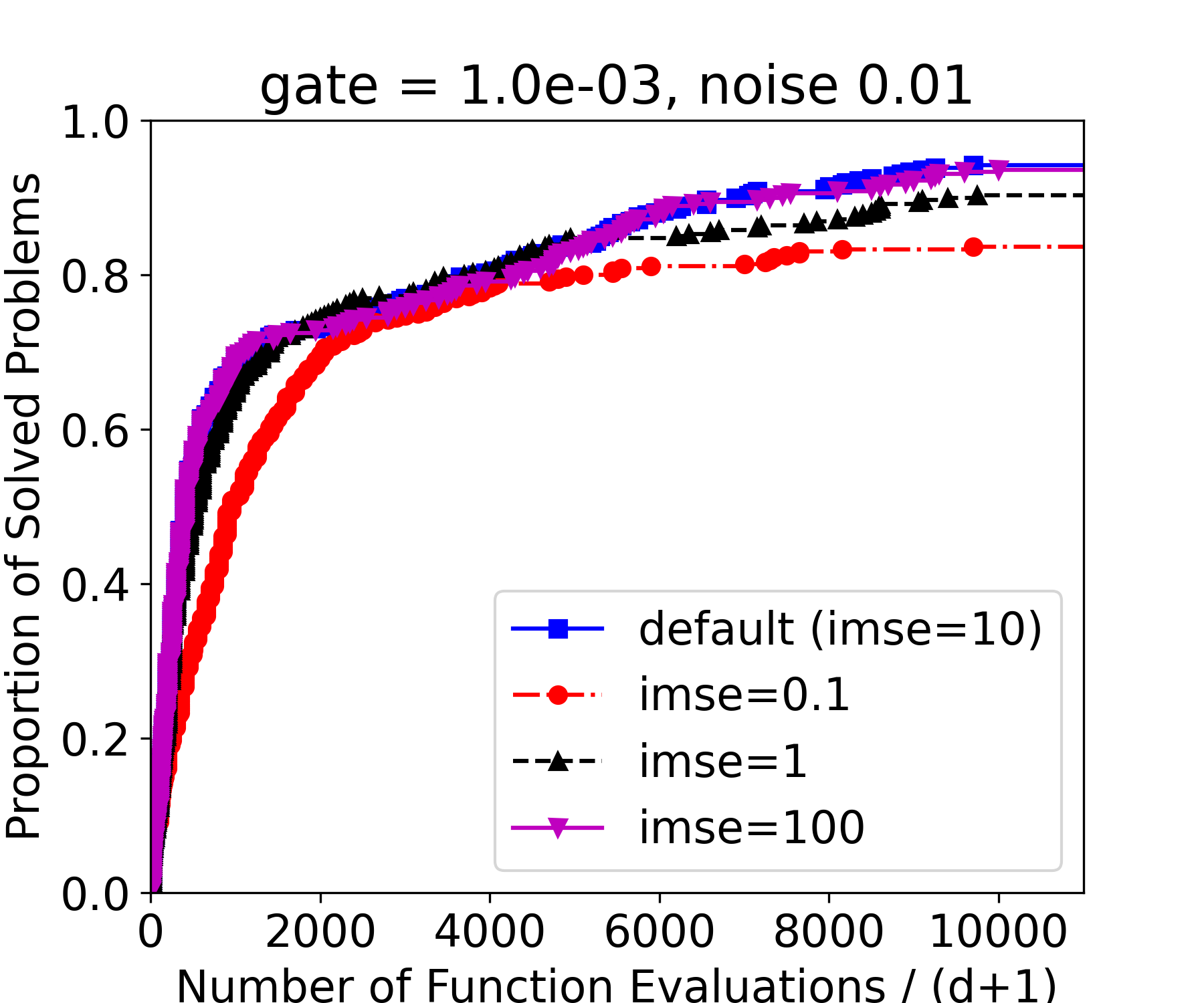}%
\includegraphics[width=0.333\textwidth, trim= 0 0 40 15, clip]{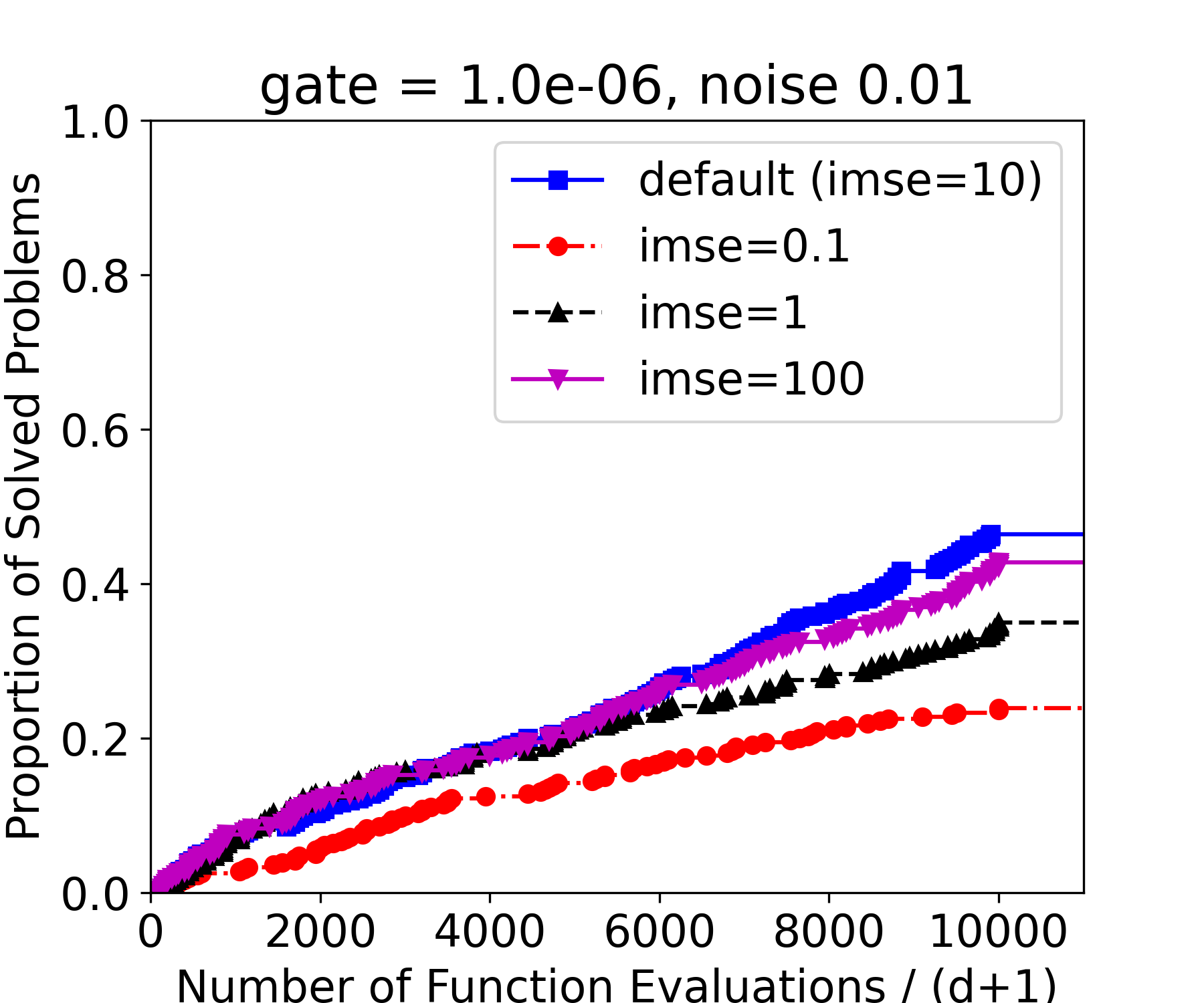}\\
\includegraphics[width=0.333\textwidth, trim= 0 0 40 15, clip]{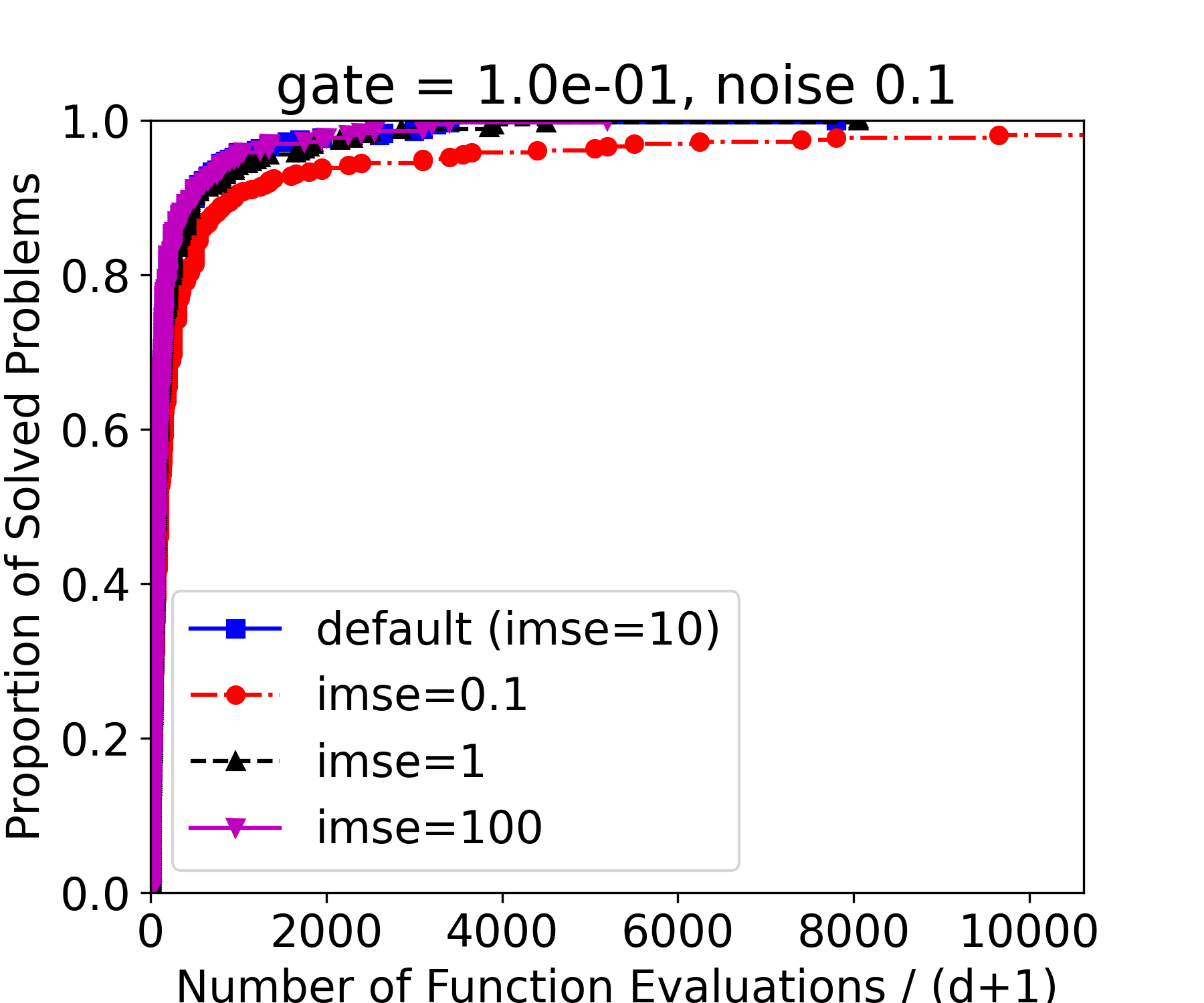}%
\includegraphics[width=0.333\textwidth, trim= 0 0 40 15, clip]{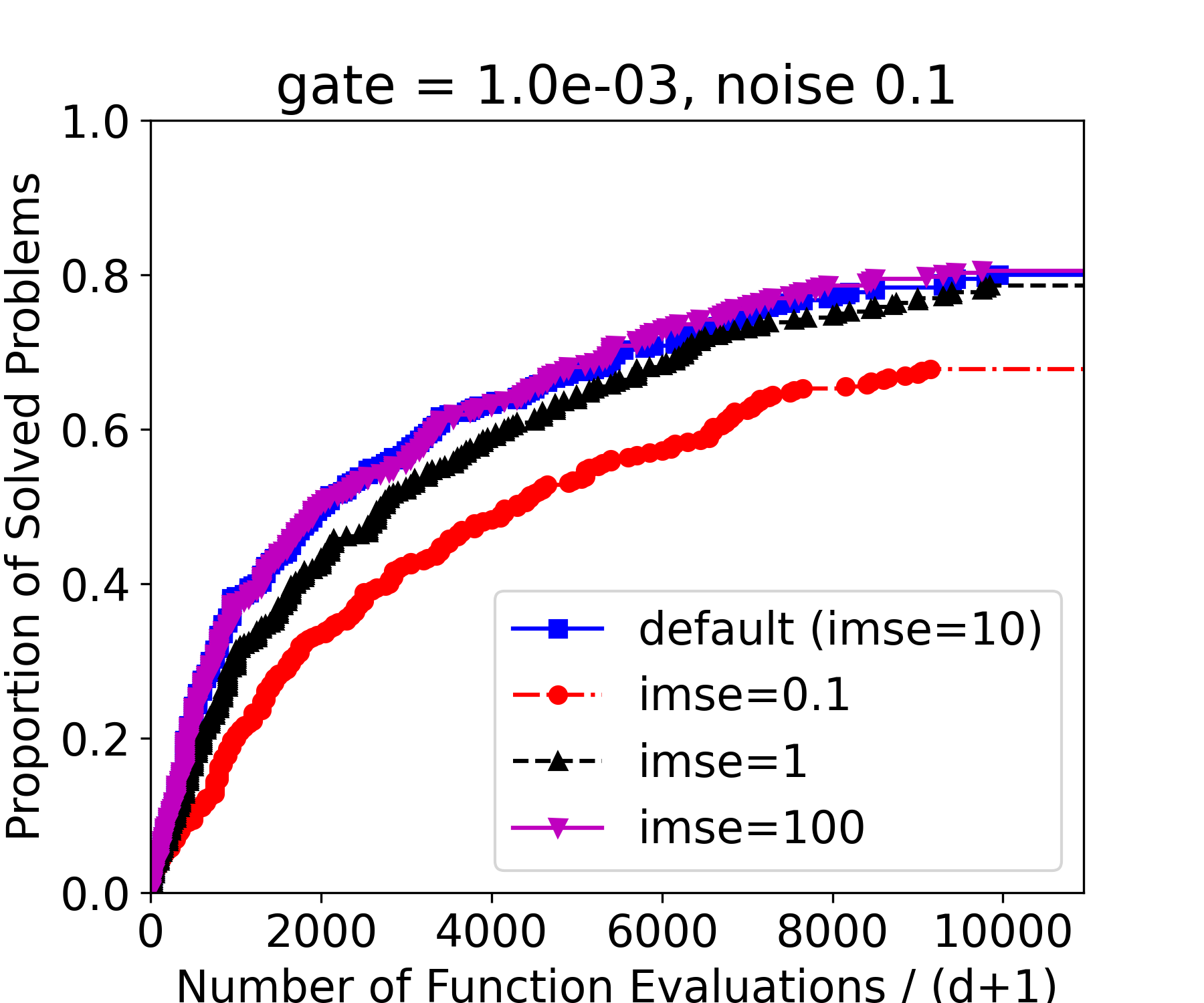}%
\includegraphics[width=0.333\textwidth, trim= 0 0 40 15, clip]{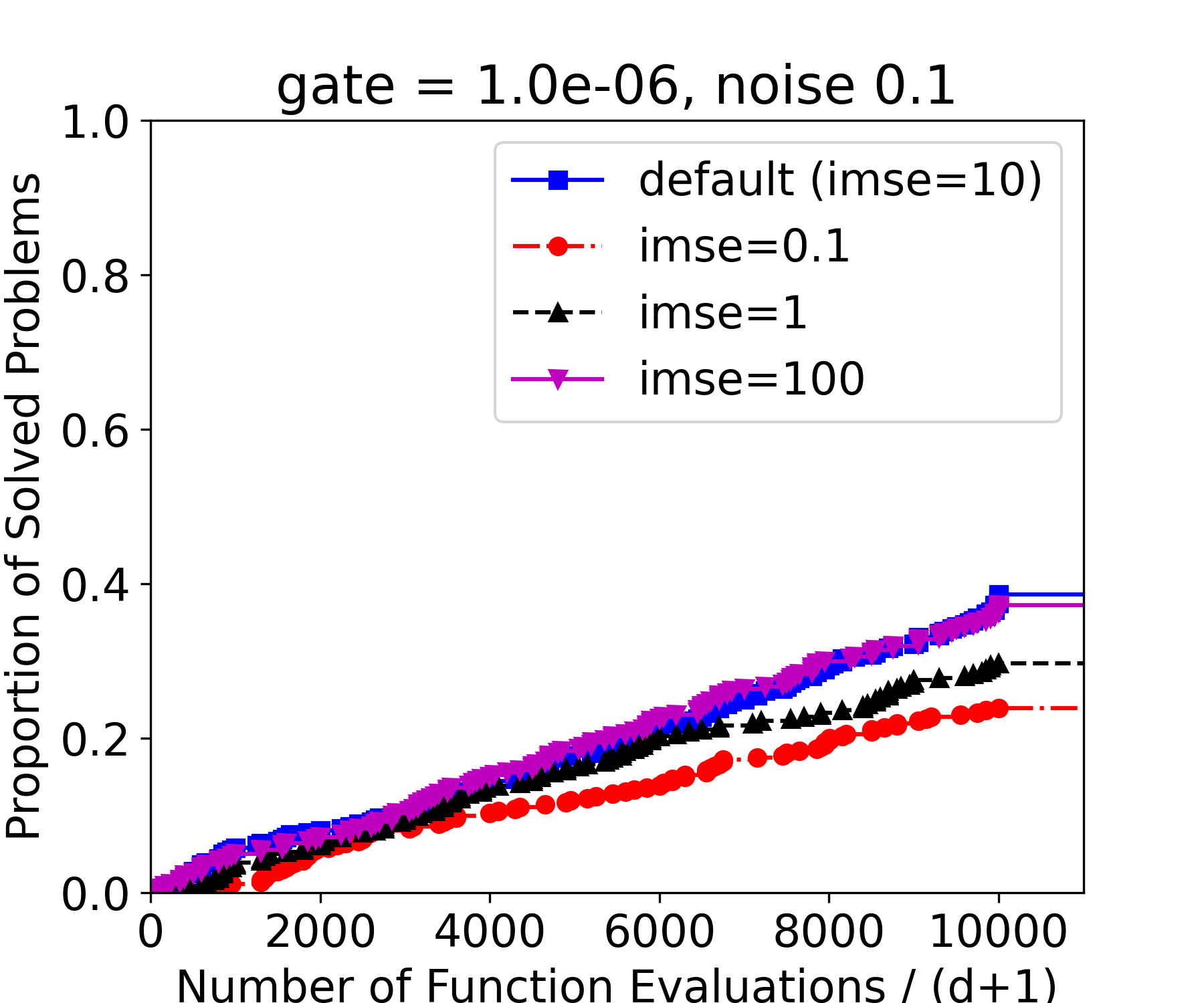}
\caption{Data profiles for varying the IMSE - variance ratio, over 40 repetitions.}
\label{fig:localbenchprofsimse}
\end{figure}

In Step 19 of \Cref{alg:OGPIT}, the variance of the mean needs to be 10 times
larger than the mean predictive variance (IMSE) to allow for decreasing the TR radius. In 
\Cref{fig:localbenchprofsimse}, we study other values. Except for the
 deterministic version where it is not used, smaller values for the ratio
 (0.1 or 1) show a degraded performance and highlight the importance of
 preventing the decrease when the noise is dominating. The largest values,
 100, gets sligthly worse results than the default value.

\subsection{Value of the minimal variance reduction}

\begin{figure}
\includegraphics[width=0.333\textwidth, trim= 0 0 40 15, clip]{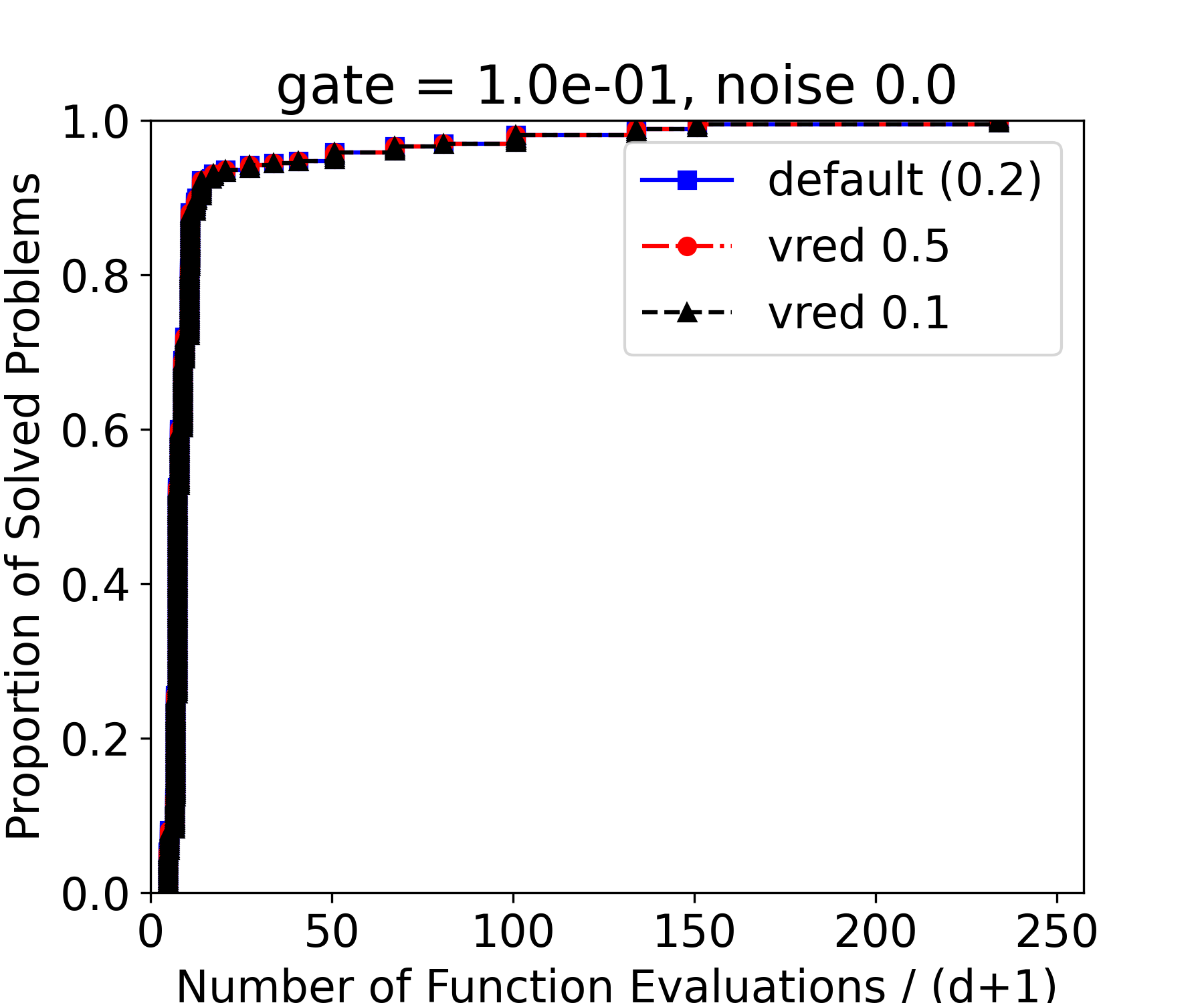}%
\includegraphics[width=0.333\textwidth, trim= 0 0 40 15, clip]{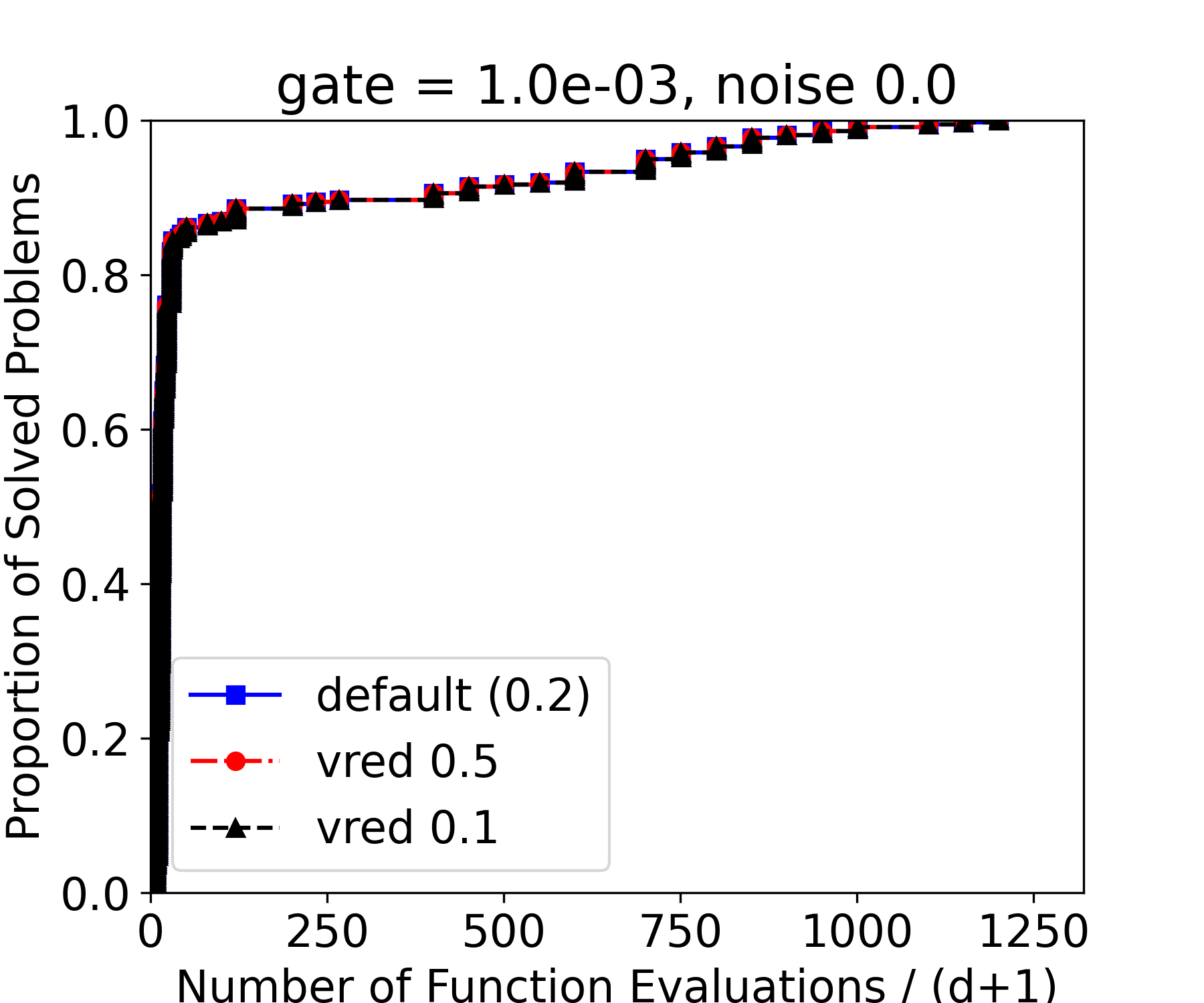}%
\includegraphics[width=0.333\textwidth, trim= 0 0 40 15, clip]{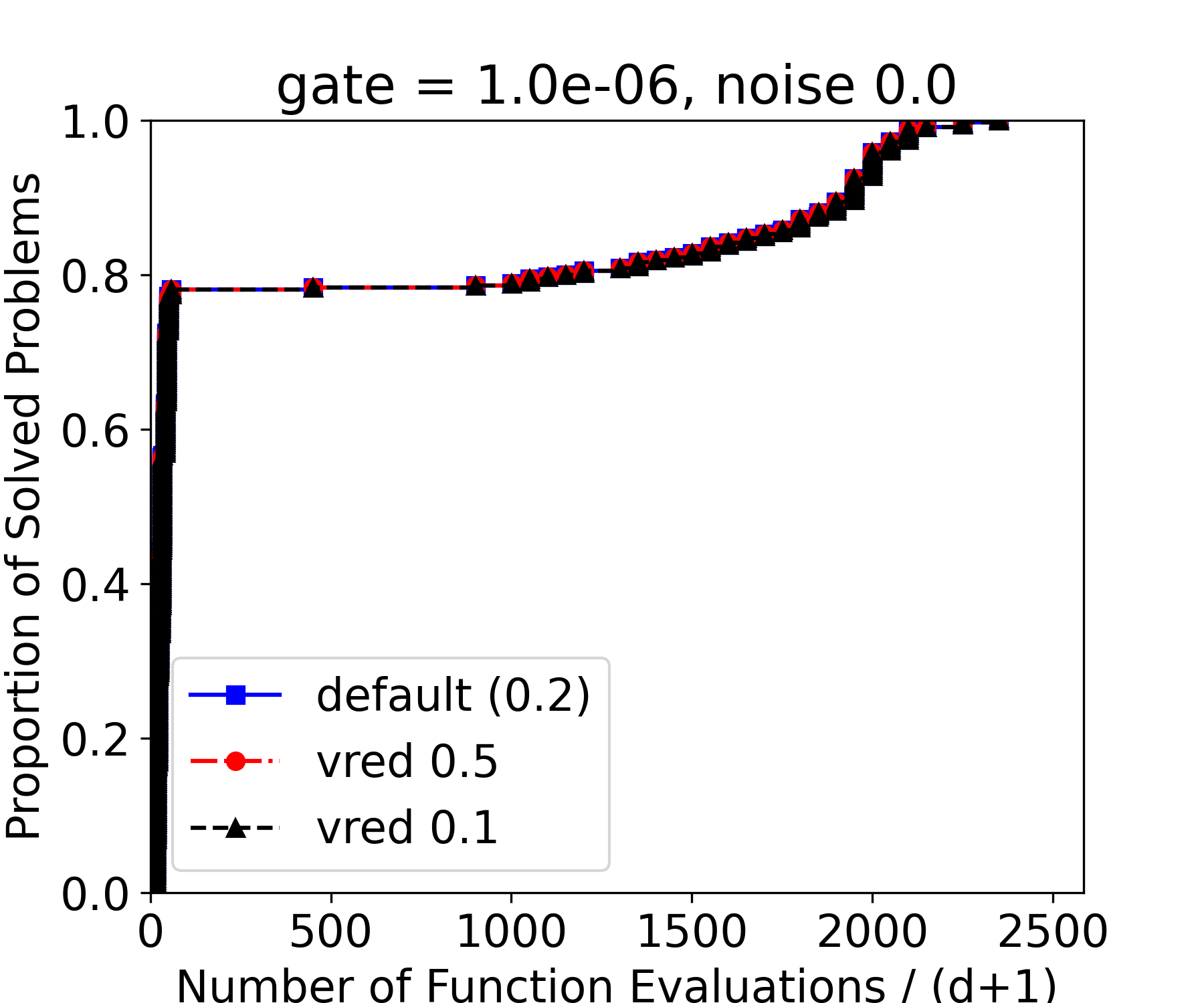}\\
\includegraphics[width=0.333\textwidth, trim= 0 0 40 15, clip]{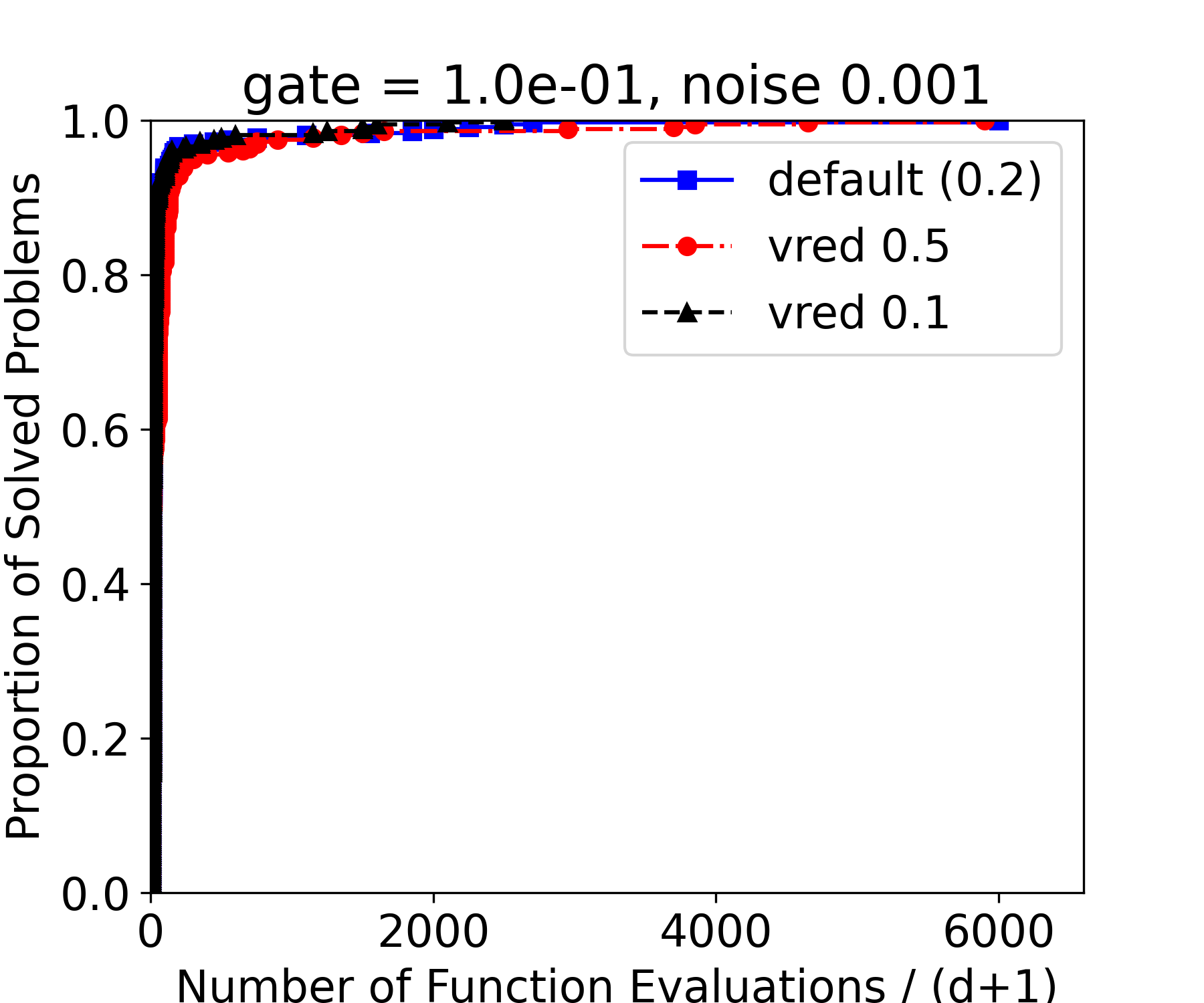}%
\includegraphics[width=0.333\textwidth, trim= 0 0 40 15, clip]{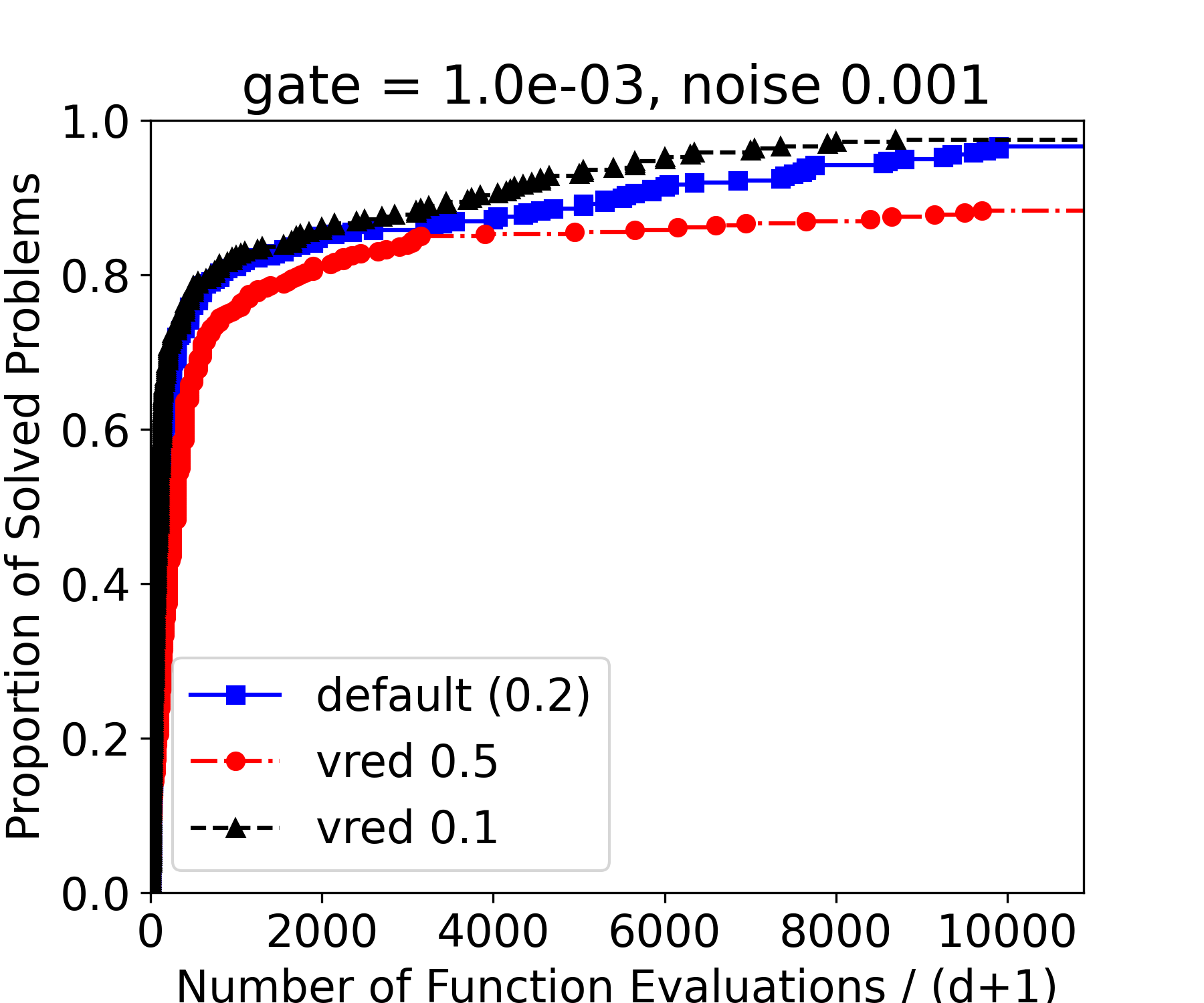}%
\includegraphics[width=0.333\textwidth, trim= 0 0 40 15, clip]{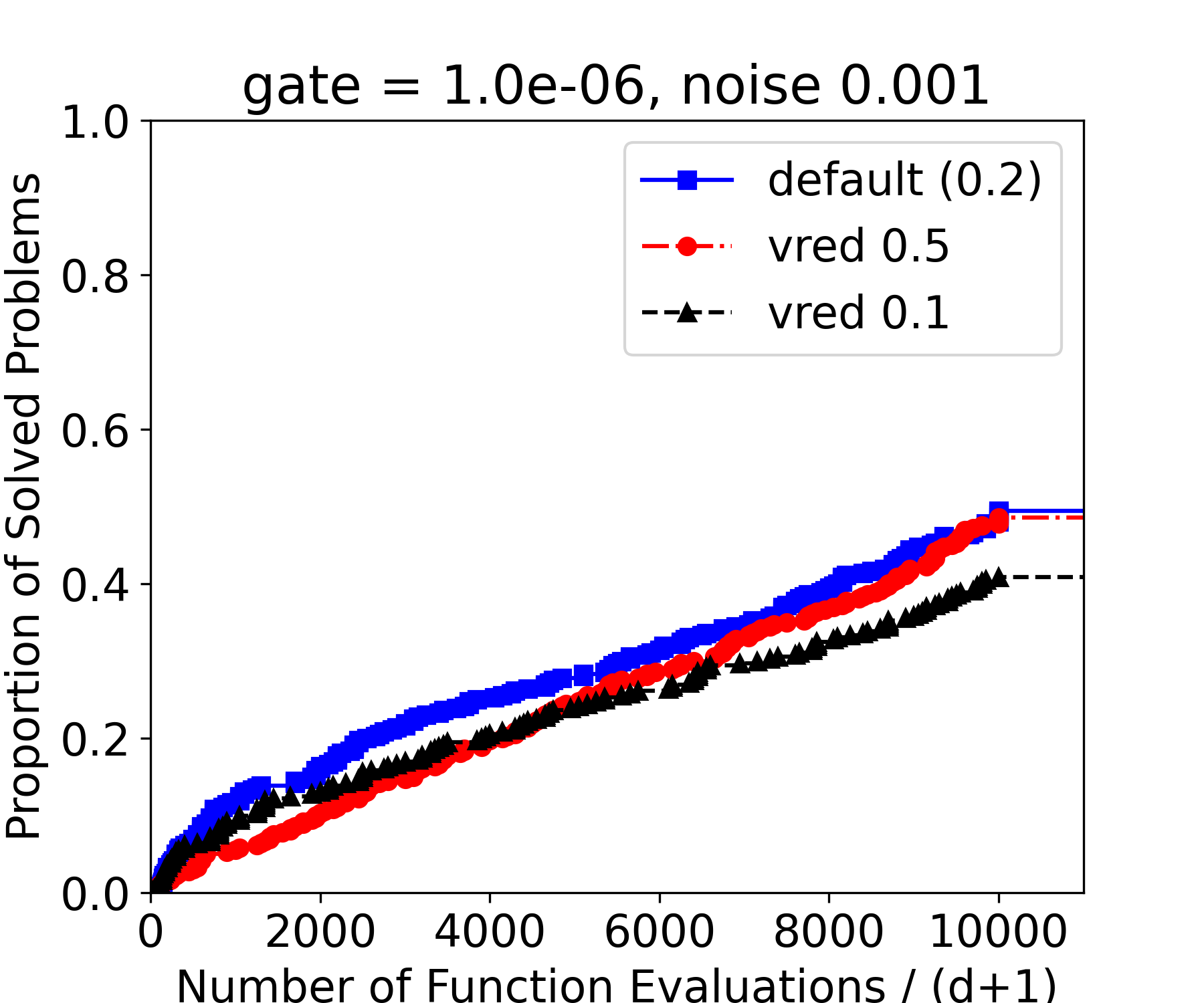}\\
\includegraphics[width=0.333\textwidth, trim= 0 0 40 15, clip]{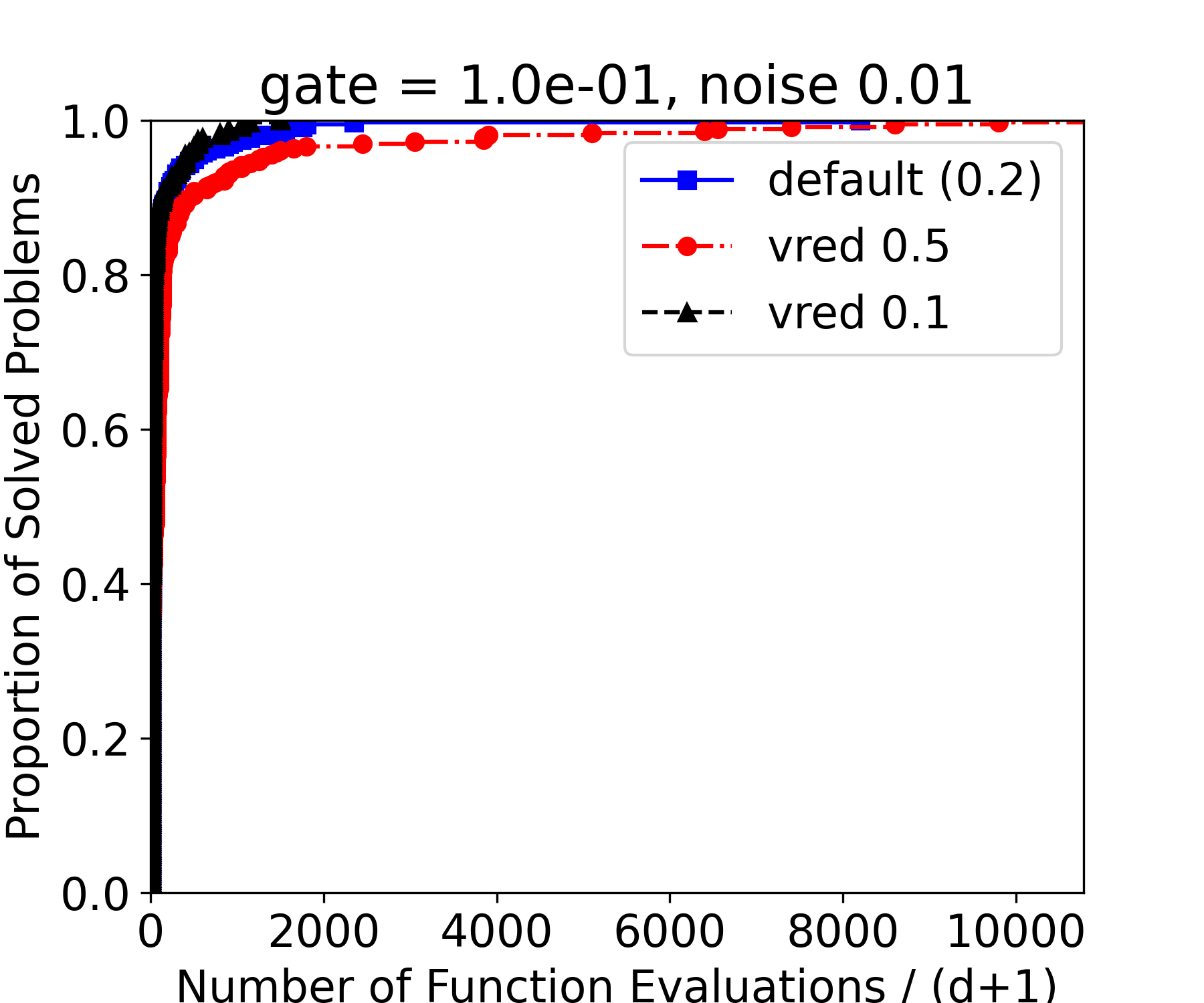}%
\includegraphics[width=0.333\textwidth, trim= 0 0 40 15, clip]{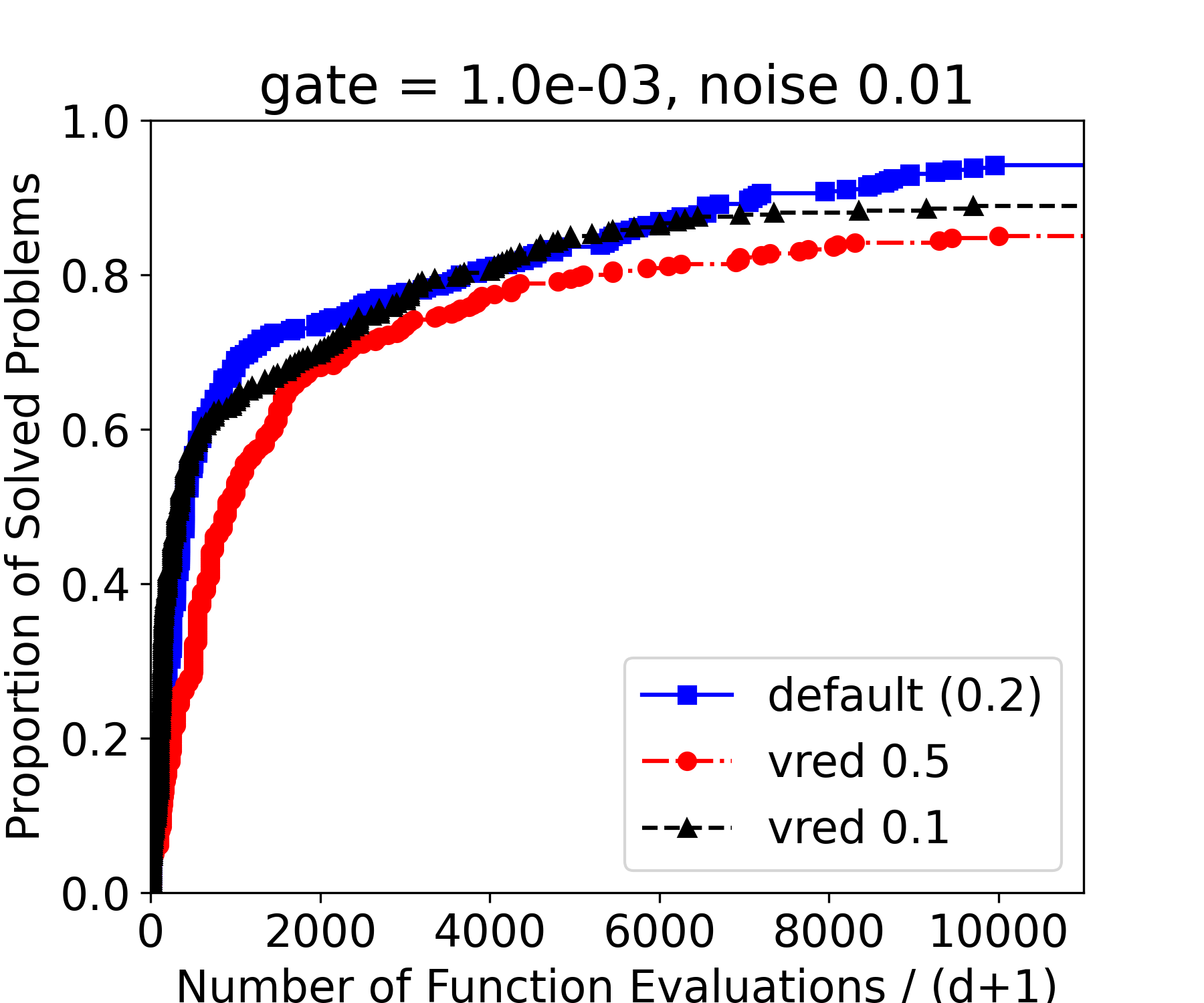}%
\includegraphics[width=0.333\textwidth, trim= 0 0 40 15, clip]{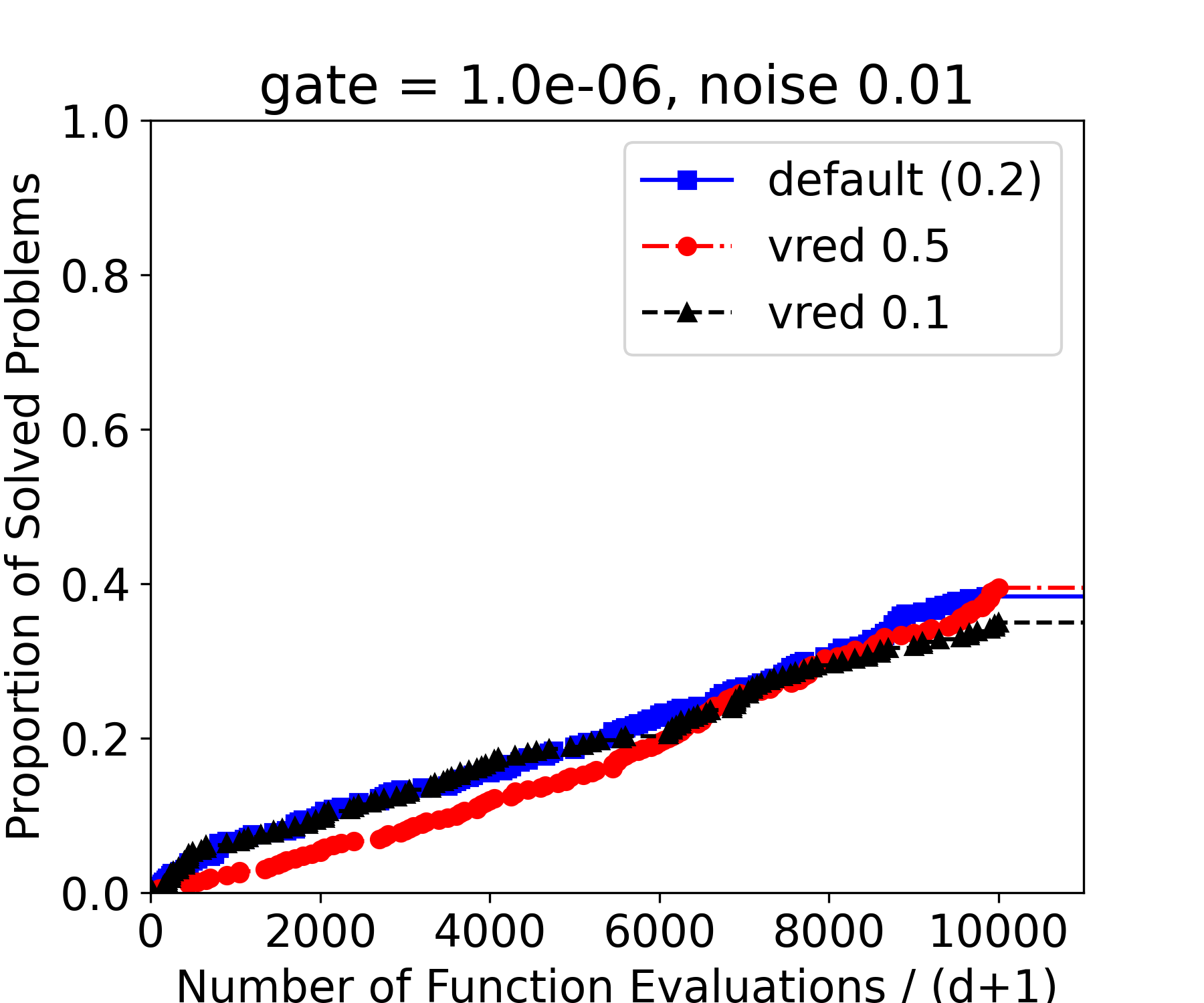}\\
\includegraphics[width=0.333\textwidth, trim= 0 0 40 15, clip]{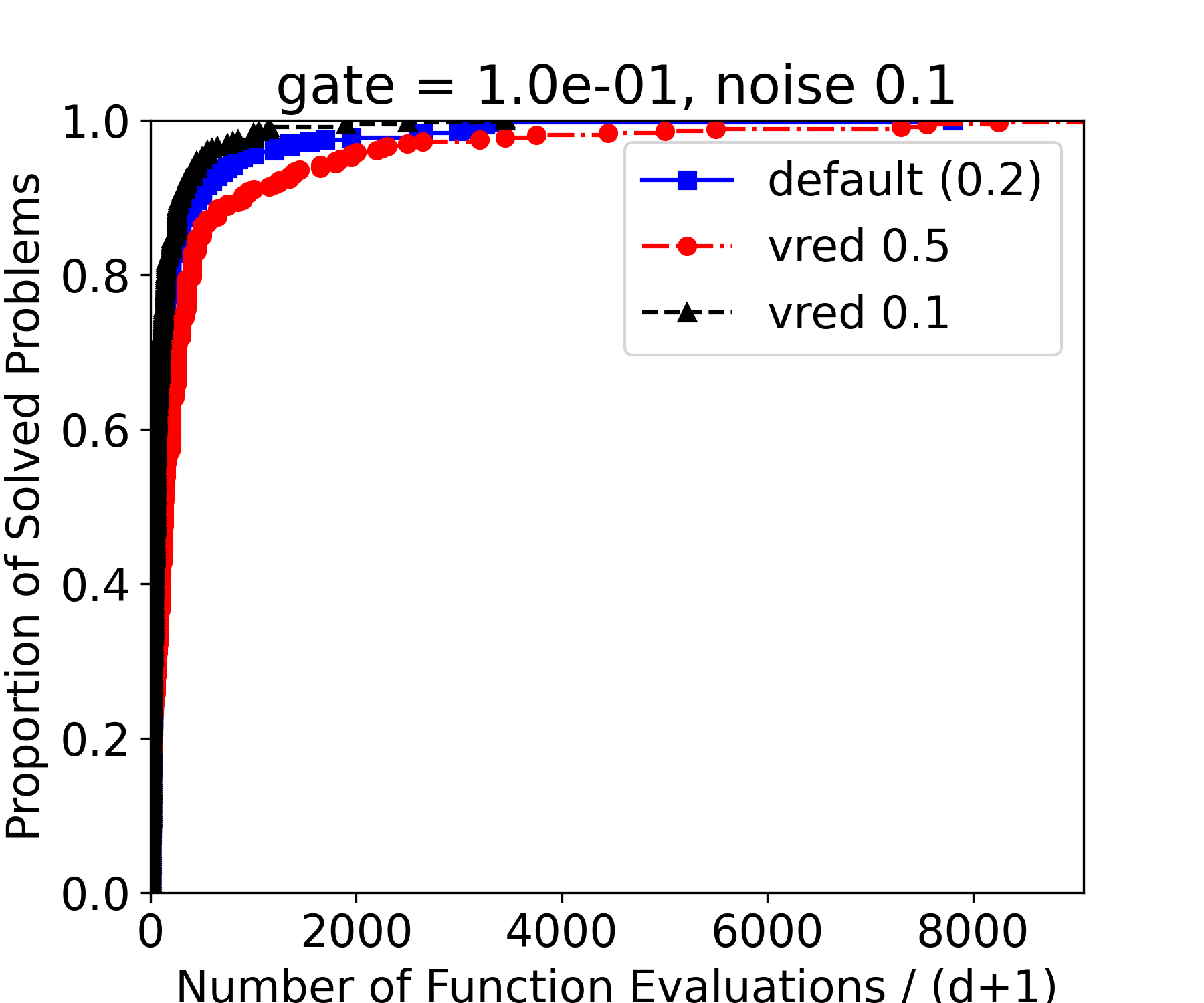}%
\includegraphics[width=0.333\textwidth, trim= 0 0 40 15, clip]{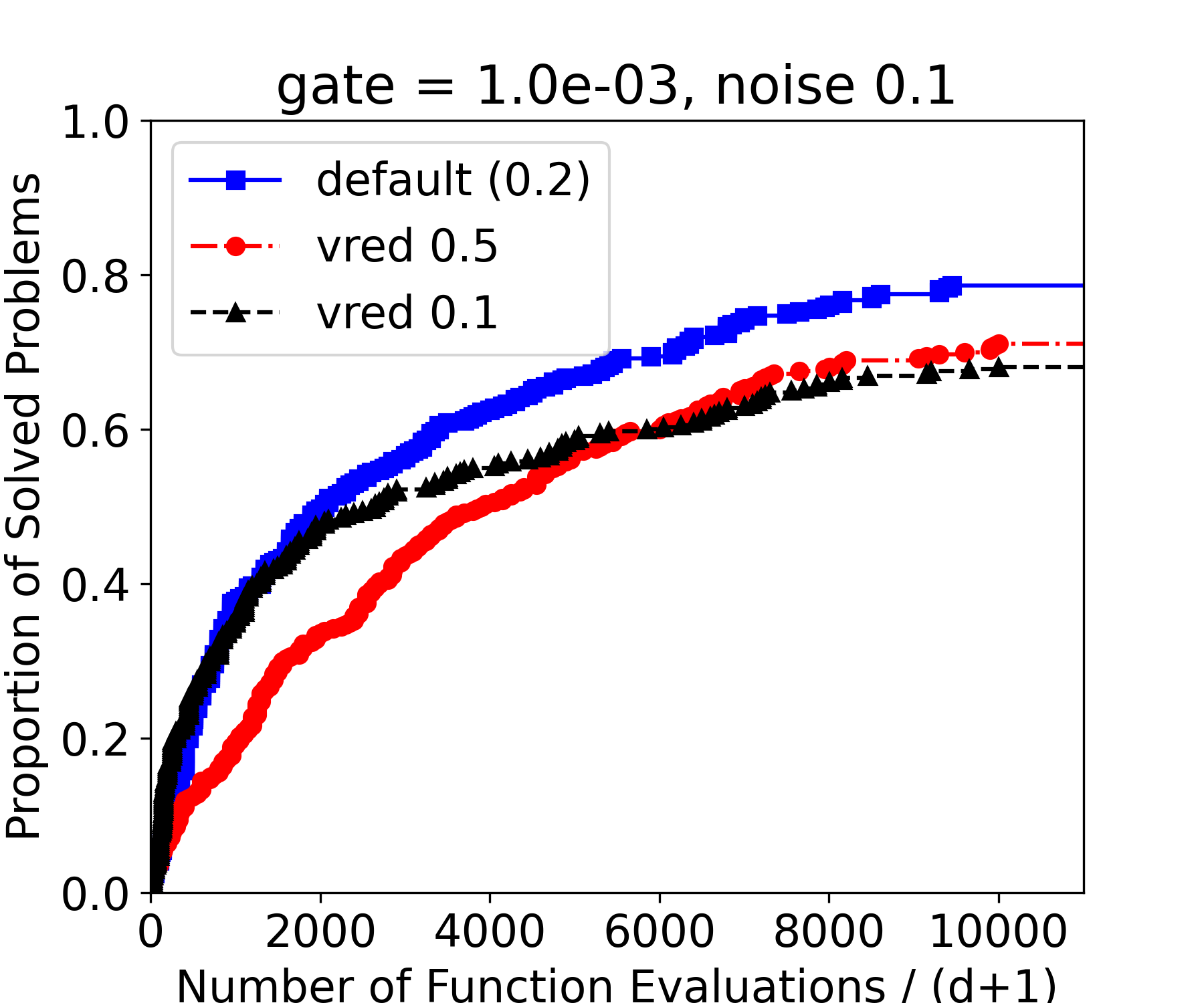}%
\includegraphics[width=0.333\textwidth, trim= 0 0 40 15, clip]{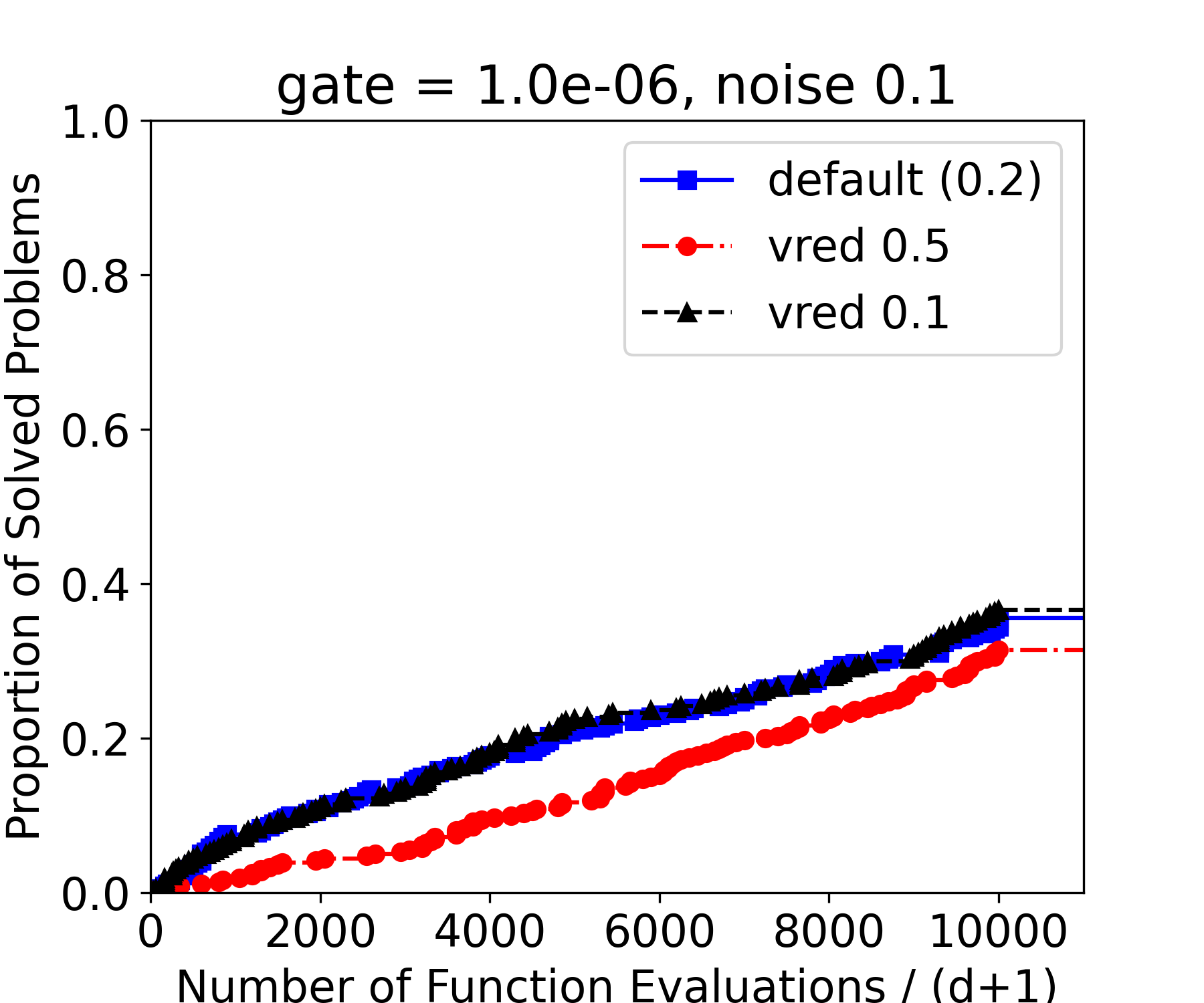}
\caption{Data profiles with varying minimal variance reduction $T_a$.}
\label{fig:localbenchvred}
\end{figure}

\Cref{fig:localbenchvred} deals with the imposed minimal variance reduction at
 the new design, $T_a$. The default value shows the best overall performance,
 while the largest value ($T_a = 0.5$) can be detrimental.

\subsection{Value of the predictive variance ratio}

\begin{figure}
\includegraphics[width=0.333\textwidth, trim= 0 0 40 15, clip]{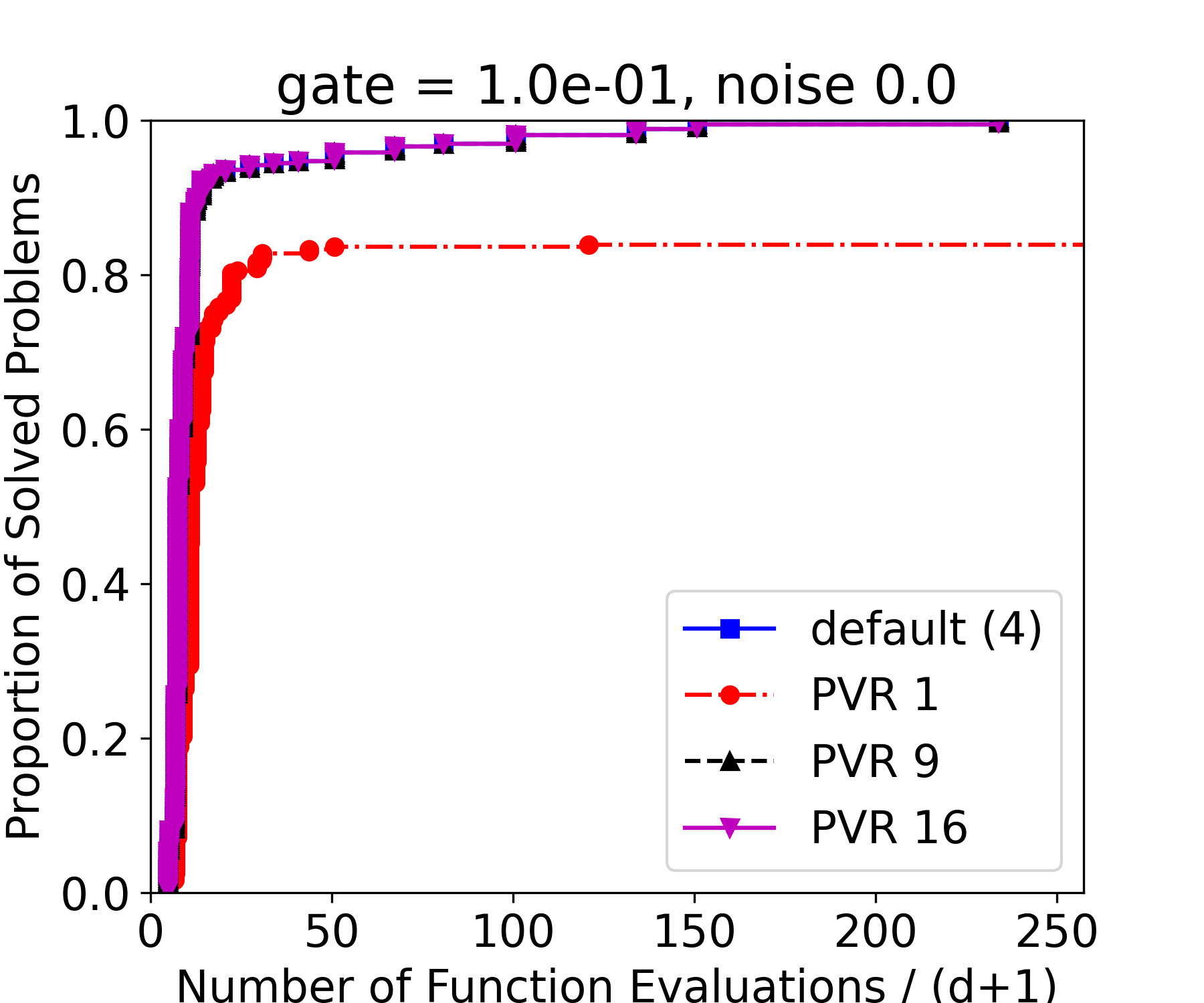}%
\includegraphics[width=0.333\textwidth, trim= 0 0 40 15, clip]{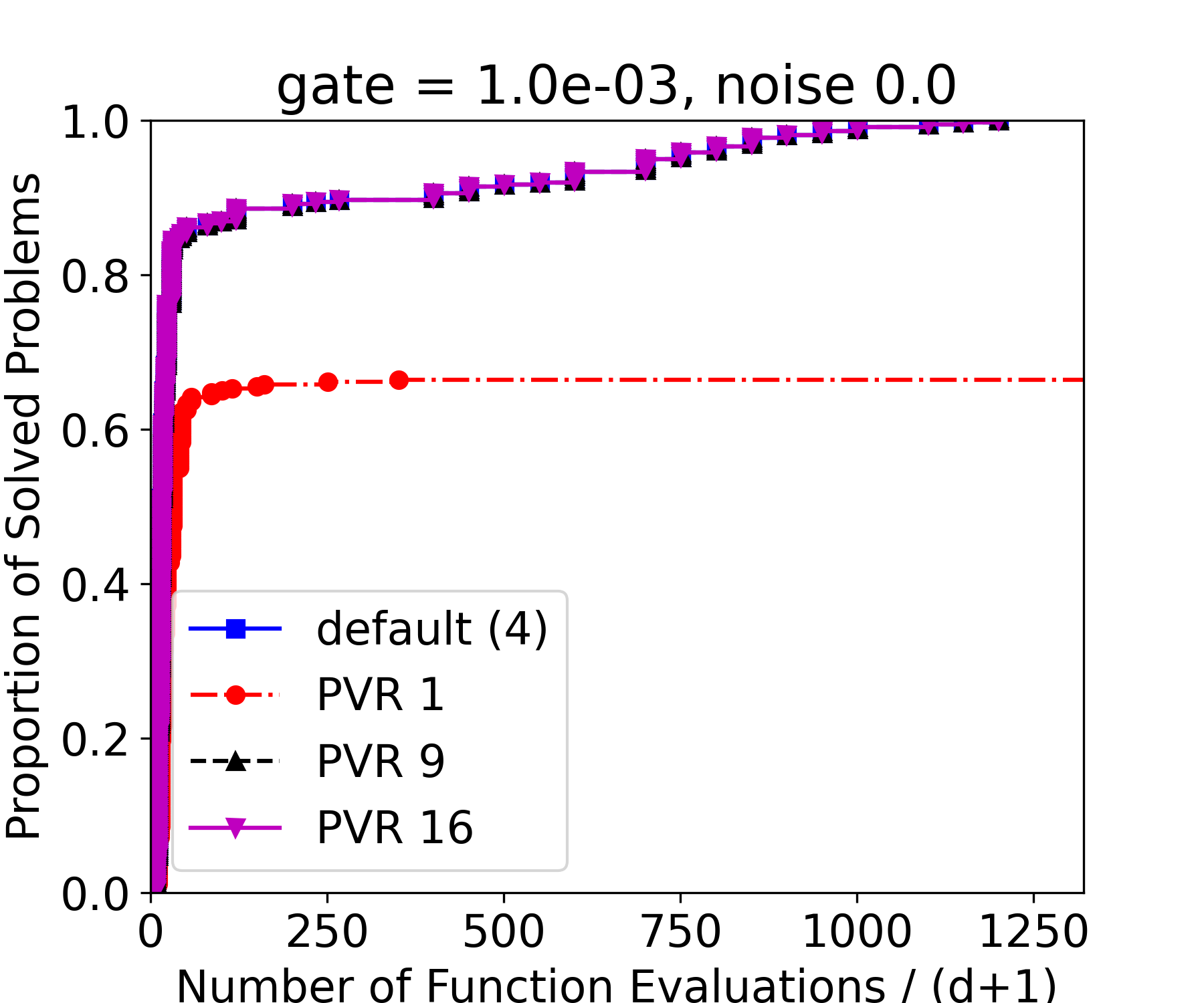}%
\includegraphics[width=0.333\textwidth, trim= 0 0 40 15, clip]{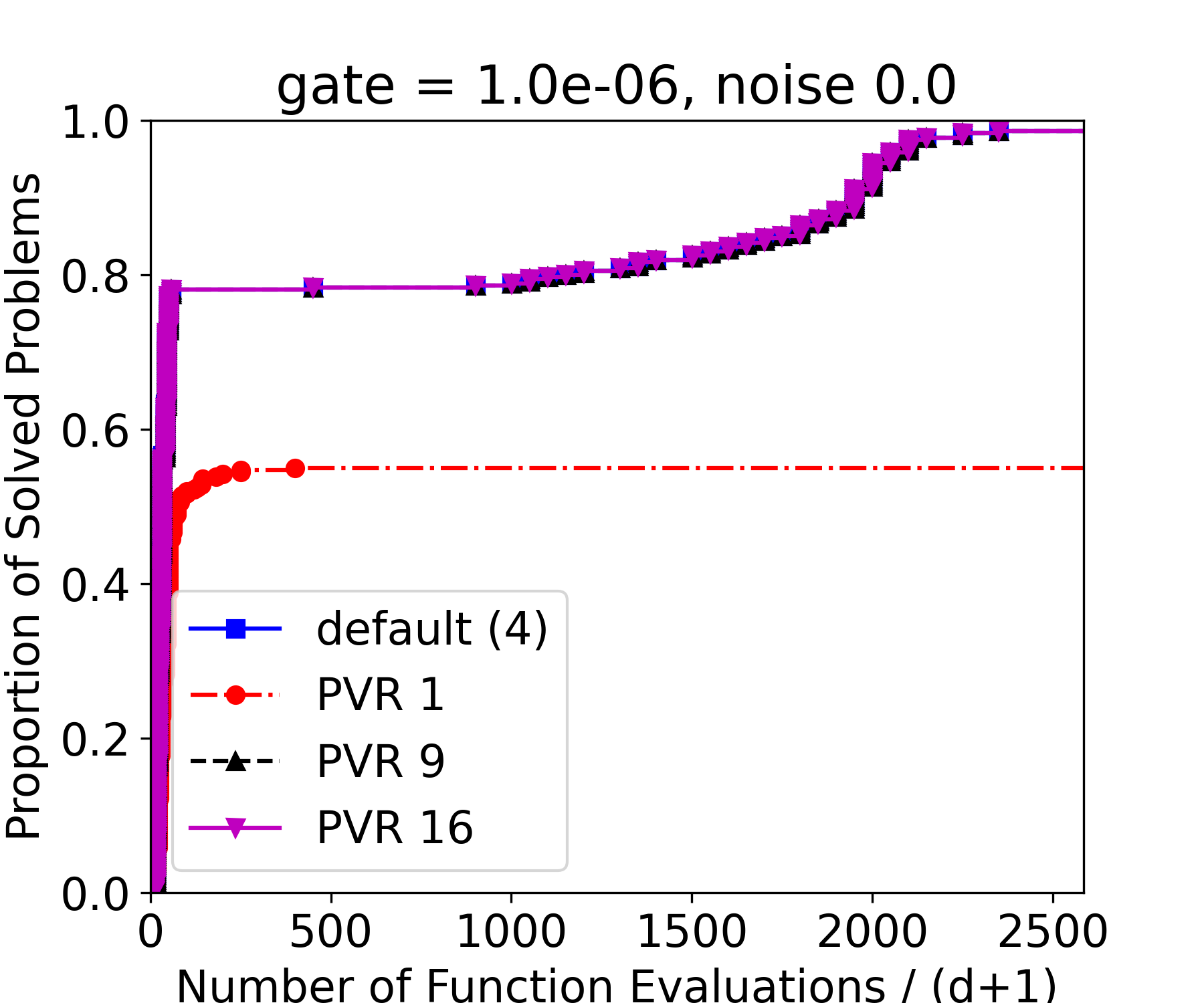}\\
\includegraphics[width=0.333\textwidth, trim= 0 0 40 15, clip]{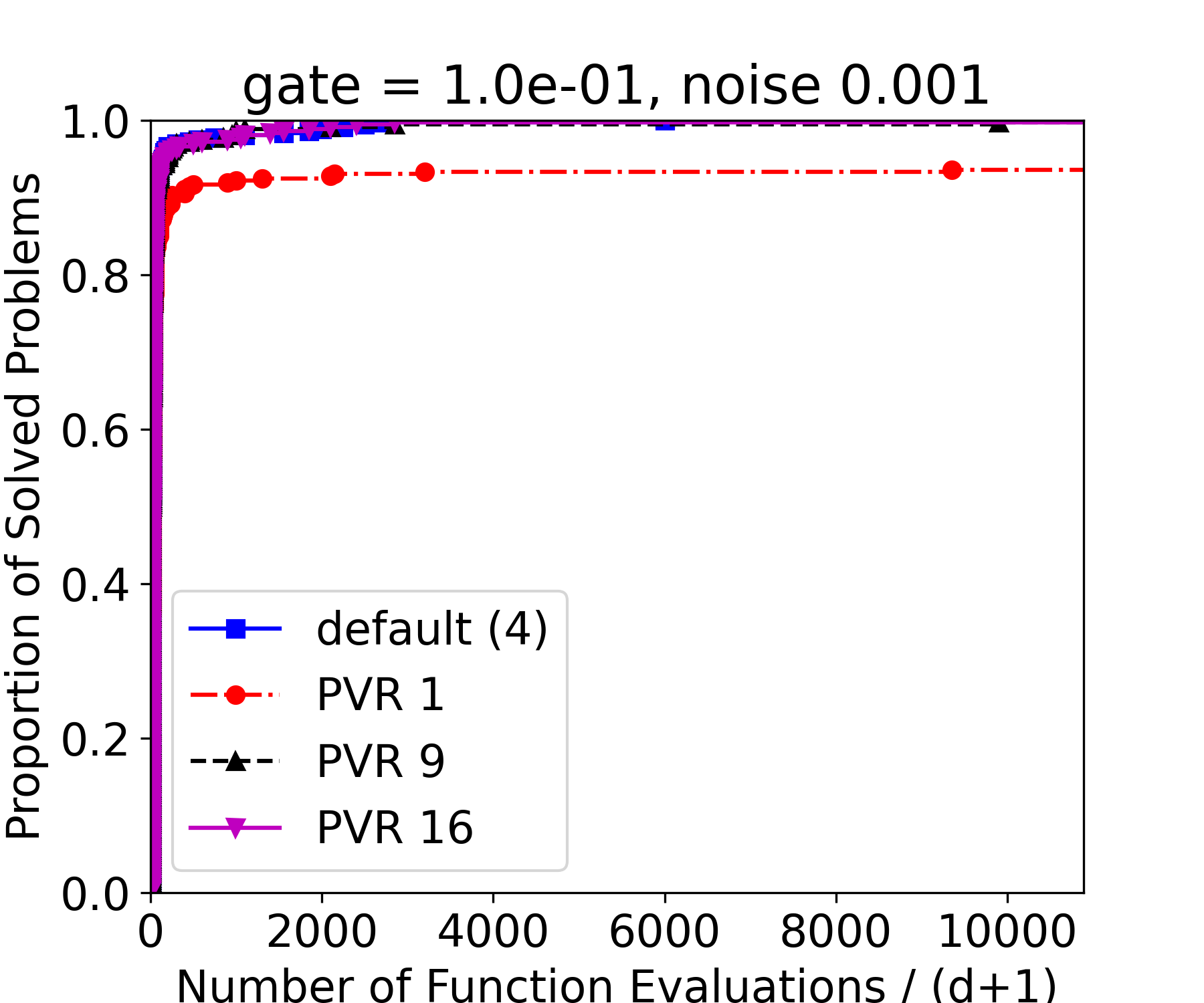}%
\includegraphics[width=0.333\textwidth, trim= 0 0 40 15, clip]{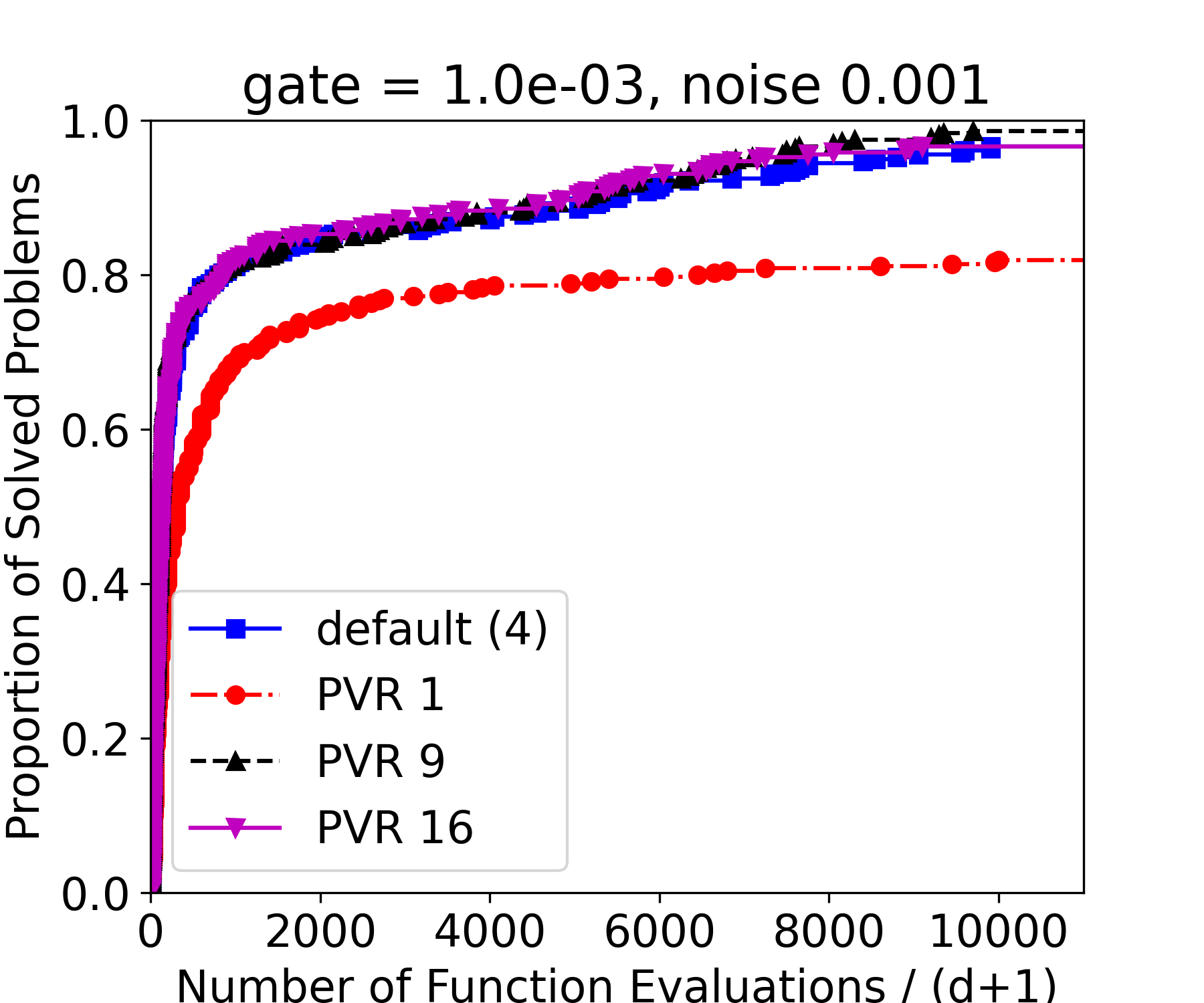}%
\includegraphics[width=0.333\textwidth, trim= 0 0 40 15, clip]{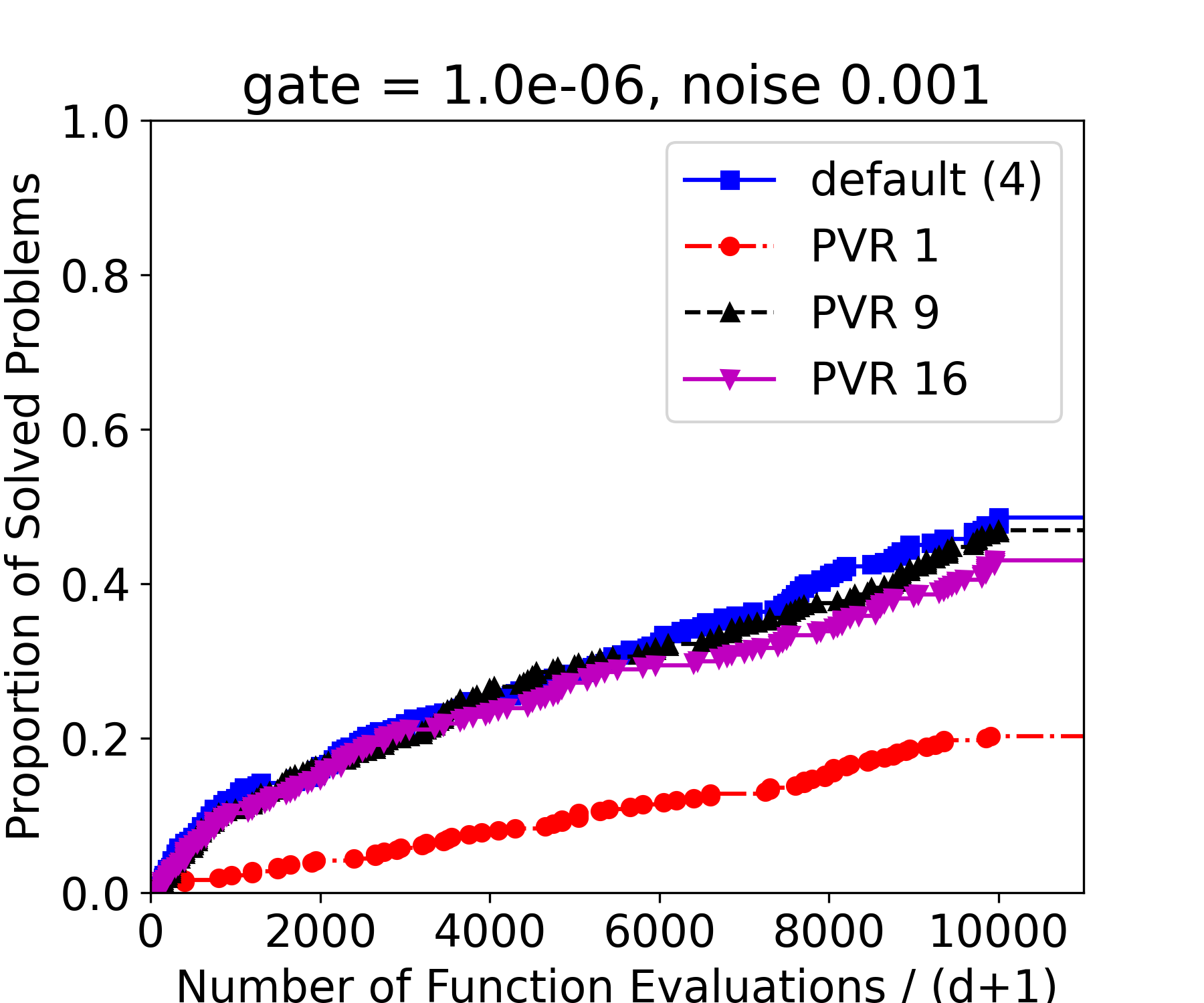}\\
\includegraphics[width=0.333\textwidth, trim= 0 0 40 15, clip]{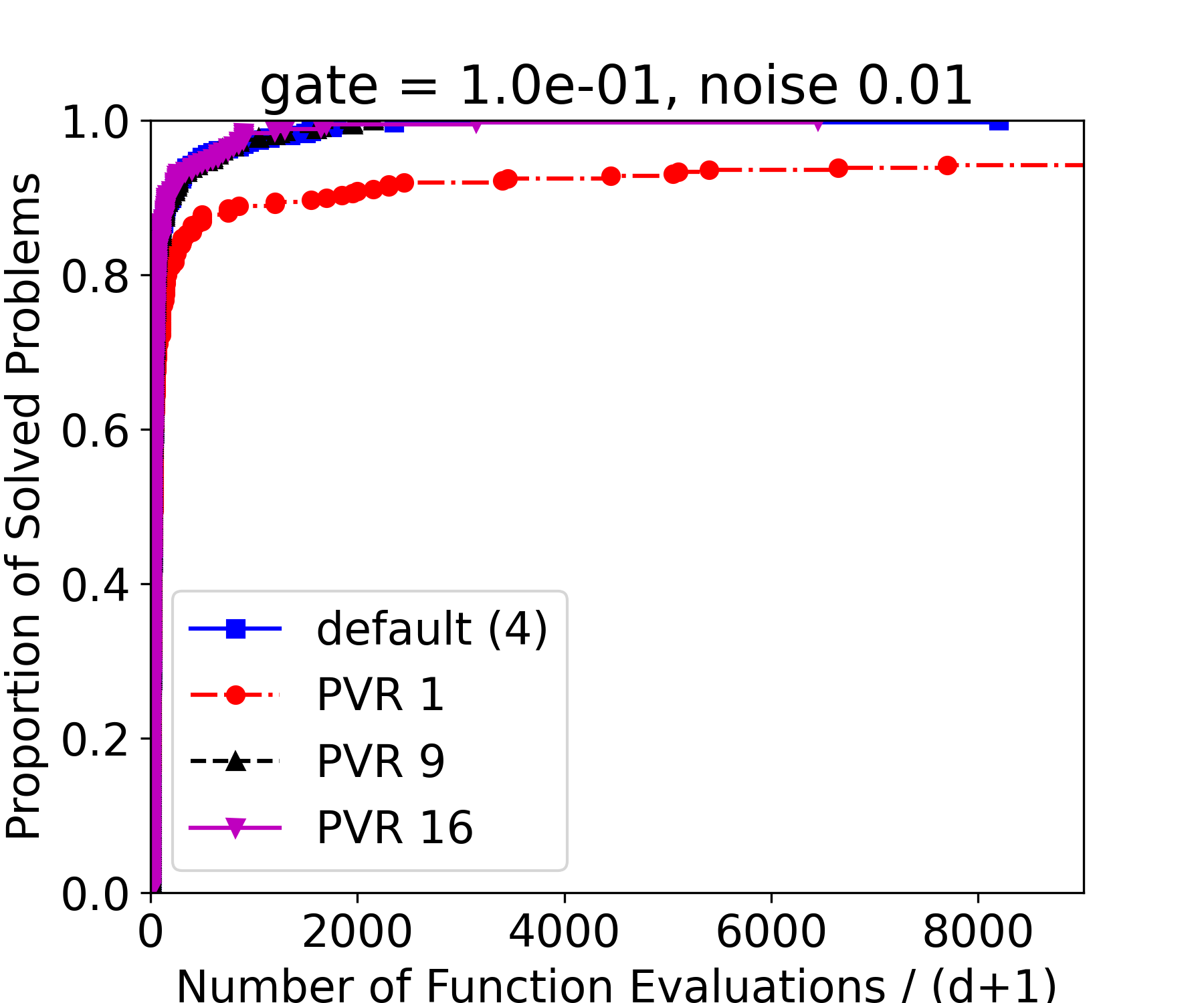}%
\includegraphics[width=0.333\textwidth, trim= 0 0 40 15, clip]{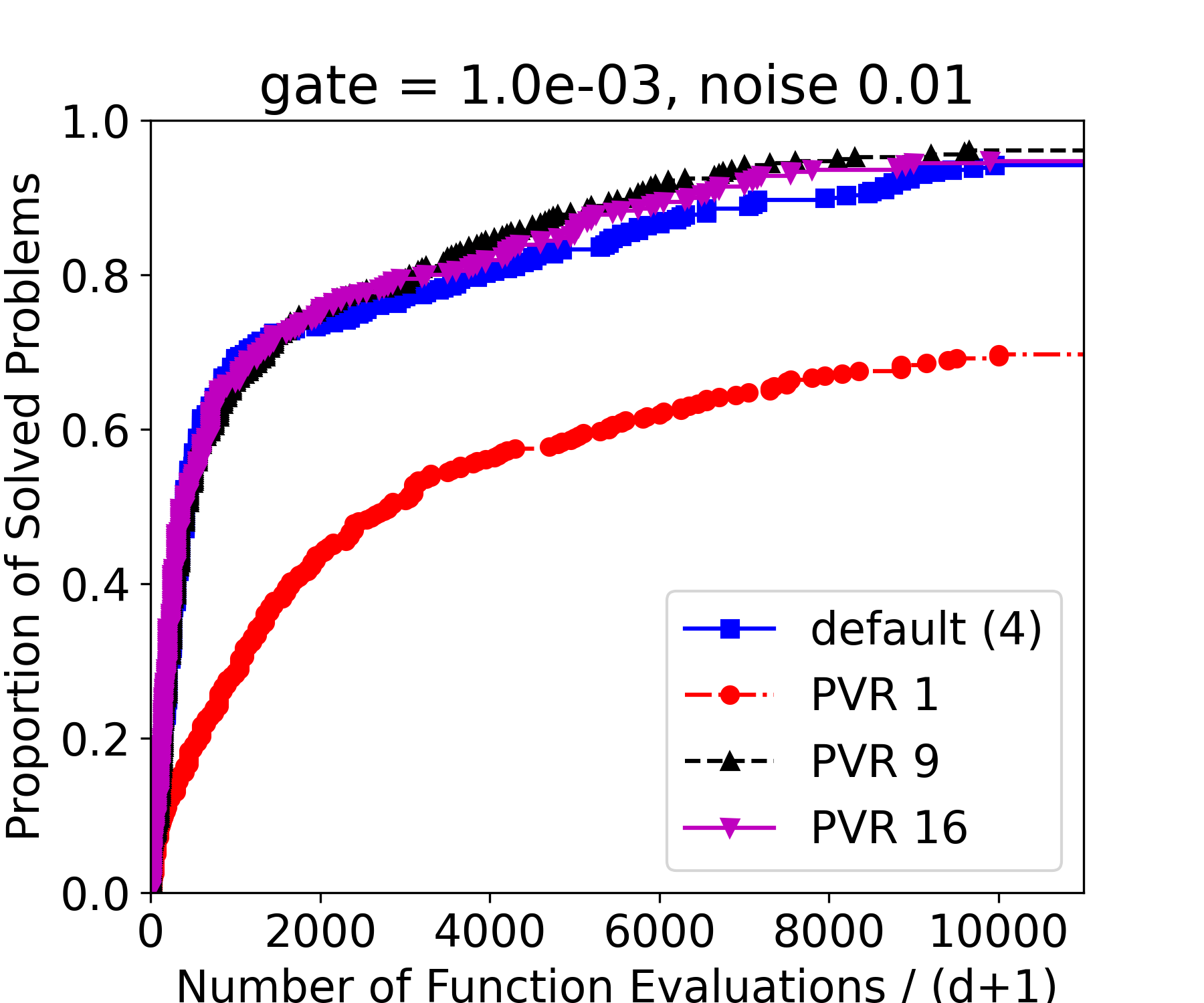}%
\includegraphics[width=0.333\textwidth, trim= 0 0 40 15, clip]{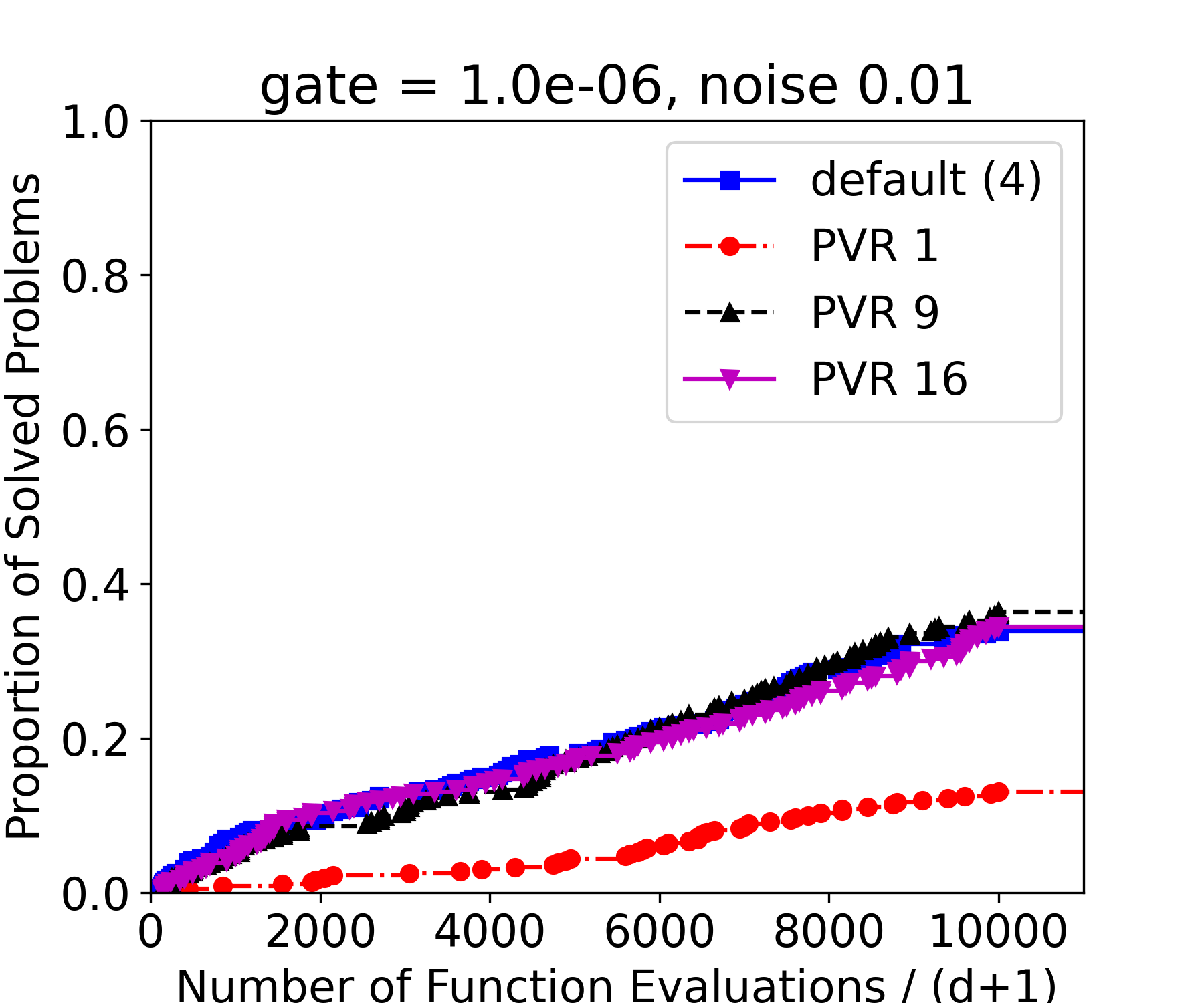}\\
\includegraphics[width=0.333\textwidth, trim= 0 0 40 15, clip]{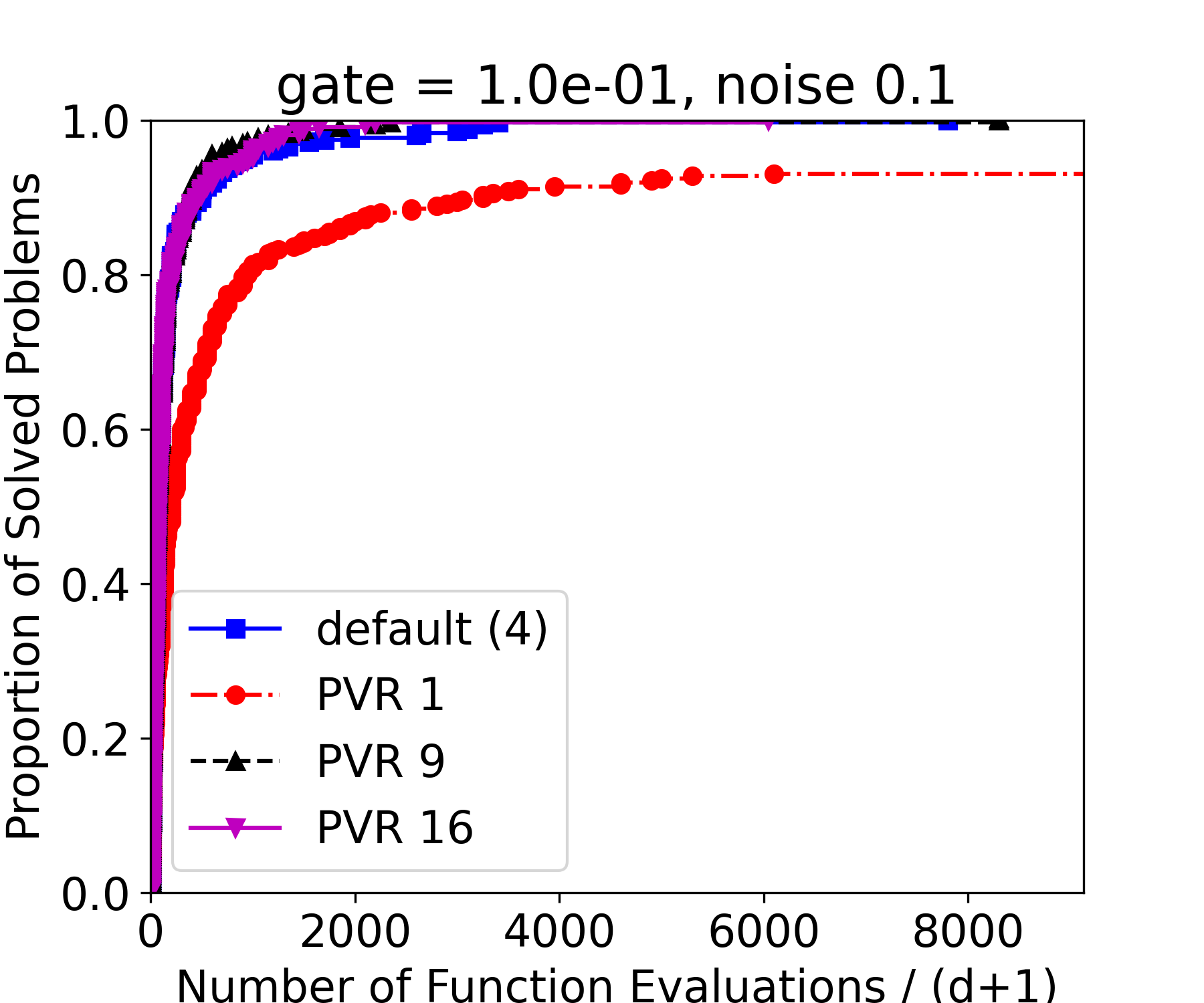}%
\includegraphics[width=0.333\textwidth, trim= 0 0 40 15, clip]{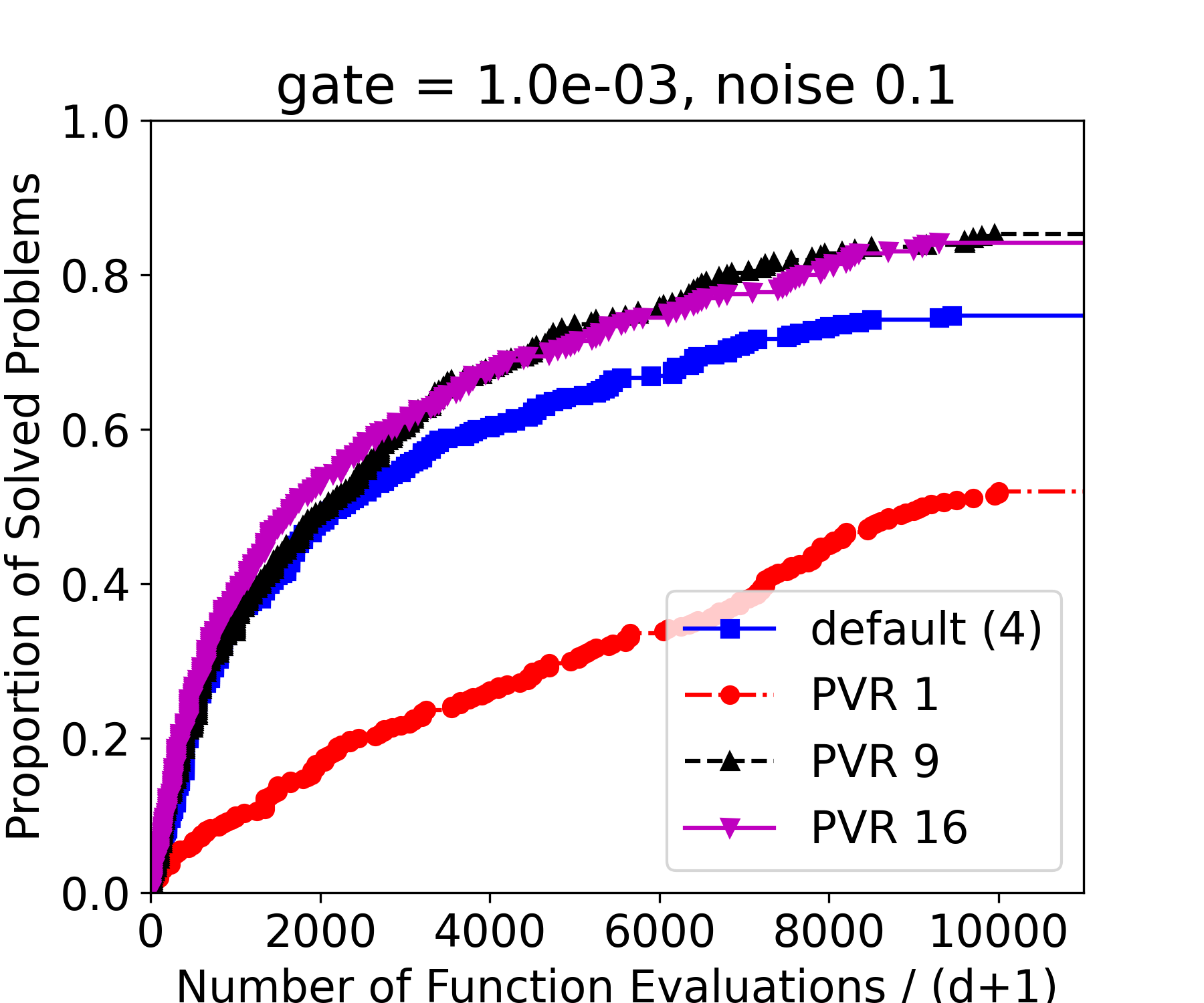}%
\includegraphics[width=0.333\textwidth, trim= 0 0 40 15, clip]{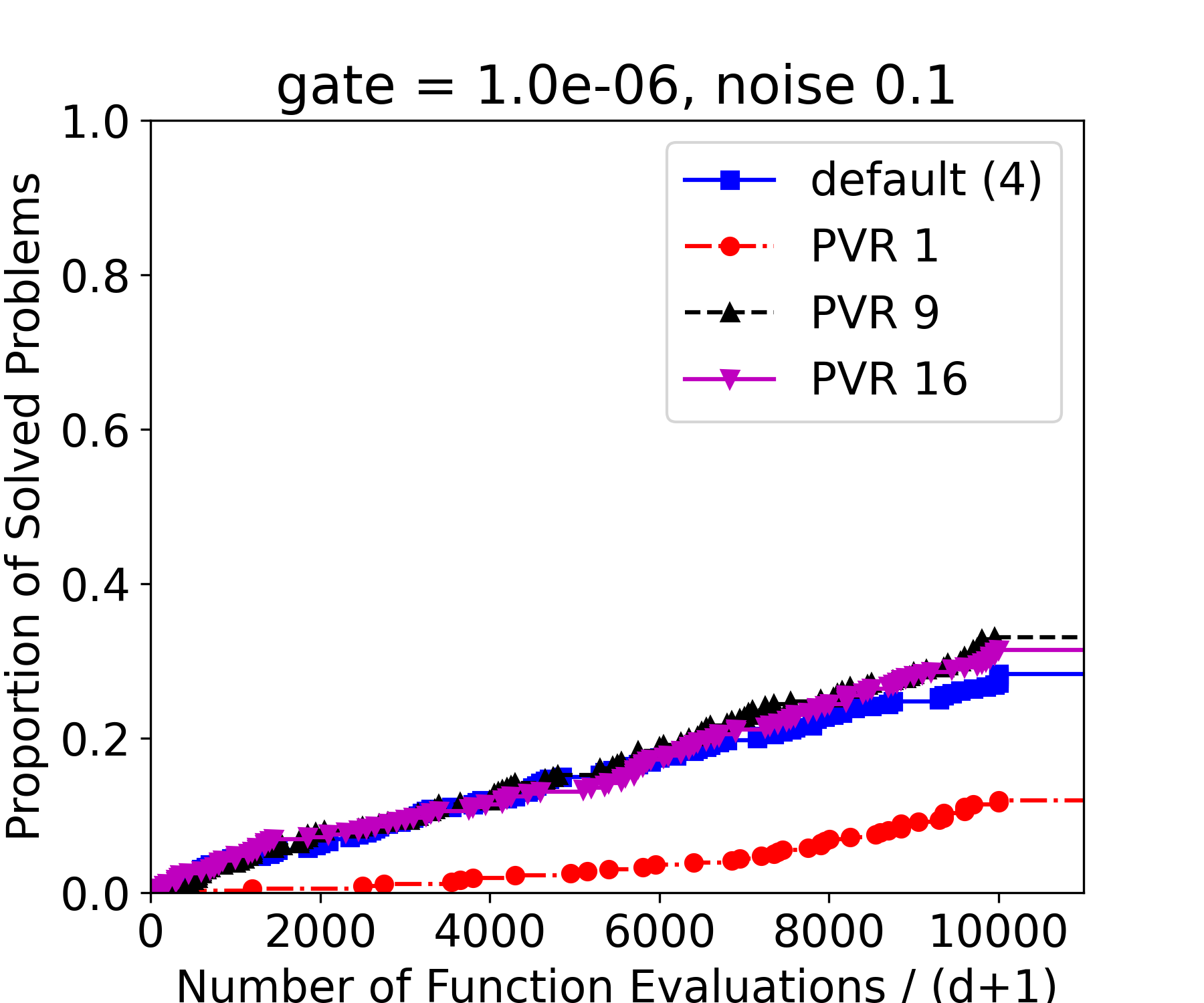}
\caption{Data profiles with varying predictive variance ratio between $\vecx_c$ and $\vecx_{n+1}$.}
\label{fig:localbenchcvred}
\end{figure}

\Cref{fig:localbenchcvred} deals with the relative ratio of predictive
 variances between the center and new point. For the largest noise variance,
 higher values $(9, 16)$ than the default (4) may be more suitable, but not too large since, overall, 9 is
 slightly dominating 16 at the latest iterations.

\section{Comparison to fixed replication budget}
\label{ap:B}

\begin{figure}
\includegraphics[width=0.333\textwidth, trim= 0 10 40 20, clip, trim= 20 35 20 20, clip]{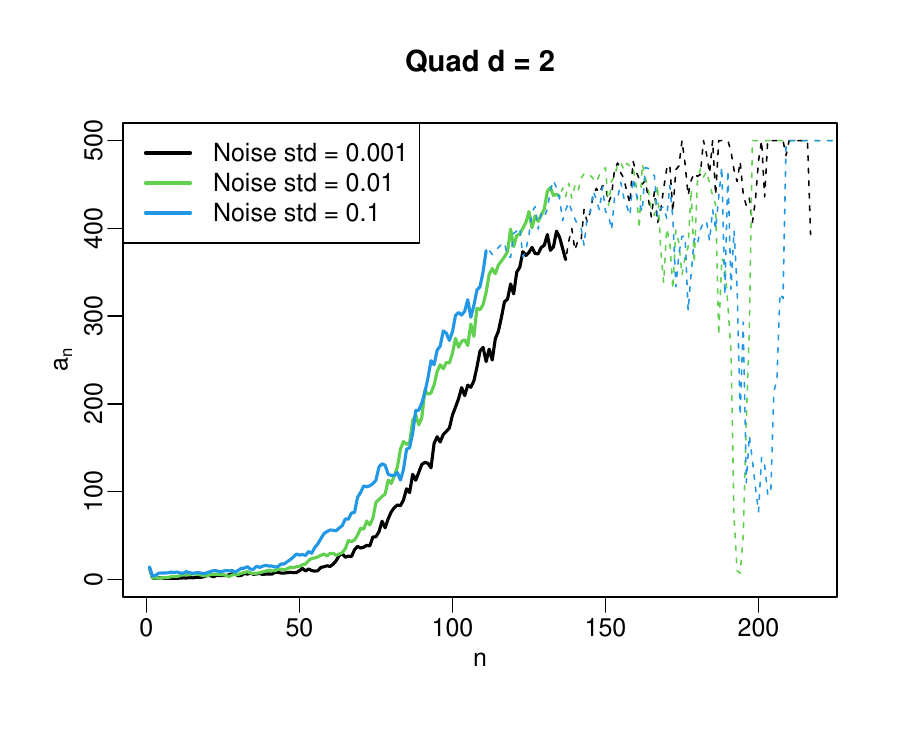}%
\includegraphics[width=0.333\textwidth, trim= 0 10 40 20, clip, trim= 20 35 20 20, clip]{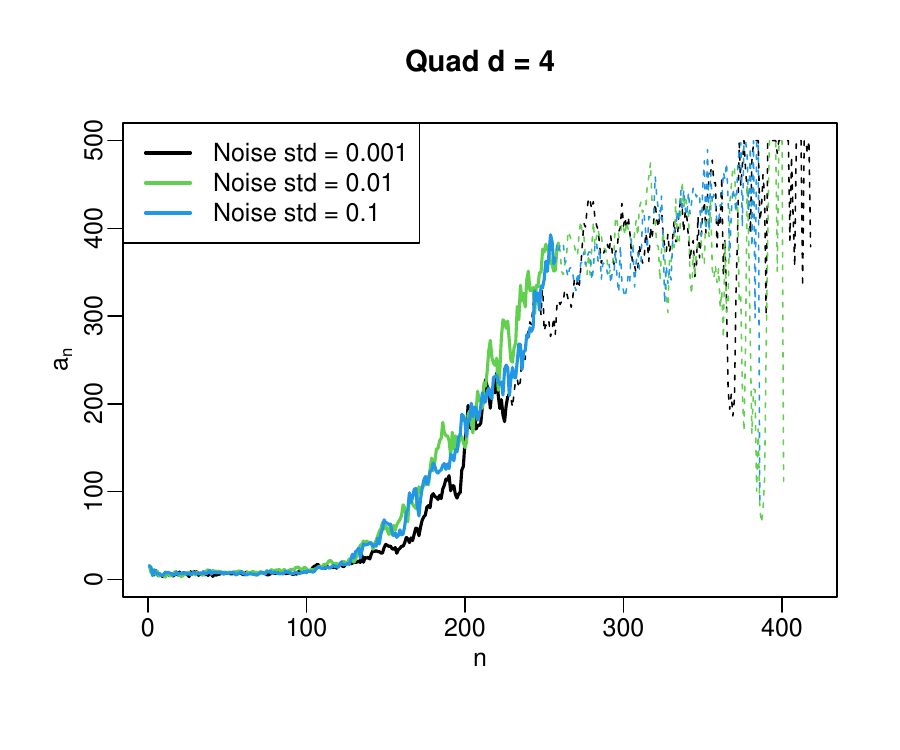}%
\includegraphics[width=0.333\textwidth, trim= 0 10 40 20, clip, trim= 20 35 20 20, clip]{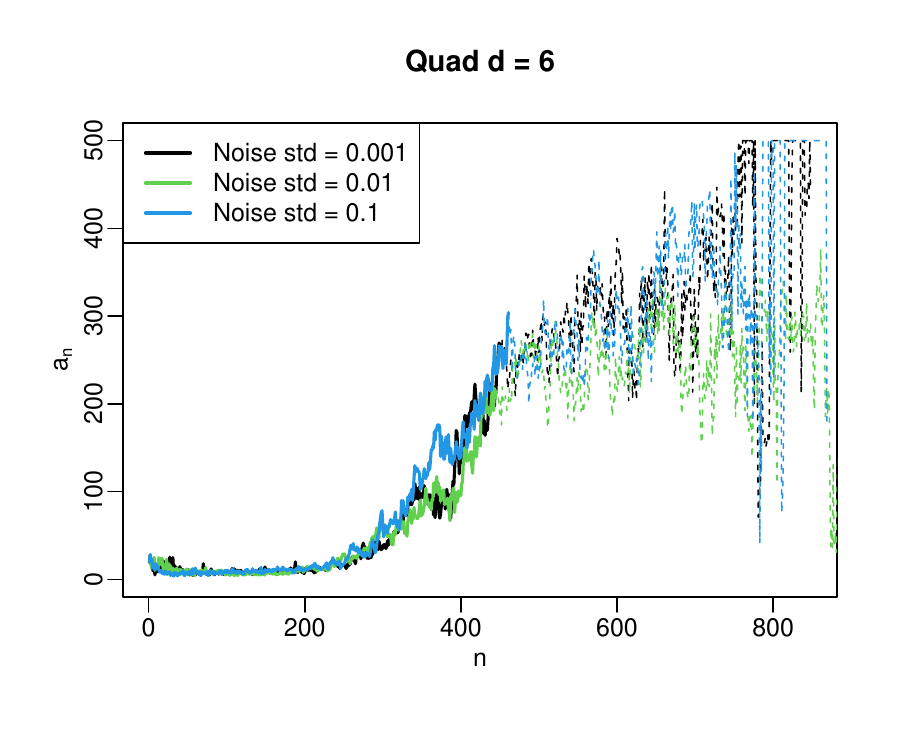}\\
\includegraphics[width=0.333\textwidth, trim= 0 10 40 20, clip, trim= 20 35 20 20, clip]{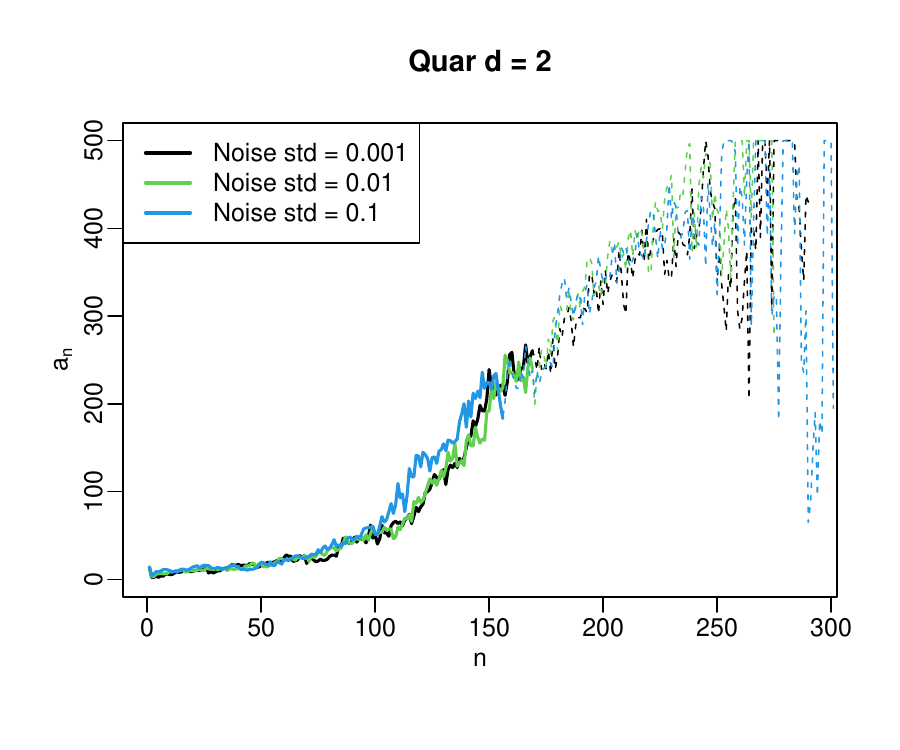}%
\includegraphics[width=0.333\textwidth, trim= 0 10 40 20, clip, trim= 20 35 20 20, clip]{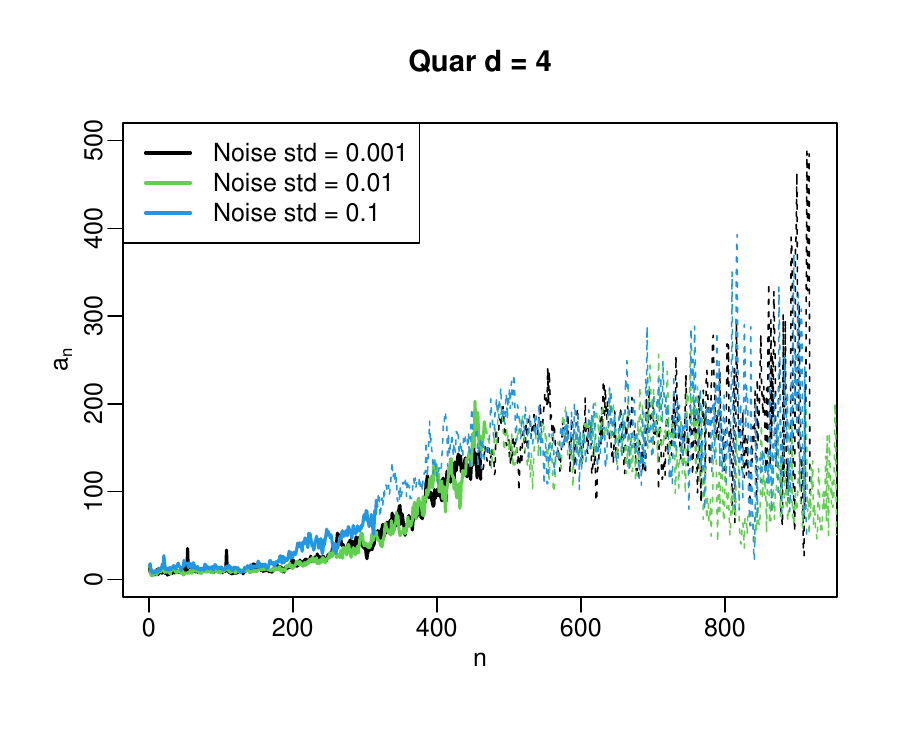}%
\includegraphics[width=0.333\textwidth, trim= 0 10 40 20, clip, trim= 20 35 20 20, clip]{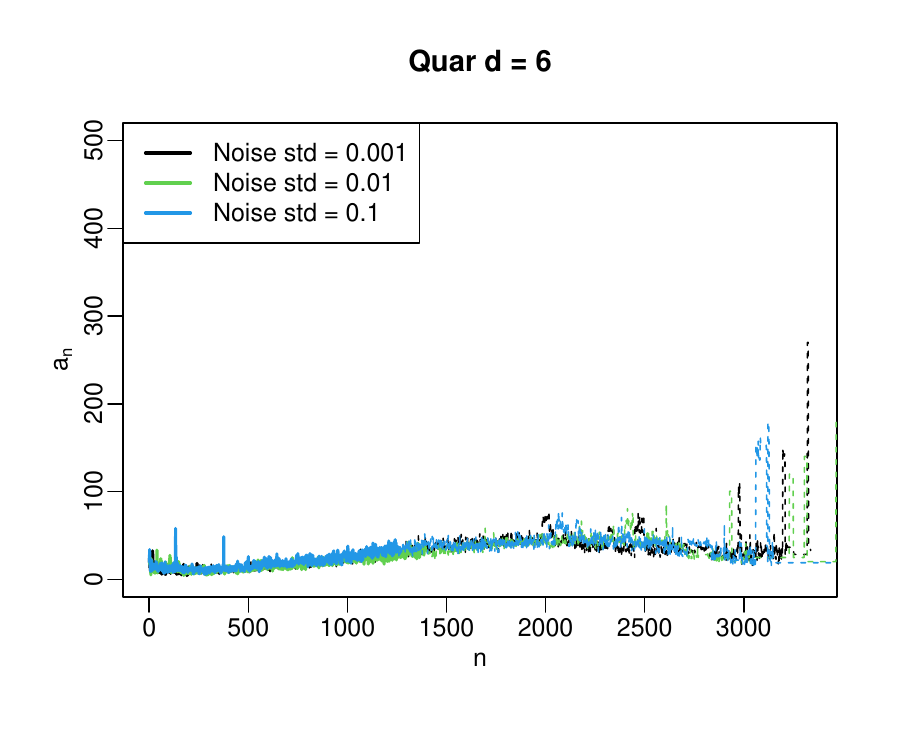}\\
\includegraphics[width=0.333\textwidth, trim= 0 10 40 20, clip, trim= 20 35 20 20, clip]{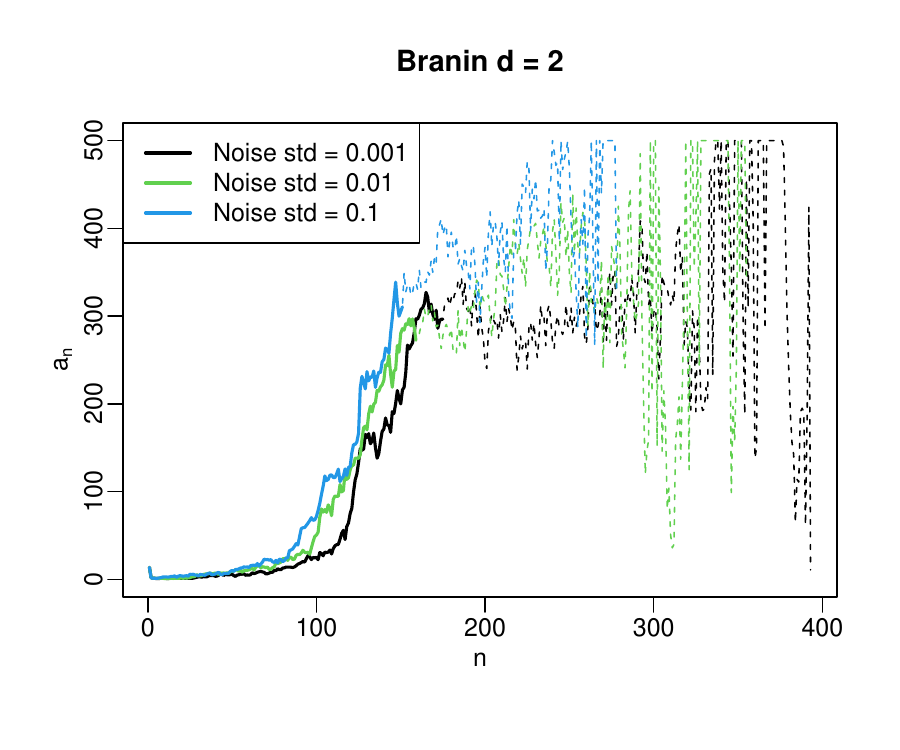}%
\includegraphics[width=0.333\textwidth, trim= 0 10 40 20, clip, trim= 20 35 20 20, clip]{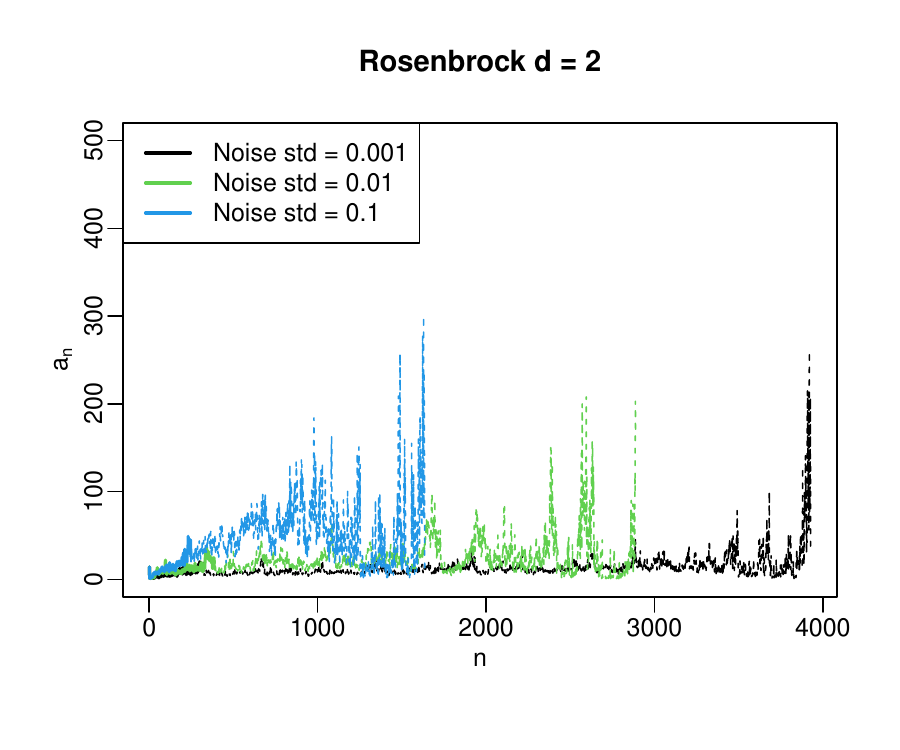}%
\includegraphics[width=0.333\textwidth, trim= 0 10 40 20, clip, trim= 20 35 20 20, clip]{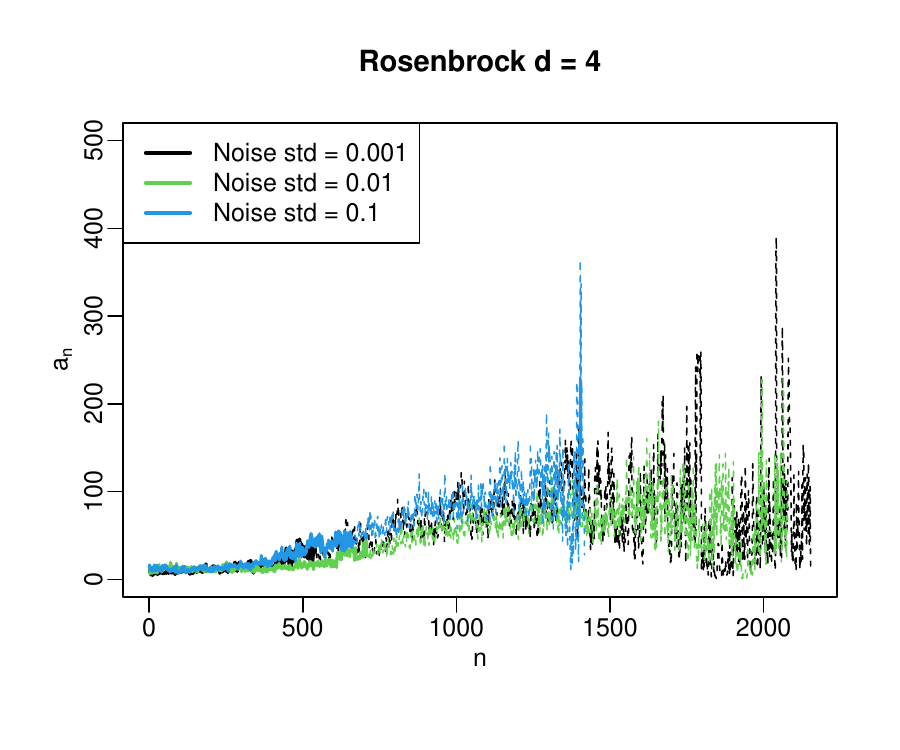}
\caption{Average number of replication at each iteration over 40 repetitions
 on the Benchmark 1 test functions. Since the total budget may be reached
 earlier if more replications are done, dashed lines indicates that less than
 40 repetitions are averaged over.}
\label{fig:repanalysis}
\end{figure}

Compared to using a fixed replication budget, we illustrate the dynamic of the
replication effort with OGPIT in \Cref{fig:repanalysis}. It depends on the
test function, but it is possible to see the effect of the dimension where
the increase in the number of replication is delayed. The effect of the noise
magnitude is visible mostly for the largest noise case. The lines are full as
long as no run finished before this iteration, then dashed lines are used. It
explains the more chaotic behavior at the latest iterations: there are less repetitions to average
over and if the budget is almost exhausted, the number of replications at the
last iteration may be small. The maximum replication budget $p_{max}$ is set
to 500, which is reached mostly on the latest iterations for a few
functions.

\begin{figure}
\includegraphics[width=0.333\textwidth, trim= 0 0 40 15, clip]{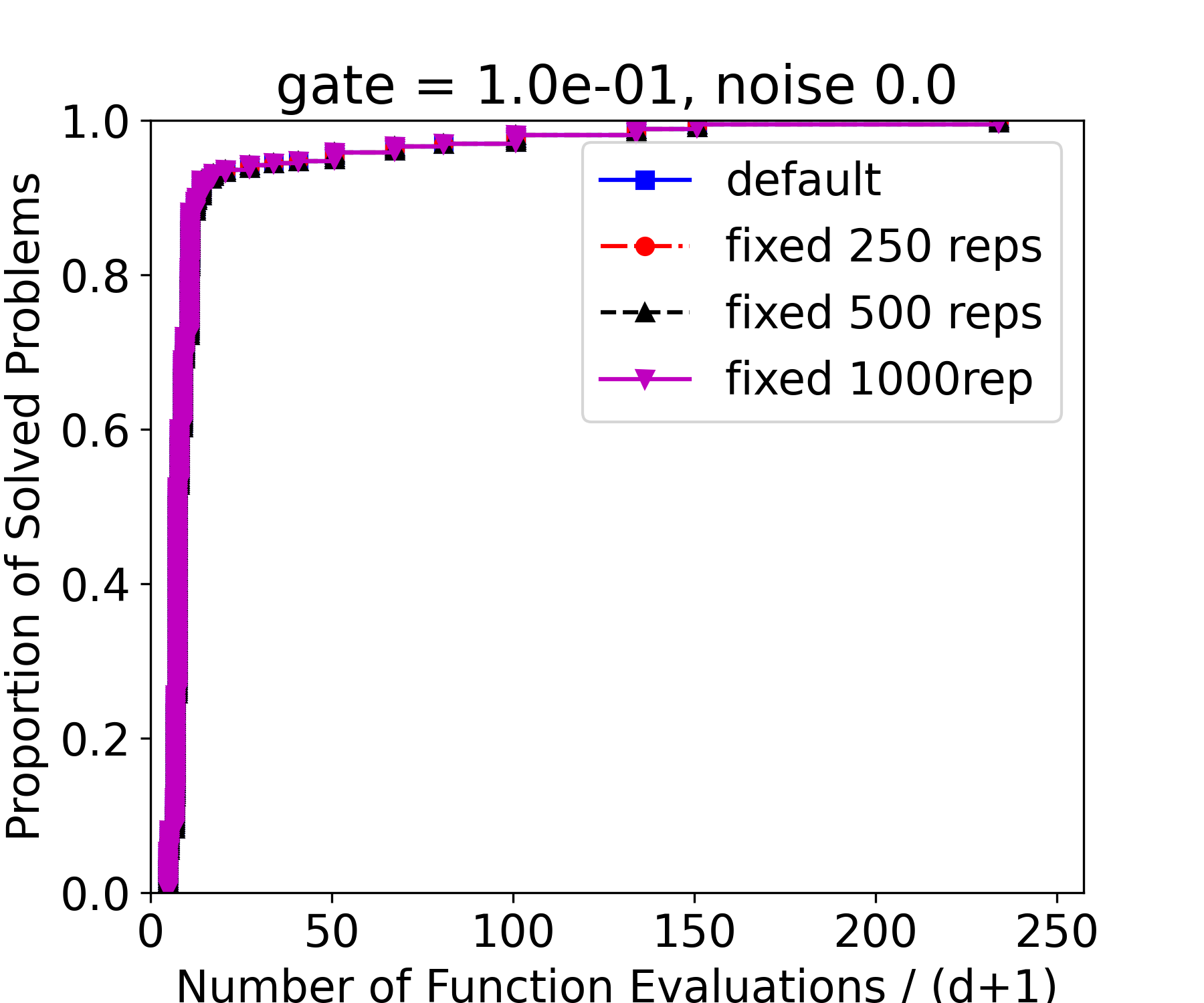}%
\includegraphics[width=0.333\textwidth, trim= 0 0 40 15, clip]{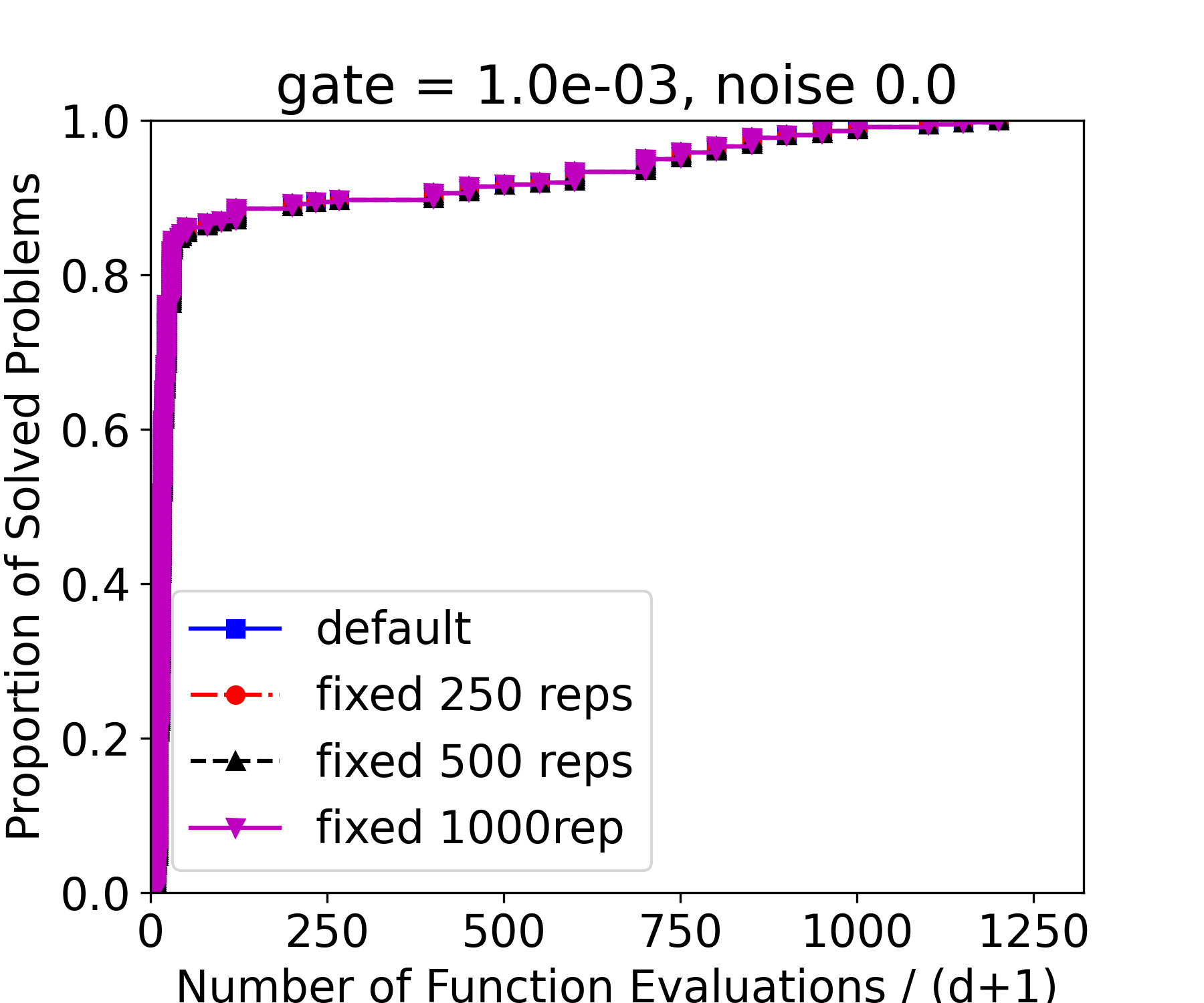}%
\includegraphics[width=0.333\textwidth, trim= 0 0 40 15, clip]{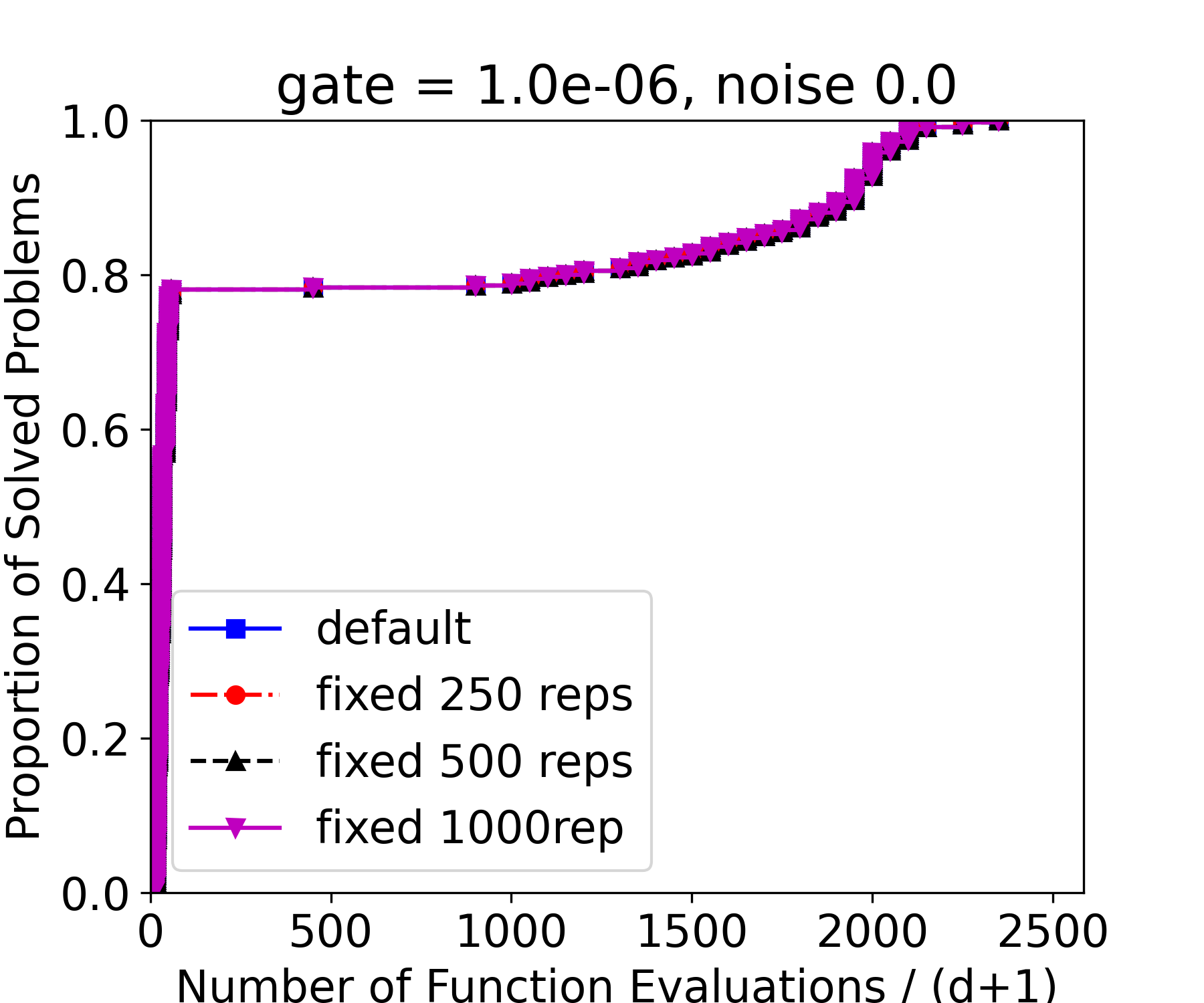}\\
\includegraphics[width=0.333\textwidth, trim= 0 0 40 15, clip]{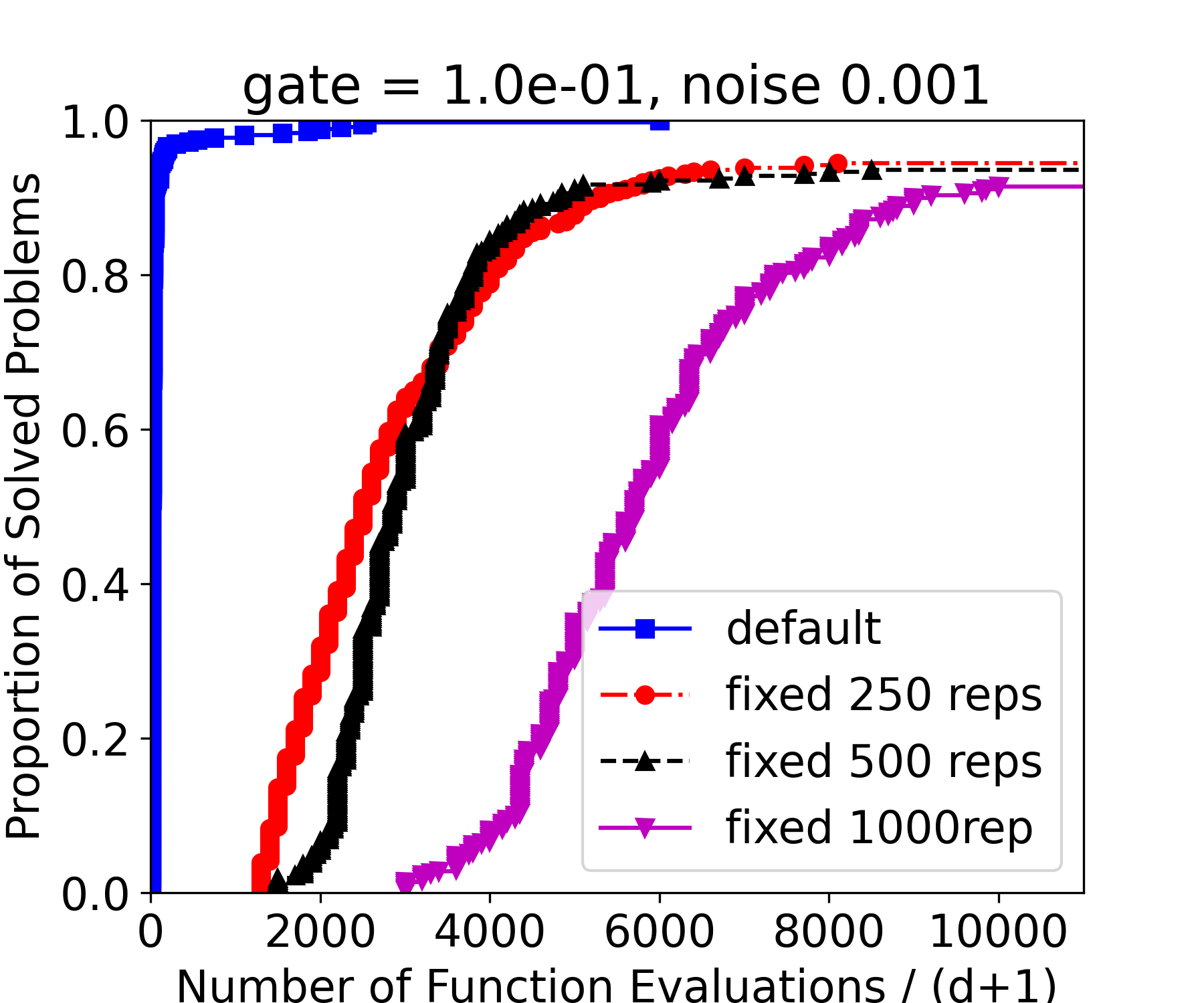}%
\includegraphics[width=0.333\textwidth, trim= 0 0 40 15, clip]{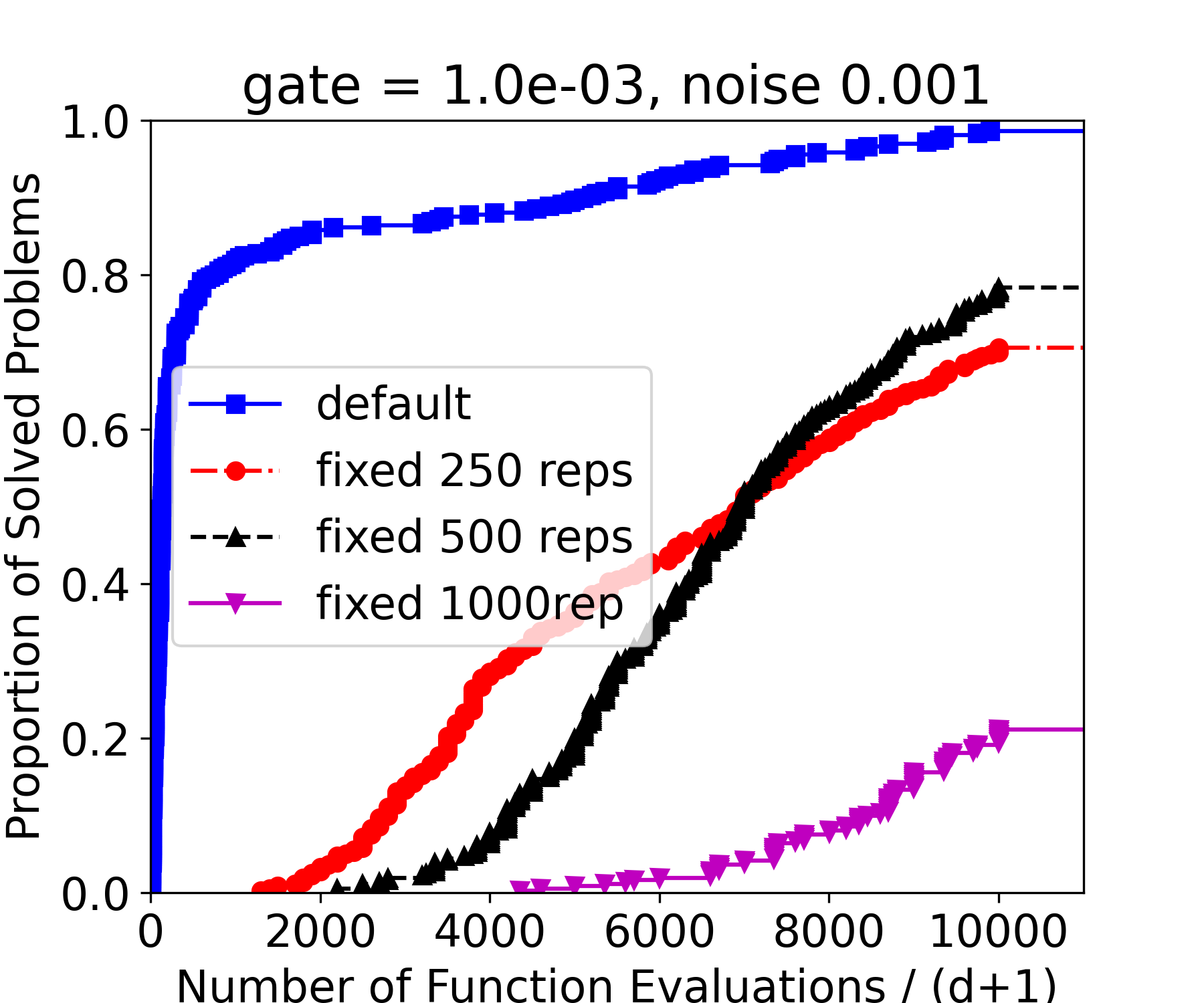}%
\includegraphics[width=0.333\textwidth, trim= 0 0 40 15, clip]{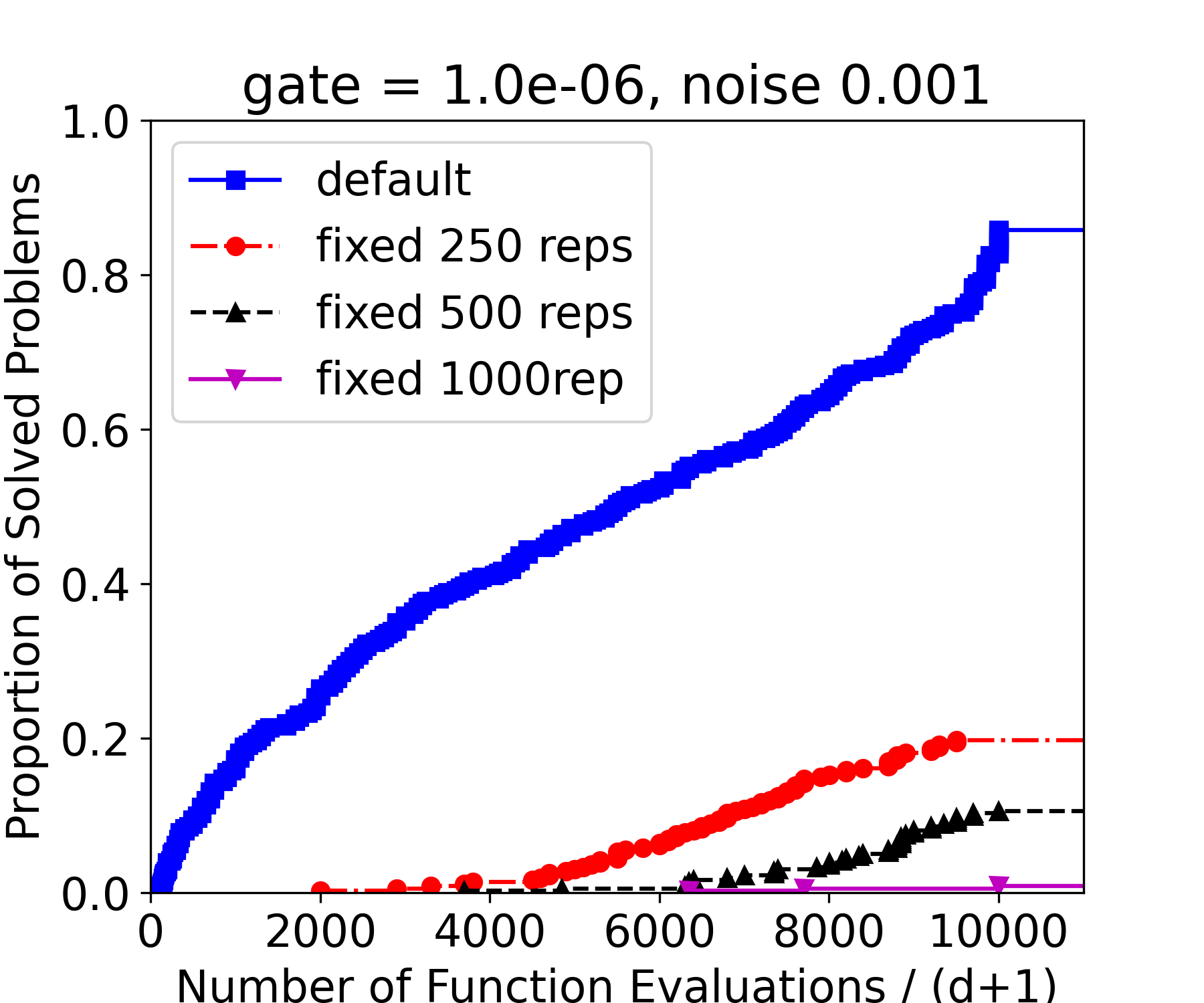}\\
\includegraphics[width=0.333\textwidth, trim= 0 0 40 15, clip]{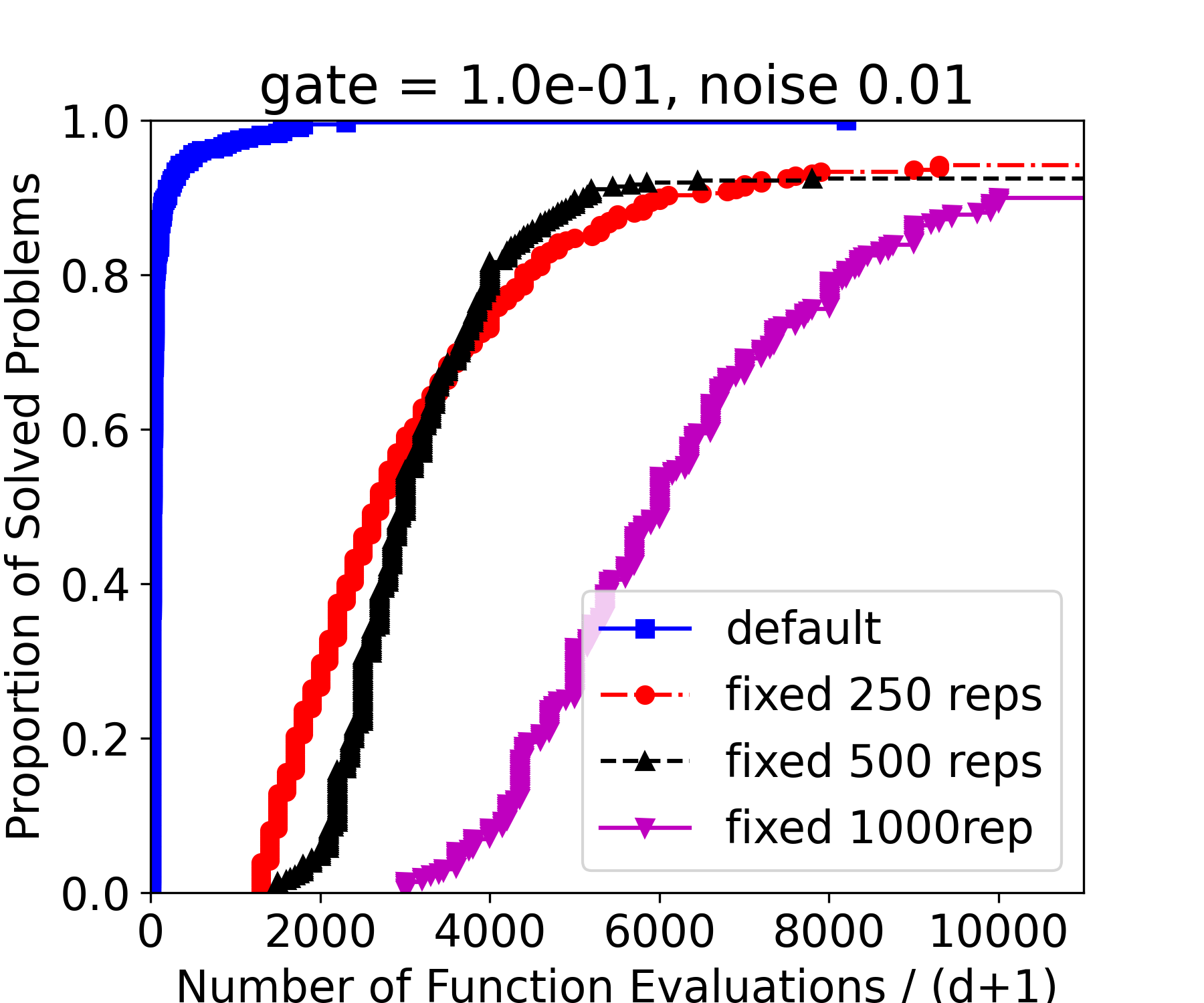}%
\includegraphics[width=0.333\textwidth, trim= 0 0 40 15, clip]{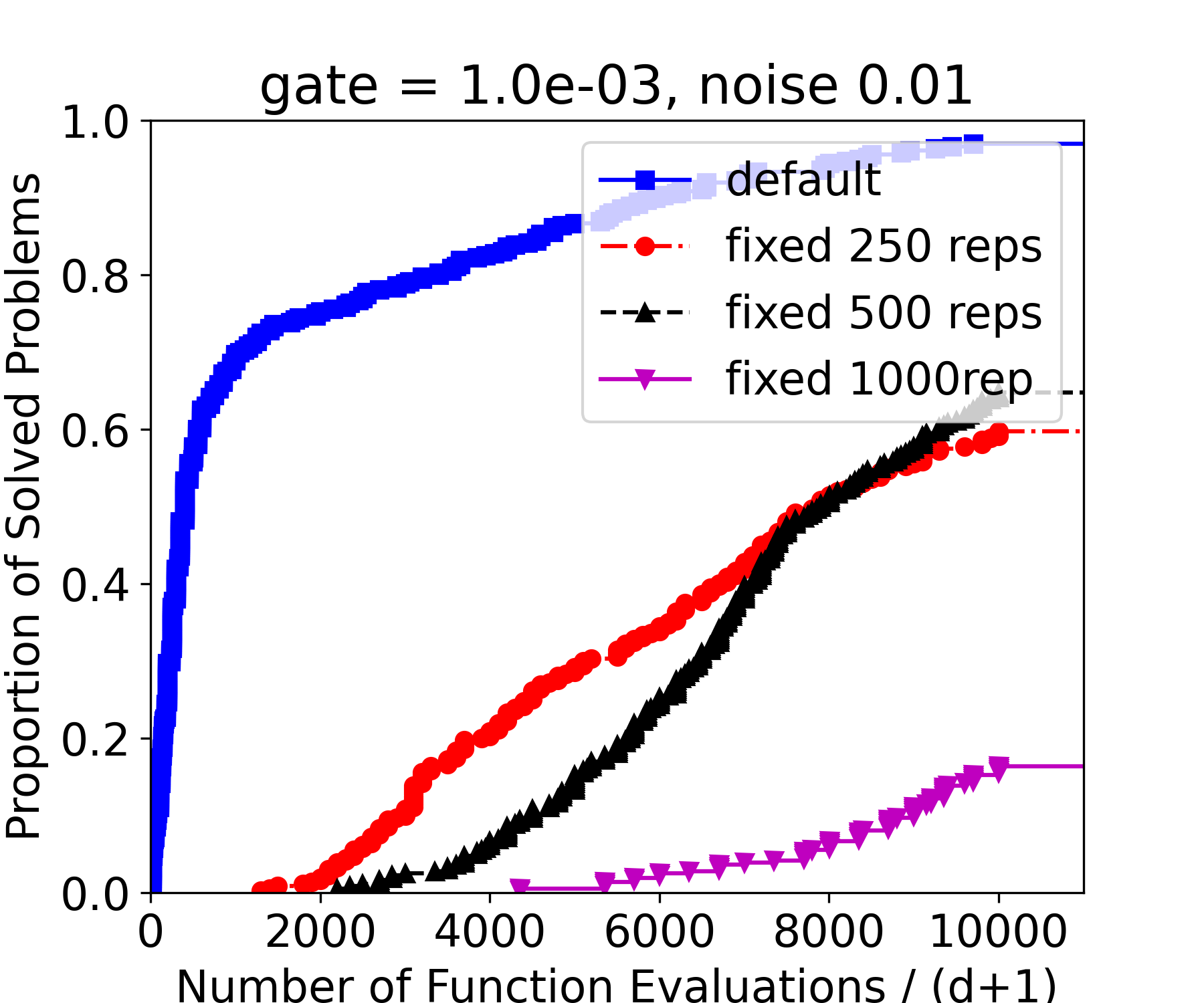}%
\includegraphics[width=0.333\textwidth, trim= 0 0 40 15, clip]{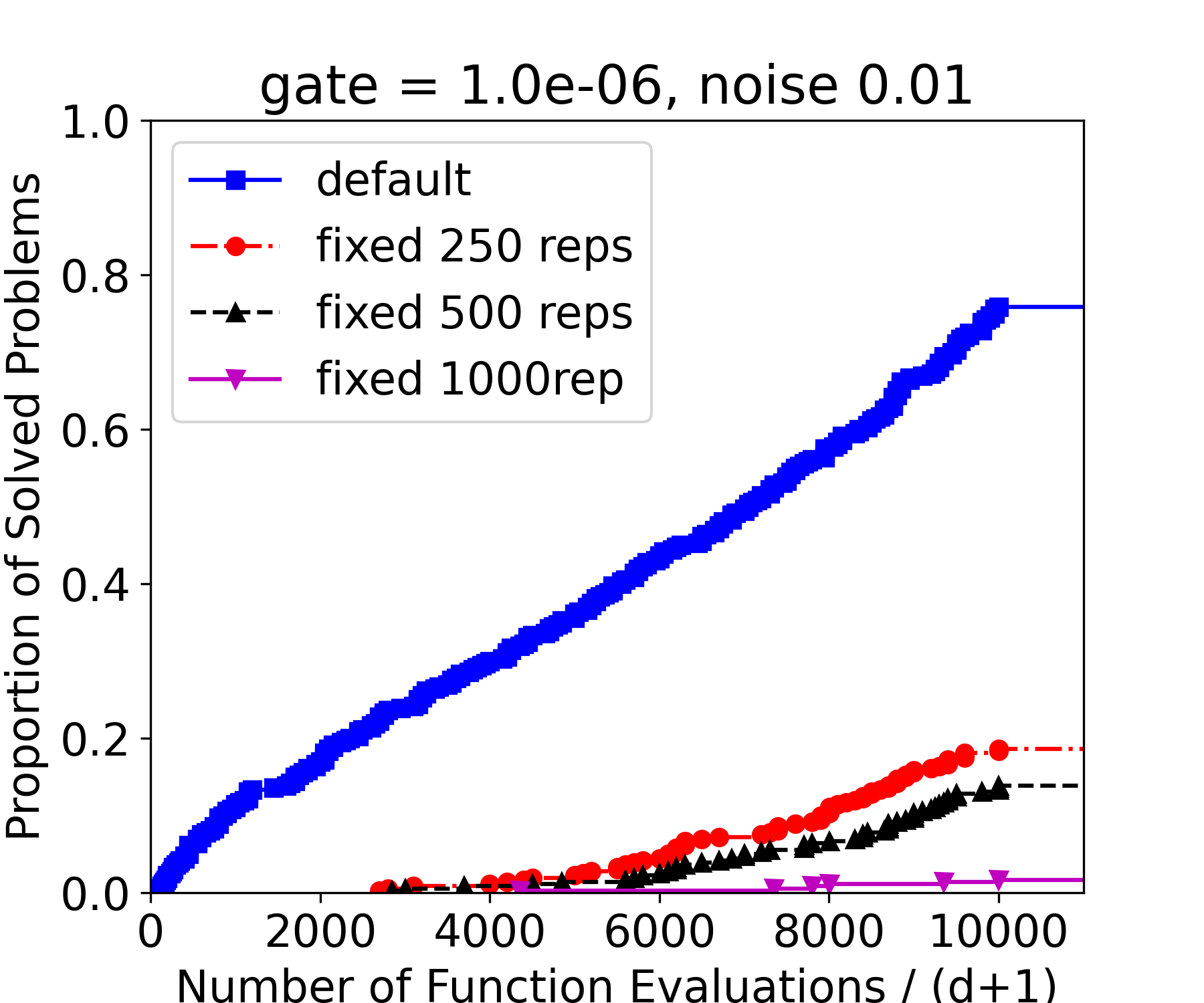}\\
\includegraphics[width=0.333\textwidth, trim= 0 0 40 15, clip]{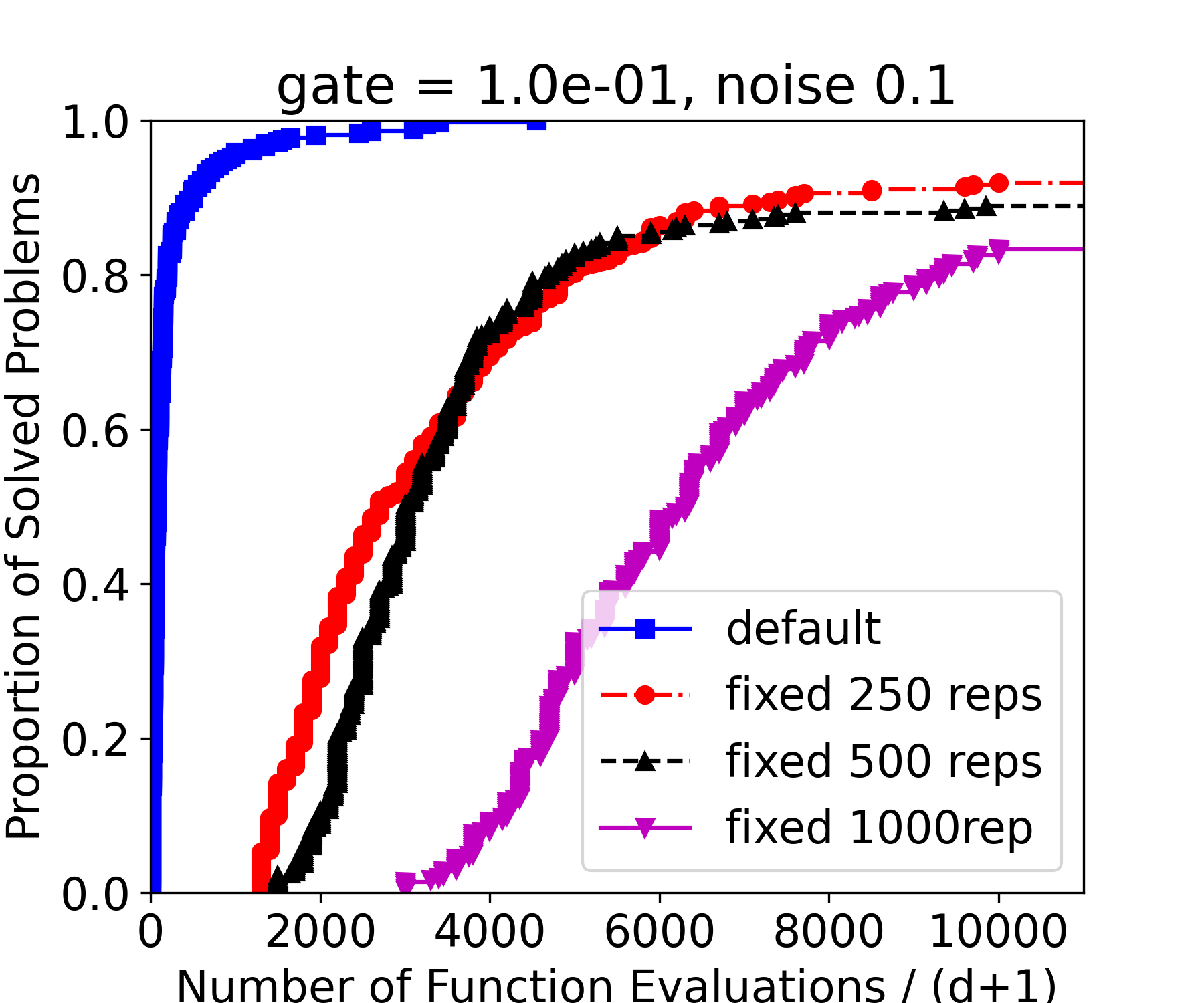}%
\includegraphics[width=0.333\textwidth, trim= 0 0 40 15, clip]{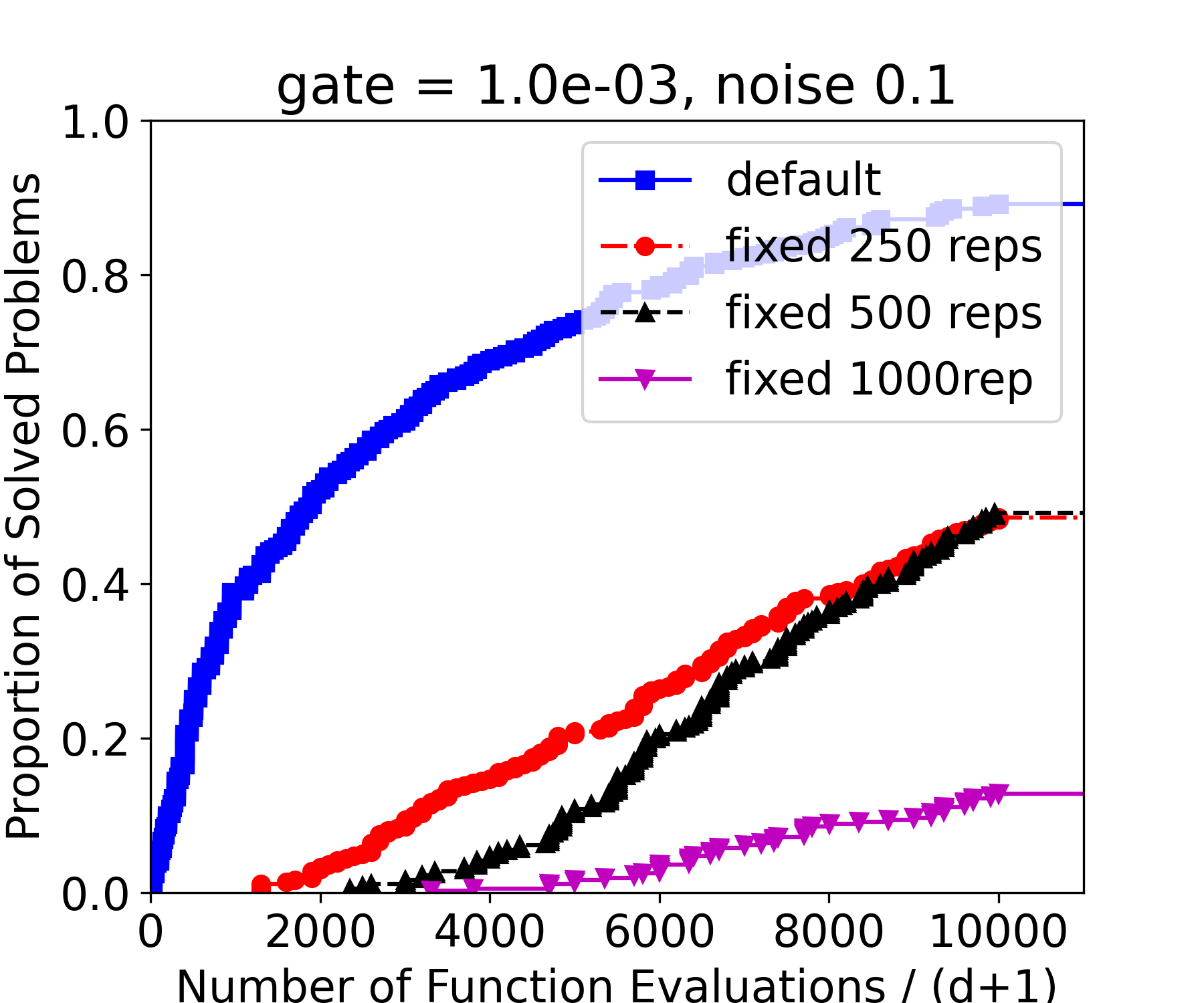}%
\includegraphics[width=0.333\textwidth, trim= 0 0 40 15, clip]{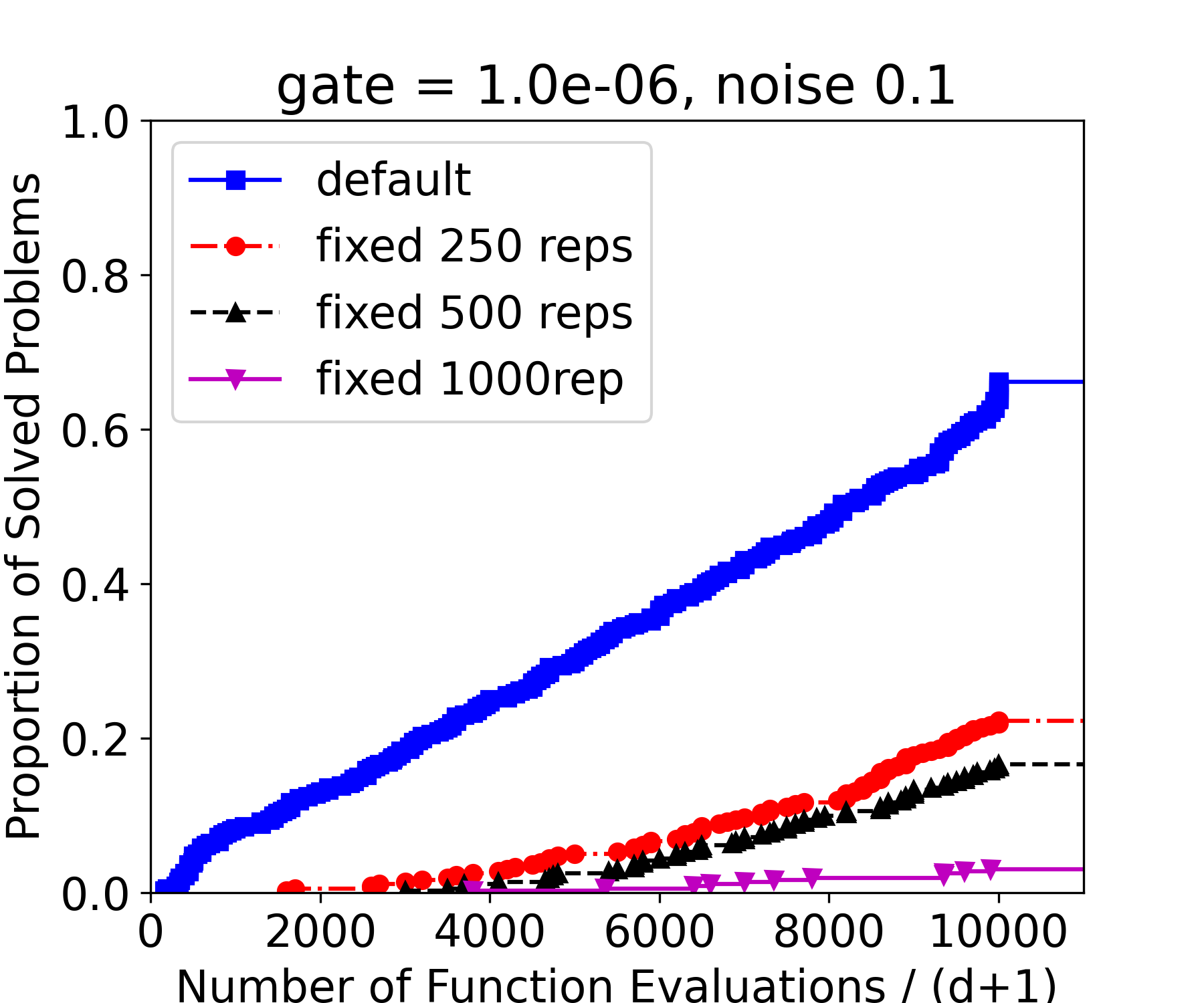}
\caption{Data profiles with varying fixed replication budgets, over 40 repetitions.}
\label{fig:fixedrep}
\end{figure}

From the previous figure, setting an appropriate fixed replication number seems difficult. In 
\Cref{fig:fixedrep}, we show that, indeed, fixing the replication budget is
 detrimental: it impacts the performance at the early iterations while more
 replicates are needed at the end (as can be seen especially in the middle
 column.)  

\newpage
\section{Starting points for SNOWPAC in Benchmark 1}

SNOWPAC starts from a single starting point rather than an initial design of experiments. These are given in \Cref{tab:startP}, in the unit hypercube.

\begin{table}[hp]
\caption{Starting points $\vecx^0$ used for SNOWPAC in Benchmark 1.}
\label{tab:startP}
\centering
\begin{tabular}{c|c|c|c}
 Problem & $2d$ Sphere & $4d$ Sphere & $6d$ Sphere\\ \hline
 $\vecx_0$ & $[0.3, 0.7]$         & $[0.3,0.7,0.2,0.8]$  & $[0.5, \cdots, 0.5]^6$\\ \hline
 Problem & $2d$ Branin    & $2d$ Rosenbrock & $4d$ Rosenbrock          \\ \hline
 $\vecx_0$ & $[0.1, 0.1]$ & $[0.2, 0.8]$  & $[0.2, 0.8, 0.2, 0.8]$
\end{tabular}
\end{table}

\end{document}